\documentclass[12pt,a4paper,leqno]{article}

\usepackage[english]{babel}
\usepackage{amscd,amsmath,amssymb,amsfonts}

\usepackage{epsf,pstricks,epsfig}
\usepackage{epic,eepic}
\usepackage{graphicx}

\newcommand{\myfig}[2]{
\begin{center}
{\epsfxsize=#2\hsize \epsfbox{#1}}\nobreak
\end{center}}

\numberwithin{equation}{section}

\newtheorem{theorem}[subsection]{Theorem}
\newtheorem{proposition}[subsection]{Proposition}
\newtheorem{lemma}[subsection]{Lemma}
\newtheorem{remark}[subsection]{Remark}
\newtheorem{corollary}[subsection]{Corollary}

\newtheorem{definition}[subsection]{Definition}

\newtheorem{claim}[subsection]{Claim}
\input xy
\xyoption{all}

\newcommand{\Spec}{\mathrm{Spec}}
\newcommand{\Proj}{\mathrm{Proj}}
\newcommand{\Dir}{\mathrm{Dir}}
\newcommand{\Char}{\mathrm{char}}
\newcommand{\gr}{\mathrm{gr}}
\newcommand{\bw}{\mathbf{w}}
\newcommand{\bv}{\mathbf{v}}
\newcommand{\fp}{\mathfrak{p}}
\newcommand{\fm}{\mathfrak{m}}
\newcommand{\Lat}{\mathrm{Lat}}

\date{18th February 2013}

\begin{document}


\title
{Canonical embedded and non-embedded resolution of singularities for excellent two-dimensional schemes}
\author{Vincent Cossart and Uwe Jannsen and Shuji Saito }
\maketitle

\def\ml{\{\{}
\def\mr{\}\}}

\def\reg{{\mathrm{reg}}}
\def\sing{{\mathrm{sing}}}
\def\qreg{{\mathrm{qreg}}}
\def\sqreg{\mathrm{sqreg}}
\def\Xreg{X_{\reg}}
\def\XBreg{X_{\cB\reg}}
\def\XnBreg{X_{n,\cB\reg}}
\def\XiBreg{X_{i,\cB\reg}}
\def\XdBdreg{X'_{\cB'\reg}}
\def\Xsing{X_{\sing}}
\def\XBsing{X_{\cB\sing}}

\def\sreg{{\mathrm{sreg}}}
\def\XBsreg{X_{\cB\sreg}}
\def\XnBsreg{X_{n,\cB\sreg}}
\def\XiBsreg{X_{i,\cB\sreg}}
\def\XdBdsreg{X'_{\cB'\sreg}}

\def\Xqreg{X_{\mathrm{qreg}}}
\def\XBqreg{X_{\cB\qreg}}
\def\XBsqreg{X_{\cB\sqreg}}
\def\XiBsqreg{(X_i)_{\cB_i\sqreg}}

\def\tu{\tilde{u}}
\def\tU{\tilde{U}}
\def\bv{{\bf v}}
\def\bw{{\bf w}}
\def\bwp{{\bf w}^+}
\def\bwm{{\bf w}^-}
\def\Lob{\Lambda^{\OB}}
\def\Lamob{\Lambda^{\OB}}
\def\inLob{in_\delta}

\def\qdiB{q'_{i,B}}
\def\qbdiB{\overline{q'_{i,B}}}
\def\Reta{R_{\eta}}
\def\fmeta{\fm_{\eta}}
\def\Rdeta{R'_{\eta}}
\def\Rdetad{R'_{\eta'}}
\def\fmdeta{\fm'_{\eta}}
\def\fmdetad{\fm'_{\eta'}}
\def\Retad{R_{\eta'}}
\def\fmetad{\fm_{\eta'}}

\def\WP{well-prepared }
\def\WPness{well-preparedness }

\def\projlim#1{\underset{#1}{\varprojlim}}

\def\hf{\hat{f}}
\def\hy{\hat{y}}
\def\hLam{\hat{\Lambda}}
\def\tFidg{\tilde{F}'_{i,\Gamma'}}

\def\aol{\overline{\alpha}}
\def\bol{\overline{\beta}}

\def\ol#1{\overline{#1}}
\def\k{\kappa}

\def\fm{{\mathfrak m}}
\def\fp{{\mathfrak p}}
\def\fq{{\mathfrak q}}
\def\fP{{\mathfrak P}}
\def\fQ{{\mathfrak Q}}
\def\fa{{\mathfrak a}}
\def\frb{{\mathfrak b}}

\def\gr{{\mathrm{gr}}}

\def\grq{{\gr}_{\fq}}
\def\grp{{\gr}_{\fp}}
\def\grm{{\gr}_{\fm}}
\def\grmd{{\gr}_{\fm'}}

\def\grqR{{\gr}_{\fq}(R)}
\def\grpR{{\gr}_{\fp}(R)}
\def\grmR{{\gr}_{\fm}(R)}
\def\grmdRd{{\gr}_{\fm'}(R')}

\def\grqRJ{{\gr}_{\fq}(R/J)}
\def\grpRJ{{\gr}_{\fp}(R/J)}
\def\grmRJ{{\gr}_{\fm}(R/J)}

\def\InqJ{In_{\fq}(J)}
\def\InpJ{In_{\fp}(J)}
\def\InmJ{In_{\fm}(J)}
\def\InmJd{In_{\fm'}(J')}
\def\InLJ{In_L(J)}
\def\InhLJ{\widehat{In}_L(J)}
\def\Inh{\widehat{In}}

\def\inq{in_{\fq}}
\def\inp{in_{\fp}}
\def\inm{in_{\fm}}
\def\inmd{in_{\fm'}}

\def\bNN{\bN^{\bN}}
\def\bNNhf{\bN^{\bN}_{hf}}
\def\bNNhp{\bN^{\bN}_{hp}}
\def\bNNsq{\bNN\times\bN}

\def\Bl{B\ell}
\def\eb{\overline{e}}
\def\qaq{\quad\text{ and }\quad}
\def\qwith{\quad\text{ with }}
\def\qfor{\quad\text{ for }}

\def\hR{\hat{R}}
\def\hJ{\hat{J}}
\def\hfm{\hat{\fm}}

\def\rmapo#1{\overset{#1}{\longrightarrow}}
\def\rmapu#1{\underset{#1}{\longrightarrow}}
\def\lmapu#1{\underset{#1}{\longleftarrow}}
\def\rmapou#1#2{\overset{#1}{\underset{#2}{\longrightarrow}}}
\def\lmapo#1{\overset{#1}{\longleftarrow}}
\def\isom{\stackrel{\sim}{\rightarrow}}

\def\fm{{\mathfrak m}}
\def\fn{{\mathfrak n}}
\def\fp{{\mathfrak p}}

\def\HS{{\mathrm{HS}}}
\def\trdeg{{\mathrm{trdeg}}}
\def\codim{{\mathrm{codim}}}
\def\Sym{{\mathrm{Sym}}}
\def\Spec{{\mathrm{Spec}}}
\def\Ker{{\mathrm{Ker}}}
\def\rank{{\mathrm{rank}}}
\def\Lat{{\mathrm{Lat}}}
\def\div{{\mathrm{div}}}

\def\cB{{\mathcal B}}
\def\cC{{\mathcal C}}
\def\cO{{\mathcal O}}
\def\cX{{\mathcal X}}
\def\cV{{\mathcal V}}
\def\cY{{\mathcal Y}}
\def\cT{{\mathcal T}}

\def\k{\kappa}

\def\bB{{\mathbb B}}
\def\bC{{\mathbb C}}
\def\bR{{\mathbb R}}
\def\bZ{{\mathbb Z}}
\def\bQ{{\mathbb Q}}
\def\bG{{\mathbb G}}
\def\bL{{\mathbb L}}
\def\bN{{\mathbb N}}
\def\bP{{\mathbb P}}
\def\bF{{\mathbb F}}
\def\bH{{\mathbb H}}
\def\bA{{\mathbb A}}

\def\qaq{\quad\text{ and }\quad}
\def\qwith{\quad\text{ with }}
\def\qfor{\quad\text{ for }}

\def\fb{\overline{f}}

\def\srp{system of regular parameters }
\def\snc{simple normal crossing divisor }

\def\BB#1{B^{(#1)}}
\def\BBd#1{{B'}^{(#1)}}
\def\BBo#1{B_0^{(#1)}}
\def\BBa#1{B_1^{(#1)}}
\def\BBq#1{B_q^{(#1)}}
\def\BBm#1{B_m^{(#1)}}
\def\BBn#1{B_n^{(#1)}}
\def\BBq#1{B_q^{(#1)}}
\def\BBu#1{B_u^{(#1)}}
\def\BBth#1{B_{\theta}^{(#1)}}
\def\BBthd#1{{B'_{\theta}}^{(#1)}}
\def\BBud#1{{B'_u}^{(#1)}}
\def\BBvd#1{{B'_v}^{(#1)}}
\def\BBv#1{B_v^{(#1)}}

\def\BZ{B_Z}
\def\BZZ#1{\BZ^{(#1)}}
\def\BZd{B_{Z'}}
\def\BZZd#1{\BZd^{(#1)}}

\def\Ix{I(x)}
\def\Iy{I(y)}
\def\Ixd{I(x')}

\def\NI{N}
\def\NIx{\NI(x)}
\def\NIxD{\NI(x)_D}
\def\NIy{\NI(y)}
\def\NIxd{\NI(x')}
\def\NIdxd{\NI'(x')}

\def\tB{\tilde{B}}
\def\Bx{B(x)}
\def\Bxd{B(x')}
\def\Bdxd{B'(x')}
\def\By{B(y)}
\def\Bdyd{B'(y')}

\def\oXx{o_X(x)}
\def\oXy{o_X(y)}

\def\coXx{{\mid O(x)\mid}}
\def\coXy{{\mid O(y)\mid}}

\def\OB{O}
\def\OBd{{\OB'}}
\def\oXdxd{o_{X'}(x')}
\def\oXy{o_{X'}(y')}
\def\OBx{O(x)}
\def\tOBx{\widetilde{O(x)}}
\def\OBy{O(y)}
\def\OBxd{O(x')}
\def\OBdxd{O'(x')}
\def\OBdyd{O'(y')}
\def\OBxq{O(x_q)}
\def\OBxX{O_X(x)}
\def\OBxq{O(x_q)}
\def\OBxX{O_X(x)}
\def\OBxXm#1{O_{X_{#1}}(x_{#1})}
\def\OBxXn#1{O_{X_{#1}}(x_{#1})}
\def\OBxdXd{O_{X'}(x')}
\def\OByX{O(x')}

\def\rOB{\overline{O}}
\def\rOBd{\overline{O'}}

\def\OBX#1{OB_X(#1)}
\def\OBXm#1#2{OB_{X_{#1}}(#2)}
\def\OBYm#1#2{OB_{Y_{#1}}(#2)}
\def\OBXnm#1#2{OB_{X_{#1}}(#2)}
\def\OBYn#1#2{OB_{Y_{#1}}(#2)}

\def\OBy{O(y)}
\def\NBx{N(x)}
\def\NBy{N(y)}

\def\NBdxd{N'(x')}

\def\tbB{\widetilde{\bB}}
\def\bBx{\bB(x)}

\def\cBx{{\mathcal B(x)}}
\def\cBy{{\mathcal B(y)}}
\def\cBdxd{{\cB'(x')}}
\def\cBdyd{{\cB'(y')}}
\def\cBI{{\cB_{in}}}

\def\cBxD{{\mathcal B}(x,D)}

\def\cOB{{O}}
\def\tOB{\widetilde{O}}
\def\cOBx{{O(x)}}
\def\cOBy{{O(y)}}
\def\cNBx{{N(x)}}

\def\cOBdxd{\cOB'(x')}
\def\cOBdyd{\cOB'(y')}
\def\cNBdxd{{N'(x')}}

\def\tcB{\widetilde{\cB}}
\def\tcOBx{\widetilde{\cOB(x)}}
\def\tcOBy{\widetilde{\cOB(y)}}
\def\tOBxd{\widetilde{\OB(x')}}
\def\tN{\widetilde{N}}

\def\HSff#1#2{H^{(#1)}_{#2}}
\def\HSf#1{H_{#1}}
\def\HSfX#1{H_X(#1)}
\def\HSfY#1{H_Y(#1)}
\def\HSfXd#1{H_{X'}(#1)}
\def\HSf#1#2{H_{#1}(#2)}
\def\HSgeqa#1{{{#1}(\geq \nu)}}
\def\HSa#1{{{#1}(\nu)}}
\def\HSb#1{{{#1}(\mu)}}
\def\HS#1#2{{#1(#2)}}
\def\HSmax#1{{{#1}_{max}}}
\def\HSP#1{{#1}_{P\mbox{-}max}}

\def\HSmaxneP#1{{#1_{P-max}^{nes}}}
\def\Pmaxne#1{{P_{#1}^{nes,max}}}

\def\HSmaxne#1{{#1}_{max}^{nes}}
\def\HSane#1{{#1}(\nu)^{nes}}
\def\HSne#1#2{{#1(#2)^{nes}}}

\def\HSfXd#1{H_{X'}(#1)}
\def\HSfY#1{H_Y(#1)}
\def\HSfW#1{H_W(#1)}
\def\HSfWd#1{H_{W'}(#1)}

\def\HSfXBB{H_{X^{B}}}
\def\HSfXB#1{H_{X^{B}}(#1)}
\def\HSfOBxX#1{H_{\OBxX}(#1)}

\def\tHSfXX{H^{O}_{X}}
\def\tHSfX#1{H^{O}_{X}(#1)}
\def\tHSf#1{H^{O}_{#1}}

\def\tHSfXXd{H^{O}_{X'}}
\def\tHSfXd#1{H^{O}_{X'}(#1)}

\def\tdHSfXd{H^{O'}_{X'}}

\def\cHSfXX{{H^O_X}}
\def\cHSfX#1{{H^O_X(#1)}}
\def\cHSfXd#1{{H^O_{X'}(#1)}}

\def\tnu{\widetilde{\nu}}
\def\tmu{\widetilde{\mu}}

\def\tHSgeqa#1{{#1}(\geq \tnu)}
\def\tHSa#1{{#1}(\tnu)}
\def\tHS#1#2{{#1}(#2)}
\def\tHSb#1{{#1}(\mu)}
\def\tHSmax#1{{#1}^{O}_{max}}
\def\tHSmaxiX{{X^{\OB}_i}_{max}}


\def\cHSmax#1{#1_{max}^O}

\def\SigmaXmax{\Sigma_X^{max}}
\def\SigmaXnmax{\Sigma_{X_n}^{max}}
\def\Sigmamax{{\Sigma^{\max}}}
\def\SigmaUmax{{\Sigma_U^{max}}}

\def\Sigmane#1{\Sigma^{nes}_{#1}}
\def\Sigmamaxne#1{\Sigma^{nes,max}_{#1}}

\def\SigmaXmaxO{\Sigma_X^{max,0}}

\def\SigmaXdmax{\Sigma_{X'}^{max}}

\def\SigmaXP{\Sigma_X^{P\mbox{-}max}}
\def\SigmaXnP{\Sigma_{X_n}^{P\mbox{-}max}}
\def\SigmamaxneP#1{{\Sigma_{#1}^{nes,P-max}}}

\def\SigmaXdmax{\Sigma_{X'}^{max}}
\def\SigmaYmax{\Sigma_Y^{max}}

\def\tSigmaX{\Sigma_X^{O}}
\def\tSigmaXn{\Sigma^{O}_{X_n}}
\def\tSigmaXmax{\Sigma_X^{O,max}}
\def\tSigmaYmax{\Sigma_Y^{O,max}}
\def\tSigmaXnmax{\Sigma_{X_n}^{O,max}}
\def\tSigmaUmax{\Sigma_U^{O,max}}
\def\tSigmamax{\Sigma^{O,max}}

\def\cSigmaX{{\Sigma^O_X}}
\def\cSigmaXmax{{\Sigma^{O,max}_X}}
\def\cSigmaXdmax{{\Sigma^{O,max}_{X'}}}
\def\cSigmadXd{{\Sigma^{O'}_{X'}}}
\def\cSigmaXd{{\Sigma^{O}_{X'}}}

\def\te{\widetilde{e}}
\def\tlam{\widetilde{\lambda}}
\def\Dir{{\mathrm{Dir}}}
\def\IDir{I{\mathrm{Dir}}}
\def\Dirob{{\mathrm{Dir}}^{O}}
\def\tDir{\widetilde{\Dir}}

\def\cDirx{{\Dir^{\OB}_x}}

\def\Cb{\overline{C}}
\def\ub{\overline{u}}

\def\Lamdob{{\Lambda'}^{\OB}}
\def\eob{e^{O}}

\def\vob{\bv^{\OB}}
\def\alphaob{\alpha^{\OB}}
\def\betaob{\beta^{\OB}}
\def\epsilonob{\epsilon^{\OB}}
\def\gammaob{\gamma^{\OB}}
\def\gammamob{\gamma^{-\OB}}
\def\gammapob{\gamma^{+\OB}}
\def\Deltaob{\Delta^{\OB}}
\def\deltaob{\delta^{\OB}}

\def\fob{f^{\OB}}
\def\fdob{{f'}^{\OB}}
\def\gob{g^{OB}}
\def\ffob#1{(f^{(#1)})^{\OB}}

\def\ff#1{f^{(#1)}}
\def\yy#1{y^{(#1)}}

\def\EP{\widetilde{V}}
\def\Bm{B_{\max}}

\def\bvob{\bv^{\OB}}

\def\tR{\tilde{R}}
\def\tJ{\tilde{J}}
\def\tfm{\tilde{\fm}}

\pagenumbering{arabic}

\setcounter{page}{1}
\setcounter{section}{-1}

\smallskip

\footnotetext[1]{
Universit\'e de Versailles,
laboratoire de math\'ematiques UMR8100,
45 avenue des \'Etats-Unis,
F78035 Versailles cedex,
France,
{\it Vincent.Cossart@math.uvsq.fr}}

\footnotetext[2]{Fakult\"at f\"ur Mathematik,
Universit\"at Regensburg,
Universit\"atsstr. 31,
93053 Regensburg,
Germany,
{\it uwe.jannsen@mathematik.uni-regensburg.de}
}


\footnotetext[3]{
Interactive Research Center of Science,
Graduate School of Science and Engineering,
Tokyo Institute of Technology,
2-12-1 Ookayama, Meguro,
Tokyo 152-8551,
Japan
{\it sshuji@msb.biglobe.ne.jp}
}


\contentsline {section}{\numberline {0}Introduction}{3}
\contentsline {section}{\numberline {1}Basic invariants for singularities}{12}
\contentsline {section}{\numberline {2}Permissible blow-ups}{30}
\contentsline {section}{\numberline {3}$B$-Permissible blow-ups - the embedded case}{40}
\contentsline {section}{\numberline {4}$B$-Permissible blow-ups - the non-embedded case}{53}
\contentsline {section}{\numberline {5}Main theorems and strategy for their proofs}{64}
\contentsline {section}{\numberline {6}$(u)$-standard bases}{86}
\contentsline {section}{\numberline {7}Characteristic polyhedra of $J\subset R$}{96}
\contentsline {section}{\numberline {8}Transformation of standard bases under blow-ups}{108}
\contentsline {section}{\numberline {9}Termination of the fundamental sequences of $B$-permissible blowups, and the case $e_x(X)=1$}{118}
\contentsline {section}{\numberline {10}Additional invariants in the case $e_x(X)=2$}{124}
\contentsline {section}{\numberline {11}Proof in the case $e_x(X)=\overline {e}_x(X)=2$, I: some key lemmas}{128}
\contentsline {section}{\numberline {12}Proof in the case $e_x(X)=\overline {e}_x(X)=2$, II: separable residue extensions}{133}
\contentsline {section}{\numberline {13}Proof in the case $e_x(X)=\overline e_x(X)=2$, III: inseparable residue extensions}{138}
\contentsline {section}{\numberline {14}Non-existence of maximal contact in dimension $2$}{151}
\contentsline {section}{\numberline {15}An alternative proof of Theorem \ref {thm.max-elim}}{159}
\contentsline {section}{\numberline {16}Functoriality, locally noetherian schemes, algebraic spaces and stacks}{162}




\newpage
\section{Introduction}

The principal aim of this paper is to show the following three theorems on the resolution of
singularities of an arbitrary reduced excellent noetherian scheme $X$ of dimension at most two.
{\it In the following, all schemes will be assumed to be noetherian, but see the end of the
introduction and section 16 for locally noetherian schemes}.

\begin{theorem}\label{thm.intro.1} (Canonical controlled resolution)
There exists a canonical finite sequence of morphisms
$$
\begin{CD}
\pi:\, X' =X_n @>>> \ldots @>>> X_1 @>>> X_0=X
\end{CD}
$$
such that $X'$ is regular and, for each $i$, $X_{i+1} \rightarrow X_i$ is the blow-up of $X_i$
in a permissible center $D_i \subset X_i$
which is contained in $(X_i)_{sing}$, the singular locus of $X_i$. This sequence is functorial
in the sense that it is compatible with automorphisms of $X$ and (Zariski or \'etale) localizations.
\end{theorem}

We note that this implies that $\pi$ is an isomorphism over $X_{reg}=X - X_{sing}$, and we
recall that a subscheme $D \subset X$ is called permissible, if $D$ is regular and $X$ is
normally flat along $D$ (see \ref{def.nf}). The compatibility with automorphisms means that
every automorphism of $X$ extends to the sequence in a unique way. The compatibility with the
localizations means that the pull-back via a localization $U \rightarrow X$ is the canonical
resolution sequence for $U$ after suppressing the morphisms which become isomorphisms over $U$.
It is well-known that Theorem \ref{thm.intro.1} implies:

\begin{theorem}\label{thm.intro.2} (Canonical embedded resolution) Let $i: X \hookrightarrow Z$ be a closed
immersion, where $Z$ is a regular excellent scheme. Then there is a canonical commutative diagram
$$
\begin{CD}
X' @>{i'}>> Z' \\
@V{\pi}VV @VV{\pi_Z}V \\
X @>{i}>> Z
\end{CD}
$$
where $X'$ and $Z'$ are regular, $i'$ is a closed immersion, and $\pi$ and $\pi_Z$ are proper
and surjective morphisms inducing isomorphisms
$\pi^{-1}(X - X_{sing}) \mathop{\longrightarrow}\limits^{\sim} X - X_{sing}$
and $\pi_Z^{-1}(Z - X_{sing}) \mathop{\longrightarrow}\limits^{\sim} Z - X_{sing}$.
Moreover, the morphisms $\pi$ and $\pi_Z$ are compatible with automorphisms of $(X,Z)$
and (Zariski or \'etale) localizations in $Z$.
\end{theorem}

In fact, starting from Theorem \ref{thm.intro.1} one gets a canonical sequence
$Z' = Z_n \rightarrow \ldots Z_1 \rightarrow Z_0 = Z$ and closed immersions
$X_i \hookrightarrow Z_i$ for all $i$, such that $Z_{i+1} \rightarrow Z_i$
is the blow-up in $D_i \subset (X_i)_{sing} \subset Z_i$ and $X_{i+1} \subseteq Z_{i+1}$ is identified
with the strict transform of $X_i$ in the blow-up $Z_{i+1} \rightarrow Z_i$. Then $Z_{i+1} \rightarrow Z_i$
is proper (in fact, projective) and surjective, and $Z_{i+1}$ is regular since $Z_i$ and $D_i$ are.

\medbreak
For several applications the following refinement is useful:

\begin{theorem}\label{thm.intro.3} (Canonical embedded resolution with boundary)
Let $i: X \hookrightarrow Z$ be a closed immersion into a regular scheme
$Z$, and let $B \subset Z$ be a simple normal crossings divisor such that no irreducible component of $X$
is contained in $B$. Then there is a canonical commutative
diagram
$$
\begin{CD}
X' @>{i'}>> Z' @. \; \supset \; B' \\
@V{\pi_X}VV @V{\pi_Z}VV \\
X @>{i}>> Z  @. \; \supset \; B
\end{CD}
$$
where $i'$ is a closed immersion, $\pi_X$ and $\pi_Z$ are projective, surjective,
and isomorphisms outside $X_{sing}\cup (X\cap B)$, and $B'=\pi_Z^{-1}(B) \cup E$, where $E$ is the
exceptional locus of $\pi_Z$ (which is a closed subscheme such that $\pi_Z$ is an isomorphism over
$Z-\pi_Z(E)$). Moreover, $X'$ and $Z'$ are regular,
$B'$ is a simple normal crossings divisor on $Z'$, and $X'$ intersects
$B'$ transversally on $Z'$. Furthermore, $\pi_X$ and $\pi_Z$ are compatible with
automorphisms of $(Z,X,B)$ and with (Zariski or \'etale) localizations in $Z$.
\end{theorem}

More precisely, we prove the existence of a commutative diagram

$$
\begin{CD}
\pi_B:\; B' = B_m @. B_{m-1} @. \ldots @. B_1 @. B_0 = B \\
 @VVV @VVV @VVV @VVV @VVV \\
\pi_Z:\; Z' = Z_m @>>> Z_{m-1} @>>> \ldots @>>> Z_1 @>>> Z_0 = Z \\
 @AAA @AAA @AAA @AAA @VVV \\
\pi\;:\; X' = X_m @>>> X_{m-1} @>>> \ldots @>>> X_1 @>>> X_0 = X
\end{CD}
$$
where the vertical morphisms are closed immersions and, for each $i$,
$X_{i+1} = \Bl_{D_i}(X_i) \rightarrow X_i$ is the blow-up of $X_i$ in a
permissible center $D_i \subset (X_i)_{sing}$, $Z_{i+1}=\Bl_{D_i}(Z_i) \rightarrow Z_i$
is the blow-up of $Z_i$ in $D_i$ (so that $Z_{i+1}$ is regular and $X_{i+1}$ is identified
with the strict transform of $X_i$ in $Z_{i+1}$), and $B_{i+1}$ is the complete transform
of $B_i$, i.e., the union of the strict transform of $B_i$ in $Z_{i+1}$ and the exceptional divisor of
the blow-up $Z_{i+1}\to Z_i$. Furthermore, $D_i$ is $B_i$-permissible, i.e.,
$D_i \subset X_i$ is permissible, and normal crossing with $B_i$ (see
Definition \ref{def.DncBZ}), which implies that $B_{i+1}$ is a simple normal crossings divisor on
$Z_{i+1}$ if this holds for $B_i$ on $Z_i$.

\bigbreak\bigbreak
In fact, the second main theorem of this paper, Theorem \ref{thm.main.3}, states a somewhat more general
version, in which $B$ can contain irreducible components of $X$. Then one can assume that
$D_i$ is not contained in the strongly $B_i$-regular locus $\XBsreg$ (see Definition \ref{def.B-reg}),
and one gets that $X'$ is normal crossing with $B$ (Definition \ref{def.DncBZ}).
This implies that $\pi$ is an isomorphism above $\XBsreg \subseteq X_{reg}$, and, in particular,
again over $X_{reg} - B$. In addition, this Theorem also treats non-reduced schemes $X$, in which
case $(X')_{red}$ is regular and normal crossing with $B$ and $X'$ is normally flat along $(X')_{red}$.

\medbreak
Moreover, we obtain a variant, in which we only consider strict transforms for the normal crossings divisor,
i.e., where $B_{i+1}$ is the strict transform of $B_i$. Then we only get the normal crossing of $X'$
(or $X'_{red}$ in the non-reduced case) with the strict transform $\tilde{B}$ of $B$ in $Z'$.

\medbreak
Theorem \ref{thm.intro.1}, i.e., the case where we do not assume any embedding for $X$,
will also be proved in a more general version: Our first main theorem, Theorem \ref{thm.main.1}, allows a
non-reduced scheme $X$ as well as a so-called boundary on $X$,
a notion which is newly introduced in this paper (see section \ref{BpbuII}).
Again this theorem comes in two versions, one with complete transforms and one with strict transforms.

\medbreak
Our approach implies that Theorem \ref{thm.main.1} implies Theorem \ref{thm.main.3}.
In particular, the canonical resolution sequence of Theorem \ref{thm.main.3} for $B=\emptyset$
and strict transforms (or of Theorem \ref{thm.intro.3} for this variant)
coincides with the intrinsic sequence for $X$ from Theorem \ref{thm.intro.1}.
Thus, the readers only interested in Theorems \ref{thm.intro.1} and \ref{thm.intro.2}
can skip sections 3 and 4 and ignore any mentioning
of boundaries/normal crossings divisors (by assuming them to be empty).

\bigbreak
We note the following corollary.

\begin{corollary}\label{cor.intro.1}
Let $Z$ be a regular excellent scheme (of any dimension), and let $X \subset Z$ be a
reduced closed subscheme of dimension at most two. Then there exists a projective surjective
morphism $\pi: Z' \longrightarrow Z$ which is an isomorphism over $Z - X$,
such that $\pi^{-1}(X)$, with the reduced subscheme structure, is a simple normal crossings
divisor on $Z'$.
\end{corollary}

In fact, applying Theorem \ref{thm.intro.3} with $B = \emptyset$, we get a projective
surjective morphism $\pi_1: Z_1 \longrightarrow Z$ with regular $Z_1$, a regular
closed subscheme $X_1 \subset Z_1$ and a simple normal crossings divisor $B_1$
on $Z_1$ such that $\pi_1$ is an isomorphism over $Z - X$ (in fact, over
$Z - (X_{sing})$), and $\pi_1^{-1}(X) = X_1 \cup B_1$. Moreover, $X_1$
and $B_1$ intersect transversally. In particular, $X_1$ is normal crossing with $B_1$
in the sense of Definition \ref{def.DncBZ}. Hence we obtain the wanted situation
by composing $\pi_1$ with $\pi_2: Z' \longrightarrow Z_1$, the blow-up of $Z_1$
in the $B_1$-permissible (regular) subscheme $X_1$, and letting $X' = \pi_2^{-1}(X_1\cup B_1)$
which is a simple normal crossings divisor, see Lemma 3.2.

\medbreak
Moreover we mention that Theorem \ref{thm.intro.3} is applied in a paper of the second
and third author \cite{JS}, to prove a conjecture of Kato and finiteness of certain
motivic cohomology groups for varieties over finite fields. This was a main motivation for
these authors to work on this subject.

\bigbreak\bigbreak
To our knowledge, none of the three theorems is known, at least not in the stated generality.
Even for $\dim(X) = 1$ we do not know a reference for these results, although
they may be well-known. For $X$ integral of dimension 1, Theorem \ref{thm.intro.1} can be found in \cite{Be}
section 4, and a proof of Theorem \ref{thm.intro.3} is written in \cite{Ja}.

\medbreak
In 1939 Zariski \cite{Za1} proved Theorem \ref{thm.intro.1} (without discussing canonicity or functoriality)
for irreducible surfaces over algebraically closed fields of characteristic zero.
Five years later, in \cite{Za3}, he proved Corollary \ref{cor.intro.1} (again without canonicity or functoriality)
for surfaces over fields of characteristic zero which are embedded in a non-singular threefold.
In 1966, in his book \cite{Ab4}, Abhyankar extended this last result to all algebraically closed fields,
making heavy use of his papers \cite{Ab3} and \cite{Ab5}.
Around the same time, Hironaka \cite{H6} sketched a shorter proof
of the same result, over all algebraically closed ground fields, based on a different method.
Recently a shorter account of Abhyankar's results was given by Cutkosky \cite{Cu2}.
For all excellent schemes of characteristic zero, i.e., whose residue fields all have characteristic
zero, and of arbitrary dimension, Theorems \ref{thm.intro.1} and \ref{thm.intro.3} were proved by
Hironaka in his famous 1964 paper \cite{H1} (Main Theorem $1^\ast$, p. 138, and Corollary 3, p. 146),
so Theorem \ref{thm.intro.2} holds as well, except that the approach is not constructive, so
it does not give canonicity or functoriality. These issues were addressed and solved in the later
literature, especially in the papers by Villamayor, see in particular \cite{Vi}, and by
Bierstone-Millman, see \cite{BM1}, by related, but different approaches. In these references,
a scheme with a fixed embedding into a regular scheme is considered, and in \cite{Vi}, the process depends on
the embedding. The last issue is remedied by a different approach in \cite{EH}. In positive characteristics,
canonicity was addressed by Abhyankar in \cite{Ab6}.

\medbreak
There are further results on a weaker type of resolution for surfaces, replacing the blowups
in regular centers by different techniques. In \cite{Za2} Zariski showed how to resolve a surface
over a not necessarily algebraically closed field of characteristic zero by so-called local
uniformization which is based on valuation-theoretic methods. Abhyankar \cite{Ab1} extended
this to all algebraically closed fields of positive characteristics, and later \cite{Ab2} extended
several of the results to more general schemes whose closed points have perfect residue fields.
In 1978 Lipman \cite{Li} gave a very simple procedure to obtain resolution of singularities for arbitrary excellent
two-dimensional schemes $X$ in the following way: There is a finite sequence $X_n \rightarrow X_n \rightarrow \ldots X_1 \rightarrow X$
of proper surjective morphisms such that $X_n$ is regular. This sequence is obtained by
alternating normalization and blowing up in finitely many isolated singular points.
But the processes of uniformization or normalization are not controlled in the sense of Theorem \ref{thm.intro.1},
i.e., not obtained by permissible blow-ups, and it is not known how to extend them to
an ambient regular scheme $Z$ like in Theorem \ref{thm.intro.2}. Neither is it
clear how to get Theorem \ref{thm.intro.3} by such a procedure. In particular, these weaker results were not sufficient for
the mentioned applications in \cite{JS}. This is even more the case for the weak resolution of
singularities proved by de Jong \cite{dJ}.

\medbreak
It remains to mention that there are some results on weak resolution of singularities
for threefolds over a field $k$ by Zariski \cite{Za3} (char$(k)$ = 0), Abhyankar \cite{Ab4}
($k$ algebraically closed of characteristic $\neq 2,3,5$ -- see also \cite{Cu2}),
and by Cossart and Piltant \cite{CP1}, \cite{CP2} ($k$ arbitrary), but this is not the topic of
the present paper.

\bigbreak\bigbreak
Our approach is roughly based on the strategy of Levi-Zariski used in \cite{Za1}, but
more precisely follows the approach (still for surfaces) given by Hironaka in the paper \cite{H6}
cited above. The general strategy is very common by now: One develops certain invariants which
measure the singularities and aims at constructing a sequence of blow-ups for which the invariants
are non-increasing, and finally decreasing, so that in the end one concludes one has reached
the regular situation. The choices for the centers of the blow-ups are made by considering
the strata where the invariants are the same. In fact, one blows up `the worst locus',
i.e., the strata where the invariants are maximal, after possibly desingularising these strata.
The main point is to show that the invariants do finally decrease. In characteristic zero this is done by
a technique introduced by Hironaka in \cite{H1}, which is now called the method of maximal
contact (see \cite{AHV} and Giraud's papers \cite{Gi2} and \cite{Gi3} for some theoretic background),
and an induction on dimension.

\medbreak
But it is known that the theory of maximal contact does not work in positive characteristic. There
are some theoretic counterexamples in \cite{Gi3}, and some explicit counterexamples for threefolds
in characteristic two by Narasimhan \cite{Na1}, see also \cite{Co2} for an interpretation in our sense.
It is not clear if the counterexamples in \cite{Na2}, for threefolds in any positive characteristic,
can be used in the same way. But in section 15 of this paper, we show
that maximal contact does not even exists for surfaces, in any characteristic, even if maximal
contact is considered in the weakest sense. Therefore the strategy of proof has to be different,
and we follow the one outlined in \cite{H6}, based on certain polyhedra (see below).
That paper only considers the case of a hypersurface, but in another paper \cite{H3} Hironaka
develops the theory of these polyhedra for ideals with several generators, in terms of certain
`standard bases' for them (which also appear in \cite{H1}). The introduction of \cite{H3} expresses
the hope that this theory of polyhedra will be useful for the resolution of singularities, at least for surfaces.
Our paper can be seen as a fulfilment of this program.

\medbreak
In his fundamental paper \cite{H1}, Hironaka uses two important invariants for measuring the singularity
at a point $x$ of an arbitrary scheme $X$. The primary is the $\nu$-invariant $\nu_x^\ast(X)\in \bNN$,
and the secondary one is the dimension $e_x(X) \in \bN$ (with $0 \leq e_x(X) \leq \dim(X)$)
of the so-called directrix $\Dir_x(X)$ of $X$ at $x$. Both only depend
on the cone $C_x(X)$ of $X$ at $x$. Hironaka proves that for a permissible
blow-up $X' \rightarrow X$ and a point $x'\in X'$ with image $x\in X$ the $\nu$-invariant is
non-increasing: $\nu_{x'}^\ast(X') \leq \nu_x^\ast(X)$. If equality holds here (one says $x'$ is near to $x$),
then the (suitably normalized) $e$-invariant is non-decreasing. So the main problem is to show that there
cannot be an infinite sequence of blow-ups with `very near' points $x' \mapsto x$ (which means that they
have the same $\nu$- and $e$-invariants).

\medbreak
To control this, Hironaka in \cite{H3} and \cite{H6} introduces a tertiary, more complex invariant, the polyhedron associated
to the singularity, which lies in $\mathbb R_{\geq 0}^e$. It depends not just on $C_x(X)$, but on the
local ring $\cO_{X,x}$ of $X$ at $x$ itself, and also on various choices: a regular local ring $R$ having
$\cO_{X,x}$ as a quotient, a system of regular parameters $y_1,\ldots,y_r,u_1,\ldots,u_e$ for $R$
such that $u_1,\ldots,u_e$ are `parameters' for the directrix $\Dir_x(X)$,
and equations $f_1,\ldots,f_m$ for $\cO_{X,x}$ as a quotient of $R$ (more precisely, a $(u)$-standard base
of $J = \ker(R \rightarrow \cO_{X,x})$\,). In the situation of Theorem \ref{thm.intro.2}, $R$ is naturally
given as $\cO_{Z,x}$, but in any case, such an $R$ always exists after
completion, and the question of ruling out an infinite sequence of very near points only depends
on the completion of $\cO_{X,x}$ as well. In the case considered in section 13,
it is not a single strictly decreasing invariant which comes out of these polyhedra, but rather the behavior
of their shape which tells in the end that an infinite sequence of very near points cannot exits.
This is sufficient for our purpose, but it might be interesting to find a strictly decreasing invariant
also in this case. In the particular situation considered in \cite{H6} (a hypersurface over
an algebraically closed field), this was done by Hironaka; see also \cite{Ha} for a variant.

\medbreak
As a counterpart to this local question, one has to consider a global strategy and the global
behavior of the invariants, to understand the choice of permissible centers and the global
improvement of regularity. Since the $\nu$-invariants are
nice for local computations, but their geometric behavior is not so nice, we use the Hilbert-Samuel
invariant $H_{\cO_{X,x}}\in \bNN$ as an alternative primary invariant here.
They were extensively studied by Bennett \cite{Be}, who proved similar non-increasing
results for permissible blow-ups, which was then somewhat improved by Singh \cite{Si1}.
Bennett also defined global Hilbert-Samuel functions $H_X: X \rightarrow \bNN$, which, however,
only work well and give nice strata in the case of so-called weakly biequidimensional excellent schemes.
We introduce a variant (Definition \ref{def.HS.schemes}) which works for arbitrary (finite dimensional)
excellent schemes. This solves a question raised by Bennett. The associated Hilbert-Samuel strata
$$
\HSa X = \{ x\in X \mid \, H_X(x) = \nu\} \quad\quad\mbox{for}\quad \nu \in \bNN\;,
$$
are then locally closed, with closures contained in $\HSgeqa X = \{x\in X \mid H_X(x)\geq \nu\}$.
In particular, $\HSa X$ is closed for maximal $\nu$ (here $\mu\geq\nu$ if $\mu(n)\geq\nu(n)$ for all $n$).

\medbreak
Although our main results are for two-dimensional schemes, the major part of this paper
is written for schemes of arbitrary dimension, in the hope that this might be useful for
further investigations. Only in part of section 5 and in sections 10 through 14
we have to exploit some specific features of the low-dimensional situation.
According to our understanding, there are mainly two
obstructions against the extension to higher-dimensional schemes: The fact that
in Theorem \ref{thmdirectrix} (which gives crucial information on the locus of near points)
one has to assume $\Char(k(x)) =0$ or $\Char(k(x))\geq \dim(X)/2+1$,
and the lack of good invariants of the polyhedra for $e >2$, or of other suitable tertiary
invariants in this case.

\medbreak
We have tried to write the paper in such a way that it is well
readable for those who are not experts in resolution of singularities (like two of us)
but want to understand some results and techniques and apply them in arithmetic or algebraic
geometry. This is also a reason why we did not use the notion of idealistic exponents \cite{H7}.
This would have given the extra burden to recall this theory, define characteristic
polyhedra of idealistic exponents, and rephrase the statements in \cite{H5}.
Equipped with this theory, the treatment of the functions defining the scheme and
the functions defining the boundary would have looked more
symmetric; on the other hand, the global algorithm clearly distinguishes these two.

\bigbreak\bigbreak
We now briefly discuss the contents of the sections. In section 1 we discuss the primary and
secondary invariants (local and global) of singularities mentioned above. In section 2 we discuss
permissible blow-ups and the behavior of the introduced invariants for these, based on
the fundamental results of Hironaka and Bennett (and Singh).

\medskip
In section 3 we study similar questions
in the setting of Theorem \ref{thm.intro.3}, i.e., in a `log-situation' $X \subset Z$ where one has a
`boundary': a normal crossings divisor $\bB$ on $Z$. We define a class of log-Hilbert-Samuel functions
$H_X^{\OB}$, depending on the choice of a `history function' $\OB: X \rightarrow \{\mbox{ subdivisors of }\bB\}$
characterizing the `old components' of $\bB$ at $x\in X$. Then $H_X^\OB(x) = (H_X(x),n)$, where $n$
is the number of old components at $x$. This gives associated log-Hilbert-Samuel strata
$$
\tHSa X = \{ x\in X \mid \,H_X^\OB(x) = \tnu \} \quad\quad\mbox{for}\quad \tnu \in \bNN\times\bN \;.
$$
For a $\bB$-permissible blow-up $X' \rightarrow X$,
we relate the two Hilbert-Samuel functions and strata, and study some transversality properties.

\medskip
In section 4 we extend this theory to the situation where we have just an excellent scheme $X$
and no embedding into a regular scheme $Z$. It turns out that one can also define the
notion of a boundary $\cB$ on $X$: it is just a tuple $(B_1,\ldots,B_r)$ (or rather
a multiset, by forgetting the ordering) of locally principal closed subschemes $B_i$ of $X$.
In the embedded situation $X \subset Z$, with a normal crossings divisor $\bB$ on $Z$,
the associated boundary $\cB_X$ on $X$ is just given by the traces of the components of $\bB$
on $X$ and we show that they carry all the information which is needed.
Moreover, this approach makes evident that the constructions and strategies defined later
are intrinsic and do not depend on the embedding. All results in section 3 can be carried over
to section 4, and there is a perfect matching (see Lemma \ref{lem.comp.emb-nemb}). We could have
started with section 4 and derived the embedded situation in section 3 as a special case,
but we felt it more illuminating to start with the familiar classical setting; moreover,
some of the results in section 4 (and later in the paper) are reduced to the embedded
situation, by passing to the local ring and completing (see Remark \ref{rem.emb-nemb},
Lemma \ref{lem.comp.emb-nemb} and the applications thereafter).

\medskip
In section 5 we state the Main Theorems \ref{thm.main.1} and \ref{thm.main.3}, corresponding to
somewhat more general versions of Theorem \ref{thm.intro.1} and \ref{thm.intro.3}, respectively, and
we explain the strategy to prove them. Based on an important theorem by Hironaka
(see \cite{H2} Theorem (1,B) and the following remark), it suffices to find
a succession of permissible blowups for which the Hilbert-Samuel invariants
decrease. Although this principle seems to be well-known, and might be obvious for surfaces, we
could not find a suitable reference and have provided a precise statement and (short) proof of this fact in any dimension
(see Corollaries \ref{cor.max-elim} and \ref{cor.max-elim2} for the case without boundary,
and Corollaries \ref{cor.Omax-elim} and \ref{cor.Omax-elim1} for the case with boundary).
The problem arising is that the set of Hilbert functions is ordered by the total (or product) order ($H \leq H'$ iff
$H(n) \leq H(n)$ for all $n$), and that with this order there are infinite decreasing sequences
in the set $\bN^\bN$ of all functions $\nu: \bN \rightarrow \bN$. This is overcome by the fact that the
subset of Hilbert functions of quotients of a fixed polynomial ring $k[x_1,\ldots,x_n]$ is a noetherian ordered set,
see Theorem \ref{thm.ordHP}.
After these preparations we define a canonical resolution sequence (see Remark \ref{rem.strategy} for the
definition of so-called $\tnu$-eliminations, and Corollaries \ref{cor.max-elim} and \ref{cor.max-elim2}
for the definition of the whole resolution sequence out of this).

\medbreak
We point out that in Remark \ref{rem.strategy}, we define these canonical resolution sequences,
i.e., an explicit strategy for resolution of singularities for any dimension. It would be interesting
to see if this strategy always works.

\medbreak
The proof of the finiteness of these resolution sequences for dimension two is reduced to two key theorems, Theorem \ref{fu.thm0} and
\ref{Thm2}, which exclude the possibility of certain infinite chains of blow-ups with near (or $\OB$-near) points.
The key theorems concern only isolated singularities and hence only the local ring of $X$ at a closed point $x$,
and they hold for $X$ of arbitrary dimension, but with the condition that the `geometric' dimension
of the directrix is $\leq 2$ (which holds for $\dim(X) = 2$). As mentioned above, for this local
situation we may assume that we are in an embedded situation.

\medbreak
As a basic tool for various considerations, we study a situation as mentioned above,
where a local ring $\cO$ (of arbitrary dimension) is a quotient $R/J$ of a regular local ring $R$.
In section 6 we discuss suitable systems $(y,u)= (y_1,\ldots ,y_n,u_1,\ldots,u_e)$
of regular parameters for $R$ and suitable families $f = (f_1,\ldots,f_m)$ of generators
for $J$. A good choice for $(y,u)$ is obtained if $u$ is admissible for $J$ (Definition \ref{def2.11})
which means that $u_1,\ldots,u_e$ are affine parameters of the directrix of $\cO$ (so that
$e$ is the $e$-invariant recalled above).
We study valuations associated to $(y,u)$ and initial forms (with respect to these valuations)
of elements in $J$ and their behavior under change of the system of parameters.
As for $f$, in the special case that $J$ is generated by one
element (case of hypersurface singularities), any choice of $f = (f_1)$ is good.
In general, some choices of $f = (f_1,\ldots,f_m)$ are better than the other. A favorable
choice is a standard basis of $J$ (Definition \ref{def0.5}) as introduced in \cite{H1}. In \cite{H3}
Hironaka introduced the more general notion of a $(u)$-standard base of $J$
which is more flexible to work with and plays an important role in our paper.

\medbreak
In section 7 we recall, in a slightly different way, Hironaka's
definition \cite{H3} of the polyhedron $\Delta(f,y,u)$ associated to a system of parameters
$(y,u)$ and a $(u)$-standard basis $f$, and the polyhedron $\Delta(J,u)$ which
is the intersection of all $\Delta(f,y,u)$ for all choices of $y$ and $f$ as above (with fixed
$u$).
We recall Hironaka's crucial result from \cite{H3} that $\Delta(f,y,u) = \Delta(J,u)$ if $u$ is
admissible and $(f,y,u)$ is what Hironaka calls well-prepared, namely normalized (Definition \ref{def.normalized})
and not solvable (Definition \ref{def.solvable}) at all vertices.
Also, there is a certain process of making a given $(f,y,u)$ normalized (by changing $f$)
and not solvable (by changing $y$) at finitely many vertices, and at all vertices, if
$R$ is complete.
One significance of this result is that it provides a natural way of
transforming a $(u)$-standard base into a standard base under the assumption
that $u$ is admissible.

\medbreak
As explained above, it is important to study permissible blow-ups
$X' \rightarrow X$ and near points $x'\in X'$ and $x\in X$.
In this situation, to a system $(f,y,u)$ at $x$ we associate certain new systems
$(f',y',u')$ at $x'$. A key result proved in section 8 is that
if $f$ is a standard base, then $f'$ is a $(u')$-standard
base. The next key result is that the chosen $u'$ is admissible.
Hence, by Hironaka's crucial result mentioned above,
we can transform $(f',y',u')$ into a system $(g',z',u')$, where $g'$ is a
standard base.

\medbreak
The Key Theorems \ref{fu.thm0} and \ref{Thm2} concern certain sequences of
permissible blowups, which arise naturally from the canonical resolution sequence.
We call them fundamental sequences of $B$-permissible blowups (Definition \ref{Def.fupb})
and fundamental units of permissible blowups (Definition \ref{Def.fupb2}) and use them
as a principal tool. These are sequences of $B$-permissible blowups
$$
X_m \rightarrow \ldots \rightarrow X_1 \rightarrow X_0 = X\,,
$$
where the first blowup is in a closed point $x\in X$ (the initial point), and where the later blowups are
in certain maximal $B$-permissible centers $C_i$, which map isomorphically onto each other, lie above $x$, and
consist of points near to $x$. For a fundamental sequence there is still a $B$-permissible
center $C_m \subset X_m$ with the same properties; for a fundamental unit there is none,
but only a chosen closed point $x_m \in X_m$ (the terminal point) which is near to $x$.
In section 9 we study some first properties of these fundamental sequences. In particular
we show a certain bound for the $\delta$-invariant of the associated polyhedra.
This suffices to show the first Key Theorem \ref{fu.thm0} (dealing with the case $e_x(X) = 1$),
but is also used in section 14.

\medbreak
For the second Key Theorem \ref{Thm2} (dealing with the case $e_x(X)=2$),
one needs some more information on the ($2$-dimensional) polyhedra, in particular,
some additional invariants. These are introduced in section 10. Then Theorem \ref{Thm2}
is proved in the next three sections. It states that there is no infinite sequence
$$
\ldots \rightarrow \cX_2 \rightarrow \cX_1 \rightarrow \cX_0
$$
of fundamental units of blow-ups such that the closed initial points and terminal
points match and are isolated in their Hilbert-Samuel strata.
After some preparations in section 11, section 12
treats the case where the residue field extension $k(x')/k(x)$ is trivial (or separable).
This is very much inspired by \cite{H6}, which however only treats the special situation
of a hypersurface in a regular threefold over an algebraically closed field and does not
contain proofs for all claims.
Then section 13 treats the case where there occur inseparable residue
field extensions $k(x')/k(x)$.
This case was basically treated in \cite{Co1} but we give
a more detailed account and fill gaps in the original proof,
with the aid of the results of section 8, and Giraud's notion of
the ridge \cite{Gi1}, \cite{Gi3} (fa\^ite in French) a notion which generalizes the directrix.

\medbreak
In section 14, we show that maximal contact does not exist for surfaces
in positive characteristic $p$. For each $p$ a counterexample is given which then works for any
field of that characteristic.

\medbreak
In section 15 we give a more algebro-geometric proof of the fact that it suffices
to show how to eliminate the maximal Hilbert-Samuel stratum.

\medbreak
Finally, in section 16 we give a re-interpretation of the functoriality we obtain
for our resolution, for arbitrary flat morphisms with regular fibers, and we apply
this to show resolution for excellent schemes which are only locally noetherian,
and for excellent algebraic stacks with atlas of dimension at most two.

\medbreak
It will be clear from the above how much we owe to all the earlier work on
resolution of singularities, in particular to the work of Hironaka which gave
the general strategy but also the important tools used in this paper.

\bigbreak
{\it Conventions and concluding remarks:} All schemes are assumed to be finite dimensional.
Regular schemes are always assumed to be locally noetherian. Recall also that excellent schemes
are by definition locally noetherian.

\medbreak
In this introduction and in sections 1 to 15, the readers should
best assume that all schemes are noetherian. At some places we write locally
noetherian to indicate that certain definitions make sense and certain results still hold
for schemes which are only locally noetherian. Resolution for such schemes is treated in section 16.

\newpage
\section{Basic invariants for singularities}

In this section we introduce some basic invariants for
singularities.

\bigskip\noindent
{\bf 1.1 Invariants of graded rings and homogeneous ideals in polynomial rings}

\bigskip
Let $k$ be a field and $S = k [X_1, \dots, X_n]$ be a polynomial
ring with $n$ variables. Let $S_\nu\subset S$ be the $k$-subspace
of the homogeneous polynomials of degree $n$ (including $0$). Fix
a homogeneous ideal $I \subseteq S$.

\bigskip
\begin{definition}\label{def0.1}
For integers $i\geq 1$ we define $\nu^{\,i}(I) \in \bN \cup \{ \infty \}$
as the supremum of the $\nu \in \bN$ satisfying the condition that
there exist homogeneous $\varphi_1, \dots, \varphi_{i - 1} \in I$
such that
$$ S_{\mu} \cap I = S_{\mu} \cap <\varphi_1, \dots, \varphi_{i -1}> \quad
\mbox{for all } \mu < \nu. $$
\end{definition}
\medskip

By definition we have $\nu^{1} (I) \leq \nu ^{2} (I) \leq \dots$.
We write
$$
\nu^*(I)=( \nu^{1} (I), \nu^{2} (I),\dots, \nu^{m}
(I),\infty,\infty, \dots)
$$
and call it the $\nu$-invariant of $I$. We have the following
result (cf. \cite{H1} Ch. III \S1, Lemma 1).

\begin{lemma}\label{lem0.2}
Let $I =<\varphi_1, \dots,\varphi_m>$
with homogeneous elements $\varphi_i$ of degree $\nu_i$ such that:
\begin{itemize}
\item[$(i)$]
$\varphi_i \not\in <\varphi_1, \dots, \varphi_{i -1}>$ for all
$i = 1, \dots, m$,
\item[$(ii)$]
$\nu_1 \leq \nu_2 \leq \dots \leq \nu_m$.
\end{itemize}
Then we have
$$ \nu^{\,i} (I) = \left\{ \begin{array}{ccc}
\nu_i & , & i \leq m, \\
\infty & , & i > m.\\
\end{array} \right. $$
\end{lemma}

\begin{definition}\label{def0.2}
Let $\varphi=(\varphi_1, \dots, \varphi_m)$ be a system of homogeneous
elements in $S$ and $I\subset S$ be the ideal that it generates.
\begin{itemize}
\item[(1)]
$\varphi$ is weakly normalized if it satisfies the condition $(i)$ of
Lemma \ref{lem0.2}.
\item[(2)]
$\varphi$ is a standard base of $I$ if it satisfies the conditions $(i)$
and $(ii)$ of Lemma \ref{lem0.2}.
\end{itemize}
\end{definition}

We have the following easy consequence of \ref{lem0.2}.

\begin{corollary}\label{cordef0.1}
Let $I\subset S$ be a homogeneous ideal and let
$\psi=(\psi_1, \dots, \psi_j)$ be a system of homogeneous elements in $I$
which is weakly normalized.
\begin{itemize}
\item[(1)]
The following conditions are equivalent:
\begin{itemize}
\item[$(ii)$]
$\deg \psi_i = \nu^{\,i} (I)$ for $i = 1, \dots, j$.
\item[$(ii')$]
For all $i = 1,\dots, j$, $\psi_i$ has minimal degree in
$I$ such that $\psi_i \not\in <\psi_1, \dots, \psi_{i -1}>$.
\end{itemize}
\item[(2)]
If the conditions of (1) are satisfied, then
$\psi$ can be extended to a standard base of $I$.
\end{itemize}
\end{corollary}
\medbreak

By the lemma a standard base of $I$ and $\nu^\ast(I)$ are obtained as follows:\\
Put $\nu_1:= \min \{ \nu \; \vert \; \exists \; \varphi\in S_{\nu} \cap I
\setminus \{ 0 \} \}$ and pick
$\varphi_1 \in S_{\nu_1} \cap I-\{0\}$.\\
Put $\nu_2:=\min \{ \nu \; \vert \; \exists \; \varphi \in (S_{\nu} \cap I)
\setminus (S_{\nu} \cap <\varphi_1>)$ and pick
$\varphi_2 \in (S_{\nu_2}\cap I)\setminus (S_{\nu_2}\cap<\varphi_1>) $. \\
Proceed until we get $I= <\varphi_1,\dots,\varphi_m>$. Then
$(\varphi_1,\dots,\varphi_m)$ is a standard base of $I$ and
$\nu^\ast(I)=(\nu_1,\dots,\nu_m,\infty,\infty,\dots)$. \medbreak

\begin{remark}\label{rem1.5}
Let $\psi_1, \dots, \psi_\ell \in I$ be homogeneous generators of
$I$ such that $\nu_1 \leq \nu_2 \leq \dots \leq \nu_\ell$, where
$\nu_i = \deg(\psi_i)$. Then the above considerations show that
$$
(\nu_1, \nu_2, \ldots, \nu_\ell, \infty, \ldots ) \leq
\nu^\ast(I)\,,
$$
because a standard base of $I$ is obtained by possibly omitting
some of the $\psi_i$ for $i\geq 2$.
\end{remark}

\medbreak In what follows, for a $k$-vector space (or a
$k$-algebra) $V$ and for a field extension $K/k$ we write
$V_K=V\otimes_k K$. From Lemma \ref{lem0.2} the following is
clear.

\begin{lemma}\label{lem.nu-baseextension}
For the ideal $I_K \subseteq S_K$ we have
$$
\nu^\ast(I) = \nu^\ast(I_K)\,.
$$
\end{lemma}

A second invariant is the directrix. By \cite{H1} Ch. II \S4, Lemma 10, we have:
\medbreak

\begin{lemma}\label{lem.dir}
Let $K/k$ be a field extension. There is a smallest $K$-subvector
space $\cT(I,K) \subseteq (S_1)_K = \oplus_{i=1}^{n} KX_i$
such that
$$
(I_K \cap K[\cT(I,K)]) \cdot S_K = I_K,
$$
where $K[\cT(I,K)] = \Sym_K(\cT(I,K))\subseteq \Sym_K((S_1)_K) = S_K$.
In other words $\cT(I,K)\subset (S_1)_K$ is the minimal $K$-subspace
such that $I_K$ is generated by elements in $K[\cT(I,K)]$. For $K=k$
we simply write $\cT(I)=\cT(I,k)$.
\end{lemma}

Recall that $C(S/I) = \Spec(S/I)$ is called the cone of the graded
ring $S/I$.

\begin{definition}\label{def.dir.graded}
For any field extension $K/k$, the closed subscheme
$$
\Dir_K(S/I) = \Dir(S_K/I_K) \subseteq C(S_K/I_K) = C(S/I)\times_kK
$$
defined by the surjection $S_K/I_K \twoheadrightarrow
S_K/\cT(I,K)S_K$ is called the directrix of $S/I$ in $C(S/I)$ over $K$. By
definition
$$
\Dir_K(S/I) \cong \Spec(\Sym_K((S_1)_K/\cT(I,K)))\,.
$$
We define
$$
e(S/I)_K = \dim(\Dir_K(S/I)) = \dim(S) - \dim_K(\cT(I,K)) = n -
\dim_K(\cT(I,K))\,,
$$
so that $\Dir_K(S/I) \cong \mathbb A_K^{e(S/I)_K}$, and simply write $e(S/I)$ for $e(S/I)_K$ with $K=k$.
\end{definition}

\begin{remark}\label{rem.dir.graded}
(a)\, By definition $\Dir_K(S/I)$ is determined by the pair
$S_K \supseteq I_K$, but indeed it has an intrinsic definition
depending on $A_K$ for $A=S/I$ only: Let $S_A
= \Sym_k(A_1)$, which is a polynomial ring over $k$.
Then the surjection $S_K \rightarrow A_K$ factors
through the canonical surjection $\alpha_{A,K}: S_{A,K} \rightarrow
A_K$, and the directrix as above is identified with the
directrix defined via $S_{A,K}$ and $\ker(\alpha_{A,K})$. In this way
$\Dir(A_K)$ is defined for any graded $k$-algebra $A$ which is
generated by elements in degree 1.

(b)\, Similarly, for any standard graded $k$-algebra $A$, i.e., any finitely generated graded
$k$-algebra which is generated by $A_1$, we may
define its intrinsic $\nu$-invariant by $\nu^\ast(A) =
\nu^\ast(\ker( \alpha_A))$, where $\alpha_A: \Sym_k(A_1)
\twoheadrightarrow A$ is the canonical epimorphism. In the
situation of Definition \ref{def0.1} we have
$$
\nu^\ast(I) = (1, \ldots, 1, \nu^1(S/I),\nu^2(S/I), \ldots )\,,
$$
with $n-t$ entries of 1 before $\nu^1(S/I)$, where
$t=\dim_k(A_1)$, and $\nu^1(S/I)>1$.

(c)\, If $X$ is a variable, then obviously $\nu^\ast(A[X]) = \nu^\ast(A)$ and
$\nu^\ast(I[X]) = \nu^\ast(I)$ in the situation of \ref{def0.1}. On the other
hand, $\Dir_K(A[X]) \cong \Dir_K(A)\times_K\mathbb A^1_K$, i.e., $e_K(A[X]) = e_K(A)+1$.
\end{remark}

\begin{lemma}\label{dimDir.graded}
Let the assumptions be as in \ref{def.dir.graded}.
\begin{itemize}
\item[(1)] $e(S/I)_K \leq \dim(S/I)$. \item[(2)] For field
extensions $k\subset K\subset L$, we have $e(S/I)_K\leq e(S/I)_L$.
\item[(3)] The equality holds if one of the following conditions
holds:
\begin{itemize}
\item[($i$)] $L/K$ is separable (not necessarily algebraic).
\item[($ii$)] $e(S/I)_K =
\dim(S/I)$.
\end{itemize}
In particular it holds if $K$ is a perfect field.
\end{itemize}
\end{lemma}
\medbreak\noindent \textbf{Proof } The inequality in (1) is trivial, and (2)
follows since $\mathcal T(I,K)\otimes L \subseteq \mathcal T(I,L)$,
This in turn implies claim (3) for condition (ii). Claim (3) for
condition (i) is proved in \cite{H1} II, Lemma 12, p. 223, for arbitrary degree of
transcendence. The case of a finite separable extension is easy: It is obviously
sufficient to consider the case where $L/K$ is finite Galois
with Galois group $G$. Then, by Hilbert's theorem 90, for any
$L$-vector space $V$ on which $G$ acts in a semi-linear way the
canonical map $V^G\otimes_KL \rightarrow V$ is an isomorphism.
This implies that $\cT(I,L)^G \otimes_K L \mathop{\rightarrow}\limits^{\sim} \cT(I,L)$ and the claim
follows. $\square$

\bigbreak
Finally we recall the Hilbert function (not to be confused with the
Hilbert polynomial) of a graded algebra.
Let $\bN$ be the set of the natural numbers (including $0$) and let
$\bNN$ be the set of the functions $\nu: \bN\to \bN$.
We endow $\bNN$ with the product order defined by:
\begin{equation}\label{eq.prod.ord}
\nu\geq\mu \Leftrightarrow \nu(n)\geq\mu(n)
\quad \text{ for any } n\in \bN.
\end{equation}
\begin{definition}\label{def.hilb.fcn}
For a finitely generated graded $k$-algebra $A$ its Hilbert function
is the element of $\bNN$ defined by
$$
H^{(0)}(A)(n) = H(A)(n) = \dim_{k}(A_n)\quad (n\in \bN).
$$
For an integer $t\geq 1$ we define $H^{(t)}(A)$ inductively by:
$$
H^{(t)}(A)(n) =\sum_{i=0}^n H^{(t-1)}(A)(i).
$$
We note
$$
H^{(t-1)}(A)(n) =H^{(t)}(A)(n)-H^{(t)}(A)(n-1) \leq H^{(t)}(n).
$$
\end{definition}

\begin{remark}\label{rem.hilb.fcn} (a)\, Obviously, for any field extension $K/k$ we have
$$
H^{(0)}(A_K) = H^{(0)}(A)\,.
$$
(b)\, For a variable $X$ we have
$$
H^{(t)}(A[X]) = H^{(t+1)}(A)\,.
$$
(c)\, For any $\nu \in \bNN$ and any $t\in \bN$ define $\nu^{(t)}$ inductively by $\nu^{(t)}(n) =
\sum_{i=0}^{n} \nu^{(t-1)}(i)$.\\
\end{remark}

\smallskip
In a certain sense, the Hilbert function measures how far $A$ is away from being a
polynomial ring:

\begin{definition}\label{def.Phi}
Define the function $\Phi = \Phi^{(0)}\in \bNN$
by $\Phi(0) = 1$ and $\Phi(n) = 0$ for $n>0$.
Define $\Phi^{(t)}$ for $t \in \bN$ inductively as above, i.e., by $\Phi^{(t)}(n)=
\sum_{i=0}^n \Phi^{(t-1)}(i)$.
\end{definition}

Then one has
$$
\Phi^{(t)}(n) = H^{(0)}(k[X_1,\ldots,X_t])(n)=\binom{n+t-1}{n}
\quad\text{for all }n\geq 0.
$$
and:

\begin{lemma}\label{lem.HS.bound1} Let $A$ be a finitely generated graded algebra of dimension
$d$ over a field $k$, which is generated by elements in degree one (i.e., is a standard graded algebra).
\begin{itemize}
\item[(a)]
Then $H^{(0)}(A) \geq \Phi^{(d)}$, and equality holds if and
only if $A \cong k[X_1,\ldots,X_d]$.
\item[(b)]
For a suitable integer $m \geq 1$ one has $m\Phi^{(d)} \geq H^{(0)}$.
\item[(c)]
If $H^{(0)}(A) = \Phi^{(N)}$ for $N\in \bN$, then $N=d$.
\end{itemize}
\end{lemma}

\medbreak\noindent \textbf{Proof } (a): We may take a base change with an extension field $K/k$,
and therefore may assume that $k$ is infinite. In this case there
is a Noether normalization $i: S = k[X_1,\ldots,X_d] \hookrightarrow
A$ such that the elements $X_1,\ldots,X_d$ are mapped to $A_1$,
the degree one part of $A$, see \cite{Ku}, Chap. II, Theorem 3.1
d). This means that $i$ is a monomorphism of graded $k$-algebras. But
then $H^{(0)}(A) \geq H^{(0)}(k[X_1,\ldots,X_d]) = \Phi^{(d)}$,
and equality holds if and only if $i$ is an isomorphism.

\smallskip
(b): Since $A$ is a finitely generated $S$-module, it also has finitely many
homogenous generators $a_1,\ldots,a_m$ of degrees $d_1,\ldots,d_m$.
This gives a surjective map of graded $S$-modules
$$
\oplus_{i=1}^m S[-d_i] \twoheadrightarrow A,
$$
where $S[m]$ is $S$ with grading $S[m]_n=S_{n+m}$, and hence
$$
H(A)(n) \leq \sum_{i=1}^m \Phi^{(d)}(n-d_i) \leq m \Phi^{(d)}(n)\,.
$$
(c) follows from (a) and (b), because $\Phi^{(d)}$ and $\Phi^{(N)}$ have different
asymptotics for $N\neq d$.
\begin{flushright}$\square$\end{flushright}

\medbreak
We shall need the following property.

\begin{theorem}\label{thm.ordHF} Let $k$ be a field and for $n \in \bN$, let $HF_n \subset \bN^\bN$ be the
set of all Hilbert functions $H(A)$ of all standard graded $k$-algebras $A$ with $H(A)(1) \leq n$.

(a) $HF_n$ is independent of $k$.

(b) $HF_n$ is a noetherian ordered set, i.e.,
\begin{itemize}
\item[(i)] $HF_n$ is {\it well-founded}, i.e., every strictly descending sequence $H_1 > H_2 > H_3 > \ldots$
in $HF_n$ is finite.
\item[(ii)] For every infinite subset $M \subset HF_n$ there are elements $H, H' \in M$ with $H < H'$.
\end{itemize}
\end{theorem}

\medbreak\noindent \textbf{Proof } If $A$ is a standard graded $k$-algebra, then $H(A)(1)\leq n$ holds if and only if
$A$ is a quotient of $S_n=k[X_1,\ldots,X_n]$, i.e., $A = k[X_1,\ldots,X_n]/I$ for a homogeneous ideal $I$.
On the other hand, it is known that $H(S_n/I) = H(S_n/I')$ for some monomial ideal $I' \subset S_n$,
i.e., an ideal generated by monomials in the variables $X_1,\ldots,X_n$ (see \cite{CLO} 6 \S 3,
where this is attributed to Macaulay; to wit, one may take for $I'$ the ideal of
leading terms for $I$, with respect to the lexicographic order on monomials, loc. cit. Proposition 9).
This proves (a). Moreover, for (b) we may assume that all considered Hilbert functions are of the form
$H(S_n/I)$ for some monomial ideal $I \subset S_n$. Thus in (b) (ii) we are led to the consideration
of an infinite set $N$ of monomial ideals in $S_n$. But then the main theorem in \cite{Macl} (Theorem 1.1)
says that there are $I, I' \in N$ with $I \subsetneqq I'$, so that $H(S_n/I) > H(S_n/I')$. For (b) (i)
we may again assume that all $H_i$ are of the form $H(S_n/I_i)$ for some monomial ideals $I_i \subset S_n$,
and by \cite{Macl} 1.1 that one finds an infinite chain
$I_{i(1)} \subsetneqq I_{i(2)} \subsetneqq I_{i(3)} \subsetneqq \ldots $ among these ideals,
necessarily strictly increasing since the sequence of the $H_{i(\nu)}$ is, which is a contradiction
because $S_n$ is noetherian.

\begin{remark} (a) For the functions $H^{(1)}(A)$ instead of the functions $H(A)$ property (b) (i)
was shown in \cite{BM2} Theorem 5.2.1. \\
(b) Our proof was modeled on \cite{AP} Corollary 3.6, which is just formulated for monomial ideals.
See also \cite{AP} Corollary 1.8 for another argument that the set of monomial
ideals in $S_n$, with respect to the reverse inclusion, is a noetherian ordered set.
Finally, in \cite{AH} it is shown, by more sophisticated methods, that even the set
$HF$ of all Hilbert functions is noetherian.
\end{remark}

\bigskip\noindent
{\bf 1.2 Invariants for local rings}

\bigskip
For any ring $R$ and any prime ideal $\fp\subset R$ we set
$$
\grpR=\underset{n\geq 0}{\bigoplus}\;\fp^n/\fp^{n+1},
$$
which is a graded algebra over $R/\fp$.

\bigskip
In what follows we assume that $R$
is a noetherian regular local ring with maximal ideal $\fm$ and residue field $k=R/\fm$.
Moreover, assume that $R/\fp$ is regular. Then we have
$$
\grpR=\Sym_{R/\fp}\big(\grp^1(R)\big), \quad \grp^1(R)=\fp/\fp^2.
$$
where $\Sym_{R/\fp}(M)$ denotes the symmetric algebra of a free
$R/\fp$-module $M$. Concretely, let $(y_1,\dots, y_r,u_1,\dots,u_e)$ be a system of regular parameters for $R$
such that $\fp=(y_1,\dots,y_r)$.
Then $\grpR$ is identified with a polynomial ring over $R/\fp$:
$$
\grpR = (R/\fp) [Y_1, \dots, Y_r]
\quad\text{where  }
Y_i = y_i \mod \fp^2 \; (1\leq i\leq r).
$$
\medbreak
Fix an ideal $J \subset R$. In case $J\subset\fp$ we set
$$
\grpRJ=\underset{n\geq 0}{\bigoplus}(\fp/J)^n/(\fp/J)^{n+1},
$$
and define an ideal $\InpJ\subset \grpR$ by the exact sequence
$$ 0 \to \InpJ \to \grpR \to \grpRJ \to 0 .$$
Note
$$
\InpJ = \bigoplus\limits_{n \geq 0} (J \cap \fp^n + \fp^{n +1})/\fp^{n+1}.
$$
\bigskip

For $f \in R$ and a prime ideal $\fp\subset R$ put
\begin{equation}\label{order}
 v_\fp(f) = \left\{ \begin{array}{ccc}
\max \{ \nu \; \vert \; f \in \fp^{\nu} \} & , & f \ne 0, \\
\infty & , & f = 0 \\
\end{array} \right.
\end{equation}
called the {\it order} of $f$ at $\fp$. For prime ideals
$\fp\subset\fq \subset R$, we have the following semi-continuity result
(cf. \cite{H1} Ch.III \S3):
\begin{equation}\label{eq.ordersc}
v_\fp(f)\leq v_\fq(f)\qfor \forall f\in R.
\end{equation}

\medbreak
Define the {\it initial form} of $f$ at $\fp$ as
$$ \inp(f) := f \mod \fp^{v_\fp(f) +1} \in \fp^{v_\fp(f)} /\fp^{v_\fp(f) +1}\;
\in \grpR.
$$
In case $J\subset \fp$ we have
$$ \InpJ = \{\inp(f)\; |\; f\in J\} .$$
\medbreak

\begin{definition}\label{def0.5}
\begin{itemize}
\item[(1)]
A system $(f_1, \dots, f_m)$ of elements in $J$ is a {\it standard base} of
$J$, if
\par\noindent
$(in_{{\mathfrak m}} (f_1), \dots, in_{{\mathfrak m}} (f_m))$ is
a standard base of $In_{{\mathfrak m}} (J)$ in the polynomial ring
$\gr_{{\mathfrak m}} (R)$.
\item[(2)]
We define $\nu^*(J,R)$ as the $\nu^*$-invariant (cf. Definition \ref{def0.1})
for $In_{{\mathfrak m}} (J)\subset \gr_{{\mathfrak m}} (R)$.
\item[(2)]
The absolute $\nu^\ast$-invariant $\nu^\ast(\cO)$ of a noetherian local ring
$\cO$ with maximal ideal $\fn$ is defined as the absolute
$\nu^\ast$-invariant (cf. Remark \ref{rem.dir.graded} (b))
$\nu^\ast(gr_{\fn}\cO)$.
\end{itemize}
\end{definition}

It is shown in \cite{H3} (2.21.d) that a standard base
$(f_1,\dots,f_m)$ of $J$ generates $J$.

\bigskip

Next we define the directrix $\Dir(\cO)$ of any noetherian local ring $\cO$.
First we introduce some basic notations. \medbreak
Let $\fn$ be the maximal ideal of $\cO$, let
$x$ be the corresponding closed point of $\Spec(\cO)$,
and let $k(x) = \cO/\fn$ be the residue field at $x$.
Define \medbreak

$T(\cO)=\Spec(\Sym_k(\fn/\fn^2))$: the Zariski tangent space of $\Spec(\cO)$ at $x$,
\medbreak

$C(\cO) = \Spec (\gr_{\fn}(\cO)) =
C(\gr_{\fn}(\cO))$: the {\it (tangent) cone}
of $\Spec(\cO)$ at $x$. \medbreak

We note that $\dim C(\cO) = \dim \cO$ and that the map
$$
\Sym_k({\fn}/{\fn}^2)
\twoheadrightarrow \gr_{\fn} (\cO)
$$
gives rise to a closed immersion
$C(\cO) \hookrightarrow T(\cO)$.
\medbreak

\begin{definition}\label{def.dir}
Let $K/k(x)$ be a field extension. Then the directrix of $\cO$ over
$K$,
$$\Dir_K(\cO) \subseteq C(\cO)\times_{k(x)}K \subseteq T(\cO)\times_{k(x)}K$$
is defined as the directrix $\Dir_K(\gr_{\fn}(\cO)) \subseteq C(\gr_{\fn}(\cO))$ of
$\gr_{\fn}(\cO)$ over $K$ (cf. Remark \ref{rem.dir.graded} (a)).
We set
$$
e(\cO)_K = \dim(\Dir_K(\cO))
$$
and simply write $e(\cO)$ for $e(\cO)_K$ with $K=k(x)$.
\end{definition}

\begin{remark}\label{rem.dir}
By definition, for $R$ regular as above and an ideal $J\subset \fm$ we have
$$
\Dir_K(R/J) = \Spec(\gr_\fm(R)_K/\cT(J,K)\gr_\fm(R)_K)
 \subseteq \Spec(\gr_\fm(R)_K/In_\fm(J)_K) = C(R/J)_K\,,
$$
where $\cT(J,K) = \cT(In_\fm(J),K) \subseteq \gr^1_{\fm} (R)_K$ is
the smallest $K$-sub vector space such that
$$
(In_\fm(J)_K \cap K[\cT(J,K)] \cdot \gr_\fm(R)_K = In_\fm(J)_K\,,
$$
i.e., such that $In_\fm(J)_K$ is generated by elements in
$K[\cT(J,K)]$. Moreover
$$
\Dir_K(R/J) \cong \Sym_K(\gr_\fm^1(R)_K/\cT(J,K))\subseteq T(R)_K\,.
$$
For $K=k$ we simply write $\cT(J)=\cT(J,k)$.
\end{remark}

\begin{lemma}\label{dimDir}
Let the assumptions be as above.
\begin{itemize}
\item[(1)]
$e(R/J)_K \leq \dim(R/J)$.
\item[(2)]
For field extensions $k\subset K\subset L$, we have
$e(R/J)_K\leq e(R/J)_L$. The equality holds if one of the
following conditions holds:
\begin{itemize}
\item[($i$)] $L/K$ is separable (not necessarily algebraic).
\item[($ii$)] $e(R/J)_K = \dim(R/J)$.
\end{itemize}
\end{itemize}
\end{lemma}
\medbreak\noindent \textbf{Proof } This follows from Lemma
\ref{dimDir.graded}, because $\dim(R/J)=\dim C(R/J)$. $\square$
\medbreak

\begin{definition}\label{def.dir3}
We define
$\eb(R/J)=e(R/J)_{\overline{k}}$ for an algebraic closure $\overline{k}$
of $k$. By Lemma \ref{dimDir} we have
$e(R/J)_K \leq \eb(R/J) \leq \dim(R/J)$ for any extension $K/k$.
\end{definition}
\medbreak

For later use we note the following immediate consequence of the construction
of a standard base below Corollary \ref{cordef0.1}.

\begin{lemma}\label{stbTJ}
Let $\cT=\cT(J)$ be as in Lemma \ref{lem.dir} . There exists a
standard base $(f_1,\dots,f_m)$ of $J$ such that $in_\fm(f_i)\in
k[\cT]$ for for all $i$.
\end{lemma}

\bigskip

Finally, we define the Hilbert functions of a noetherian local ring $\cO$ with
maximal ideal $\fm$ and residue field $F$ as
those of the associated graded ring:
$$
\HSff t \cO(n) = H^{(t)}(\gr_\fm(\cO))\,.
$$
Explicitly, the Hilbert function is the element of
$\bNN$ defined by
$$
\HSff 0 \cO(n) =\dim_{F}(\fm^n/\fm^{n+1})\quad (n\in \bN).
$$
For an integer $t\geq 1$ we define $\HSff t \cO $ inductively by:
$$
\HSff t \cO(n) =\sum_{i=0}^n \HSff {t-1} \cO(i).
$$
In particular, $\HSff 1 \cO(n)$ is the length of the $\cO$-module
$\cO/\fm^{n+1}$ and $\HSff 1 \cO$ is called the Hilbert-Samuel function of $\cO$.
We will sometimes call all functions of the form $\HSff n \cO$ Hilbert-Samuel functions.
The Hilbert function measures how far away $\cO$ is from being a regular ring:

\begin{lemma}\label{lem.HS.bound2}
Let $\cO$ be a noetherian local ring of dimension $d$, and define $\Phi^{(t)}$ as in
Definition \ref{def.Phi}. Then $H_\cO^{(0)}
\geq \Phi^{(d)}$, and equality holds if and only if $\cO$ is
regular.
\end{lemma}

\medbreak\noindent \textbf{Proof } (see also \cite{Be} Theorem (2) and
\cite{Si2}, property (4) on p. 46) Since $\dim(\cO)
= \dim(gr_\fm \cO)$, where $\fm$ is the maximal ideal of $\cO$,
and since $\cO$ is regular if and only if $gr_\fm \cO \cong
k[X_1,\ldots,X_d]$ where $k=\cO/\fm$, this follows from Lemma \ref{lem.HS.bound1}.

\bigbreak\noindent
For later purpose, we note the following facts.

\begin{lemma}\label{HSquotient}
Let the assumptions be as above and let $\overline{\cO}=\cO/\fa$ be
a quotient ring of $\cO$. Then $\HSff t \cO \geq \HSff t {\overline{\cO}}$
and the equality holds if and only if $\cO=\overline{\cO}$.
\end{lemma}
\medbreak
\textbf{Proof  }
Let $\overline{\fm}$ be the maximal ideal of $\overline{\cO}$.
The inequality holds since the natural maps
$\fm^{n+1}\to  \overline{\fm}^{n+1}$ are surjective.
Assume $\HSff s \cO = \HSff s {\overline{\cO}}$ for some $s\geq 0$.
By Definition \ref{def.hilb.fcn} it implies
$\HSff t \cO = \HSff t {\overline{\cO}}$ for all $t\geq 0$,
in particular for $t=1$. This implies that the natural maps
$\pi_n:\cO/\fm^{n+1}\to \overline{\cO}/\overline{\fm}^{n+1}$
are isomorphisms for all $n\geq 0$. Noting
$\Ker(\pi_n)\simeq \fa/\fa\cap\fm^{n+1}$, we get
$\displaystyle{\fa\subset \underset{n\geq 0}{\cap}\fm^{n+1}=(0)}$
and hence $\cO=\overline{\cO}$.
$\square$
\medbreak

\begin{lemma}\label{lem.HS.comparison}
Let $\cO$ and $\cO'$ be noetherian local rings.
\begin{itemize}
\item[(a)]
For all non-negative integers $a$ and $e$ one has
$$
\dim \cO \geq e   \quad\Longleftrightarrow\quad H_\cO^{(a)} \geq \Phi^{(e+a)}\,.
$$
\item[(b)]
For all non-negative integers $a$ and $b$ one has
$$
H_\cO^{(a)} \geq H_{\cO'}^{(b)} \quad \Longrightarrow \quad \dim \cO + a \geq \dim \cO' + b.
$$
\end{itemize}
\end{lemma}
\medbreak\noindent \textbf{Proof } Let $d = \dim \cO$. (a):
If $d \geq e$, then we get $H_\cO^{(a)}\geq \Phi^{(d+a)} \geq \Phi^{(e+a)}$, by Lemma \ref{lem.HS.bound2}.
Conversely, assume $H_\cO^{(a)}\geq \Phi^{(e+a)}$. Then form Lemma \ref{lem.HS.bound1} we get
$$
m \Phi^{(d+a)} \geq H_\cO^{(a)}\geq \Phi^{(e+a)}
$$
for some integer $m \geq 1$. If $d<e$ this is a contradiction, because of the asymptotic behavior of
$\Phi^{(t)}$. Hence $d\geq e$.

For (b) let $d' = \dim \cO'$. If $H_\cO^{(a)} \geq H_{\cO'}^{(b)}$ then
$$
H_\cO^{(a)} \geq H_{\cO'}^{(b)} \geq \Phi^{(d'+b)}\,
$$
and by (a) we have $d \geq d' + b -a$. (Note: If $d'+b-a<0$, the statement is empty; if $e=d'+b-a\geq 0$,
then we can apply (a).)

\bigskip\noindent
{\bf 1.3 Invariants for schemes}

\bigskip
Let $X$ be a locally noetherian scheme.

\begin{definition}\label{def.invschemes}
For any point $x \in X$ define
$$
 \nu^\ast_x(X) = \nu^\ast(\cO_{X,x})\qaq \Dir_x(X) = \Dir(\cO_{X,x})
$$
and
$$
e_x(X)= e(\cO_{X,x}) = \dim_{k(x)}(\Dir_x(X)),\quad \eb_x(X)=
\eb(\cO_{X,x}),\quad e_x(X)_K= e(\cO_{X,x})_K,
$$
where $K/k(x)$ is a field extension. If $X \subseteq Z$ is a
closed subscheme of a (fixed) regular excellent scheme $Z$, define
$$
\nu^\ast_x(X,Z) = \nu^\ast(J_x,\cO_{Z,x})\,,
$$
where $\cO_{X,x} = \cO_{Z,x}/J_x$. We also define
$$
\IDir_x(X)\subset \gr_{\fm_x}(\cO_{Z,x})
$$
to be the ideal defining
$\Dir_x(X)\subset T_x(Z)=\Spec(\gr_{\fm_x}(\cO_{Z,x}))$, where
$\fm_x$ is the maximal ideal of $\cO_{Z,x}$.
\end{definition}

We note that always
$$
\Dir_x(X)\subseteq C_x(X)\subseteq T_x(X) \;\;(\;\; \subseteq T_x(Z) \mbox{   for   } X \subset Z \mbox{   as above   }\;)   ,
$$
where $C_x(X)=C(\cO_{X,x})$ is called the tangent cone of $X$
at $x$ and $T_x(X)=T(\cO_{X,x})$ is called the Zariski
tangent space of $X$ at $x$ (similarly for $Z$).

\begin{lemma}\label{completion.nu.2}
Let $X$ be a locally noetherian scheme.
\begin{itemize}
\item[(1)] Let $\pi: X'\to X$ be a morphism with $X'$ locally noetherian, and let $x'\in X'$
be a point lying over $x\in X$. Assume that $\pi$ is quasi-\'etale
at $x'$ in the sense of Bennett \cite{Be} (1.4), i.e., that $\cO_{X,x} \rightarrow \cO_{X',x'}$
is flat and $\fm_x\cO_{X',x'} = \fm_{x'}$ where $\fm_x \subset \cO_{X,x}$ and
$\fm_{x'}\subset \cO_{X',x'}$ are the respective maximal ideals. (In particular,
this holds if $\pi$ is \'etale.)
Then there is a canonical isomorphism
\begin{equation}\label{dir.etale.2}
C_{x'}(X') \cong C_x(X)\times_{k(x)}k(x')
\end{equation}
so that $\nu^\ast_{x'}(X') = \nu^\ast_x(X)$. If $k(x')/k(x)$ is separable, then
$$
 \Dir_{x'}(X') \cong \Dir_x(X)\times_{k(x)}k(x'), \quad\mbox{and hence}\quad e_{x'}(X') = e_x(X) .
$$
Consider in addition that there is a cartesian diagram
$$
\begin{CD}
X' @>{i}>> Z' \\
@V{\pi}VV @VV{\pi_Z}V \\
X @>{i}>> Z
\end{CD}
$$
where $i$ and $i'$ are closed immersions, $Z$ and $Z'$ are regular excellent schemes
and $\pi_Z$ is quasi-\'etale at $x'$. Then
$$
T_{x'}(Z) \cong T_x(Z)\times_{k(x)}k(x') \qaq \nu^\ast_{x'}(X',Z')
= \nu^\ast_x(X,Z)\,.
$$
\item[(2)]
Let $\pi: X' \rightarrow X$ be a morphism and let $x \in X$ and $x' \in X$ with $\pi(x') = x$.
Assume that $X'$ is locally noetherian, that $\pi$ is flat and that the fibre $X'_x=X'\times_Xx$ of $\pi$
over $x$ is regular at $x'$ (e.g., assume that $\pi$ is smooth around $x'$).
Then there are non-canonical isomorphisms
\begin{equation}\label{eq.dir.smooth}
C_{x'}(X') \cong C_x(X)\times_{k(x)}C_{x'}(X'_x) \cong C_x(X)\times_{k(x)}\mathbb A_{k(x)}^d\,,
\end{equation}
where $d = \dim(\cO_{X'_x,x'}) = \codim_{X'_x}(x')$. Hence
$$
\Dir_{x'}(X') \cong \Dir_x(X)\times_{k(x)}\mathbb A_{k(x')}^d, \quad
e_{x'}(X') = e_x(X) + d \qaq \nu^\ast_{x'}(X') = \nu^\ast_x(X)\,.
$$

\item[(3)] Let $X_0$ be an excellent scheme, and let $f: X \to X_0$ be a morphism with $X$ locally noetherian.
Let $x_0 \in X_0$ and
$$
\hat{X_0} = \Spec(\hat{\cO}_{X_0,x_0}), \quad \hat{X} = X \times_{X_0} \hat{X_0}\,,
$$
where $\hat{\cO}_{X_0,x_0}$ is the completion of ${\cO}_{X_0,x_0}$. Then
for any $x \in X$ and any $\hat{x} \in \hat{X}$
lying over $x$, there are non-canonical isomorphisms
\begin{equation}\label{eq.dir.completion.2}
C_{\hat{x}}(\hat{X}) \cong C_{x}(X)\times_{k(x)} C_{\hat{x}}({\hat{X}}_{x}) \cong C_x(X)\times_{k(x)}\mathbb A_{k(\hat{x})}^d\,,
\end{equation}
where $\hat{X}_x = \hat{X}\times_{X}x$ is the fibre over $x$ for the morphism $\pi: \hat{X} \rightarrow X$ and
$d = \dim(\cO_{{\hat{X}}_x,\hat{x}})=\codim_{\hat{X}_x}(\hat{x})$. Hence
$$
\Dir_{\hat{x}}(\hat{X}) \cong \Dir_x(X)\times_{k(x)}\mathbb A_{k(\hat{x})}^d, \quad
e_{\hat{x}}(\hat{x}) = e_x(X) + d \qaq \nu^\ast_{\hat{x}}(\hat{X}) =
\nu^\ast_x(X)\,.
$$
Assume further that there is a commutative diagram
\begin{equation}\label{eq.Z-diagram}
\begin{CD}
X @>{i}>> Z\\
@V{f}VV @VV{g}V \\
X_0  @>{i_0}>> Z_0
\end{CD}
\end{equation}
where $Z_0$ is a regular excellent scheme, $Z$ is a regular locally noetherian scheme,
and $i_0$ and $i$ are closed immersions.
Denote
$$
\hat{Z_0}=\Spec(\hat{\cO}_{Z_0,x_0}),\quad \hat{Z}= Z\times_{Z_0}\hat{Z_0},
$$
where $\hat{\cO}_{Z_0,x_0}$ is the completion of ${\cO}_{Z_0,x_0}$, so
that $\hat{X}= X\times_{Z_0}\hat{Z}=X\times_{Z}\hat{Z}$
can be regarded as a closed subscheme of $\hat{Z}$. Then $\hat{Z}$ is regular, and
$$
\nu^\ast_{\hat{x}}(\hat{X},\hat{Z}) = \nu^\ast_x(X,Z) \,.
$$
\end{itemize}
\end{lemma}
\medbreak\noindent \textbf{Proof } (1): It suffices to show
(\ref{dir.etale.2}). Let $(A,\fm_A) \rightarrow (B,\fm_B)$ be a flat
local morphism of local noetherian rings, with $\fm_A B = \fm_B$.
Then we have isomorphisms
\begin{equation}\label{etale.iso.2}
\fm_A^n\otimes_A B \mathop{\longrightarrow}\limits^{\sim} \fm_B^n
\end{equation}
for all $n\geq 0$. In fact, this holds for $n=0$, and, by
induction and flatness of $B$ over $A$, the injection $\fm_A^n
\hookrightarrow \fm_A^{n-1}$ induces an injection
$$
\fm_A^n\otimes_A B \hookrightarrow \fm_A^{n-1}\otimes_A B
\mathop{\rightarrow}\limits^{\sim} \fm_B^{n-1}\,,
$$
whose image is $\fm_A^n B = \fm_B^n$.
From \eqref{etale.iso.2} we now deduce isomorphisms for all $n$
$$
(\fm_A^n/\fm_A^{n+1})\otimes_{k_A}k_B \cong
(\fm_A^n/\fm_A^{n+1})\otimes_A B
\mathop{\longrightarrow}\limits^{\sim} \fm_B^n/\fm_B^{n+1}\,,
$$
where $k_A= A/\fm_A$ and $k_B = B/\fm_B$, and hence the claim
(\ref{dir.etale.2}).

(2): As for \eqref{eq.dir.completion.2}, consider the local rings $A= \cO_{X,x}$ with maximal
ideal $\fm$ and residue field $k=A/\fm$, and $A' = \cO_{X',x'}$, with maximal
ideal $\fm'$ and residue field $k' = A'/\fm'$.
By assumption, $\pi$ is flat and has regular fibers, and hence the same is true for the local morphism
$A \rightarrow A'$ since it is obtained from $\pi$ by localization in $X$ and $X'$.
$A\otimes_{\cO_{X,x}} \hat{\cO}_{X,x}$.
Hence, by \cite{Si2}, Lemma (2.2), the closed subscheme $\Spec(A'/\fm A') \hookrightarrow \Spec(A')$ is permissible,
i.e. it is regular and $\gr_{\fm A'}(A')$ is flat over $A'/\fm A'$.
By \cite{H1} Ch. II, p. 184, Proposition 1, we get a non-canonical isomorphism
$$
\gr_{\fm'}(A') \cong (\gr_{\fm A'}(A')\otimes_{A'}k')\otimes_{k'} \gr_{\fm'/\fm A'}(A'/\fm A')\,.
$$
On the other hand, by flatness of $A'$ over $A$ we get canonical isomorphisms
$$
\fm^n A'/\fm^{n+1} A' \cong (\fm^n/\fm^{n+1})\otimes_AA'
$$
for all $n$. Hence we get an isomorphism $\gr_{\fm A'}(A')\otimes_{A'}k' \cong \gr_\fm(A)\otimes_AA'$
and the above isomorphism becomes
$$
\gr_{\fm'}(A') \cong  \gr_\fm(A)\otimes_k \gr_{\fm'/\fm A'}(A'/\fm A')\,,
$$
which is exactly the first isomorphism in \eqref{eq.dir.completion.2}.
Since $A'/\fm A' = \cO_{X'_x,x'}$
is regular of dimension $d$, we have an ismorphism $\gr_{\fm'/\fm A'}(A'/\fm A') \cong k'[T_1,\ldots,T_d]$,
which gives the second isomorphism in \eqref{eq.dir.completion.2}.

The statements for the directrix and the $\nu^\ast$-invariant of $X'$ at $x'$
now follow from Remark \ref{rem.dir.graded} (c).

(3): The claims not involving $Z_0$ and $Z$ follow by applying (2) to the morphism
$$
\tilde{\pi}: \hat{X} = X\times_{X_0}\Spec(\cO_{X_0,x_0})\times_{\Spec(\cO_{X,x})}\hat{X_0}\rightarrow X\times_{X_0}\Spec(\cO_{X_0,x_0})\,.
$$
In fact, since $X_0$ is excellent, the morphism $\pi_0: \cO_{X_0,x_0} \rightarrow \hat{\cO}_{X_0,x_0}$ is flat,
with geometrically regular fibers, and so the same is true for $\tilde{\pi}$, which is a base change of $\pi_0$.
Now consider the diagram \eqref{eq.Z-diagram}. The above, applied to $Z \rightarrow Z_0$, gives isomorphisms
$$
C_{\hat{x}}(\hat{Z}) \cong C_x(Z)\times_{k(x)}C_{\hat{x}}({\hat{Z}}_{x}) \cong C_x(Z)\times_{k(x)}\mathbb A_{k(\hat{x})}^d \cong \mathbb A_{k(\hat{x})}^{N+d}
$$
where $N = \dim(\cO_{Z,x})$, because $Z$ is regular and $\hat{Z}_x \cong \hat{X}_x$.
By Lemmas \ref{lem.HS.bound2} and \ref{lem.HS.bound1} this
implies that $\hat{Z}$ is regular at $x$. The final equality in (3) follows from the isomorphism
$T_{\hat{x}}(\hat{Z}) \cong C_{\hat{x}}(\hat{Z})$ and Remark \ref{rem.dir.graded} (c).
$\square$

\bigskip\noindent
Next we introduce Hilbert-Samuel functions for excellent schemes.
Recall that an excellent scheme $X$ is catenary so that
for any irreducible closed subschemes
$Y\subset Z$ of $X$, all maximal chains of irreducible closed subschemes
$$Y=Y_0\subset Y_1\subset \cdots\subset Y_r=Z$$
have the same finite length $r$ denoted by $\codim_Z(Y)$.
For any irreducible closed subschemes $Y\subset Z\subset W$ of $X$ we have
$$\codim_W(Y)=\codim_W(Z) + \codim_Z(Y).$$

\begin{definition}\label{def.HS.schemes}
Let $X$ be a locally noetherian catenary scheme (e.g., an excellent scheme),
and fix an integer $N \geq \dim X$. (Recall that all schemes are assumed to be finite
dimensional.)
\begin{itemize}
\item[(1)]
For $x \in X$ let $I(x)$ be the set of irreducible components $Z$ of $X$ with $x\in Z$.
\item[(2)]
Define the function $\phi_X := \phi_X^N: X \rightarrow \bN$ by $\phi_X(x) = N - \psi_X(x)$,
where
$$
\psi_X(x) = \mbox{min } \{ \codim_Z(x) \mid Z \in I(x) \}\,.
$$
\item[(3)]
Define $H_X := H_X^N:  X \rightarrow \bNN$ as follows. For $x\in X$ let
$$
H_X(x) = H_{\cO_{X,x}}^{(\phi_X(x))} \in \bNN\,.
$$
\item[(4)]
For $\nu\in \bNN$ we define
$$
\HSgeqa X:=\{x\in X|\; \HSf X x \geq \nu\}
\quad\text{ and }\quad
\HSa X:=\{x\in X|\; \HSf X x =\nu\}.
$$
The set $X(\nu)$ is called the Hilbert-Samuel stratum for $\nu$.
\end{itemize}
\end{definition}

By sending $Z$ to its generic point $\eta$, the set $I(x)$ can be identified with the set
of generic points (and hence the set of irreducible components) of the local ring
$\cO_{X,x}$. Therefore $\psi_X(x)$ depends only on $\cO_{X,x}$, and $\phi_X(x)$ and $\HSf X x$
depend only on $N$ and $\cO_{X,x}$.

\begin{remark}\label{rem.HSfunct}
(a) The choice of $N$ does not matter. For sake of definiteness, we could have taken
$N = \dim X$, and the readers are invited to do so, if they prefer.
But when dealing with two different schemes $X$ and $X'$, the
difference $\dim X - \dim{X'}$ would always appear. Note that we can even have
$\dim U < \dim X$ for an open subscheme $U \subseteq X$. In the applications,
there will usually be a common bound for the dimensions of the schemes considered,
which we can take for $N$.

(b) A more sophisticated way would be to consider the whole array
$$
\underline{H}_X := (H_X^0, H_X^1, H_X^2,\ldots )
$$
where, for a function $\varphi: \bN \rightarrow \bN$ and possibly negative
$m \in \bZ$, $\varphi^{(m)}: \bN \rightarrow \bN$ is defined inductively by
the formulae
$$
\varphi^{(0)} = \varphi \quad,\quad \varphi^{(m+1)}(n) = \sum_{i=0}^{n} \varphi^{(m)}(i) \quad,
\quad \varphi^{(m-1)}(n) = \varphi^{(m)}(n) - \varphi^{(m)}(n-1)\;.
$$
Then $(\varphi^{(m_1)})^{(m_2)} = \varphi^{(m_1+m_2)}$ for $m_1, m_2 \in \bZ$,
and for a second function $\psi: \bN \rightarrow \bN$ one has
$$
\varphi \leq \psi \; \Longrightarrow \; \varphi^{(1)} \leq \psi^{(1)}\,,
$$
(in the product order \eqref{eq.prod.ord}), but the converse does not hold in general.
With these definitions we could use the ordering
$$
\varphi \preceq \psi \; :\Longleftrightarrow \; \varphi^{(m)} \leq \psi^{(m)} \quad \mbox{ for all }m \gg 0
$$
(not to be confused with the ordering $\varphi(n) \leq \psi(n)$ for $n\gg 0$ ) in all
applications below, and no choice of $N$ appears.
Note that $\varphi = \psi \Longleftrightarrow \varphi^{(m)} = \psi^{(m)}$,
for all $m\in\bZ$.
\end{remark}

In the rest of the paper, a choice of $N$ will be assumed and often suppressed in the notations.
We shall need the following semi-continuity property of $\phi_X$.

\begin{lemma}\label{lem.phi.sc}
\begin{itemize}
\item[(1)]
For $x,y \in X$ with $x \in \overline{\{y\}}$ one has $I(y) \subseteq I(x)$ and
$$
\phi_X(y) \leq \phi_X(x) + \codim_{\overline{\{y\}}}(x)\,.
$$
\item[(2)]
For $y \in X$, there is a non-empty open subset $U \subseteq \overline{\{y\}}$
such that
$$
I(x) = I(y) \qaq \phi_X(y) = \phi_X(x) + \codim_{\overline{\{y\}}}(x)
$$
for all $x\in U$.
\end{itemize}
\end{lemma}

\medbreak\noindent \textbf{Proof } (1): The inclusion $I(y) \subseteq I(x)$ is clear, and for
$Z \in I(y)$ one has
\begin{equation}\label{eq.codim}
\codim_Z(x) = \codim_Z(y) + \codim_{\overline{\{y\}}}(x)\,.
\end{equation}
Thus $\psi_X(x) \leq \psi_X(y) + \codim_{\overline{\{y\}}}(x)$, and
the result follows.
\medbreak
(2): Let $Z_1,\ldots,Z_r$ be the irreducible components of $X$ which do not
contain $y$. Then we may take $U = \overline{\{y\}}\cap (X \setminus \bigcup_{i=1}^r\,Z_i)$.
In fact, if $x \in U$, then $I(y) = I(x)$, and from \eqref{eq.codim} we get
$$
\psi_X(x) = \psi_X(y) + \codim_{\overline{\{y\}}}(x)\,.
$$

\bigskip
Now we study the Hilbert-Samuel function $H_X$. The analogue of Lemma \ref{lem.HS.bound1} (a)
and Lemma \ref{lem.HS.bound2} is the following, where $N$ is as in \ref{def.HS.schemes}.

\begin{lemma}\label{lem.HS.bound3}
For $x \in X$ one has $H_X(x) \geq \Phi^{(N)}$,
and equality holds if and only if $x$ is a regular point.
\end{lemma}

\medbreak\noindent \textbf{Proof } We have
$$
H_X(x) = H_{\cO_{X,x}}^{(N-\psi_X(x))}\geq \Phi^{(N-\psi_X(x)+\codim_X(x))}\geq \Phi^{(N)}\,.
$$
Here the first inequality follows from Lemma \ref{lem.HS.bound2}, and
the second inequality holds because $\codim_X(x) \geq \psi_X(x)$. If $H_X(x) = \Phi^{(N)}$,
then all inequalities are equalities, and hence, again by Lemma \ref{lem.HS.bound2},
$x$ is a regular point. Conversely, if $x$ is regular, then there is only one
irreducible component of $X$ on which $x$ lies, and hence $\codim_X(x) = \psi_X(x)$,
so that the second inequality is an equality. Moreover, by Lemma \ref{lem.HS.bound2},
the first inequality is an equality.

\begin{remark}\label{HSregular}
In particular, $X$ is regular if and only if
$\HSa X=\emptyset$ for all $\nu\in \bNN$ except for $\nu=\nu_X^{reg}$,
where $\nu_X^{reg}=\HSf X x$ for a regular point $x$ of $X$ which is
independent of the choice of a regular point, viz., equal to $\Phi^{(N)}$.
\end{remark}

We have the following important upper semi-continuity of the Hilbert-Samuel function.

\begin{theorem}\label{HSf.usc}
Let $X$ be a locally noetherian catenary scheme.
\begin{itemize}
\item[(1)]
If $x\in X$ is a specialization of $y\in X$, i.e., $x \in \overline{\{y\}}$, then
$H_X(x) \geq H_X(y)$.
\item[(2)]
For any $y\in X$, there is a dense open subset $U$ of $\overline{\{y\}}$
such that
$H_X(x) = H_X(y)$ for all $x\in U$.
\item[(3)]
The function $H_X$ is upper semi-continuous, i.e., for any $\nu\in \bNN$,
$\HSgeqa X$ is closed in $X$.

\end{itemize}
\end{theorem}

\medbreak\noindent \textbf{Proof } (1): We have
$$
H_X(x) = H_{\cO_{X,x}}^{(\phi_X(x))}
\geq H_{\cO_{X,y}}^{(\phi_X(x)+\codim_{\overline{\{y\}}}(x))}
\geq H_{\cO_{X,y}}^{(\phi_X(y))} = H_X(y)\,.
$$
Here the first inequality holds by results of Bennett (\cite{Be}, Theorem (2)),
as improved by Singh (\cite{Si1}, see p. 202, remark after Theorem 1),
and the second holds by Lemma \ref{lem.phi.sc} (1).

\smallskip
(2): First of all, there is a non-empty open set $U_1 \subseteq \overline{\{y\}}$
such that $\overline{\{y\}}\subseteq X$ is permissible (see \ref{def.nf}) at each $x \in U_1$ (\cite{Be}
Ch. 0, p. 41, (5.2)). Then
$$
H_{\cO_{X,x}}^{(0)} = H_{\cO_{X,y}}^{(\codim_{\overline{\{y\}}}(x))}
$$
for all $x\in U_1$, see \cite{Be} Ch. 0, p. 33, (2.1.2). On the other hand, by Lemma
\ref{lem.phi.sc} (2) there is a non-empty open subset $U_2\subseteq \overline{\{y\}}$
such that $\phi_X(y) = \phi_X(x) + \codim_{\overline{\{y\}}}(x)$ for $x\in U_2$.
Thus, for $x\in U = U_1\cap U_2$ we have
$$
H_X(x) =  H_{\cO_{X,x}}^{(\phi_X(x))} = H_{\cO_{X,y}}^{(\phi_X(x) + \codim_{\overline{\{y\}}}(x))} =
H_{\cO_{X,y}}^{(\phi_X(y))} = H_X(y)\,.
$$

\smallskip
By the following lemma, (3) is equivalent to the conjunction of (1) and (2).

\bigbreak
Let $X$ be a locally noetherian topological space which is Zariski, i.e., in which every closed irreducible subset
admits a generic point. (For example, let $X$ be a locally noetherian scheme.) Recall that a map
$H: X \longrightarrow G$ into an ordered abelian group $(G,\leq)$ is called upper semi-continuous
if the set
$$
X_{\geq\nu} := X^H_{\geq \nu} := \{ x \in X \mid H(x) \geq \nu \}
$$
is closed for all $\nu \in G$. We note that this property is compatible with
restriction to any topological subspace. In particular, if $X$ is a scheme, $H$ is
also upper semi-continuous after restricting it to a subscheme or an arbitrary localization.

\begin{lemma}\label{lem.up.sc}
(a) The map $H$ is upper semi-continuous if and only if the following holds.
\begin{itemize}
\item[(1)]
If $x,y\in X$ with $x \in \overline{\{y\}}$, then $H(x) \geq H(y)$.
\item[(2)]
For all $y\in X$ there is a dense open subset $U \subset \overline{\{y\}}$ such that
$H(x) = H(y)$ for all $x \in U$.
\end{itemize}

\smallskip (b)
Assume $H$ is upper semi-continuous. Then the set $X_{\nu} = \{x\in X \,\mid\, H(x) = \nu \}$
is locally closed for each $\nu\in G$, and its closure is contained in $X_{\geq\nu}$. In particular,
$X_\nu$ is closed if $\nu$ is a maximal element in $G$. Moreover the set
$X_{max} = \{ x\in X\,\mid\,x\in X_\nu\mbox{ for some maximal }\nu\in G\,\}$ is closed.

\smallskip
(c) If $X$ is noetherian, then $H$ takes only finitely many values.
\end{lemma}

{\bf Proof } We may restrict to the case where $X$ is noetherian by taking an open covering $(U_\alpha)$ by noetherian
subspaces. In fact, a subset of $X$ is closed (resp. locally closed) if and only if this holds for the intersection
with each $U_\alpha$. Moreover, in (1), if $x\in U_\alpha$, then $y \in U_\alpha$. In (2) we may take $U \subset U_\alpha$
where $y \in U_\alpha$. If $X$ is noetherian, then, by \cite{Be}, Ch. III, Lemma (1.1), $X_{\geq \nu}$ is closed
if and only if the following conditions hold.
\begin{itemize}
\item[(1')]
If $y \in X_{\geq \nu}$, then every $x \in \overline{\{y\}}$ is in $X_{\geq \nu}$.
\item[(2')]
If $y \notin X_{\geq \nu}$, then $(X-X_{\geq \nu})\cap\overline{\{y\}}$ is open in $\overline{\{y\}}$.
\end{itemize}
Obviously, (1) above is equivalent to (1') for all $\nu$. Assuming (1), and (2') for all $\nu$,
we get (2) by taking the following set for $\mu=H(y)$:
$$
U = X(\mu)\cap \overline{\{y\}} = \cup_{\nu>\mu} \,(X - X_{\geq \nu})\cap \overline{\{y\}}\,.
$$

Conversely assume (1) and (2), and let $\nu \in G$. By noetherian induction we show that
$X_{\geq \nu}$ is closed. Assume $\emptyset \neq Y\subset X$ is
a minimal closed subset on which this is wrong. Since $Y$ has only finitely many irreducible components,
$Y$ must be irreducible. Let $\eta$ be the generic point of $Y$, and let $H(\eta) = \mu$.
By (2) there is a dense open subset $U \subset Z$ such that $H(x) = \mu$ for all $x\in U$.
If $\nu=\mu$, then $Y_{\geq \nu} = Y$ by (1), so it is closed.
If $\nu\neq \mu$, then $Y_{\geq \nu}\subset A$, hence is closed in $A$ (minimality of $Y$), hence in $Y$ --
contradiction!

\smallskip
Similarly, we show that $H$ only takes finitely many values on (noetherian) $X$, which shows (c).
Assume that $Y\neq \emptyset$ is a minimal closed subset on which this is wrong; again $Y$ is necessarily irreducible.
Let $\eta$ be the generic point of $Y$, and let $U$ be as in (1), for $\eta$. By minimality of $Y$,
$H$ takes only finitely many values on the closed set $A = Y - U$ which is a contradiction.
The first claim in (b) follows from the equality
$$
X_\nu = X_{\geq \nu} - \cup_{\mu > \nu}\; X_\mu = X_{\geq \nu} - \cup_{\mu > \nu}\; X_{\geq \mu}\,,
$$
and the other claims follow easily (recall we may assume (c)).

\begin{definition}\label{Sigmamax}
For $X$ locally noetherian catenary, let
$\Sigma_X:=\{\HSf {X} x |\; x\in X\} \subset \bNN\,,$
and let $\SigmaXmax$ be the set of the maximal elements in
$\Sigma_X$. The set
\begin{equation}\label{HSlocul.def}
\HSmax X =\underset{\nu\in \SigmaXmax}{\cup}\HSa X
\end{equation}
is called the Hilbert-Samuel locus of $X$ (note that \eqref{HSlocul.def} is a disjoint union).
\end{definition}

By definition $\HSa X\not=\emptyset$ if and only if $\nu\in \Sigma_X$.
By Lemma \ref{lem.up.sc} we have

\begin{lemma}\label{lem.Sigmamax} Let $X$ be a locally noetherian catenary scheme.
\begin{itemize}
\item[(a)]
For each $\nu\in \bN^\bN$ the set $\HSa X$ is locally closed in $X$
and its closure in $X$ is contained in the closed set
$$
\HSgeqa X = \underset{\mu\geq\nu}{\cup} \HSb X.
$$
In particular $\HSa X$ is closed in $X$ if $\nu\in \SigmaXmax$, and $\HSmax X$
is closed.
\item[(b)]
If $X$ is noetherian, then $\Sigma_X$ is finite.
\end{itemize}
\end{lemma}

In the following we will regard the sets $\HSa X$, $\HSgeqa X$ and $\HSmax X$ as (locally closed) subschemes
of $X$, endowed with the reduced subscheme structure.

\begin{lemma}\label{completion.HS}
Let $X$ be a locally noetherian scheme.
\begin{itemize}
\item[(1)]
Let $\pi: X' \rightarrow X$ be a morphism with $X'$ locally noetherian and let $x \in X$ and $x' \in X'$
with $\pi(x') = x$. Assume that $\pi$ is flat and that the fibre $X'_x=X'\times_Xx$ of $\pi$
above $x$ is regular at $x'$ (e.g., assume that $\pi$ is smooth around $x'$, or that $\pi$ is
quasi-\'etale at $x'$ in the sense of Lemma \ref{completion.nu.2} (1)). Then one has
\begin{equation}\label{eq.HS.smooth}
H_{\cO_{X',x'}}^{(0)} = H_{\cO_{X,x}}^{(d)} \quad\mbox{   and   }\quad
\psi_{X'}(x') = \psi_{X}(x) + d\,,
\end{equation}
where $d = \dim(\cO_{X'_x,x'}) = \codim_{X'_x}(x')$. Hence, if $X$ and $X'$ are catenary,
then $H_{X'}(x') = H_X(x)$. In particular, if $\pi$ is
flat with regular fibers (e.g., assume that $\pi$ is \'etale or smooth), then
\begin{equation}\label{eq.HS.smooth2}
\pi^{-1}(\HSa X) = \HS {X'} {\nu} \cong \HSa X\times_XX'\quad\quad\mbox{for
all $\nu \in \bN^\bN$}\,.
\end{equation}
Moreover, $X$ is regular at $x$ if and only if $X'$ is regular at $x'$.

\item[(2)]
Let $X_0$ be an excellent scheme, and let $f:X \to X_0$ be a morphism which is locally of finite type.
Let $x_0\in X_0$ and
$$
\hat{X_0} = \Spec(\hat{\cO}_{X_0,x_0}), \quad
\pi: \hat{X} = X \times_{X_0} \hat{X_0} \to X\,,
$$
where $\hat{\cO}_{X_0,x_0}$ is the completion of ${\cO}_{X_0,x_0}$. Then
for any $x \in X$ and any $\hat{x} \in \hat{X}$ lying over $x$ we have
\begin{equation}\label{eq.HS.completion}
H_{\cO_{\hat{X},\hat{x}}}^{(0)} = H_{\cO_{X,x}}^{(d)} \quad\mbox{   and   }\quad
\psi_{\hat{X}}(\hat{x}) = \psi_{X}(x) + d\,,
\end{equation}
where
$$
d=\codim_{\hat{X}_x}(\hat{x}), \quad \hat{X}_x = \hat{X}\times_{X}y.
$$
In particular, $H_{\hat{X}}(\hat{x}) = H_{X}(x)$, and
\begin{equation}\label{eq.HS.completion2}
\pi^{-1}(\HSa {X}) = \HS {\hat{X}} {\nu} \cong \HSa {X} \times_X\hat{X}
\quad\quad \mbox{for all $\nu \in \bNN$}\,.
\end{equation}
\end{itemize}
\end{lemma}

For the proof we shall use the following two lemmas.

\begin{lemma}\label{lem.flat}
Let $f:W\to Z$ be a flat morphism of locally noetherian schemes and let $w\in W$ and
$z=f(w)\in Z$. Then we have
$$
\codim_W(w)=\codim_Z(z)+\dim\cO_{W_z,w},
\quad W_z=W\times_Z z.
$$
In particular $w$ is a generic point of $W$ if and only if $f(w)$
is a generic point of $Z$ and $\codim_{W_z}(w)=0$.
\end{lemma}
\medbreak\textbf{Proof}
See \cite{EGAIV}, 2, (6.1.2).

\begin{lemma}\label{lem.qf}
$f:W\to Z$ be a quasi-finite morphism of excellent schemes and
let $w\in W$ and $z=f(w)\in Z$.
Assume that $W$ and $Z$ are irreducible. Then we have
$$
\codim_W(w)=\codim_Z(z).
$$
\end{lemma}
\medbreak\textbf{Proof}
See \cite{EGAIV}, 2, (5.6.5).

\bigbreak\noindent \textbf{Proof of Lemma \ref{completion.HS}}
(1): The first equality in \eqref{eq.HS.smooth} follows from
Lemma \ref{completion.nu.2} \eqref{eq.dir.completion.2}.
We show the second equality. We may assume that $X$ is reduced.
Since $\pi:X'\to X$ is flat, Lemma \ref{lem.flat} implies that
if $\eta'$ is a generic point of $X'$ with $x'\in \overline{\{\eta'\}}$,
then $\eta=\pi(\eta')$ is a generic point of $X$ with $x\in \overline{\{\eta\}}$.
Moreover, if $\xi$ is a generic point of $X$ with
$x\in Z:=\overline{\{\xi\}}$, then there exists a generic point
$\xi'$ of $X'$ such that $x'\in \overline{\{\xi'\}}$.
Indeed one can take a generic point $\xi'$ of $\pi^{-1}(Z)$ such that $x'\in \overline{\{\xi'\}}$.
Then Lemma \ref{lem.flat} applied to the flat morphism $\pi_Z:X'\times_X Z\to Z$
implies that $\pi(\xi')=\xi$ and $\xi'$ is of codimension
$0$ in $\pi_Z^{-1}(\xi)=\pi^{-1}(\xi)$ and hence $\xi'$ is a generic point of $X'$.
This shows that we may consider the case where $X$ is integral, and it suffices show the following.

\begin{claim}\label{claim.eq.HS.completion}
Assume $X$ is integral. Let $W$ be an irreducible component of $X'$ containing $x'$. Then
$$
\codim_{W}(x')=\codim_{X}(x) + d,
\quad d=\codim_{X'_x}(x').
$$
\end{claim}

\medskip
Since the question is local at $x$, we may assume $X=\Spec(A)$ for the local noetherian ring
$A = \cO_{X,x}$.
If $A$ is normal, then by \cite{EGAIV}, 2, (6.5.4) the local ring $B=\cO_{X',x'}$ is
normal as well, because the fibers of $\pi$ are regular and hence normal. Therefore
both $A$ and $B$ are integral, we have $\psi_X(x) = \dim(\cO_{X,x})$ and $\psi_{X'}(x') = \dim(\cO_{X',x'})$,
and the claim follows from Lemma \ref{lem.flat}.

In general let $\tilde{A}$ be the normalization of $A$, let $\tilde{X} = \Spec(\tilde{A})$
(the normalization of $X$) and consider the cartesian diagram
$$
\begin{CD}
\tilde{X}' @>{g'}>> X' \\
@V{\tilde{\pi}}VV @V{\pi}VV  \\
\tilde{X} @>{g}>> X
\end{CD}
$$
in which the vertical morphisms are flat.
Since $X$ is excellent, the morphism $g$ is finite and surjective, and so is its base change $g'$.
We claim that a point $\xi \in \tilde{X}'$ is generic if and only if
$\eta'= g'(\xi)$ is a generic point of $X'$. Indeed, by Lemma
\ref{lem.flat}, $\xi$ is generic if and only if it maps to the generic point $\tilde{\eta}$ of $\tilde{X}$
and is of codimension zero in the fibre $\tilde{\pi}^{-1}(\tilde{\eta})$.
Since the fibres over $\tilde{\eta}$ and the generic point $\eta$ of $X$ are the same,
this is the case if and only if $\eta'$ satisfies the corresponding properties for $\pi$,
i.e., if $\eta'$ is a generic point of $X'$.
Now let $W$ be an irreducible component of $X'$ containing $x'$.
By the last claim together with the surjectivity of $g'$,
there is an irreducible component $\tilde{W}$ of $\tilde{X}'$ dominating $W$.
Since $g'$ is finite, it is a closed map, and so we have even $g'(\tilde{W}) = W$.
Thus there is a point $z'\in \tilde{W}$ with $g'(z') = x'$.
If $z = \tilde{\pi}(z')$, then $g(z) = y$ and
\begin{equation}\label{eq.10}
\codim_{\tilde{X}'_z}(z') = \codim_{X'_x}(x') = d\,,
\end{equation}
where
$\tilde{X}'_z = X'_x\times_x z$.
We now conclude
$$
\codim_W(x') = \codim_{W'}(z') = \codim_{\tilde{X}}(z) + d = \codim_{X}(x) + d \,,
$$
where the first (resp. the third) equality follows from Lemma \ref{lem.qf}
applied to the finite morphism $\tilde{W}\to W$ (resp. $\tilde{X}\to X$) and
the second equality follows from \eqref{eq.10} and the first case of the proof
noting that $\tilde{X}$ is normal. This shows Claim \ref{claim.eq.HS.completion}
and completes the proof of \eqref{eq.HS.smooth}. The next claims,
including the first equality in \eqref{eq.HS.smooth2}, are now obvious. To see the
isomorphism of schemes in \eqref{eq.HS.smooth2}, note the following. For any morphism
$f: S' \rightarrow S$ of schemes and any subscheme $T \subset S$, $T\times_SS'$ is
identified with a subscheme of $S'$, whose underlying topological space is homeomorphic
with $f^{-1}(T)$. So we have to show that $\HSa X\times_XX'$ is reduced. But this follows
from \cite{EGAIV}, 2, (3.3.5) and the flatness of $\pi$. The last statement in (1)
now follows from Remark \ref{HSregular}. (See also \cite{EGAIV}, 2, (6.5.2) for another proof.)

Finally, (2) follows from applying (1) to the morphism $\hat{X}' \rightarrow X'\times_X\Spec(\cO_{X,x})$
(compare the proof of Lemma \ref{completion.nu.2}, deduction of (3) from (2)).

\begin{remark}\label{rem.comp.Bennett}
Bennett \cite{Be} defined global Hilbert-Samuel functions by $H_{X,x}^{(i)} = H_{\cO_{X,x}}^{(i+d(x))}$,

\vspace{-1mm}
where $d(x)=\dim(\overline{\{x\}})$, and showed that these functions have good properties for
so-called weakly biequidimensional excellent schemes. By looking at the generic points and the closed
points one easily sees that this function coincides with our function $H_X(x)$ (for $N = \dim X$)
if and only if $X$ is biequidimensional.
\end{remark}

\newpage
\section{Permissible blow-ups}
\bigskip

We discuss some fundamental results on the behavior of the $\nu^*$-invariants,
the $e$-invariants, and the Hilbert-Samuel functions under permissible blow-ups.

\medbreak

Let $X$ be a locally noetherian scheme and let $D\subset X$ be a closed reduced
subscheme. Let $I_D\subset \cO_X$ the ideal sheaf of $D$ in $X$
and $\cO_D=\cO_X/I_D$. Put
$$
\gr_{I_D}(\cO_X) =\underset{t\geq 0}{\bigoplus} I_D^t/I_D^{t+1}.
$$

\begin{definition}\label{def.nf}
\begin{itemize}
\item[(1)]
$X$ is normally flat along $D$ at $x\in D$ if the stalk
$\gr_{I_D}(\cO_X)_x$ of $\gr_{I_D}(\cO_X)$ at $x$ is a flat $\cO_{D,x}$-module.
$X$ is normally flat along $D$ if $X$ is normally flat along $D$
at all points of $D$, i.e., if $\gr_{I_D}(\cO_X)$ is a flat $\cO_D$-module.
\item[(2)]
$D\subset X$ is permissible at $x\in D$ if $D$ is regular at $x$, if $X$
is normally flat along $D$ at $x$, and if $D$ contains no irreducible component
of $X$ containing $x$.
$D\subset X$ is permissible if $D$ is permissible at all points of $D$.
\item[(3)]
The blowup $\pi_D:\Bl_D(X) \to X$ in a permissible center $D\subset X$,
is called a permissible blowup.
\end{itemize}
\end{definition}
\medbreak

For a closed subscheme $D\subset X$, the normal cone of $D\subset X$ is
defined as:
$$C_D(X)=\Spec(\gr_{I_D}(\cO_X)) \to D .$$

\begin{theorem}\label{nf.thm1}
\begin{itemize}
\item[(1)]
There is a dense open subset $U\subset D$ such that $X$ is normally flat
along $D$ at all $x\in U$.
\item[(2)] The following conditions are equivalent:
\begin{itemize}
\item[$(i)$]
$X$ is normally flat along $D$ at $x\in D$.
\item[$(ii)$]
$T_x(D)\subset \Dir_x(X)$ and the natural map $C_x(X) \to C_D(X)\times_D x$
induces an isomorphism
$C_x(X)/T_x(D)\isom C_D(X)\times_D x$,
where $T_x(D)$ acts on $C_x(X)$ by the addition in $T_x(X)$.
\end{itemize}
Assume that in addition that $X$ is a closed subscheme of
a regular locally noetherian scheme $Z$. Let $x\in D$ and set $R=\cO_{Z,x}$ with the
maximal ideal $\fm$ and $k=R/\fm$. Let $J\subset R$ (resp. $\fp\subset R$) be
the ideal defining $X\subset Z$ (resp. $D\subset Z$). Then the following conditions
are equivalent to the conditions (i) and (ii) above.
\begin{itemize}
\item[$(iii)$]
Let $u:\grpR\otimes_{R/\fp} k \to \grmR$ be the natural map.
Then $\InmJ$ is generated in $\grmR$ by $u(\InpJ)$.
\item[$(iv)$]
There exists a standard base
$f=(f_1,\dots,f_m)$ of $J$ such that $v_\fm(f_i)=v_\fp(f_i)$ for all
$i=1,\dots,m$ (cf. \eqref{order}) .
\end{itemize}
\item[(3)]
Let $\pi: X' \rightarrow X$ be a morphism, with $X'$ locally noetherian, and let
$D' = X'\times_XD$, regarded as a closed subscheme of $X'$.
Let $x' \in D'$ and $x=\pi(x') \in D$, and assume that $\pi$ is flat at $x'$ and
that the fiber $X'_x$ over $x$ is regular at $x'$. Then $X'$ is
normally flat along $D'$ at $x'$ if and only if $X$ is normally flat along $D$ at $x$.
Moreover, $D'$ is regular (resp. permissible) at $x'$ if and only $D$ is regular (resp.
permissible) at $x$.
\end{itemize}
\end{theorem}
\medbreak
\textbf{Proof } (1) follows from \cite{H1} Ch. I Theorem 1 on page 188
and (2) from \cite{Gi1}, II \S2, Theorem 2.2 and 2.2.3 on page II-13 to II-15.
(3): Since $\cO_{X,x} \rightarrow \cO_{X',x'}$ is flat, one has
$$
\gr_{D'}(\cO_{X'})_{x'} \cong \gr_D(\cO_X)_x\otimes_{\cO_{X,x}}\cO_{X',x'}\cong \gr_D(\cO_X)_x\otimes_{\cO_{D,x}}\cO_{D',x'}\,.
$$
Thus the first claim follows from the following general fact: If
$A \rightarrow B$ is a flat morphism of local rings (hence faithfully flat), and $M$ is
an $A$-module, then $M$ is flat over $A$ if and only if $M\otimes_AB$ is flat over $B$.
Since $D' \rightarrow D$ is flat and its fiber over $x$ is the same as for $\pi$, hence regular,
the next claim (on regularity) follows from \ref{completion.HS} (1) (or \cite{EGAIV}, 2, (6.5.2)).
The last claim (on nowhere density) follows from the flatness of $\pi$ by the arguments used at the beginning of the proof of Lemma
\ref{completion.HS}: If $\eta'$ is a generic point of $X'$ with $x' \in \overline{\{\eta'\}} \subset D' = \pi^{-1}(D)$,
then $\eta = \pi(\eta')$ is a generic point of $X$ with $x \in \overline{\{\eta\}} \subset D$.
Conversely, if start with $\eta$ as in the latter situation, there is a generic point $\eta'$ of $X'$
mapping to $\eta$ with $x' \in \overline{\{\eta'\}}$. But then $\eta' \in D'$.
\bigskip

There is a numerical criterion for normal flatness
due to Bennett,
which we carry over to our setting:
Let the assumption be as in the beginning of this section.
Assume in addition that $X$ is catenary (e.g., let $X$ be excellent)
and let $\HSfX x$ be as in Definition \ref{def.HS.schemes}.

\begin{theorem}\label{HSf.perm}
Assume that $D$ is regular. Let $x\in D$, and let $y$ be
the generic point of the component of $D$ which contains $x$.
Then the following conditions are equivalent.
\begin{itemize}
\item[(1)]
$X$ is normally flat along $D$ at $x$.
\item[(2)]
$H_{\cO_{X,x}}^{(0)} = H_{\cO_{X,y}}^{(\codim_{Y}(x))}$,
where $Y = \overline{\{y\}}$, the closure of $y$ in $X$.
\item[(3)]
$\HSfX x=\HSfX y$.
\end{itemize}
\end{theorem}
\textbf{\bf Proof }
The equivalence of (1) and (2) was proved by Bennett (\cite{Be}, Theorem (3)).
The rest is a special case of the following lemma.

\begin{lemma}\label{lem.HSf-spec}
Let $x, y \in X$ with $x \in \overline{\{y\}}$. Then the following are equivalent.
\begin{itemize}
\item[(1)]
$H_{\cO_{X,x}}^{(0)} = H_{\cO_{X,y}}^{(\codim_{\overline{\{y\}}}(x))}$.
\item[(2)]
$\HSfX x=\HSfX y$.
\end{itemize}
If these conditions hold, then $I(x) = I(y)$ and $\psi_X(x) = \psi_X(y) + \codim_Y(x)$.
\end{lemma}

\textbf{\bf Proof } By the definition of the considered functions, for the equivalence of
(1) and (2) it suffices to show that either of (1) and (2) implies
\begin{equation}\label{eq.1.7}
\phi_X(y) = \phi_X(x) + \codim_Y(x)\,,
\end{equation}
where $Y = \overline{\{y\}}$. Assume (2). By Lemma \ref{lem.HS.comparison} we have
$$
\dim \cO_{X,x} + \phi_X(x) = \dim \cO_{X,y} + \phi_X(y).
$$
On the other hand,
$$
\dim \cO_{X,x} = \codim_X(x) \geq \codim_Y(x) + \codim_X(y) = \codim_Y(x) +
\dim \cO_{X,y}\,.
$$
Thus we get
\begin{equation*}
\phi_X(y) = \dim \cO_{X,x} +\phi_X(x) -\dim \cO_{X,y}\geq
\phi_X(x) + \codim_Y(x)\,,
\end{equation*}
which implies \eqref{eq.1.7} by Lemma \ref{lem.phi.sc}.

\bigskip
Next assume (1). Let $\cO = \cO_{X,x}$ and
let $\fp \subset \cO$ be the prime ideal corresponding to $y$. It suffices to show that
$\fp$ contains all minimal prime ideals of $\cO$. In fact, this means that $y$
is contained in all irreducible components of $X$ which contain $x$, i.e., that $I(y) = I(x)$.
Since, for any irreducible $Z \in I(y)$ we have
$$
\codim_Z(x) = \codim_Z(y) + \codim_Y(x)\,,
$$
we deduce the equality $\psi_X(x) = \psi_X(y) + \codim_Y(x)$ and hence
\eqref{eq.1.7}. At the same time we have proved the last claims of the lemma. As for the
claim on $\cO$ and $\fp$ let
$$
(0) = \fP_1\cap \ldots \cap \fP_r
$$
be a reduced primary decomposition of the zero ideal of $\cO$, and let
$\fp_i = \mbox{Rad}(\fP_i)$ be the prime ideal
associated to $\fP_i$. Then the set $\{\fp_1,\ldots,\fp_r\}$ contains all
minimal prime ideals of $\cO$, and it suffices to show that $\fp$ contains
all ideals $\fP_i$ (Note that an ideal $\fa$ is contained in $\fp$ if and
only if $\mbox{Rad}(\fa)$ is contained in $\fp$).
Assume the contrary. We may assume $\fP_1 \nsubseteq \fp$.
Put $\cO'=\cO/\fQ$ with $\fQ=\fP_2\cap \cdots \cap \fP_r$ and let
$\cO_\fp$ (resp. $\cO'_\fp$) be the localization of $\cO$ (resp. $\cO'$) at
$\fp$ (resp. $\fp\cO'$). Then
$\cO_\fp=\cO'_\fp$ and $\cO/\fp=\cO'/\fp\cO'$ and
$$
H_{\cO}^{(0)} \geq H_{\cO'}^{(0)}\geq H_{\cO'_\fp}^{(d)}=
H_{\cO_\fp}^{(d)}\; ,
$$
where $d=\codim_{Y}(x)=\dim \cO/\fp=\dim \cO'/\fp\cO'$.
The first inequality follows from Lemma \ref{HSquotient} and the second
inequality follows from \cite{Be}, Theorem (2). Hence (2) implies
$H_{\cO}^{(1)} = H_{\cO'}^{(1)}$. By Lemma \ref{HSquotient} this implies
$\cO=\cO'$ so that $\fQ=0$, contradicting the assumption that the primary
decomposition is reduced.
$\square$
\medbreak

The above criterion is complemented by the following observation.

\begin{lemma}\label{lem.constant.HSf}
Let $X$ be connected. If there is an irreducible component $Z\subseteq X$
such that $\HSfX x=\HSfX y$ for all $x,y\in Z$ (i.e., $Z \subseteq X(\nu)$ for some $\nu \in \bNN$),
then $Z=X$.
\end{lemma}
\par\noindent
\textbf{Proof }
We have to show that $X$ is irreducible. Assume not. Then there exists an $x\in X$ which is contained
in two different irreducible components. Let $\cO = \cO_{X,x}$ be the local ring of $X$ at
$x$, let
\begin{equation}\label{eq.constant.HSf1}
<0> \; = \; \fP_1\cap \ldots \cap \fP_r
\end{equation}
be a reduced primary decomposition of the zero ideal of $\cO$, and let
$\fp_i = \mbox{Rad}(\fP_i)$ be the prime ideal
associated to $\fP_i$. By assumption we have $r\neq 1$.
Then there is an $i$ such that the trace of $Z$ in $\cO$ is given by $\fp_i$.
Let $\eta_i \in \Spec(\cO/\fP_i) \subset \Spec(\cO) \subset X$ be the
generic point of $Z$ (corresponding to $\fp_i$). By assumption we have
$$
H_{\cO_{X,\eta_i}}^{(N)} = \HSfX {\eta_i} =  \HSfX x = H_{\cO}^{(N-\psi_X(x))}
$$
where $N$ is as in Definition \ref{def.HS.schemes}, because $\psi_X(\eta_i)=0$. On the other hand, we have
\begin{equation}\label{eq.constant.HSf2}
H_{(\cO/{\fP_i})_{\fp_i}}^{(N)}
\leq H_{\cO/{\fP_i}}^{(N-c)} \leq H_{\cO}^{(N-c)} \leq H_{\cO}^{(N-\psi_X(x))}\,,
\end{equation}
where $c := \dim(\cO/{\fP_i}) = \dim(\cO/{\fp_i}) = \codim_{\overline{\{\eta_i\}}}(x)=\codim_Z(x)$.
Here the first inequality holds by the results of Bennett/Singh recalled in the proof
of Theorem \ref{HSf.usc} (1), the second inequality follows from Lemma \ref{HSquotient},
and the last inequality holds since $\psi_X(x) \leq c $.
Now, since $\fp_i$ is a minimal prime ideal and \eqref{eq.constant.HSf1} is reduced,
we have an isomorphism
$$
\cO_{X,\eta_i} = \cO_{\fp_i} \cong (\cO/{\fP_i})_{\fp_i}\,.
$$
Therefore we have equalities in \eqref{eq.constant.HSf2}. By the other direction of
Lemma \ref{HSquotient} we conclude that $\cO = \cO/{\fP_i}$, i.e., $\fP_i=0$,
which is a contradiction if $r\neq 1$.$\square$

\bigskip
We now prove a semi-continuity property for $e_x(X)$, the dimension of the directrix
at $x$.

\begin{theorem}\label{nf.dir}
Let $X$ be an excellent scheme, or a scheme embeddable into a regular scheme,
and let $x,y\in X$ such that $D=\overline{\{y\}}$ is permissible and $x\in D$. Then
$$
e_y(X) \leq e_x(X) - \dim(\cO_{D,y})\,.
$$
\end{theorem}

The question only depends on the local ring $\cO_{X,x}$, and by Lemma \ref{completion.nu.2}
(3) we may assume that we consider the spectrum of a complete local ring
and the closed point $x$ in it. Moreover, by the Cohen structure theory of complete local noetherian
rings (see \cite{EGAIV}, 1, (0.19.8.8)) every such ring is a quotient of a (complete) regular local ring.
Therefore Theorem \ref{nf.dir} is implied by the following result.

\begin{theorem}\label{nf.dir2}
Let $R$ be a regular local ring with maximal ideal $\fm$.
Let $J\subset R$ be an ideal. Let $\fp\subset R$ be a prime ideal such that
$J\subset \fp$ and that $\Spec(R/\fp)\subset \Spec(R/J)$ is permissible.
Let $R_\fp$ be the localization of $R$ at $\fp$ and $J_\fp=J R_\fp$. Then
$$
e(R_\fp/J_\fp)\leq e(R/J) -\dim(R/\fp).
$$
\end{theorem}
\par\noindent
\textbf{Proof }
Set $A=R/\fp$ and let $K$ be the quotient field of $A$.
Set $M_K=M\otimes_A K$ for an $A$-module $M$.
By definition there exists a $K$-subspace
$V=\IDir(R_\fp/J_\fp)\subset \grp^1(R)_K$ of dimension $s$ such that
$\dim_K(V)=\dim(R_\fp)-e(R_\fp/J_\fp)$ and
\begin{equation}\label{nf.dir.eq1}
\InpJ_K =(K[V]\cap\InpJ_K)\cdot\grpR_K.
\end{equation}

\begin{lemma}\label{nf.dir.lem1}
Assume that there exist free $A$-modules $T,S\subset \grp^1(R)$ such that
$\grp^1(R)=T\oplus S$ and $V=T_K$. Let
$u:\grpR\otimes_{R/\fp} k \to \grmR$ be the natural map.
Then $u(T)\supset \IDir(R/J)\cap \grm^1(R)$.
\end{lemma}
\medbreak

Note that the assumption of the lemma is satisfied if $\dim(A)=1$,
by the theory of elementary divisors.
Theorem \ref{nf.dir} is a consequence of the conclusion of
Lemma \ref{nf.dir.lem1} by noting
$$\dim_k(u(T))=\dim_k(T\otimes_A k)=\dim_K(V)=\dim(R_\fp)-e(R_\fp/J_\fp)$$
so that the lemma finishes the proof of the theorem in case $\dim(R/\fp)=1$.
\medbreak

We show Lemma \ref{nf.dir.lem1}. The assumption of the lemma implies
\begin{equation}\label{nf.dir.eq2}
\grpR =\Sym_A(\grp^1(R))=A[T]\otimes_A A[S],
\end{equation}
where $A[T]=\Sym_A(T)$ (resp. $A[S]=\Sym_A(S)$) is the sub $A$-algebra
of $\grpR$ generated by $T$ (resp. $S$).

\begin{claim}\label{nf.dir.claim1}
$$
\InpJ =(A[T]\cap\InpJ)\cdot\grpR.
$$
\end{claim}
\medbreak
By Theorem \ref{nf.thm1}(2)$(iii)$ the claim implies that
$\InmJ$ is generated by $u(A[T])\cap \InmJ$, which implies
Lemma \ref{nf.dir.lem1}. Thus it suffices to show the claim.
Note
\begin{equation}\label{nf.dir.eq3}
(K[V]\cap \InpJ_K)\cap \grpR=
(A[T]\cap \InpJ)_K \cap \grpR= A[T]\cap \InpJ.
\end{equation}
Indeed \eqref{nf.dir.eq2} implies
$K[V]\cap \grpR=A[T]$ and the flatness of $\grpRJ=\grpR/\InpJ$ implies
$\InpJ_K\cap \grpR=\InpJ$.
\medbreak

Take any $\phi\in \InpJ$. By \eqref{nf.dir.eq1} and \eqref{nf.dir.eq3}
there exists $c\in A$ such that
$$
c \phi=\underset{1\leq i\leq m}{\sum} a_i\psi_i,\quad
\psi_i\in A[T]\cap \InpJ,\; a_i\in \grpR,
$$
Choosing a basis $Z_1,\dots,Z_r$ of the $A$-module $S$,
\eqref{nf.dir.eq2} allows us to identify $\grpR$ with the polynomial ring
$A[T][Z]=A[T][Z_1,\dots,Z_r]$ over $A[T]$.
Then, expanding
$$
\phi=\underset{B\in \bZ^r_{\geq 0}}{\sum} Z^B \phi_B,\quad
a_i=\underset{B\in \bZ^r_{\geq 0}}{\sum} Z^B a_{i,B},
\quad\text{ with  } \phi_B\in A[T],\; a_{i,B}\in A[T],
$$
we get
$$
c \phi_B=\underset{1\leq i\leq m}{\sum} a_{i,B} \psi_i
\quad\text{for any } B\in \bZ^r_{\geq 0}.
$$
Since $\underset{1\leq i\leq m}{\sum} a_{i,B} \psi_i\in A[T]\cap\InpJ$,
this implies $\phi_B\in A[T]\cap\InpJ$ by \eqref{nf.dir.eq3} so that
$\phi\in (A[T]\cap\InpJ)\cdot\grpR$.
This completes the proof of the claim.
\medbreak

To finish the proof of Theorem \ref{nf.dir}, it suffices to reduce it to
the case $\dim(R/\fp)=1$ as remarked below Lemma \ref{nf.dir.lem1}.
Assume $\dim(R/\fp)> 1$ and take a prime ideal $\fq\supset \fp$ such that
$R/\fq$ is regular of dimension $1$.
Let $R_\fq$ be the localization of $R$ at $\fq$ and $J_\fq=J R_\fq$.
Noting
$$
\gr_{\fp R_\fq}(R_\fq/J_\fq)\simeq \grpRJ\otimes_{R/\fp} R_\fq/\fp R_\fq,
$$
the assumption implies that
$\Spec(R_\fq/\fp R_\fq)\subset \Spec(R_\fq/J R_\fq)$ is permissible.
By the induction on $\dim(R)$, we have
$$
e(R_\fp/J_\fp)\leq e(R_\fq/J_\fq) -\dim(R_\fq/\fp R_\fq)
=e(R_\fq/J_\fq) -(\dim(R/\fp)-1)
$$
(Note that any regular local ring is catenary \cite{EGAIV}, 1, (0.16.5.12)). Hence we are reduced to show
$$
e(R_\fq/J_\fq)\leq e(R/J) -\dim(R/\fq)=e(R_\fq/J_\fq)\leq e(R/J) -1.
$$
This completes the proof of Theorem \ref{nf.dir}.

\bigskip

Bennett and Hironaka proved results about the behavior of Hilbert-Samuel functions
in permissible blow-ups which are fundamental in resolution of singularities.
We recall these results (as well as some improvements by Singh) and carry them over to our setting.

\begin{theorem}\label{thm.pbu.inv}
Let $X$ be an excellent scheme, or a scheme which is embeddable in a regular scheme.
Let $D \subset X$ be a permissible closed subscheme, and let
$$
\pi_X : X'=B\ell_D(X) \to X\quad \text{and }\quad
$$
be the blowup with center $D$.
Take any points $x\in D$ and $x'\in \pi_X^{-1}(x)$ and let
$\delta=\delta_{x'/x}:= \trdeg_{k(x)}(k(x'))$. Then:
\begin{itemize}
\item[(1)]
$\; H_{\cO_{X',x'}}^{(\delta)} \leq H_{\cO_{X,x}}^{(0)}$ and $\;\phi_{X'}(x') \leq \phi_X(x) + \delta$
and $\;\HSfXd  {x'}\leq \HSfX {x}$.
\item[(2)]
$\; \HSfXd  {x'} = \HSfX {x} \quad \Leftrightarrow \quad
H_{\cO_{X',x'}}^{(\delta)} = H_{\cO_{X,x}}^{(0)}$.
\item[(3)]
If the equalities in (2) hold, then $\phi_{X'}(x') = \phi_X(x) + \delta$,
the morphism $\cO_{X,x} \rightarrow \cO_{X',x'}$ is injective, and $I(x')=\{Z'\mid Z\in I(x)\}$
where $Z'$ denotes the strict transform of $Z\in I(x)$.
\item[(4)]
If the equalities in (2) hold, then, for any field extension $K/k(x')$ one has
$$
e_{x'}(X')_K \leq e_x(X)_K-\delta_{x'/x}.
$$
\end{itemize}
Assume in addition that $X \hookrightarrow Z$ is a closed immersion into a regular
scheme $Z$, and let
$$
\pi_Z : Z'=B\ell_D(Z) \to Z
$$
be the blowup with center $D$. Then:
\begin{itemize}
\item[(5)]
$\; \nu^*_{x'}(X',Z') \leq \nu^*_x(X,Z)$ in the lexicographic order.
\item[(6)]
$\; H_{\cO_{X',x'}}^{(\delta)} = H_{\cO_{X,x}}^{(0)}
\quad \Leftrightarrow \quad \nu^*_{x'}(X',Z') = \nu^*_x(X,Z)$.
\end{itemize}
\end{theorem}

 \textbf{Proof}
In a slightly weaker form, viz., $H_{\cO_{X',x'}}^{(1+\delta)} \leq H_{\cO_{X,x}}^{(1)}$,
the first inequality in (1) was proved by Bennett (\cite{Be} Theorem (2)), and
Hironaka gave a simplified proof (\cite{H4} Theorem I). In the stronger form above
it was proved by Singh (\cite{Si1}, Remark after Theorem 1).
For the second inequality, since $\dim(X)=\dim(X')$, it suffices to show
\begin{equation}\label{eq.HS.blowup5.5}
\psi_{X}(x) \leq \psi_{X'}(x') + \delta.
\end{equation}
Note that $X$ is universally catenary by assumption (because excellent schemes and closed subschemes
of regular schemes are universally catenary (by definition and by \cite{EGAIV}, 2, (5.6.4), respectively).
Let $Y_1,\dots,Y_r$ be the irreducible components of $X$ and let $Y'_i$
be the strict transform of $Y_i$ under $\pi_X$. Then $Y_1',\dots,Y'_r$ are
the irreducible components of $X'$. If $x' \in Y'_i$, then $x\in Y_i$ and
\cite{EGAIV}, 2, (5.6.1) implies (note that $\cO_{Y_i,x}$ is universally catenary since $X$ is)
\begin{equation}\label{eq.HS.blowup6}
\codim_{Y_i}(x)=\codim_{Y'_i}(x') + \delta.
\end{equation}
\eqref{eq.HS.blowup5.5} follows immediately from this.
\medbreak

The last inequality in (1) now follows:
\begin{equation}\label{eq.HS.blowup7}
H_{X'}(x') = H_{\cO_{X',x'}}^{(\phi_{X'}(x'))} \leq H_{\cO_{X',x'}}^{(\phi_X(x)+\delta)} \leq H_{\cO_{X,x}}^{(\phi_X(x))} = H_X(x).
\end{equation}
Claims (5) and (6) were proved by Hironaka in \cite{H4} Theorems II and III, hence
it remains to show (2), (3) and (4).

As for (2), assume $\HSfXd  {x'} = \HSfX {x}$, i.e., that equality holds everywhere in \eqref{eq.HS.blowup7}.
To show $H_{\cO_{X',x'}}^{(\delta)} = H_{\cO_{X,x}}^{(0)}$, it suffices to show that
\begin{equation}\label{eq.HS.blowup8}
\phi_{X'}(x') = \phi_X(x) + \delta\,.
\end{equation}
Let $d= \dim(\cO_{X,x})$ and $d' = \dim(\cO_{X',x'})$.
By Lemma \ref{lem.HS.comparison}(b) the assumption implies
\begin{equation}\label{eq.HS.near2}
d' + \phi_{X'}(x') = d + \phi_X(x)\,.
\end{equation}
On the other hand, by Lemma \ref{lem.HS.comparison}(b), the inequality
$H_{\cO_{X',x'}}^{(\delta)} \leq H_{\cO_{X,x}}^{(0)}$ from (1)
implies
\begin{equation}\label{eq.HS.blowup4}
d'+ \delta \leq d\,.
\end{equation}
From \eqref{eq.HS.near2} and \eqref{eq.HS.blowup4} we deduce
$\phi_X(x) + \delta \leq \phi_{X'}(x')\,,$
which implies \eqref{eq.HS.blowup8}, in view of (1).
\medbreak

Conversely assume $H_{\cO_{X',x'}}^{(\delta)} = H_{\cO_{X,x}}^{(0)}$. To show $\HSfXd  {x'} = \HSfX {x}$,
it again suffices to show \eqref{eq.HS.blowup8}.
We have to show
\begin{equation}\label{eq.HS.near3}
\psi_X(x) = \psi_{X'}(x') + \delta\,.
\end{equation}

In view of \eqref{eq.HS.blowup6} and with the notations there, it suffices to show
that $x' \in Y_i'$ if $x \in Y_i$, i.e., the third claim of (3). For this, by the lemma below, it suffices to show
the injectivity of $\cO = \cO_{X,x} \rightarrow \cO_{X',x'}=\cO'$, i.e., the second claim of (3).

\begin{lemma}
Let $f: A \to B$ be a homomorphism of noetherian rings.
Assume that there is a minimal prime ideal of $A$, which does not lie
in the image of $\Spec(B)\to \Spec(A)$. Then $f$ is not injective.
\end{lemma}

{\bf Proof} Let $(0) = \fP_1\cap \ldots \cap \fP_r$ be a primary decomposition of
the zero ideal in $B$, and let $\fp_i = \mbox{Rad}(\fP_i)$, which is a prime ideal.
Assume $A \rightarrow B$ is injective. Then, with $\fQ_i = \fP_i\cap A$, $(0) = \fQ_1,\cap\ldots\cap\fQ_r$
is a primary decomposition in $A$. Therefore $\{\fq_,\ldots,\fq_r\}$, with $\fq_i=\mbox{Rad}(\fQ_i) = \fp_i\cap A$,
is the set of all associated prime ideals of $A$ and hence contains all minimal prime ideals of $A$,
see \cite{Ku} VI Theorems 2.18 and 2.9. This contradicts the assumption.

\medbreak
Before we go on, we note the following consequence of Lemmas \ref{completion.nu.2} and \ref{completion.HS}:
For all claims of Theorem \ref{thm.pbu.inv}, we may, via base change, assume that $X=\Spec(\cO_{X,x})$,
and that $\cO_{X,x}$ is a quotient of a regular local ring $R$. In fact, the last property holds either by
assumption, or $X$ is excellent and we may base change with the completion $\hat{\cO}$ of $\cO=\cO_{X,x}$ which then is a
quotient of a regular ring by Cohen's structure theorem. It also suffices to check the injectivity of $\cO\rightarrow \cO'$ after
this base change, because $\cO \rightarrow \hat{\cO}$ is faithfully flat.

\medbreak
Now we use the results of Hironaka in \cite{H4}. First consider the case $k(x') = k(x)$, where $\delta=0$.
Let $\fp \subset \cO = \cO_{X,x}$ be the prime ideal corresponding
to $D$, let $A = gr_\fp(\cO_{X,x})\otimes_{\cO/\fp}k(x)$, and let $F = \pi_X^{-1}(x) \subset X'$ be the scheme
theoretic fibre of $\pi_X: X' \rightarrow X$ over $x$. Then $F = \Proj(A)$, and by \cite{H4} (4.1) we have inequalities
\begin{equation}\label{eq.HS.blowup10}
H_{\cO_{X',x'}}^{(1+\delta)} \leq H_{\cO_{F,x'}}^{(2+\delta+s)} \leq H_{C_{X,D,x}}^{(1+s)} = H_{\cO_{X,x}}^{(1)}\,,
\end{equation}
where $s=\dim \cO_{D,x}$ and $C_{X,D,x} = \Spec(A)$. By our assumption, these are all equalities.
Moreover, if $\cO = R/J$ for a regular local ring $R$ with maximal ideal $\fm$, then
there is a system of regular parameters $(x_0,\ldots,x_r,y_1,\ldots,y_s)$ of $R$ such that
the ideal $P$ of $D$ in $R$ is generated by $x_0,\ldots,x_r$. If $X_0,\ldots,X_r,Y_1,\ldots,Y_s$
are the initial forms of $x_0,\ldots,x_r,y_1,\ldots,y_s$ (with respect to $\fm$), then the fiber
$E$ over $x$ in the blowup $Z'$ of $Z=\Spec(R)$ in the center $D$ is isomorphic to $\Proj(k[X.])$ for
$k=k(x)$ and the polynomial ring $k[X.]=k[X_0,\ldots,X_r]$, and there is a homogeneous ideal
$I \subset k[X.]$ such that $A=k[X.]/I$ and $F\subset E$ identifies with the canonical
immersion $\Proj(k[X.]/I) \subset \Proj(k[X.])$. We may assume that $x' \in F \subset E$
lies in the standard open subset $D_+(X_0) \subset \Proj(k[X.]) = E$.
If furthermore $k(x') = k(x)=k$ as we assume, then $\delta=0$, and by equality in the middle
of \eqref{eq.HS.blowup10} and the proof of lemma 8 in \cite{H4} there is a graded $k$-algebra
$B$ and an isomorphism of graded $k[X_0]$-algebras $A \cong B[X_0]$.

\bigskip
On the other hand, let $\fn$ be the maximal ideal of $\cO$, and let
$z_1,\ldots,z_s \in \fn$ be elements whose images $Z_1,\ldots,Z_s \in \fn/\fn^2$
form a basis of $\fn/\fp+\fn^2$. Then one has an isomorphism of graded algebras
$$
A[Z_1,\ldots,Z_s] \cong gr_\fn(\cO)
$$
induced by the canonical map $A \rightarrow gr_\fn(\cO)$, because $\fp \subset \cO$ is permissible
(\cite{H1} II 1. Proposition 1). This shows that the image of $X_0$ in $\fn/\fn^2$,
is not zero and not a zero divisor in $gr_\fn(\cO)$. Therefore the image of $x_0 \in R$ in $\cO$
is not zero and not a zero divisor in $\cO$.

\medbreak
Now we claim that every element in the kernel of $\cO \rightarrow \cO'$ is annihilated by a power of $x_0$,
which then gives a contradiction if this kernel is non-zero.
Let $R'= \cO_{Z',x'}$ be the monoidal transform of $R$ with center $P$ corresponding to $x'\in Z'$,
where $\fp= P/J$, and let $J'\subset R'$ be the strict transform of $J$, so that $\cO':= \cO_{X',x'}=R'/J'$.
Then, since
$$
R' = R[\mbox{\large $\frac{x_1}{x_0},\ldots,\frac{x_r}{x_0}$}]_{<\frac{x_1}{x_0},\ldots,\frac{x_r}{x_0}>}
\quad\mbox{and}\quad PR'=x_0R'\,,
$$
it follows from \cite{H1} III Lemma 6, p. 216 that there are generators $f_1,\ldots,f_m$
of $J$ and natural numbers $n_1,\ldots,n_m$ such that $J'$ is generated by
$f_1/x_0^{n_1},\ldots,f_m/x_0^{n_m}$. Evidently this implies
that every element in the kernel of $\cO \rightarrow \cO'$ is annihilated by a power of $x_0$.

\medbreak
Now consider the case that the residue field extension $k(x')/k(x)$ is arbitrary. We reduce to the residually
rational case ($k(x')=k(x)$)
by the same technique as in \cite{H4}. As there, one may replace $X$ by $\Spec \cO_{X,x}$,
and consider a cartesian diagram
$$
\begin{CD}
x' @.\quad\quad  X' @<{i'}<< \tilde X' @. \quad\quad  \tilde x' \\
@VVV \quad\quad @V{\pi_X}VV @VV{\pi_{\tilde X}}V \quad \quad @VVV \\
x @.\quad\quad  X  @<{i}<< \tilde X @. \quad \quad  \tilde x \,,
\end{CD}
$$
where $i$ is a faithfully flat monogenic map which is either finite or the projection $\tilde X = \bA^1_X \rightarrow X$,
and $\tilde f$ is the blow-up of $\tilde D = i^{-1}(D)$, which is again permissible.
Moreover, $\tilde x \in \tilde X$ is the generic point of $i^{-1}(x)$ such that $k(\tilde x)$ is
a monogenic field extension of $k(x)$, and there is a point $\tilde x' \in \tilde X'$
which maps to $x'\in X'$ and $\tilde x \in \tilde X$ and satisfies $k(\tilde x') = k(x')$.
Furthermore one has the inequalities
$$
H_{\cO_{X',x'}}^{(1+\delta)} \leq H_{\cO_{\tilde X',\tilde x'}}^{(1+\tilde \delta)} \leq H_{\cO_{\tilde X,\tilde x}}^{(1)} = H_{\cO_{X,x}}^{(1)}\,,
$$
where $\tilde \delta=tr.deg(k(\tilde x')/k(\tilde x))$ ($=\delta$ if $k(\tilde x)/k(x)$ is algebraic,
and $\delta-1$ otherwise). By our assumption all inequalities
become in fact equalities, and by induction on the number of generators of $k(x')$ over
$k(x)$ (starting with the residually rational case proved above), we may assume that
$\cO_{\tilde X,\tilde x} \rightarrow \cO_{\tilde X',\tilde x'}$ is injective. Since
$\cO_{X,x} \rightarrow \cO_{\tilde X,\tilde x}$ is injective, we obtained the injectivity of
 $\cO_{X,x} \rightarrow \cO_{X',x'}$. This completes the proof of (2), and while doing it, we
 also proved the claims in (3).

\medskip
Finally we prove (4). Still under the assumption that $X$ is embedded in a regular scheme $Z$,
Hironaka proved in \cite{H2} Theorem (1,A), that the equality
$\nu^*_{x'}(X',Z') = \nu^*_x(X,Z)$ implies the inequality in (4).
Together with (2) and (6) this implies (4) and finishes the proof of Theorem \ref{thm.pbu.inv}.
$\square$

\begin{corollary}\label{cor.pbu.inv}
For $\nu\in \SigmaXmax$ (cf. Definition \ref{Sigmamax}),
either $\nu\not\in \Sigma_{X'}$ or $\nu\in \SigmaXdmax$. We have
$$
\HSa {X'} \subset \pi_X^{-1}(\HSa X)\qaq
\pi_X^{-1}(\HSa X) \subset \underset{\mu\leq\nu}{\bigcup} \HSb {X'}.
$$
\end{corollary}

\begin{definition}\label{def.nearpoint}
Let the assumption be as in Theorem \ref{thm.pbu.inv} and put $k'=k(x')$.
\begin{itemize}
\item[(1)]
$x'\in \pi_X^{-1}(x)$ is near to $x$ if $\HSf {X'} {x'} =\HSf X x$.
\item[(2)]
$x'$ is very near to $x$ if it is near to $x$ and
$e_{x'}(X')+\delta_{x'/x} =e_x(X)_{k'}=e_x(X)$.
\end{itemize}
\end{definition}

We recall another result of Hironaka, as improved by Mizutani, which plays a crucial role in this paper.
Let again $X$ be an excellent scheme or a scheme embeddable in a regular scheme.
Let $D \subset X$ be a permissible closed subscheme, and let
$$
\pi_X : X'=B\ell_D(X) \to X\quad \text{and }\quad
$$
be the blowup with center $D$. Take any points $x\in D$ and $x'\in \pi_X^{-1}(x)$.

\begin{theorem}\label{thmdirectrix}
Assume that $x'$ is near to $x$.
Assume further that $\Char(k(x)) = 0$, or $\Char(k(x)) \geq \dim (X)/2 + 1$,
where $k(x)$ is the residue field of $x$.
Then
$$
x'\in {\mathbb P}(\Dir_x(X)/T_x(D))\subset \pi_X^{-1}(x)\,,
$$
where ${\mathbb P}(V)$ is the projective space associated to a vector space $V$.
\end{theorem}

{\bf Proof } First we note that the inclusion above is induced by the inclusion of
cones (i.e., spectra of graded algebras)
$\Dir_x(X)/T_x(D) \subset C_x(X)/T_x(D)$ and the isomorphism $C_x(X)/T_x(D) \cong C_D(X)_x$
of cones from Theorem \ref{nf.thm1} (2)(ii). More precisely, it is induced by applying the
$\Proj$-construction to the surjection of graded $k(x)$-algebras
$$
A_D = gr_{\fn_x}(\cO_{X,x})/\Sym(T_x(D)) \twoheadrightarrow \Sym(Dir_x(X))/\Sym(T_x(D))=B_D
$$
where we identify affine spaces with the associated vector spaces and note that $\Proj(A_D)=\pi_X^{-1}(x)$.
Since the claim is local, we may pass to the local ring $\cO=\cO_{X,x}$ of $X$ at $x$. Further,
as in the proof of Theorem \ref{thm.pbu.inv}, we may assume that $X$ is embedded into a regular scheme $Z$.
If $\pi_Z : Z'=B\ell_D(Z) \to Z$ denotes the blowup of $Z$ in $D$, we have a further
inclusion
$$
{\mathbb P}(C_x(X)/T_x(D)) = \pi_X^{-1}(x) \subset {\mathbb P}(T_x(Z)/T_x(D))= \pi_Z^{-1}(x)\,.
$$
Therefore the claim that $x' \in  {\mathbb P}(\Dir_x(X)/T_x(D))$ follows from \cite{H4}, Theorem IV,
\cite{H5}, Theorem 2, and \cite{Miz}. In fact, by the second reference there is a certain canonical
subgroup scheme $B_{{\mathbf P},x'} \subset {\mathbf V} = T_x(Z)/T_x(D)$ just depending on
$x' \in {\mathbf P}= {\mathbb P}({\mathbf V})$, which has the following properties. It is
defined by homogeneous equations in the coordinates of ${\mathbf V}$, hence is a subcone of
${\mathbf V}$, and the associated subspace ${\mathbb P}(B_{\mathbf P,x'})$
contains $x'$. Moreover, by the first reference it is a vector subspace of ${\mathbf V}$
if $\Char(k(x)) = 0$, or $\Char(k(x)) = p > 0$ with $p \geq \dim(B_{{\mathbf P},x'})$, and
in \cite{Miz} this was improved to the (sharp) bound $p \geq \dim(B_{{\mathbf P},x'})/2 + 1$.
On the other hand, by the first reference, the action of $B_{{\mathbf P},x'}$ on ${\mathbf V}$
respects $C_x(X)/T_x(D)$ if $x'$ is near to $x$. Since $0 \in C_x(X)$, we conclude that $B_{\mathbf P,x'}$
is contained in $C_x(X)/T_x(D)$, and hence has dimension at most $d=\dim(X)-\dim(D)\leq \dim(X)$.
Therefore, by the assumption $p \geq \dim(X)/2 + 1$,  $B_{\mathbf P,x'}$
is a vector subspace of $C_x(X)/T_x(D)$ and is thus contained in the biggest such subspace - which
is $\Dir_x(X)/T_x(D)$. Therefore $x' \in {\mathbb P}(B_{{\mathbf P},x'}) \subset {\mathbb P}(\Dir_x(X)/T_x(D))$.

\medbreak

\begin{lemma}\label{Lem1}
Consider $\HSa X$ for $\nu\in\SigmaXmax$.
Let $\pi: X'=\Bl_D(X)\to X$ be the blowup with permissible center $D$
contained in $\HSa X$.
Let $Y\subset\HSa X$ be an irreducible closed subset which contains $D$
as a proper subset.
Then:
\begin{itemize}
\item[(1)]
$Y'\subset \HSa {X'}$, where $Y'\subset X'$ be the proper transform.
\item[(2)]
Assume $\Char(k(x)) = 0$, or $\Char(k(x)) \geq \dim (X)/2+1$. Then we have
$e_x(X)\geq 1$ for any $x\in D$ and $e_{x'}(X')\geq 1$ for any
$x'\in \pi^{-1}(D)\cap Y'$.
\end{itemize}
\end{lemma}
\medbreak\noindent
\textbf{Proof}
Let $\eta$ (resp. $\eta'$) be the generic point of $Y$ (resp. $Y'$).
Take points $x\in D$ and $x'\in \pi^{-1}(x)\cap Y'$.
Then we have
$$
\HSfX x \geq \HSfXd {x'} \geq \HSfXd  {\eta'} = \HSfX {\eta},
$$
where the first inequality follows from Theorem \ref{thm.pbu.inv}(1), the
second from theorem \ref{HSf.usc}, and the last equality follow from the fact
$\cO_{X,\eta}\cong \cO_{X',\eta'}$. Since $Y\subset \HSa X$, we
have $\HSf X x =\HSf X \eta = \nu$ so that the above inequalities are
equalities. This implies $Y' \subset \HSa {X'}$, which proves (1).
Next we show (2). If $e_x(X)=0$ for $x\in D$, then there is no point of
$\Bl_x(X)$ which is near to $x$ by Theorem \ref{thmdirectrix}.
Thus (1) implies $e_x(X)\geq 1$. To show $e_{x'}(X')\geq 1$ for
$x'\in \pi^{-1}(D)\cap Y'$, let $W\subset X'$ be the closure of $x'$ in $X'$.
By assumption $W$ is a proper closed subset of $Y'$. Since $e_{x'}(X')$ is a
local invariant of $\cO_{X',x'}$, we may localize $X'$ at $x'$ to assume $W$
is regular. Then, by Theorem \ref{HSf.perm}, $W'\subset X'$ is permissible.
Now the assertion follows from the previous assertion applied to
$\Bl_W(X') \to X'$.
$\square$
\medbreak

\newpage
\section{$B$-Permissible blow-ups - the embedded case}\label{BpbuI}
\bigskip

Let $Z$ be a regular scheme and let $\bB \subset Z$ be a simple normal crossing divisor on $Z$.
For each $x\in Z$, let $\bB(x)$ be the subdivisor of $B$ which is the union
of the irreducible components of $\bB$ containing $x$.

\begin{definition}\label{def.DncBZ}
Let $D\subset Z$ be a regular subscheme and $x\in D$. We say $D$ is normal
crossing (n.c.) with $\bB$ at $x$ if there exists a system $z_1,\dots,z_d$
of regular parameters of $R:=\cO_{Z,x}$ satisfying the following conditions:
\begin{itemize}
\item[(1)]
$D\times_Z\Spec(R)=\Spec(R/\langle z_1,\dots,z_r \rangle)$ for some
$1\leq r\leq d$.
\item[(2)]
$\bBx\times_Z\Spec(R)=\Spec(R/ \langle \underset{j\in J}{\prod} z_j \rangle)$
for some $J\subset \{1,\dots,d\}$.
\end{itemize}
We say $D$ is transversal with $\bB$ at $x$ if in addition the parameters can
be chosen such that $J \subset \{r+1,\ldots,d\}$.
We say $D$ is n.c. (resp. transversal) with $\bB$ if $D$ is n.c. (resp. transversal)
with $\bB$ at any point $x\in D$.
\end{definition}

Note that $D$ is n.c. with $\bB$ if and only if
$D$ is transversal with the intersection of any set of irreducible components of $\bB$
which do not contain $D$.
\medbreak

Let $D\subset Z$ be n.c. with $\bB$. Consider $Z'= \Bl_D(Z) \rmapo{\pi_Z} Z$.
Let $\tbB=\Bl_{D\times_Z \bB}(\bB)\subset Z'$ be the strict transform of $\bB$ in $Z'$,
let $E:=\pi_Z^{-1}(D)$ be the exceptional divisor, and let $\bB' = \tbB \cup E$ be the complete transform of $\bB$ in $Z'$.
We easily see the following:

\begin{lemma}\label{lem.Bperm0}
$\tbB$ and $\bB'$ are simple normal crossing divisors on $Z'$.
\end{lemma}
\medbreak

For the following, and for the comparison with the next section, it will be more
convenient to consider the set $\cB = C(\bB)$ of irreducible components
of $\bB$.

\begin{definition}\label{def.ncbnd}
A simple normal crossings boundary on $Z$ is a set $\cB = \{B_1,\ldots,B_n\}$ of
regular divisors on $Z$ such that the associated divisor $\div(\cB) = B_1\cup\ldots\cup B_n$
is a (simple) normal crossings divisor. For $x\in Z$ let $\cBx = \{ B\in \cB \mid x \in B\}$.
Often the elements of $\cB$ are also called components of $\cB$.
\end{definition}

An equivalent condition is that the $B_i$ intersect transversally, i.e., that for each
subset $\{i_1,\ldots,i_r\}\subset \{1,\ldots,n\}$ the intersection
$B_{i_1}\times_Z\ldots \times_Z B_{i_r}$ is regular of pure codimension $r$ in $Z$.
The associations
\begin{equation}\label{eq.corr.ncd-ncb}
\bB \mapsto \cB = C(\bB) \quad\quad , \quad\quad \cB \mapsto \bB = \div(\cB)
\end{equation}
give mutual inverse bijections between the set of simple normal crossing (s.n.c.)
divisors on $Z$ and the set of simple normal crossing (s.n.c.) boundaries on $Z$,
and we will now use the second language. Via \eqref{eq.corr.ncd-ncb} the above definitions
correspond to the following in the setting of boundaries.

\begin{definition}\label{def.nc-tr.bnd}
(a) A regular subscheme $D\subset Z$
is transversal with a s.n.c. boundary $\cB$ at $x$ if for $\cBx = \{B_1,\ldots,B_r\}$ it
intersects all multiple intersections $B_{i_1}\times_Z\ldots\times_Z B_{i_s}$ transversally
at $x$, and $D$ is normal crossing (n.c.) with $\cB$ at $x$ if it is transversal with
$\cBx - \cBxD$ at $x$, where $\cBxD = \{ B \in \cBx \mid D \subset B\}$. $D$ is transversal
(resp. n.c.) with $\cB$, if it is transversal (resp. n.c.) with $\cB$ at all points $x\in D$.
\smallskip

\pagebreak

(b) If $D$ is n.c. with $\cB$, and $Z'=\Bl_D(Z) \rmapo{\pi_Z} Z$ is the blow-up of $D$, then
the strict and the complete transform of $\cB$ are defined as
$$
\tcB := \{ \tB \mid B \in \cB\} \quad\quad , \quad\quad \cB' :=  \tcB\,\cup\,\{E\}\,
$$
respectively, where $\tB=\Bl_{D\times_Z B}(B)$ is the strict transform of $B$ in $Z$, and
$E = \pi_Z^{-1}(D)$ is the exceptional divisor.
(Note that $\tcB$ and $\cB'$ are s.n.c. boundaries on $Z'$ by Lemma \ref{lem.Bperm0}.)
\end{definition}

In the following, consider a regular scheme $Z$ and a simple normal crossing boundary $\cB$ on $Z$.
Moreover let $X \subset Z$ be a closed subscheme.

\begin{definition}\label{def.Bperm}
Let $D\subset X$ be a regular closed subscheme and $x\in D$.
We say $D$ is $\cB$-permissible at $x$ if $D\subset X$ is permissible
at $x$ and $D$ is n.c. with $\cB$ at $x$. We say
$D\subset X$ is $\cB$-permissible if $D\subset X$ is permissible at all $x\in D$.
\end{definition}

\begin{definition}\label{def.histfcnZ}
A history function for $\cB$ on $X$ is a function
\begin{equation}\label{eq.OI}
\OB: X \to \{\text{subsets of $\cB$}\}\;;\; x \to \OBx,
\end{equation}
which satisfies the following conditions:
\begin{itemize}
\item[(O1)]
For any $x\in X$, $\OBx\subset \cBx$.
\item[(O2)]
For any $x,y\in X$ such that $x\in \overline{\{y\}}$ and $\HSfX x=\HSfX y$,
we have $\OBy\subset \OBx$.
\item[(O3)]
For any $y\in X$, there exists a non-empty open subset $U\subset \overline{\{y\}}$
such that $\OBx=\OBy$ for all $x\in U$ such that $\HSfX x=\HSfX y$.
\end{itemize}
For such a function, we put for $x\in X$,
$$
\NBx= \cBx - \OBx.
$$
A component of $\cB$ is called old (resp. new) for $x$ if it is a component of
$\OBx$ (resp. \cNBx).
\end{definition}

A basic example of a history function for $B$ on $X$ is given by the following:

\begin{lemma}\label{lem.histfcnBx}
The function $\cOBx =\cBx$ ($x\in X$), is a history function for
$\cB$ on $X$. In fact it satisfies \ref{def.histfcnZ} (O2) and
(O3) without the condition $\HSfX x=\HSfX y$.
\end{lemma}

\textbf{Proof } Left to the readers.
\medbreak

Define the Hilbert-Samuel function of $(X,O)$ as:
$$
\cHSfXX : X \to \bNNsq \;;\; x\to (\HSfX x, \coXx),
$$
where $\coXx$ is the cardinality of $\cOBx$. We endow $\bNNsq$
with the lexicographic order:
$$
(\nu,\mu)\geq (\nu',\mu')\Leftrightarrow \nu >\nu'\text{ or }
\nu=\nu'\text{ and }\mu\geq \mu'.
$$
Theorem \ref{HSf.usc} and the conditions in \ref{def.histfcnZ}
immediately imply the following:

\begin{theorem}\label{thm.BHSf.usc}
Let the assumption be as above.
\begin{itemize}
\item[(1)]
If $x\in X$ is a specialization of $y\in X$, then
$\cHSfX x \geq \cHSfX y$.
\item[(2)] For any $y\in X$, there is a
dense open subset $U$ of $\overline{\{y\}}$ such that $\cHSfX
y=\cHSfX x$ for any $x\in U$.
\end{itemize}
\end{theorem}

In other words (see Lemma \ref{lem.up.sc} (a)), the function $\cHSfXX$
is upper semi-continuous on $X$. Let
$$
\cSigmaX:=\{\cHSfX x|\; x\in X\} \subset \bNN\times\bN\,,
$$
and let $\cSigmaXmax$ be the set of of the maximal elements in $\cSigmaX$.

\begin{definition}\label{def.BHSlocus}
(1) For $\tnu \in \cSigmaX$ we define
$$
X(\tnu) = X^\OB(\tnu) = \{ x \in X \mid \cHSfX x = \tnu \}
$$
and
$$
\tHSgeqa X = X^\OB(\geq \tnu) = \{ x \in X \mid \cHSfX x \geq \tnu \}\,,
$$
and we call
\begin{equation}\label{eq.BHSlocus}
\cHSmax X =\underset{\tnu\in \cSigmaXmax}{\cup} X(\tnu)\,.
\end{equation}
the $\cOB$-Hilbert-Samuel locus of $X$.

\smallskip
(2) We define
$$
\cDirx(X):=
\Dir_x(X)\cap \underset{B\in \OBx}{\bigcap} T_x(B) \; \subset T_x(Z).
$$
$$
\eob_x(X)=\dim_{k(x)}(\cDirx(X)).
$$
\end{definition}

By Lemma \ref{lem.up.sc}, $X(\tnu)$ is locally closed, with closure contained in $\tHSgeqa X$.
Moreover, $X(\tnu)$ is closed for $\tnu\in \cSigmaXmax$, the union in \eqref{eq.BHSlocus} is disjoint,
$\cHSmax X$ is a closed subset of $X$, and $\cSigmaX$ is finite if $X$ is noetherian.
Theorems \ref{nf.thm1} and \ref{HSf.perm} imply the following:

\begin{theorem}\label{thm.OBHSf.perm}
Let $D\subset X$ be a regular closed subscheme and $x\in D$.
Then the following conditions are equivalent:
\begin{itemize}
\item[(1)]
$D\subset X$ is permissible at $x$ and there is an open neighborhood
$U$ of $x$ in $Z$ such that $D\cap U\subset B$ for every $B \in \OBx$.
\item[(2)]
$\tHSfX x =\tHSfX y$ for any $y\in D$ such that $x$ is a specialization of $y$.

\end{itemize}
Under the above conditions, we have
\begin{equation}\label{BDir}
T_x(D)\subset \cDirx(X).
\end{equation}
\end{theorem}

\begin{definition}\label{def.Operm}
A closed subscheme $D \subset X$ is called $\cOB$-permissible at $x$,
if it satisfies the equivalent conditions in Theorem
\ref{thm.OBHSf.perm}, and it is called $\cOB$-permissible, if it is
$\cOB$-permissible at all $x\in D$.
\end{definition}

\begin{remark}\label{rem.Operm}
Note that $D\subset X$ is $\cB$-permissible at $x$ if and only if
$D$ is $\cOB$-permissible at $x$ and n.c with $\cNBx$ at $x$.
\end{remark}

\begin{theorem}\label{nf.Bdir}
Let the assumption be as in
Theorem \ref{thm.OBHSf.perm} and assume that $D$ is irreducible.
Assume:
\begin{itemize}
\item[(1)]
$D\subset X$ is $\cOB$-permissible at $x$.
\item[(2)]
$e_\eta(X) = e_x(X) -\dim(\cO_{D,x})$ (cf. Theorem \ref{nf.dir}),
\end{itemize}
where $\eta$ is the generic point of $D$. Then we have
$$
\eob_\eta(X)\leq \eob_x(X) -\dim(\cO_{D,x}).
$$
\end{theorem}
\par\noindent
\textbf{Proof }
First we claim that (1) and (2) hold after replacing $x$ by any
point $y\in D$ such that $x\in \overline{\{y\}}$ and $\overline{\{y\}}$ is
regular at $x$. Indeed the claim follows from the inequalities:
$$
\begin{aligned}
& \tHSfX \eta\leq \tHSfX y \leq \tHSfX x,\\
& \begin{aligned}
e_\eta(X)
&\leq e_y(X)-\dim(\cO_{D,y})\\
&\leq \big(e_x(X)-\dim(\cO_{\overline{\{y\}},x})\big)-\dim(\cO_{D,y})
= e_x(X)-\dim(\cO_{D,x}),\\
\end{aligned}
\end{aligned}
$$
which follows from Theorems \ref{thm.BHSf.usc} and \ref{nf.dir}.
By the claim we are reduced to the case where $\dim(\cO_{D,x})=1$
by the same argument as in the proof of Theorem \ref{nf.dir}.
Let $R=\cO_{Z,x}$ and let $J\subset \fp\subset R$ be the ideals defining
$X\subset Z$ and $D\subset Z$, respectively. Let $R_\fp$ is the localization
of $R$ at $\fp$ and $J_\fp=J R_\fp$.
By Lemma \ref{nf.dir.lem1} (2) implies that there exists a part of a system
of regular parameters $y=(y_1,\dots,y_r)$ of $R$
such that $y\subset \fp$ and that
$$
\begin{aligned}
&I\Dir(R_\fp/J_\fp)=\langle in_{\fp}(y_1),\dots,in_{\fp}(y_r)\rangle
\;\subset \gr_{\fp}(R_\fp),\\
&I\Dir(R/J)=\langle in_{\fm}(y_1),\dots,in_{\fm}(y_r)\rangle
\;\subset \gr_{\fm}(R),\\
\end{aligned}
$$
We can take $\theta_1,\dots,\theta_s\in R$ such that
$(y_1,\dots,y_r,\theta_1,\dots,\theta_s)$ is a part of a system of regular parameters of $R$ and that there exists an irreducible component $B_i$
of $\OBx$ for each $i=1,\dots,s$ such that
$B_i\times_Z\Spec(R)=\Spec(R/\langle \theta_i\rangle)$ and
$$
I\Dirob(R/J)=\langle in_{\fm}(y_1),\dots,in_{\fm}(y_r),
in_{\fm}(\theta_1),\dots,in_{\fm}(\theta_s) \rangle ,
$$
where $I\Dirob(R/J)\subset \grmR$ is the ideal defining
$\Dirob_x(X)\subset T_x(Z)$. Now (1) implies $\OBx=\OB(\eta)$ so that
$D\subset B_i$ and $\theta_i\in \fp$ for all $i=1,\dots,s$.
Hence $(y_1,\dots,y_r,\theta_1,\dots,\theta_s)$ is a part of a system of regular parameters of $R_{\fp}$ and
$$
I\Dirob(R_{\fp}/J_{\fp}) \supset \langle in_{\fp}(y_1),\dots,in_{\fp}(y_r),
in_{\fp}(\theta_1),\dots,in_{\fp}(\theta_s) \rangle ,
$$
which implies the conclusion of Theorem \ref{nf.Bdir}.
$\square$
\bigskip

Let $D\subset X$ be a $\cB$-permissible closed subscheme. Consider the diagram
\begin{equation}\label{eq.Bpermbu1}
\begin{CD}
X' =\;@. \Bl_D(X) @>{i'}>> \Bl_D(Z) @.\; =  Z' \\
@.  @V{\pi_X}VV @VV{\pi_Z}V   @. \\
@.  X  @>{i}>> Z @.    \,,
\end{CD}
\end{equation}
and let $\cB'$ and $\tcB$ be the complete and strict transform of $\cB$ in $Z'$, respectively.
For a given history function $\OBx$ ($x\in X$), we define
functions $\OBd,\, \tOB \,: X'\to \{\text{subsets of $\cB'$}\}$ as follows:
Let $x'\in X'$ and $x=\pi_X(x')\in X$. Then define
\begin{equation}\label{def.OIxd}
\OBdxd=\left.\left\{\gathered
\tOBx\cap\cBdxd \\
 \cBdxd \\
\endgathered\right.\quad
\begin{aligned}
&\text{if $\HSfXd {x'}=\HSfX x$}\\
&\text{otherwise},
\end{aligned}\right.
\end{equation}
where $\tOBx$ is the strict transform of $\OBx$ in $Z'$, and
\begin{equation}\label{def.tOIxd}
\tOB(x')=\left.\left\{\gathered
\tOBx\cap\cBdxd \\
 \tcB(x') \\
\endgathered\right.\quad
\begin{aligned}
&\text{if $\HSfXd {x'}=\HSfX x$}\\
&\text{otherwise}.
\end{aligned}\right.
\end{equation}
Note that $\OBdxd = \tOB(x') = \tOBx\cap\tcB(x')$ if $x'$ is near to $x$.

\begin{lemma}\label{lem.Bperm1}
The functions
$x' \to \OBdxd,\,x' \to \tOB(x')$
are history functions for $\cB'$ on $X$.
\end{lemma}
\medbreak

The proof of Lemma \ref{lem.Bperm1} will be given later.

\begin{definition}\label{def.compl.transf.cBO}
We call $(\cB',\OB')$ and $(\tcB,\tOB)$ the complete and strict transform of $(\cB,\OB)$ in $Z'$, respectively.
\end{definition}

\begin{theorem}\label{thm.Bpbu.inv}
Take points $x\in D$ and $x'\in \pi_X^{-1}(x)$.
Then $H^{\tOB}_{X'}(x') \leq H^{O'}_{X'}(x') \leq \tHSfX x$. In particular we have
$$
\pi_X^{-1}(X^{\OBd}(\tnu)) \subset \underset{\tmu\leq\tnu}{\cup} {X'}^{\OBd}(\tmu)
\qfor \tnu\in \SigmaXmax,
$$
and the same holds for $\tOB$ in place of $\OBd$.
\end{theorem}
\textbf{Proof } This follows immediately from Theorem \ref{thm.pbu.inv}, \eqref{def.OIxd} and \eqref{def.tOIxd}.
\bigskip

In the following we mostly use the complete transform $(\cB',\OB')$ and,
for ease of notation, we often write $\cHSfXd {x'}$ and $\cSigmaXd$ instead of
$H^{O'}_{X'}(x')$ and $\Sigma^{O'}_{X'}$, similarly for $\cSigmaXdmax$  etc., because
everything just depends on $\cOB$.

\begin{definition}\label{def.Bnearpoint}
We say that $x'\in \pi_X^{-1}(x)$ is $\OB$-near to $x$ if the following equivalent conditions hold:
\begin{itemize}
\item[(1)]
$\tHSfXd {x'} =\tHSfX {x}$ \quad \quad($\;\Leftrightarrow\; H^{\OBd}_{X'}(x') =\tHSfX {x} \;\Leftrightarrow \;
H^{\tOB}_{X'}(x') = \tHSfX {x}\;$).
\item[(2)]
$x'$ is near to $x$ and contained in the strict transforms of all $B \in \OBx$.
\end{itemize}
Call $x'$ very $\OB$-near to $x$ if $x'$ is $\OB$-near and very near to $x$
and $\eob_{x'}(X')= \eob_x(X)-\delta_{x'/x}$.
\end{definition}

\smallskip
The following result is an immediate consequence of Theorem
\ref{thmdirectrix} and Definition \ref{eq.BHSlocus} (2).

\begin{theorem}\label{thmBdirectrix}
Assume that $x'\in X'$ is $\OB$-near to $x=\pi_Z(x')\in X$.
Assume further that $\Char(k(x)) = 0$, or $\Char(k(x)) \geq \dim (X)/2+1$,
where $k(x)$ is the residue field of $x$.
Then
$$
 x'\in {\mathbb P}(\cDirx(X)/T_x(D))\subset
{\mathbb P}(T_x(Z)/T_x(D))=\pi_Z^{-1}(x)
$$
\end{theorem}

\bigskip\noindent
\textbf{Proof of Lemma \ref{lem.Bperm1}}:
Take $y',x'\in X'$ such that $x'\in \overline{\{y'\}}$ and
$\HSfXd {x'}=\HSfXd {y'}$. We want to show $\OBdyd\subset \OBdxd$.
Put $x=\pi_X(x'),y=\pi_X(y') \in X$.
We have $x\in \overline{\{y\}}$. By Theorems \ref{HSf.usc} and \ref{thm.pbu.inv}
$$
\HSfXd {y'}\leq \HSfX y \leq \HSfX x \geq \HSfXd {x'}=\HSfXd {y'}.
$$
First assume $\HSfX x \geq \HSfXd {x'}$, which implies
$\HSfXd {y'} =\HSfX y = \HSfX x$. By \ref{def.histfcnZ} (O2) and \eqref{def.OIxd} we get
$$
\OBdyd=\widetilde{\OBy}\cap\cBdyd \subset \widetilde{\OBx}\cap\cBdxd = \OBdxd.
$$
Next assume $\HSfX x > \HSfXd {x'}$. Then, by Lemma \ref{lem.histfcnBx} and
\eqref{def.OIxd}, we get
$$
\cOBdyd \subset \cBdyd \subset \cBdxd = \cOBdxd.
$$
Next we show that for $y'\in X'$, there exists a non-empty open subset
$U'\subset \overline{\{y'\}}$ such that $\HSfXd {x'}=\HSfXd {y'}$ and $\cOBdxd=\cOBdyd$
for all $x'\in U'$. Put $y=\pi_X(y')$.
By Lemma \ref{lem.histfcnBx} and \ref{def.histfcnZ} (O3), there exists a non-empty open subset
$U\subset \overline{\{y\}}$ such that $\HSfX {x}=\HSfX {y}$, $\cBx=\cBy$ and $\cOBx=\cOBy$ for all $x\in U$.
By Theorem \ref{HSf.usc}(3) and Lemma \ref{lem.histfcnBx}, there exists a non-empty open subset
$U'\subset \overline{\{y'\}}\cap \pi_X^{-1}(U)$
such that $\HSfXd {x'}=\HSfXd {y'}$ and $\cBdxd=\cBdyd$ for all $x'\in U'$.
We now show $U'$ satisfies the desired property.
Take $x'\in U'$ and put $x=\pi_X(x')$. By the assumption we have
$\HSfXd {x'}=\HSfXd {y'}$ and $\HSfX {x}=\HSfX {y}$.
\medbreak

First assume $\HSfXd {x'}=\HSfX {x}$, which implies $\HSfXd {y'}=\HSfX {y}$.
By \eqref{def.OIxd} we get
$$
\OBdyd=\tcOBy\cap\cBdyd = \tcOBx\cap\cBdxd =\OBdxd.
$$
Next assume
$\HSfXd {x'}<\HSfX {x}$, which implies $\HSfXd {y'}<\HSfX {y}$.
By \eqref{def.OIxd} we get
$$
\OBdyd=\cBdyd=\cBdxd=\OBdxd.
$$
This completes the proof of Lemma \ref{lem.Bperm1} for $(\cB',\OB')$. The proof for $(\tcB,\tOB)$
is similar. $\square$

\bigskip
For $x'\in X'$ let $\cNBdxd = \cBdxd - \cOBdxd$ be the set of the new components of $(\cB',\OB')$ for $x'$.
If $x'$ is near to $x = \pi_X(x') \in D$, i.e., $\HSfXd {x'}=\HSfX x$, then
\begin{equation}\label{eq.NBdxd}
\cNBdxd= (\widetilde{\cNBx}\cap \cBdxd)\, \cup \{E\} \quad \mbox{ with }\quad
E=\pi_Z^{-1}(D)
\end{equation}
where $\widetilde{\cNBx}$ is the strict transform of $\cNBx$ in $X'$. If $x'$ is not near to $x$,
then $\cNBdxd = \emptyset$.
Similarly, $\tN(x') = \tcB(x') - \tOB(x') = \widetilde{\cNBx}\cap\tcB(x') \subset \cNBdxd$ if $x'$
is near to $x$, and $\tN(x') = \emptyset$, otherwise.
We study the transversality of $\cNBdxd$ with a certain regular subscheme of $E$.

\begin{definition}\label{def.NBncDir}
For a $k(x)$-linear subspace $T\subset T_x(Z)$, we say that $T$ is transversal with
$\cNBx$ (notation: $T \pitchfork \cNBx$) if
$$
\dim_{k(x)}\big(T \cap \underset{B \in \cNBx}{\cap} T_x(B)\big)=
\dim_{k(x)}(T)-\mid \cNBx\mid.
$$
\end{definition}

\begin{lemma}\label{lem.NBncDir1}
Let $\pi_Z:Z'=\Bl_D(Z)\to Z$ be as in \eqref{eq.Bpermbu1}.
Assume $D=x$ and $T \pitchfork \NBx$. Then the closed subscheme
$\bP(T)\subset E=\bP(T_x(Z))$ is n.c. with $\NBdxd$ and $\tN(x')$ at each $x'\in \pi_X^{-1}(x)$.
\end{lemma}
\textbf{Proof }
Let $R=\cO_{Z,x}$ with the maximal ideal $\fm$.
For each $B\in \NIx$, take $h_B\in R$ such that
$B\times_Z \Spec(R)=\Spec(R/\langle h_B\rangle)$.
Put $H_B=\inm(h_B)\in \grm^1(R)$.
In view of \eqref{eq.NBdxd} the lemma follows from the following facts:
The ideal $\langle H_B \rangle\subset\grmR$ defines the subschemes
$$
T_x(B)\subset T_x(Z)=\Spec(\grmR)\qaq
E\times_{Z'}\tB\subset E=\Proj(\grmR),
$$
where $\tB$ is the strict transform of $B$ in $Z'$.
$\square$

\begin{lemma}\label{lem.NBncDir2}
Let $\pi_Z:Z'=\Bl_D(Z)\to Z$ be as in \eqref{eq.Bpermbu1}.
Assume $T \pitchfork \NBx$ and
$T_x(D)\subset T$ and $\dim_{k(x)}(T/T_x(D))=1$. Consider
$$
\{x'\}=\bP(T/T_x(D))\subset \bP(T_x(Z)/T_x(D))=\pi_Z^{-1}(x).
$$
Let $D'\subset E$ be any closed subscheme such that $x' \in D'$ and $\pi_Z$ induces
an isomorphism $D'\simeq D$. Then $D'$ is n.c. with $\NBdxd$ and $\tN(x')$ at $x'$.
\end{lemma}
\medbreak\noindent
\textbf{Proof }
It suffices to consider $\NBdxd$.
For $B \in \cBxD$ we have $T_x(D)\subset T_x(B)$ so that the assumptions of
the lemma imply $T_x(B)\cap T= T_x(D)$.
By the argument of the last part of the proof of Lemma \ref{lem.NBncDir1},
this implies $\{x'\}=\bP(T/T_x(D)) \nsubseteq \tB$. Thus we are reduced to showing
$D'$ is n.c. with $\NBdxd \cap (\cB(x) - \cBxD)'(x')$,
which follows from:

\begin{lemma}\label{lem.NBncDir3}
Let $Z$ be a regular scheme and $D,W\subset Z$ be regular closed subschemes
such that $D$ and $W$ intersect transversally.
Let $\pi:Z'=\Bl_D(Z)\to Z$ and let $\widetilde{W}$ be the strict transform of $W$ in $Z'$.
Suppose $D'\subset E:=\pi^{-1}(D)$ is a closed subscheme such that
$\pi$ induces an isomorphism $D'\isom D$. Then
$D'$ and $\widetilde{W}$ intersect transversally.
\end{lemma}
\medbreak\noindent
\textbf{Proof }
By definition, $\widetilde{W}=\Bl_{D\cap W}(W)$. The transversality of $D$ and $W$
implies $\Bl_{D\cap W}(W)\simeq \Bl_D(Z)\times_Z W=Z'\times_Z W$.
Thus
$$
E \times_{Z'} \widetilde{W} \simeq E\times_{Z'} (Z'\times_Z W) =E\times_Z W=
E\times_D (D\times_Z W).
$$
Hence we get
$$
D'\times_{Z'} \widetilde{W} = D'\times_E (E\times_{Z'} \widetilde{W}) \simeq
D'\times_E (E\times_D (D\times_Z W)) =D'\times_D (D\times_Z W)
\simeq W\times_Z D,
$$
where the last isomorphism follows from the assumption $D'\isom D$.
This completes the proof of the lemma.
$\square$
\bigskip

\begin{theorem}\label{thm.NBncDir}
Let $\pi_Z:Z'=\Bl_D(Z)\to Z$ be as in \eqref{eq.Bpermbu1}.
Take $x\in X$ and $x'\in \pi_X^{-1}(x)$. Assume
$\Char(k(x))=0$, or $\Char(k(x))\geq \dim(X)/2+1$.
\begin{itemize}
\item[(1)]
If $x'$ is $\OB$-near and very near to $x$, then
$\eob_{x'}(X')\leq \eob_x(X)-\delta_{x'/x}$.
\item[(2)]
Assume $x'$ is very $\OB$-near and $\NBx\pitchfork \Dirob_x(X)$.
Then $\NBdxd\pitchfork \Dirob_{x'}(X')$.
\end{itemize}
\end{theorem}
\textbf{Proof }
We first show (1).
Assume that $x'$ is $\OB$-near and very near to $x$.
For the sake of the later proof of (2),
we also assume $\NBx\pitchfork \Dirob_x(X)$.
By doing this, we do not lose generality for the proof of (1)
since we may take $\NBx=\emptyset$.
Put $R=\cO_{Z,x}$ (resp. $R'=\cO_{Z',x'}$) with the maximal ideal
$\fm$ (resp. $\fm'$) and $k=k(x)=R/\fm$ (resp. $k'=k(x')=R'/\fm'$).
By the assumptions and Theorem \ref{thmBdirectrix} there exists a system of regular parameters of $R$
$$
(y_1,\dots,y_r,\theta_1,\dots,\theta_q, u_1,\dots,u_a,u_{a+1},\dots,u_{a+b},
v_1,\dots,v_s,v_{s+1},\dots,v_{s+t})
$$
satisfying the following conditions: Fixing the identification:
$$
\grm(R)=k[Y,\Theta,U,V]=
k[Y_1,\dots,Y_r,\Theta_1,\dots,\Theta_q,U_1,\dots,U_{a+b},V_1,\dots,V_{s+t}],
$$
$$
(Y_i=\inm(y_i),\; \Theta_i=\inm(\theta_i),\; U_i=\inm(u_i),\; V_i=\inm(v_i)\;\in \gr^1_\fm(R))
$$
we have
\begin{itemize}
\item[$(i)$]
$D\times_Z \Spec(R)=
\Spec(R/\langle y_1,\dots,y_r,, u_1,\dots,u_{a+b} \rangle)$.
\item[$(ii)$]
$\IDir_x(X)=\langle Y_1,\dots,Y_r \rangle \;\subset\grm(R)$
(cf. Definition \ref{def.invschemes}),
\item[$(iii)$]
For $1\leq i\leq q$, there exists $B \in\OB(x)$ such that
$B\times_Z\Spec(R)=\BBth i$,
where $\BBth i:=\Spec(R/\langle\theta_i\rangle)$, and we have
$$
\cDirx(X)=\Dir_x(X)\cap \underset{1\leq i\leq q}{\cap}T_x(\BBth i).
$$
\item[$(iv)$]
$\NBx \times_Z \Spec(R)= \underset{1\leq i\leq b}{\bigcup} \BBu {a+i} \cup
\underset{1\leq j\leq t}{\bigcup} \BBv {s+j}.$
\end{itemize}
Here
$\BBu {i}=\Spec(R/\langle u_i\rangle)$,
$\BBv {j}=\Spec(R/\langle v_j\rangle)$.
Let $\BBthd i$ (resp. $\BBud i$, resp. $\BBvd i$) be the strict transform of
$\BBth i$ (resp. $\BBu i$, resp. $\BBv i$) in $\Spec(R')$. Let
$$
\Xi=\{i\in [a+1,a+b]\;|\; x'\in \BBud i \}.
$$
By Theorem \ref{thm.Bpbu.inv} there exists
$i_0 \in [1, a+b] -\Xi$ such that
$$
(y'_1,\dots,y'_r, \theta'_1,\dots,\theta'_q,w, u'_j\;(j\in \Xi),v_1,\dots,v_{s+t}),
$$
where $w=u_{i_0}$, $y'_i=y_i/w$, $\theta'_i=\theta/w$, $u'_j=u_j/w$,
is a part of a system of regular parameters of $R'$ so that
the polynomial ring:
$$
k'[Y'_1,\dots,Y'_r,\Theta'_1,\dots,\Theta'_q,W,U'_j\; (j\in \Xi),V_1,\dots,V_{s+t}],
$$
where $Y'_i=in_{\fm'}(y'_i),\; \Theta'_i=in_{\fm'}(\theta'_i),\; W=in_{\fm'}(w),\;U'_j=in_{\fm'}(u'_j)$,
is a subring of $\gr_{\fm'}(R')$. It also implies
$$
\NBdxd:=\underset{i\in \NIdxd}{\bigcup} \BBd i
=E\; \cup \; \underset{i\in \Xi}{\bigcup} \BBud i\cup \;
\underset{1\leq j\leq t}{\bigcup} \BBvd {i+s},
$$
where $E=\pi_Z^{-1}(D)$ and $E\times_Z\Spec(R)=\Spec(R/\langle w \rangle)$.
Note
$$
\text{$T_{x'}(E)\subset T_{x'}(Z')=\Spec(\gr_{\fm'}(R'))$ is defined by
$\langle W \rangle\subset \gr_{\fm'}(R')$.}
$$
Moreover
$$
\begin{aligned}
&\text{$T_{x'}(\BBthd j)\subset T_{x'}(Z')$
is defined by $\langle \Theta_j \rangle\subset\gr_{\fm'}(R')$.}\\
&\text{$T_{x'}(\BBud i)\subset T_{x'}(Z')$
is defined by $\langle U'_i \rangle\subset\gr_{\fm'}(R')$ for $i\in \Xi$.}\\
&\text{$T_{x'}(\BBvd j)\subset T_{x'}(Z')$
is defined by $\langle V_j \rangle\subset\gr_{\fm'}(R')$.}\\
\end{aligned}
$$
On the other hand, by Theorem \ref{thm3.1} the assumption that $x'$ is
very near to $x$ implies that there exist $\lambda_1,\dots,\lambda_r\in k'$
such that
\begin{equation*}
\IDir_{x'}(X')=\langle Y'_1+\lambda_1 W,\dots,Y'_r+\lambda_r W \rangle
\;\subset \grmd(R').
\end{equation*}
so that
\begin{equation}\label{thm.NBncDir.claim.eq1}
I\Dirob_{x'}(X') \supset
\langle
Y'_1+\lambda_1 W,\dots,Y'_r+\lambda_r W,\Theta_1,\dots,\Theta_q \rangle.
\end{equation}
This clearly implies the assertion of (1).
If $x'$ is very $\OB$-near, the inclusion in \eqref{thm.NBncDir.claim.eq1}
is equality and then it implies
$\NBdxd\pitchfork \Dirob_{x'}(X')$.
Thus the proof of Theorem \ref{thm.NBncDir} is complete.
$\square$
\bigskip

\begin{corollary}\label{cor.NBncDir}
Let $\pi_Z:Z'=\Bl_D(Z)\to Z$ be as in \eqref{eq.Bpermbu1} and take closed points
$x\in D$ and $x'\in \pi_X^{-1}(x)$
such that $x'$ is $\OB$-near to $x$.
Assume $\Char(k(x))=0$, or $\Char(k(x))\geq \dim(X)/2+1$.
Assume further that there is an integer $e\geq 0$ for which the following
hold:
\begin{itemize}
\item[(1)]
$e_x(X)_{k(x')}\leq e$, and either $e\leq 2$ or $k(x')$ is separable
over $k(x)$.
\item[(2)]
$\NBx\pitchfork \Dirob_x(X)$ or $\eob_x(X)\leq e-1$.
\end{itemize}
Then $\NBdxd\pitchfork \Dirob_{x'}(X')$ or $\eob_{x'}(X')\leq e-1$.
\end{corollary}
\textbf{Proof }
We claim

\medskip
{\it If $\eob_{x'}(X')\geq e$, then $e_x(X)=e=e_{x'}(X')$ and
$\eob_x(X)=e=\eob_{x'}(X')$, so that $x'$ is very $\OB$-near to $x$.}

\medskip
First we show the first equality which implies $x'$ is very near to $x$.
Indeed the assumption implies by Lemma \ref{dimDir} and
Theorem \ref{thm.pbu.inv}
$$
e \leq \eob_{x'}(X')\leq e_{x'}(X') \leq e_x(X)_{k(x')}\leq e.
$$
Hence $e_{x'}(X') = e_x(X)_{k(x')} = e$.
It remains to show $e_x(X)=e$. If $k(x')$ is separable over $k(x)$,
this follows from Lemma \ref{dimDir} (2).
Assume $e\leq 2$ and $e_x(X)<e_x(X)_{k(x')}=2$.
Then Theorem \ref{thmdirectrix} implies that $k(x')=k(x)$ so that
$e_x(X)=e_x(X)_{k(x')}$, which is a contradiction.
Since $x'$ is very near to $x$, Theorem \ref{thm.NBncDir} (1) implies
$$
e\leq \eob_{x'}(X)\leq \eob_x(X)\leq e_x(X)\leq e,
$$
which shows the second equality and the claim is proved.
\medbreak

By the claim, if $\eob_x(X)\leq e-1$, we must have $\eob_{x'}(X')\leq e-1$.
Hence it suffices to show $\NBdxd\pitchfork \Dirob_{x'}(X')$ assuming
$\NBx\pitchfork \Dirob_x(X)$ and $\eob_{x'}(X')\geq e$.
By the claim the second assumption implies that $x'$ is very $\OB$-near
to $x$. Therefore the assertion follows from Theorem \ref{thm.NBncDir} (2).
$\square$
\bigskip

\begin{definition}\label{def.adm.bnd}
Call $(\cB,\OB)$ admissible at $x\in X$, if $N(x)\pitchfork T_x(X)$, and call $(\cB,\OB)$ admissible if it is admissible at all $x\in X$.
\end{definition}

We note that admissibility of $(\cB,\OB)$ at $x$ implies
$\cBI(x) \subseteq \OB(x)$, where $\cBI(x)$ is defined as follows.

\begin{definition}\label{def.iness.bnd}
Call $B\in \cB$ inessential at $x\in X$, if it contains all
irreducible components of $X$ which contain $x$. Let
$$
\cBI(x) = \{ B \in \cB \mid Z \subseteq B \mbox{ for all } Z \in I(x) \}
$$
be the set of inessential boundary components at $x$, where $I(x)$ is the set
of the irreducible components of $X$ containing $x$.
\end{definition}

\begin{definition}\label{def.tnureg}
Call $x\in X$ $\OB$-regular (or $X$ $\OB$-regular at $x$), if
\begin{equation}\label{eq.Oreg}
\tHSfX x = (\nu_X^{reg},|\cBI(x)|)\,,
\end{equation}
where $\nu_X^{reg}$ is as in Remark \ref{HSregular}.
Call $X$ $\OB$-regular, if it is $\OB$-regular at all $x\in X$.
\end{definition}

\begin{lemma}\label{lem.tnureg}
If $(\cB,\OB)$ is admissible at $x$ and $X$ is $\OB$-regular at $x$, then $X$ is regular and normal crossing with $\cB$ at $x$.
\end{lemma}

\textbf{Proof } The first claim follows from Lemma \ref{lem.HS.bound3}. Since $x$ is regular,
the assumption $N(x) \pitchfork T_x(X)$ means that $N(x)$ is transversal to $X$ at $x$. On the other hand,
for the (unique) connected component $W$ on which $x$ lies we have $W \subset B$ for all
$B \in \cBI(x) = \OB(x)$ where the last equality holds by assumption. Thus $X$ is n.c. with $\cB$ at $x$.
$\square$

\bigskip
We have the following transition property.

\begin{lemma}\label{lem.NBncTx}
Let $\pi_X: X' = \Bl_D(X) \rightarrow X$ be as in \eqref{eq.Bpermbu1}.
If $(\cB,\OB)$ admissible at $x\in X$, then $(\cB',\OB')$ and $(\tcB,\tOB)$ are admissible at
any $x'\in \pi_X^{-1}(x')$.
\end{lemma}
\textbf{Proof }
The proof is somewhat similar to that of Theorem \ref{thm.NBncDir}:
If $x'$ is not near to $x$, then $N'(x')$ is empty by definition. Therefore we
may consider the case where $x'$ is near to $x$.
Look at the surjection $R = \cO_{Z,x} \rightarrow \cO_{X,x}$ with kernel $J$, and let
$R = \cO_{Z,x} \rightarrow \cO_{X,x}$ be the corresponding surjection for the local rings of the
blowups at $x'$, with kernel $J'$. Then there is a system of regular parameters
for $R$
$$
(f_1,\ldots,f_m,u_1,\ldots,u_{a+b},v_1,\ldots,v_{r+s})
$$
satisfying the following conditions:
\begin{itemize}
\item[(i)]
$J$ has a standard basis $(f_1,\ldots,f_m,f_{m+1},\ldots,f_n)$ with $f_1,\ldots,f_m \in \fn-\fn^2$ and $f_{m+1},\ldots,f_n\in \fn^2$
for the maximal ideal $\fn\subset R$, so that the initial forms of $f_1,\ldots,f_m$ define $T_x(X)$
inside $T_x(Z)$.
\item[(ii)]
$D\times_Z\Spec(R) = \Spec(R/\langle f_1,\ldots,f_m,u_1,\ldots,u_{a+b}\rangle)$.
\item[(iii)]
$N(x)\times_Z\Spec(R)$ is given by $\div(u_{a+1}),\ldots,\div(u_{a+b}),\div(v_{s+1}),\ldots,\div(v_{s+t})$.
\end{itemize}
Let
$$
\Xi=\{i\in [a+1,a+b]\;|\; x'\in \div(u_i)' \}.
$$
Then there exists
$i_0 \in [1, a+b] -\Xi$ such that
$$
(f'_1,\dots,f'_m, w , u'_j\;(j\in \Xi),v_1,\dots,v_{s+t}),
$$
where $w=u_{i_0}$, $f'_i=f_i/w$, $u'_j=u_j/w$,
is a part of a system of regular parameters of $R'$. Since $x'$ is near to $x$,
we have $H_{\cO_{X',x'}}^{(\delta)}= H_{\cO_{X,x}}^{(0)}$ by Theorem \ref{thm.pbu.inv},
where $\delta=\trdeg_{k(x)}(k(x')))$.
Evaluating at $1$, we get $\dim T_{x'}(X')+\delta = \dim T_x(X)$.
Similarly we get $\dim T_{x'}(Z')+\delta= \dim T_x(Z)$, and hence
$\dim T_{x'}(Z') - \dim T_{x'}(X')= \dim T_x(Z)- \dim T_x(X)$.
It follows that the initial forms of $f'_1,\ldots,f'_m$ already define $T_{x'}(X')$ inside $T_{x'}(Z')$.
This shows that $N'(x') \pitchfork T_{x'}(X')$, because $N'(x')$ is defined
by $w, u'_j \; (j \in \Xi), v_{s+1},\ldots,v_{s+t}$.
$\square$

\bigskip
Finally we give the following functoriality for the objects we have introduced above.

\begin{lemma}\label{lem.Operm}
Consider a cartesian diagram
$$
\begin{CD}
X^\ast @>{i^\ast}>> Z^\ast \\
@V{\varphi}VV @VV{\phi}V \\
X @>{i}>> Z\,,
\end{CD}
$$
in which $Z$ is a regular scheme, $i$ (and hence $i^\ast$) is a closed immersion,
and $\phi$ (and hence $\varphi$) is a flat morphism with regular fibers. Let $\bB \subset Z$ be a simple normal
crossing divisor, let $O$ be a history function for $\bB$ (i.e., for the associated boundary $\cB$),
and let $\bB^\ast=\bB\times_ZZ^\ast$.
\begin{itemize}
\item[(1)]
$Z^\ast$ is regular, and $\bB^\ast$ is a simple normal crossing divisor on $Z^\ast$.
\item[(2)]
Let $\cB^\ast$ be the boundary associated to $\bB^\ast$, let $x^\ast\in X^\ast$ and $x=\varphi(x^\ast)$.
The map $\phi_x^{x^\ast}: \cB(x) \rightarrow \cB^\ast(x^\ast)$ which maps
$B\in \cB(x)$ to $\phi^{-1}(B)\cap\bB^\ast(x^\ast)$, the unique component of $\cB^\ast$ containing $x^\ast$,
is a bijection; its inverse maps $B^\ast\in\cB^\ast(x^\ast)$
to $\phi(B^\ast)$ (with the reduced subscheme structure).
\item[(3)]
The function
$$
\OB^\ast:= \phi^{-1}\OB: X^\ast \rightarrow \{ \,\mbox{subsets of}\, \cB^\ast\}\quad;\quad x^\ast \mapsto \phi_x^{x^\ast}(O(\varphi(x^\ast)))
$$
is a history function for $\cB^\ast$ on $X^\ast$, and one has
\begin{equation}\label{eq.Operm}
H_{X^\ast}^{\OB^\ast}(x^\ast) = {\tHSfX x}
\end{equation}
for any $x^\ast \in X^\ast$ and $x=\varphi(x^\ast)\in X$. In particular, for any
$\tnu \in \tSigmaX$ we have
\begin{equation}\label{eq.Operm2}
\phi^{-1}(\tHS X {\tnu}) = \tHS {X^\ast} {\tnu} \cong {\HS X \tnu}\times_XX^\ast\,.
\end{equation}
\item[(4)]
Let $D$ be a closed subscheme of $X$,
and let $D^\ast = D\times_XX^\ast = D\times_ZZ^\ast$, regarded as a closed subscheme of $X^\ast$.
Let $x^\ast\in D^\ast$ and $x=\varphi(x^\ast)\in D$. Then $D^\ast$ is transversal with $\cB^\ast$ (resp. normal crossing
with $\cB^\ast$, resp. $\cB^\ast$-permissible, resp. $\OB^\ast$-permissible) at $x^\ast$ if and only if
$D$ is transversal with $\cB$ (resp. normal crossing with $\cB$, resp. $\cB$-permissible,
resp. $\OB$-permissible) at $x$.
\item[(5)]
There are unique morphisms $\varphi'$ and $\phi'$ making the diagrams
$$
\begin{aligned}
\begin{CD}
(X^\ast)' = \Bl_{D^\ast}(X^\ast) @>{\varphi'}>> \Bl_D(X) = X' \\
@V{\pi_{X^\ast}}VV @VV{\pi_X}V \\
X^\ast @>{\varphi}>> X
\end{CD}
\end{aligned}
\quad\quad
\begin{aligned}
\begin{CD}
(Z^\ast)' = \Bl_{D^\ast}(Z^\ast) @>{\phi'}>> \Bl_D(Z) = Z' \\
@V{\pi_{Z^\ast}}VV @VV{\pi_Z}V \\
Z^\ast @>{\phi}>> Z
\end{CD}
\end{aligned}
$$
commutative. Moreover, the diagrams are cartesian, and the morphism $i:X \hookrightarrow Z$ induces a
morphism between the diagrams.
\item[(6)]
The diagrams in (5) identify
$(\bB^\ast)^\sim$ with $(\bB^\sim)\times_{Z'}(Z^\ast)'$ and $(\bB^\ast)'$ with
$\bB'\times_{Z'}(Z^\ast)'$, as well as
$(\OB^\ast)^\sim$ with $(\OB^\sim)^\ast:=(\phi')^{-1}(\OB^\sim)$ and $(\OB^\ast)'$ with $(\OB')^\ast:=(\phi')^{-1}(\OB')$.
\end{itemize}
\end{lemma}

\textbf{Proof } (1): In view of the remarks after Definition \ref{def.ncbnd},
it suffices to show: If $Y \subset Z$ is a regular closed subscheme which is of
pure codimension $c$, then $Y^\ast:=Y\times_ZZ^\ast \subset Z^\ast$ is regular and of pure codimension $c$
as well. But the first property is clear from Lemma \ref{completion.HS} (1), and the second one
follows from Lemma \ref{lem.flat}: If $\eta^\ast$ is a generic point of $Y^\ast$, then its image $\eta$
in $Z$ is a generic point of $Y$, and the codimension of $\eta^\ast$ in the fiber over $\eta$
is zero. The latter fiber is the same for $Z^\ast \rightarrow Z$ and $Y^\ast \rightarrow Z$.
Thus Lemma \ref{lem.flat} implies that $\codim_{Z^\ast}(\eta^\ast) = \codim_Z(\eta)$.

(2): In particular, for $B \in \cBx$, its preimage $\phi^{-1}(B)= B\times_ZZ^\ast$ is regular.
Since $\bB^\ast \rightarrow \bB$ is flat, it follows as well that each generic point
of $\bB^\ast$ maps to a generic point of $\bB$. Thus $\phi^{-1}(B)$ is the disjoint union
of those irreducible components of $\bB^\ast$ whose generic points lie above $B$.
There is exactly one component which contains $x^\ast$; this is $\phi_x^{x^\ast}(B)$.

(3): Condition (O1) of Definition \ref{def.histfcnZ} holds by construction. As for (O2)
let $x^\ast,y^\ast\in X^\ast$ with $x^\ast\in\overline{\{y^\ast\}}$ and $H_{X^\ast}(x^\ast) = H_{X^\ast}(y^\ast)$.
Let $x=\varphi(x^\ast)$ and $y=\varphi(y^\ast)$. Then $x\in\overline{\{y\}}$, and $\HSfX x=\HSfX y$
by Lemma \ref{completion.HS} (1). So we have $\OB(y) \subset \OB(x)$ by property (O2) for $\OB$.
If $B^\ast\in \OB^\ast(y^\ast)$, then $y^\ast\in B^\ast$ and $\phi(B^\ast) \in \OB(y) \subset \OB(x)$.
Hence $x^\ast\in B^\ast$ and $B^\ast\in \OB^\ast(x^\ast)$, which shows (O2) for $\OB^\ast$.
Now we show (O3) for $\OB^\ast$.
If $y^\ast\in X^\ast$, and $y=\varphi(y^\ast)$, there exists an open subset $U \subset \overline{\{y\}}$
such that $\OB(x) = \OB(y)$ for all $x\in U$ with $\HSfX x=\HSfX y$. Now $\varphi$ maps
$\overline{\{y^\ast\}}$ into $\overline{\{y\}}$, and we let $U^\ast$ be the preimage of $U$ in
$\overline{\{y^\ast\}}$. Now let $x^\ast\in U^\ast$ with $H_{X^\ast}(x^\ast) = H_{X^\ast}(y^\ast)$. Then $x=\varphi(x^\ast)\in U$,
and $\HSfX x=\HSfX y$ so that $\OB(x) = \OB(y)$. By the above we already know that
$\OB^\ast(y^\ast) \subset \OB^\ast(x^\ast)$. On the other hand, if $B^\ast\in \OB^\ast(x^\ast)$,
and $B=\phi(B^\ast)$, then $\phi^{-1}(B)$ is a disjoint union of irreducible components, and
$B^\ast$ is the unique component containing $x^\ast$. Since $x\in U$, we have $\OB(x)=\OB(y)$,
so that $y\in B$. Hence there is a unique component $B^{\ast\ast}$ of $\pi^{-1}(B)$ containing $y^\ast$.
Since $x^\ast\in\overline{\{y^\ast\}}$, it must equal to $B^\ast$. Hence $B^\ast\in\OB^\ast(y^\ast)$,
and we have shown $\OB^\ast(x^\ast) = \OB^\ast(y^\ast)$ and hence (O3) for $\OB^\ast$.
The equation \eqref{eq.Operm} follows from Lemma \ref{completion.HS} (1) and the
bijection between $\OB^\ast(x^\ast)$ and $\OB(x)$, and the equality in \eqref{eq.Operm2}
follows from this. The isomorphism in \eqref{eq.Operm2} follows as in the proof of
\eqref{eq.HS.smooth2}.

(4): This follows via the same arguments as used for (1).

(5): The first claim for $X$ follows from the universal property of blow-ups.
In fact, if $I \subset \cO_X$ is the ideal sheaf of $D \subset X$, then $\varphi^\ast(I)$ is
the ideal sheaf of $D^\ast \subset X^\ast$, and this coincides with the image ideal sheaf
$\varphi^{-1}\cO_{X^\ast}$ by the flatness of $\varphi$. On the other hand, since the ideal
sheaf of $D^\ast$ is $\varphi^{-1}\cO_{X^\ast}$, it follows from \cite{EGAII} (3.5.3)
that the left diagram is cartesian. Letting $X = Z$ we get the same for $\phi'$.
Finally, the universal property of blow-ups for the closed immersions $X \hookrightarrow Z$
and $X^\ast \hookrightarrow Z^\ast$ give the last claim.

(6): The first claim follows by by applying (5) to $\bB^\ast \rightarrow \bB$ and $D\times_Z\bB \subset \bB$
(for $(\bB)^\sim=\Bl^{D\times_\bB}(\bB)$), and by taking the base change of the first diagram in (5) with $D \hookrightarrow X$
(to treat the exceptional divisor in $(\bB)'$).
For the remaining claim let $y\in (X^\ast)'$ with images $x^\ast, x'$ and $x$ in $X^\ast, X'$ and $X$, respectively. Then
$H_{X^\ast}^{\OB^\ast}(x^\ast) = {\tHSfX x}$ and $H_{(X')^\ast}^{(\OB')^\ast}(y) = {H_{X'}(x')}$
by (3), and $\phi$ and $\phi'$ induce bijections $\OB^\ast(x^\ast) \cong \OB(x)$ and $(\OB')^\ast(y)\cong \OB'(x')$.
The second claim follows from this and the first claim, which also gives a bijection
$\widetilde{\OB^\ast(x^\ast)}\cong \widetilde{\OB(x)}^\ast$.

\newpage
\section{$B$-Permissible blow-ups - the non-embedded case}\label{BpbuII}
\bigskip

Let $X$ be a locally noetherian scheme. We start with the following
definition.

\begin{definition}\label{def.bndX}
A boundary on $X$ is a multiset $\cB = \ml B_1,\ldots,B_r \mr$ of locally
principal closed subschemes of $X$.
\end{definition}

Recall that multisets are `sets with multiplicities'; more precisely a multiset
of $r$ elements is an $r$-tuple in which one forgets the ordering. One can also think of sets
in which an element can appear several times. This then makes clear how one can define
elements, cardinalities, inclusions, intersections and unions of multisets.
Note also that the locally principal subschemes need not be divisors; e.g.,
they could be $X$ itself. Both this and
the use of multisets is convenient for questions of functoriality, see below.

\medskip
In the following, let $X$ be a locally noetherian scheme and let $\cB=
\ml B_1,\ldots,B_n \mr$ be a boundary on $X$. Sometimes, we also call
the elements of $\cB$ components of $\cB$, although they are neither
irreducible nor connected in general. For each $x\in X$, let
$\cBx \subset \cB$ be the submultiset given by the components containing $x$. We
note that this definition is compatible with arbitrary
localization in $X$. For any morphism $f: Y  \rightarrow X$ we
have a pull-back
\begin{equation}\label{eq.pullback.bnd}
\cB_Y := \cB\times_XY := \ml B_Y := B\times_XY \mid B \in \cB \mr \,.
\end{equation}
We also write $f^{-1}(\cB)$. Note that, even if we start with a true set of locally principal divisors on $X$,
the pull-back will be a multiset if there are $B_i \neq B_j$ in $\cB$ with $(B_i)_Y = (B_j)_Y$.
Also, some $(B_i)_Y$ might not be a divisor.
For $x\in X$ we let $\cB_x= f^{-1}(\cB)$ with $f:\Spec(\cO_{X,x})\to X$,
which is a boundary on $X_x=\Spec(\cO_{X,x})$.

\begin{definition}\label{def.Dncbnd}
Let $D\subset X$ be a regular subscheme and let $x\in D$. We say
$D$ is transversal with $\cB$ at $x$ if for each submultiset
$\ml B_{i_1},\ldots,B_{i_r}\mr \subseteq \cBx$, the scheme-theoretic
intersection $D\times_X B_{i_1}\times_X B_{i_2}\times_X\ldots
\times_X B_{i_r}$ is regular and of codimension $r$ in $D$ at $x$.
(So this can only hold if $\cBx$ is a true set.) We say
$D$ is normal crossing (n.c.) with $\cB$ at $x$ if $D$ is transversal
with $\cBx - \cBxD$ at $x$ where
$$
\cBxD = \{ B \in \cBx \mid D \subset B \}\,.
$$
We say that $D$ is transversal (resp. normal crossing) with $\cB$, if
$D$ is transversal (resp. n.c.) with $\cB$ at every $x \in D$.
\end{definition}

\begin{remark}\label{rem.Dncbnd}
Obviously, $D$ is transversal (resp. normal crossing) with $\cB$ at $x$ if and
only if the pull-back $\cB(x)_D$ (resp. $(\cBx - \cBxD)_D$) is a true set
and a simple normal crossings boundary on the regular scheme $D$ (see
Definition \ref{def.ncbnd}), i.e., defines a divisor with normal crossings
on $D$.
\end{remark}

\begin{definition}\label{def.BpermX}
Let $D\subset X$ be a regular closed subscheme and $x\in D$. We
say $D\subset X$ is $\cB$-permissible at $x$ if $D\subset X$ is
permissible at $x$ and $D$ is n.c. with $\cB$ at $x$. We say
$D\subset X$ is $\cB$-permissible if $D\subset X$ is
$\cB$-permissible at all $x\in D$.
\end{definition}

\medbreak Let $D\subset X$ be any locally integral closed subscheme and let $B$ be a
locally principal (closed) subscheme of $X$. We now define a canonical locally principal subscheme
$B'$ on $X'= \Bl_D(X) \rmapo{\pi_X} X$, the blow-up of $X$ in $D$.
Locally we have $X=\Spec(A)$ for a ring $A$,
$D$ is given by a prime ideal $\fp$, and $B$ is given by a principal ideal $fA$, $f\in A$.
In this situation, $\Bl_D(X) = \Proj(A(\fp))$ for the graded $A$-algebra
$A(\fp)  = \oplus_{n\geq 0} \fp^n$. Define the homogenous element
\begin{equation}\label{def.eq.Bd}
f^h =\left.\left\{\gathered
f \in A(\fp)_0 = A \\
f \in A(\fp)_1 = \fp \\
\endgathered\right.\quad
\begin{aligned}
&\text{if $f \notin \fp$}\\
&\text{if $f \in \fp$},
\end{aligned}\right.
\end{equation}
Then the graded principal ideal $A(\fp,B):= f^hA(\fp)$ only depends on $B$
and not on the equation $f$, because it does not change
if $f$ is multiplied by a unit. Thus
\begin{equation}\label{def.Bd}
B' := \Proj(A(\fp)/A(\fp,B)) \subset \Proj(A(\fp)) = X'
\end{equation}
gives a well-defined locally principal subscheme, which is a divisor if $B$ is.

\medbreak
To show that the definition glues on a general $X$ and gives a well-defined locally principal subscheme $B'$
on $X'$, we have to show that the above is compatible with localization. So let $g\in A$ and consider
$\fp_g \subset A_g$. If $f\in \fp$ then $f \in \fp_g$ and there is nothing to show. The same holds
if $f \notin \fp_g$, so that $f \notin \fp$. So we have to consider the case where $f \notin \fp$
but $f\in \fp_g$. Then $f = a/g^m$ for some $a \in \fp$ and some $m > 0$, so that $g^nf \in \fp$
for some $n > 0$, and hence $g\in \fp$ since $f\notin \fp$ and $\fp$ is a prime ideal.
Consider any $h\in\fp$ and the associated chart $D_+(h) \subset \Proj(A(\fp))$,
which is the spectrum of the subring $A(\fp)_{(h)} \subset A_h$ which is the union of the
sets
$$
\frac{\fp^n}{h^n}
$$
for all $n\geq 0$. By definition (and the assumption $f \notin \fp$), the trace of $B'$ on this
chart is defined by the ideal given by the union of the sets
$$
\frac{f\fp^n}{h^n},
$$
and the pull-back to $A_g$ is defined by the union of the sets
\begin{equation}\label{eq.eq1}
\frac{f\fp^n}{h^ng^s}
\end{equation}
for all $n, s \geq 0$. On the other hand, the pull-back of $D_+(h)$ to $A_g$ is also the spectrum of
the ring which is the union of the sets
$$
\frac{\fp_g^n}{h^n}
$$
for $n>0$, and by definition \eqref{def.eq.Bd} (and the assumption $f \in \fp_g$), the ideal
of $B'$ for this ring is the union of the sets
\begin{equation}\label{eq.eq2}
\frac{f\fp_g^{n-1}}{h^n} = \bigcup_{s\geq 0}\,\frac{f\fp^{n-1}}{h^ng^s} = \bigcup_{s\geq 0}\,\frac{fg\fp^{n-1}}{h^ng^{s+1}}
\end{equation}
for all $n\geq 0$. But the union of the sets \eqref{eq.eq1} and \eqref{eq.eq2} is the same, because $g\in \fp$.

\begin{definition}\label{def.princ.transf}
The locally principal subscheme $B'$ defined above is called the principal strict transform of $B$ in $X'$.
\end{definition}

In the following, we will always use these principal strict transforms and will call them
simply transforms.

\begin{remark}\label{rem.transf.div}
(a) There is always a commutative diagram of natural proper morphisms
$$
\xymatrix{\tB \ar[dr]_{\pi_{\tB}} \ar[rr]^{i} & & B' \ar[dl]^{\pi_{B'}}\\& B}
$$
where $\tB = \Bl_{D\times_XB}(B)$ is the scheme-theoretic strict transform of $B$ in $X'$.
All morphisms are isomorphisms over $B\backslash D$, and $i: \tB \rightarrow B'$ is a closed immersion.
However, it is not in general an isomorphism, and $\tB$ need not be a
locally principal subscheme. In fact, with the notations above, $\tB$ is locally
given by the graded ideal
\begin{equation}\label{eq.cAB}
\widetilde{A(\fp,B)} := \oplus_n\, fA\cap\fp^n \quad\supset\quad A(\fp,B) =
\left.\left\{\gathered
\oplus_n\,f\fp^n \\
\oplus_n \, f\fp^{n-1} \\
\endgathered\right.\quad
\begin{aligned}
&\text{if $f \notin \fp$}\\
&\text{if $f \in \fp$},
\end{aligned}\right.
\end{equation}
and the indicated inclusion need not be an equality (or give an isomorphism after
taking $\Proj$).
\smallbreak
(b) If $X$ and $D$ are regular, and $B$ is a regular divisor, then $\tB = B'$. In fact,
with the notations above, $\fp$ is a regular prime ideal,
and locally we have $fA\cap\fp^n = f\fp^{n-v_\fp(f)}$ for the discrete valuation $v_\fp$ associated
to $\fp$, because $v_\fp(fa) = v_\fp(f)+v_\fp(a)$.
Moreover, $v_\fp(f)\in \{0,1\}$ by assumption.
\smallbreak
(c) If $i: X \hookrightarrow Z$ is a closed immersion into a regular scheme $Z$ and $\bB$ is a
simple normal crossings divisor on $Z$, then the set $\cB = C(\bB) = \{B_1,\ldots, B_r\}$ of
irreducible components is a simple normal crossings boundary on $Z$. In particular, it is a boundary
in the sense of Definition \ref{def.bndX}, and $\cB_X = i^{-1}(\cB)$,
its pull-back to $X$, is a boundary on $X$ (which may be a multiset).
This construction connects the present section with the previous one.
(See also Lemma \ref{lem.comp.emb-nemb} below.)
\end{remark}

\medbreak
Now let $\cB = \ml B_1,\ldots,B_n \mr$ be a boundary on $X$.
Let $B_i'$ be the (principal strict) transform of $B_i$ in $X'$, $i =
1,\ldots,n$, and let $E = D \times_X X'$ be the exceptional
divisor.
\begin{definition}\label{def.transf.bnd}
Call $\widetilde{\cB} = \ml B_1',\ldots,B_n' \mr$ the
strict transform and $ \cB' = \ml B_1',\ldots, B_n', E \mr$
the complete transform of $\cB$ in $X'$.
\end{definition}

We note that $E$ is always a locally principal divisor, so that
$\widetilde{\cB}$ and $\cB'$ are boundaries on $X$. Moreover, we note the
following useful functorialities.

\begin{lemma}\label{lem.funct.transf.bnd}
(a) Let $Y\hookrightarrow X$ be a closed immersion, and assume that $D\subset Y$
is a nowhere dense, locally integral closed subscheme. Then for the closed immersion
$Y'=\Bl_D(Y)\hookrightarrow \Bl_D(X) = X'$ one has
$$
(\cB_Y)' = (\cB')_{Y'} \qaq \widetilde{\cB_Y} = (\widetilde{\cB})_{Y'}\,.
$$

\smallskip
(b) Let $\varphi: X^\ast \rightarrow X$ be a flat morphism with regular fibers,
let $D \subset X$ be a closed subscheme,
and let $D^\ast = D\times_XX^\ast$, regarded as a closed subscheme of $X^\ast$.
Then $D^\ast$ is regular if and only $D$ is regular, and in this case one has
$$
(\cB_{X^\ast})' = (\cB')_{(X^\ast)'} \qaq \widetilde{\cB_{X^\ast}} = (\widetilde{\cB})_{(X^\ast)'}\,,
$$
where $(X^\ast)'=\Bl_{D^\ast}(X^\ast)\rightarrow \Bl_D(X) = X'$ is the canonical morphism.
\end{lemma}

{\bf Proof } The question is local on $X$, so we may assume that $X = \Spec(A)$ is affine,
and take up the notations of \eqref{def.eq.Bd}.

\smallskip
(a): Here $Y=\Spec(A/\frb)$ for an ideal $\frb\subset\fp\subset A$, and for
$B = \Spec(A/fA) \in \cB$, $B\times_XY\in \cB_Y$ is given by $(fA+\frb)/\frb$. Thus the case distinction in
\eqref{def.eq.Bd} is the same for $f$ and $\fb=f+\frb$, and the claim follows
from the equality $\fb(\fp/\frb)^n = (f\fp^n+\frb)/\frb$ in the first case and
the equality $\fb(\fp/\frb)^{n-1} = (f\fp^{n-1}+\frb)/\frb$ in the second case.
For the exceptional divisor $E_X = D\times_XX'$ on $X'$ one has
$E_X\times_{X'}Y' = D\times_XY' = D\times_YY' = E_Y$, the exceptional divisor on $Y'$.

\smallskip
(b): Since $D^\ast \rightarrow D$ is flat with regular fibers, the first claim follows from
\ref{completion.HS} (1). For the following we may assume that $X^\ast = \Spec(A^\ast)$
for a flat $A$-algebra $A^\ast$. Then for $B = \Spec(A/fA)\in\cB$, $B_{X^\ast}$ equals
$\Spec(A^\ast/fA^\ast)$, and $D^\ast$ is defined by the ideal $\fp A^\ast$. Now $A/\fp$ is
integral and the map $A/\fp  \rightarrow A^\ast\otimes_A A/\fp = A^\ast/\fp A^\ast$ is flat,
hence this map is also injective.
Thus the case distinction in \eqref{def.eq.Bd} is the same for $f$ and $\fp$ on the one hand,
and for $f$ and $\fp A^\ast$ on the other. Hence $A(\fp A^\ast,B_{X^\ast}) = A(\fp,B)\otimes_A A^\ast$,
and hence $(B_{X^\ast})' = (B')_{(X^\ast)'}$ as claimed.
$\square$

\begin{definition}\label{def.histfcn.cB}
(1) A history function for a boundary $\cB$ on $X$ is a function
\begin{equation}\label{eq.OIBX}
\cOB: X \to \{\text{submultisets of $\cB$}\}\;;\; x \to \cOBx,
\end{equation}
which satisfies the following conditions:
\begin{itemize}
\item[(O1)]
For any $x\in X$, $\cOBx\subset \cBx$.

\item[(O2)] For any $x,y\in X$ such that $x\in \overline{\{y\}}$
and $\HSfX x=\HSfX y$, we have $\cOBy\subset \cOBx$.

\item[(O3)] For any $y\in X$, there exists a non-empty open subset
$U\subset \overline{\{y\}}$ such that $\cOBx=\cOBy$ for all $x\in
U$ such that $\HSfX x=\HSfX y$.
\end{itemize}
For such a function, we put for $x\in X$,
$$
\cNBx=\cBx - \cOBx.
$$
A divisor $B \in \cB$ is called old (resp. new) for $\cOB$ at $x$
if it is a component of $\cOBx$ (resp. \cNBx).

\smallskip
(2) A boundary with history on $X$ is a pair $(\cB,\cOB)$, where
$\cB$ is a boundary on $X$ and $\cOB$ is a history function for $\cB$.
\end{definition}

\medbreak
A basic example of a history function for $B$ on $X$ is given by the following:

\begin{lemma}\label{lem.histfcn.cB(x)}
The function $\cOBx =\cBx$ ($x\in X$), is a history function for
$\cB$ on $X$. In fact it satisfies \ref{def.histfcn.cB} (O2) and
(O3) without the condition $\HSfX x=\HSfX y$.

\end{lemma}

\textbf{Proof } Left to the readers.
$\square$
\medbreak

Define a function:
$$
\cHSfXX : X \to \bNNsq \;;\; x\to (\HSfX x, \coXx),
$$
where $\coXx$ is the cardinality of $\cOBx$. We endow $\bNNsq$
with the lexicographic order:
$$
(\nu,\mu)\geq (\nu',\mu')\Leftrightarrow \nu >\nu'\text{ or }
\nu=\nu'\text{ and }\mu\geq \mu'.
$$
The conditions in \ref{def.histfcn.cB} and Theorem \ref{HSf.usc}
immediately imply the following:

\begin{theorem}\label{thm.cBHSf.usc}
Let the assumption be as above.
\begin{itemize}
\item[(1)] If $x\in X$ is a specialization of $y\in X$, then
$\cHSfX x \geq \cHSfX y$. \item[(2)] For any $y\in X$, there is a
dense open subset $U$ of $\overline{\{y\}}$ such that $\cHSfX
y=\cHSfX x$ for any $x\in U$.
\end{itemize}
\end{theorem}

In other words (see Lemma \ref{lem.up.sc}), the function $\cHSfXX$
is upper semi-continuous on $X$. By noetherian induction
Theorem \ref{thm.cBHSf.usc} implies
$$
\cSigmaX:=\{\cHSfX x|\; x\in X\} \subset \bNN\times\bN
$$
is finite. We define $\cSigmaXmax$ to be the set of of the maximal
elements in $\cSigmaX$.

\begin{definition}\label{def.cBHSlocus}
(1) For $\tnu \in \cSigmaX$ we define
$$
X(\tnu) = X^\OB(\tnu) = \{ x \in X \mid \cHSfX x = \tnu \}
$$
and
$$
\tHSgeqa X = X^\OB(\geq \tnu) = \{ x \in X \mid \cHSfX x \geq \tnu \}\,,
$$
and we call
\begin{equation}\label{eq.cBHSlocus}
\cHSmax X =\underset{\tnu\in \cSigmaXmax}{\cup} X(\tnu)\,.
\end{equation}
the $\cOB$-Hilbert-Samuel locus of $X$.

\smallskip
(2) We define
$$
\cDirx(X):=
\Dir_x(X)\cap \underset{B\in \OBx}{\bigcap} T_x(B) \; \subset T_x(Z).
$$
$$
\eob_x(X)=\dim_{k(x)}(\cDirx(X)).
$$
\end{definition}

By Theorem \ref{thm.cBHSf.usc} and Lemma \ref{lem.up.sc},
$X(\tnu)$ is locally closed, with closure contained in $X(\geq
\tnu)$. In particular, $X(\tnu)$ is closed for $\tnu\in
\cSigmaXmax$, the union in Definition \ref{eq.cBHSlocus} is disjoint and
$\cHSmax X$ is a closed subset of $X$. Theorems \ref{nf.thm1} and
\ref{HSf.perm} imply the following:

\begin{theorem}\label{thm.cOBHSf.perm}
Let $D\subset X$ be a regular closed subscheme and $x\in D$.
Then the following conditions are equivalent:
\begin{itemize}
\item[(1)] $D\subset X$ is permissible at $x$ and there is an open
neighborhood $U$ of $x$ in $Z$ such that $D\cap U\subset B$ for every $B \in \cOBx$.
\item[(2)] $\cHSfX x =\cHSfX y$ for any $y\in D$ such that $x$ is
a specialization of $y$.
\end{itemize}
Under the above conditions, we have
\begin{equation}\label{eq.cBDir}
T_x(D)\subset \cDirx(X).
\end{equation}
\end{theorem}

\begin{definition}\label{def.cOperm}
A closed subscheme $D \subset X$ is called $\cOB$-permissible at $x$,
if it satisfies the equivalent conditions in Theorem
\ref{thm.cOBHSf.perm}, and it is called $\cOB$-permissible, if it is
$\cOB$-permissible at every $x\in D$.
\end{definition}

\begin{remark}\label{rem.cOperm}
Note that $D\subset X$ is $\cB$-permissible at $x$ if and only if
$D$ is $\cOB$-permissible and n.c with $\cNBx$ at $x$.
\end{remark}
\bigskip

Let $D\subset X$ be a $\cB$-permissible closed subscheme. Consider
the blow-up
\begin{equation}\label{eq.cBpermbu1}
X' =  \Bl_D(X) \rmapo{\pi_X}  X
\end{equation}
of $X$ in $D$, and let $\cB'$ be the complete transform of $\cB$ in $X'$.
For a given history function $\cOB$ for $\cB$ on
$X$, we define functions $\OBd,\, \tOB \,: X'\to \{\text{subsets of $\cB'$}\}$ as follows:
Let $x'\in X'$ and $x=\pi_X(x')\in X$. Then define
\begin{equation}\label{def.cOBd}
\cOBdxd =\left.\left\{\gathered
\tcOBx\cap \cBdxd \\
 \cBdxd \\
\endgathered\right.\quad
\begin{aligned}
&\text{if $\HSfXd {x'}=\HSfX x$}\\
&\text{otherwise},
\end{aligned}\right.
\end{equation}
where $\tcOBx$ is the strict transform of $\cOBx$ in $X'$ and
\begin{equation}\label{def.tOB}
\tOB(x')=\left.\left\{\gathered
\tOBx\cap\cBdxd \\
 \tcB(x') \\
\endgathered\right.\quad
\begin{aligned}
&\text{if $\HSfXd {x'}=\HSfX x$}\\
&\text{otherwise}.
\end{aligned}\right.
\end{equation}
Note that $\OBdxd = \tOB(x') = \tOBx\cap\tcB(x')$ if $x'$ is near to $x$.

\medbreak
The proof of the following lemma is identical with that of Lemma \ref{lem.Bperm1}.

\begin{lemma}\label{lem.cBperm1}
The function $X'\to \{\text{subsets of $\cB'$}\}\;;\; x' \to
\cOBdxd$ is a history function.
\end{lemma}
\medbreak

\begin{definition}\label{def.ccomplete.transform}
We call $(\cB',\cOB')$ and $(\tcB,\tOB)$ the complete and strict transform
of $(\cB,\cOB)$ in $X'$, respectively.
\end{definition}

\medbreak
Results from the previous section (embedded case) have
their companions in the non-embedded situation.
We start with the following theorem analogous to Theorem \ref{thm.Bpbu.inv}.

\begin{theorem}\label{thm.cBpbu.inv}
Take points $x\in D$ and $x'\in \pi_X^{-1}(x)$.
Then $H^{\tOB}_{X'}(x') \leq H^{O'}_{X'}(x') \leq \tHSfX x$. In particular we have
$$
\pi_X^{-1}(X^{\OBd}(\tnu)) \subset \underset{\tmu\leq\tnu}{\cup} {X'}^{\OBd}(\tmu)
\qfor \tnu\in \SigmaXmax,
$$
and the same holds for $\tOB$ in place of $\OBd$.
\end{theorem}
\textbf{Proof }
This follows immediately from Theorem \ref{thm.pbu.inv}, \eqref{def.cOBd} and \eqref{def.tOB}.
$\square$

\medbreak
In the following we mostly use the complete transform $(\cB',\OB')$ and,
for ease of notation, we often write $\cHSfXd {x'}$ and $\cSigmaXd$ instead of
$H^{O'}_{X'}(x')$ and $\Sigma^{O'}_{X'}$, similarly for $\cSigmaXdmax$  etc., because everything just depends on $\cOB$.

\begin{definition}\label{def.cBnearpoint}
We say that $x'\in \pi_X^{-1}(x)$ is $\cOB$-near to $x$ if the following equivalent conditions
hold:
\begin{itemize}
\item[(1)]
$\cHSfXd {x'} =\cHSfX {x}$\quad \quad($\;\Leftrightarrow\; H^{\OBd}_{X'}(x') =\tHSfX {x} \;\Leftrightarrow \;
H^{\tOB}_{X'}(x') = \tHSfX {x}\;$).
\item[(2)]
$x'$ is near to $x$ and contained in the strict transforms of
$B$ for all $B \in \cOBx$.
\end{itemize}
Call $x'$ is very $\OB$-near to $x$ if $x'$ is $\OB$-near and very near to $x$
and $\eob_{x'}(X')= \eob_x(X)-\delta_{x'/x}$.
\end{definition}

\bigskip
The following result, the non-embedded analogue of Theorem \ref{thmBdirectrix},
is an immediate consequence of Theorem
\ref{thmdirectrix} and Definition \ref{def.cBHSlocus}(2).

\begin{theorem}\label{thm.cBdirectrix}
Assume that $x'\in X'$ is $\cOB$-near to $x= \pi(x') \in X$.
Assume further that $\Char(k(x)) = 0$, or $\Char(k(x)) \geq \dim (X)/2+1$,
where $k(x)$ is the residue field of $x$.
Then
$$
 x'\in {\mathbb P}(\cDirx(X)/T_x(D))\subset
{\mathbb P}(C_x(X)/T_x(D))=\pi_X^{-1}(x)
$$
\end{theorem}

Next we show the compatibility of the notions for the non-embedded case in the present section
with the corresponding notions for the embedded situation in the previous section.

\begin{lemma}\label{lem.comp.emb-nemb}
Let $i: X \hookrightarrow Z$ be an embedding into a regular excellent scheme $Z$,
let $\cB$ be a simple normal crossings boundary on $Z$, and let $\cB_X=i^{-1}(\cB)$
be its pull-back to $X$.

\smallskip
 (1) For a closed regular subscheme $D\subset X$ and $x\in D$, $D$ is transversal (resp. normal
crossing) with $\cB$ at $x$ in the sense of Definition \ref{def.DncBZ} if
and only if it is transversal (resp. normal crossing) with $\cB_X$ in the sense of
Definition \ref{def.Dncbnd} (which is an intrinsic condition on $(X,\cB_X)$).

\smallskip
 (2) Let $\OB$ be a history function for $\cB$ in the sense of Definition \ref{def.histfcnZ}, and define the function
$$
\OB_X: X \rightarrow \{\mbox{submultisets of }\cB_X\} \quad,\quad \OB_X(x) = \ml B_X \mid B \in \OB(x) \mr\,.
$$
Then $\OB_X$ is a history function for $\cB_X$ in the sense of Definition \ref{def.histfcn.cB}, and one has
$$
H^{\OB_X}_X(x) = H^{\OB}_X(x) \qaq \Dir^{\OB_X}_x(X) = \Dirob_x(X)
$$
(and hence $e^{\OB_X}_x(X) = \eob_x(X)$) for all $x\in X$.
Also, for a $k(x)$-linear subspace $T \subset T_x(X)$, the two notions for the transversality
$T \pitchfork N(x)$ (Definition \ref{def.NBncDir} for $\cB$ and Definition \ref{def.cNBncDir} for $\cB_X$)
 are equivalent.

\smallskip
 (3) A regular closed subscheme $D \subset X$ is $\cB_X$-permissible in the sense of Definition \ref{def.BpermX}
 if and only if it is $\cB$-permissible in the sense of Definition \ref{def.Bperm}.
Moreover, it is $\OB_X$-permissible in the sense of Definition \ref{def.cOperm}
if and only if it is $\OB$-permissible in the sense of Definition \ref{def.Operm}.

\smallskip
 (4) Let $D \subset X$ be $\cB$-permissible, let $\pi_X: X' = \Bl_D(X) \rightarrow X$ and $\pi_Z: Z' = \Bl_D(Z) \rightarrow Z$
be the respective blowups in $D$ and $i': X' \hookrightarrow Z'$ the closed immersion. Moreover let
$\OB$ a history function for $\cB$. Then we have the equalities
$$
((\cB_X)',(\OB_X)') = ((\cB')_{X'},(\OB')_{X'}) \qaq (\widetilde{\cB_X},\widetilde{\OB_X}) = ((\tcB)_{X'},(\tOB)_{X'})
$$
for the complete transforms and strict transforms, respectively.
\end{lemma}

{\bf Proof } The claims in (1), (2) and (3) easily follow from the definitions. For the claim on
the directrix in (2) note that $T_x(B_X) = T_x(B)\cap T_x(X)$ (in $T_x(Z)$).
The claim in (4) follows from Lemma \ref{lem.funct.transf.bnd}.
$\square$

\medbreak
Results in the non-embedded case which depend only on the local ring at a point
(of the base scheme) can often be reduced to the embedded case. This relies on
Lemma \ref{lem.comp.emb-nemb} and the following two observations.

\begin{remark}\label{rem.emb-nemb}
Let $X$ be an excellent scheme, let $\cB$ be a boundary on $X$, and let $x\in X$.
Assume a property concerning $(X,\cB,x)$ can be shown by passing to the local ring $\cO = \cO_{X,x}$,
and its completion $\hat{\cO}$.
Then the following construction is useful. The ring $\hat{\cO}$ is the quotient of
a regular excellent ring $R$. Let $\cB(x) = \{\{ B_1,\ldots,B_r\}\}$,
and let $f_1,\ldots,f_r$ be the local functions defining them in $\cO_{X,x}$
(so we can have $f_i = f_j$ for $i\neq j$).
Then we get a surjection
$$
R[X_1,\ldots,X_r] \twoheadrightarrow \cO_{X,x}
$$
mapping $X_i$ to $f_i$, and the functions $X_i$ define a simple normal crossings boundary on
$Z = \Spec(R[X_1,\ldots,X_r])$, such that $\cB(x)$ is its pull-back under $\Spec(\cO_{X,x}) \hookrightarrow
 Z$. We may thus assume that $X$ can be embedded in a regular excellent
scheme $Z$ with simple normal crossings boundary $\cB_Z$, and that $\cB$ is the pull-back of
$\cB_Z$ to $X$.
\end{remark}

\bigskip
Now we apply Remark \ref{rem.emb-nemb} and Lemma \ref{lem.comp.emb-nemb}.

\begin{theorem}\label{nf.cBdir}
Let $D \subseteq X$ be an irreducible subscheme.
Assume:
\begin{itemize}
\item[(1)]
$D\subset X$ is $\cOB$-permissible at $x$.
\item[(2)]
$e_\eta(X) = e_x(X) -\dim(\cO_{D,x})$ (cf. Theorem \ref{nf.dir}),
\end{itemize}
where $\eta$ is the generic point of $D$. Then we have
$$
\eob_\eta(X)\leq \eob_x(X) -\dim(\cO_{D,x}).
$$
\end{theorem}

\textbf{Proof }
The question is local around $x$, and we may pass to $\cO_{X,x}$ and then to its completion,
since $X$ is excellent. By \ref{rem.emb-nemb} and \ref{lem.comp.emb-nemb} we may assume
that we are in an embedded situation.
Thus the claim follows from the corresponding result in the embedded case (Theorem \ref{nf.Bdir}).
$\square$

\medbreak
Let $\pi_X:X'=\Bl_D(X)\to X$ be as in \eqref{eq.cBpermbu1}.
For $x'\in X'$ let $\cNBdxd = \cBdxd - \cOBdxd$ be the set of the new components of $(\cB',\OB')$ for $x'$.
If $x'$ is near to $x = \pi_X(x') \in D$, i.e., $\HSfXd {x'}=\HSfX x$, then
\begin{equation}\label{eq.cNBdxd}
\cNBdxd= (\widetilde{\cNBx}\cap \cBdxd)\, \cup \{E\} \quad \mbox{ with }\quad
E=\pi_Z^{-1}(D)
\end{equation}
where $\widetilde{\cNBx}$ is the strict transform of $\cNBx$ in $X'$. If $x'$ is not near to $x$,
then $\cNBdxd = \emptyset$.
Similarly, $\tN(x') = \tcB(x') - \tOB(x') = \widetilde{\cNBx}\cap\tcB(x') \subset \cNBdxd$ if $x'$
is near to $x$, and $\tN(x') = \emptyset$, otherwise.
We study the transversality of $\cNBdxd$ with a certain regular subscheme of $E$.

\begin{definition}\label{def.cNBncDir}
For a $k(x)$-linear subspace $T\subset T_x(X)$, we say that $T$ is transversal with
$\cNBx$ (notation: $T \pitchfork \cNBx$) if
$$
\dim_{k(x)}\big(T \cap \underset{B \in \cNBx}{\cap} T_x(B)\big)=
\dim_{k(x)}(T)-\mid \cNBx\mid.
$$
\end{definition}

\begin{lemma}\label{lem.cNBncDir1}
Assume $D=x$, $T \subset \Dir_x(X)$ and $T \pitchfork \cNBx$. Then the closed subscheme
$\bP(T)\subset E_x=\bP(C_x(X))$ is n.c. with $\NBdxd$ (and hence also $\tN(x')$) at $x'$.
\end{lemma}

\textbf{Proof }
In the same way as above, we this claim follows from the corresponding result
in the embedded case (Lemma \ref{lem.NBncDir1}).
$\square$

\begin{lemma}\label{lem.cNBncDir2}
Let $\pi_X:X'=\Bl_D(X)\to X$ be as \eqref{eq.cBpermbu1}.
Assume $T \pitchfork \NBx$ and
$T_x(D)\subset T \subset \Dir_x(X)$ and $\dim_{k(x)}(T/T_x(D))=1$. Consider
$$
{x'}=\bP(T/T_x(D))\subset \bP(C_x(X)/T_x(D))=\pi_X^{-1}(x).
$$
Let $D'\subset E = \pi_X^{-1}(D)$ be any closed subscheme such that $x'\in D'$ and $\pi_X$ induces
an isomorphism $D'\simeq D$. Then $D'$ is n.c. with $\cNBdxd$ at $x'$.
\end{lemma}

\medbreak\noindent
\textbf{Proof } In the same way as above, this follows from the
corresponding result in the embedded case (Lemma \ref{lem.NBncDir2}).
$\square$

\begin{theorem}\label{thm.cNBncDir}
Let $\pi_X:X'=\Bl_D(X)\to X$ be as in \eqref{eq.cBpermbu1}.
Take $x\in X$ and $x'\in \pi_X^{-1}(x)$. Assume $\Char(k(x))=0$, or $\Char(k(x))\geq \dim(X)/2+1$.
\begin{itemize}
\item[(1)]
If $x'$ is $\OB$-near and very near to $x$, then
$\eob_{x'}(X')\leq \eob_x(X)-\delta_{x'/x}$.
\item[(2)]
Assume $x'$ is very $\OB$-near and $\NBx\pitchfork \Dirob_x(X)$.
Then $\NBdxd\pitchfork \Dirob_{x'}(X')$.
\end{itemize}
\end{theorem}
\textbf{Proof } Reduction to the embedded case (Theorem \ref{thm.NBncDir}).
$\square$
\bigskip

\begin{corollary}\label{cor.cNBncDir}
Let $\pi_X:X'=\Bl_D(X)\to X$ be as in \eqref{eq.cBpermbu1} and take closed points
$x\in D$ and $x'\in \pi_X^{-1}(x)$ such that $x'$ is $\OB$-near to $x$.
Assume $\Char(k(x))=0$, or $\Char(k(x))\geq \dim(X)/2+1$.
Assume further that there is an integer $e\geq 0$ for which the following
hold:
\begin{itemize}
\item[(1)]
$e_x(X)_{k(x')}\leq e$, and either $e\leq 2$ or $k(x')$ is separable
over $k(x)$.
\item[(2)]
$\NBx\pitchfork \Dirob_x(X)$ or $\eob_x(X)\leq e-1$.
\end{itemize}
Then $\NBdxd\pitchfork \Dirob_{x'}(X')$ or $\eob_{x'}(X')\leq e-1$.
\end{corollary}

This follows from Theorem \ref{thm.cNBncDir} like Corollary \ref{cor.NBncDir}
follows from Theorem \ref{thm.NBncDir}.

\begin{definition}\label{def.cadm.bnd}
Call $(\cB,\OB)$ admissible at $x\in X$, if $N(x)\pitchfork T_x(X)$, and call $(\cB,\OB)$ admissible if it is admissible at all $x\in X$.
\end{definition}

We note that admissibility of $(\cB,\OB)$ at $x$ implies
$\cBI(x) \subseteq \OB(x)$, where $\cBI(x)$ is defined as follows.

\begin{definition}\label{def.ciness.bnd}
Call $B\in \cB$ inessential at $x\in X$, if it contains all
irreducible components of $X$ which contain $x$. Let
$$
\cBI(x) = \{ B \in \cB \mid Z \subseteq B \mbox{ for all } Z \in I(x) \}
$$
be the set of inessential boundary components at $x$, where $I(x)$ is the set
of the irreducible components of $X$ containing $x$.
\end{definition}

\begin{definition}\label{def.ctnureg}
Call $x\in X$ $\OB$-regular (or $X$ $\OB$-regular at $x$), if
\begin{equation}\label{eq.cOreg}
\tHSfX x = (\nu_X^{reg},|\cBI(x)|)\,.
\end{equation}
Call $X$ $\OB$-regular, if it is $\OB$-regular at all $x\in X$.
\end{definition}

\begin{lemma}\label{lem.ctnureg}
If $(\cB,\OB)$ is admissible at $x$ and $X$ is $\OB$-regular at $x$, then $X$ is regular and normal crossing with $\cB$ at $x$.
\end{lemma}
\textbf{Proof } The proof is identical with that of Lemma \ref{lem.tnureg}.
$\square$
\bigskip

\begin{lemma}\label{lem.cNBncTx}
Let $\pi: X' = \Bl_D(X) \rightarrow X$ be as \eqref{eq.cBpermbu1}.
If $(\cB,\OB)$ admissible at $x\in X$, then $(\cB',\OB')$ and $(\tcB,\tOB)$ are admissible at
any $x'\in X'$ with $x=\pi(x')$.
\end{lemma}
\textbf{Proof } In the same way as in the proof of Lemma \ref{nf.cBdir},
this is reduced to the embedded case (Theorem \ref{lem.NBncTx}).
$\square$
\medbreak

Later we shall need the following comparison for a closed immersion.

\begin{lemma}\label{lem.comp.incl}
Let $i: Y \hookrightarrow X$ be a closed immersion of excellent schemes,
let $\cB$ be a boundary on $X$, and let $\cB_Y=i^{-1}(\cB)$
be its pull-back to $Y$.

\smallskip
 (1) For a closed regular subscheme $D\subset Y$ and $x\in D$, $D$ is transversal (resp. normal
crossing) with $\cB$ at $x$ if
and only if it is transversal (resp. normal crossing) with $\cB_Y$.

\smallskip
 (2) Let $\OB$ be a history function for $\cB$, and define the function
$$
\OB_Y: Y \rightarrow \{\mbox{submultisets of }\cB_Y\} \quad,\quad \OB_Y(x) = \ml B_Y \mid B \in \OB(x) \mr\,.
$$
Then $\OB_Y$ is a history function for $\cB_Y$. If $x\in Y$ and $(\cB_Y,\OB_Y)$ is admissible at $x$,
then $(\cB,\OB)$ is admissible at $x$. (The converse does not always hold.)

\smallskip
 (3) Let $D \subset Y$ be a regular closed subscheme which is permissible for $Y$ and $X$.
Then $D$ is $\cB_Y$-permissible if and only if it is $\cB$-permissible.
Moreover, it is $\OB_Y$-permissible if and only if it is $\OB$-permissible.

\smallskip
 (4) Let $D \subset Y$ be $\cB_Y$-permissible and $\cB_X$-permissible,
let $\pi_Y: Y' = \Bl_D(Y) \rightarrow Y$ and $\pi_X: X' = \Bl_D(X) \rightarrow X$
be the respective blowups in $D$ and $i': Y' \hookrightarrow X'$ the closed immersion. Moreover let
$\OB$ a history function for $\cB$. Then we have the equality
$$
((\cB_Y)',(\OB_Y)') = ((\cB')_{Y'},(\OB')_{Y'}) \qaq (\widetilde{\cB_Y},\widetilde{\OB_Y}) = ((\tcB)_{Y'},(\tOB)_{Y'})
$$
for the complete and strict transforms, respectively.
\end{lemma}

The proofs are along the same lines as for Lemma \ref{lem.comp.emb-nemb}.
For (2) note that $T_x(Y) \subset T_x(X)$ and that for subspaces
$T_1 \subset T_2 \subset T_x(X)$ one has $N(x) \pitchfork T_1 \Longrightarrow N(x) \pitchfork T_2$.

\begin{remark}\label{rem.comp.emb-nemb}
Since a regular subscheme of a regular scheme is always permissible, Lemma \ref{lem.comp.emb-nemb}
can be seen as a special case of Lemma \ref{lem.comp.incl}.
\end{remark}

The following is the non-embedded analogue of Lemma \ref{lem.Operm}.

\begin{lemma}\label{lem.Operm.nemb}
Let $\varphi: X^\ast \rightarrow X$ be a flat morphism with regular fibers of locally noetherian schemes,
let $\cB$ be a boundary on $X$, and let $\OB$ be a history function for $\cB$ on $X$.

\smallskip
(1) Let $\cB^\ast = \varphi^{-1}(\cB) = \cB_{X^\ast}$ be the pull-back to $X^\ast$, let $x^\ast \in X^\ast$ and $x = \varphi(x^\ast)$.
The map $\varphi_x^{x^\ast}: \cB(x) \rightarrow \cB^\ast(x^\ast)$ which maps $B\in \cB(x)$ to
$B^\ast = B\times_X X^\ast$ is a bijection.

(2) The function
$$
\OB^\ast:= \varphi^{-1}(\OB): X^\ast \rightarrow \{ \,\mbox{subsets of}\, \cB^\ast\}\quad;\quad x^\ast \mapsto \varphi_{\varphi(x^\ast)}^{x^\ast}(O(\varphi(x^\ast)))
$$
is a history function for $\cB^\ast$ on $X^\ast$, and one has
\begin{equation}\label{eq.nemb.Operm}
H_{X^\ast}^{\OB^\ast}(x^\ast) = {\tHSfX x}
\end{equation}
for any $x^\ast \in X^\ast$ and $x=\varphi(x^\ast)\in X$. In particular, for any
$\tnu \in \tSigmaX$ we have
\begin{equation}\label{eq.nemb.Operm2}
\varphi^{-1}(\tHS X {\tnu}) = \tHS {X^\ast} {\tnu} \cong {\HS X \tnu}\times_XX^\ast\,.
\end{equation}

(3) Let $D$ be a closed subscheme of $X$,
and let $D^\ast = D\times_XX^\ast = D\times_ZZ^\ast$, regarded as a closed subscheme of $X^\ast$.
Let $x^\ast\in D^\ast$ and $x=\varphi(x^\ast)\in D$. Then $D^\ast$ is transversal with $\cB^\ast$ (resp. normal crossing
with $\cB^\ast$, resp. $\cB^\ast$-permissible, resp. $\OB^\ast$-permissible) at $x^\ast$ if and only if
$D$ is transversal with $\cB$ (resp. normal crossing with $\cB$, resp. $\cB$-permissible,
resp. $\OB$-permissible) at $x$.

(4) There is a unique morphism $\varphi'$ making the following diagram
$$
\begin{aligned}
\begin{CD}
(X^\ast)' = \Bl_{D^\ast}(X^\ast) @>{\varphi'}>> \Bl_D(X) = X' \\
@V{\pi_{X^\ast}}VV @VV{\pi_X}V \\
X^\ast @>{\varphi}>> X
\end{CD}
\end{aligned}
$$
commutative, where $\pi_X$ and $\pi_{X^\ast}$ is the structural morphisms.
Moreover the diagram is cartesian, and it identifies
$(\cB^\ast)^\sim$ with $(\cB^\sim)^\ast:=(\varphi')^{-1}(\cB^\sim)$ and
$(\cB^\ast)'$ with $(\cB')^\ast:=(\varphi')^{-1}(\cB')$,
as well as $(\OB^\ast)^\sim$ with $(\OB^\sim)^\ast:=(\varphi')^{-1}(\OB^\sim)$
and $(\OB^\ast)'$ with $(\OB')^\ast:=(\varphi')^{-1}(\OB')$.
\end{lemma}

\textbf{Proof } (1) is trivial, and in (2), the conditions for a history function
are easily checked. Since we have a bijection between $\OB^\ast(x^\ast)$ and $\OB(x)$,
equation \eqref{eq.nemb.Operm} follows from Lemma \ref{completion.HS} (1).
By Remark \ref{rem.Dncbnd}, the first three cases in (3) follow from Lemma \ref{lem.Operm} (1).
The last claim follows from \eqref{eq.nemb.Operm}.

In (4) the first two claims are shown like in the proof of Lemma \ref{lem.Operm} (5).
The last four claims follow by applying the first two claims to every $B \in \cB$.

\newpage
\section{Main theorems and strategy for their proofs}\label{sec:strategy}
\def\OBn{O_n}
\def\tHSmaxXn{(X^O_n)_{\max}}
\def\tHSmaxXm{(X^O_m)_{\max}}
\def\tHSmaxXi{(X^{(i)})^O_{\max}}

\bigskip

We will treat the following two situations in a parallel way:

\medskip
(E)  ({\it embedded case}) $X$ is an excellent noetherian scheme, $i: X \hookrightarrow Z$
is a closed immersion into an excellent regular noetherian scheme $Z$, and $\cB$ is a simple
normal crossings boundary on $Z$ (Definition \ref{def.ncbnd}).

\medskip
(NE) ({\it non-embedded case}) $X$ is an excellent noetherian scheme, and $\cB$ is a boundary
on $X$ (Definition \ref{def.bndX}).

\medskip
In this section, all schemes will be assumed to be noetherian.

\begin{definition}\label{def.B-reg}
(1) A point $x\in X$ is called
$\cB$-regular if $X$ is regular at $x$ (i.e.
$\cO_{X,x}$ is regular) and normal crossing with $\cB$ at $x$.
Call $X$ $\cB$-regular, if every $x\in X$ is
$\cB$-regular, i.e., if $X$ is regular and $\cB$ is normal crossing with $X$.

\smallskip
(2) Call $x$ strongly $\cB$-regular if $X$ is regular at $x$ and for every $B\in \cB(x)$,
$B$ contains the (unique) irreducible component on which $x$ lies. (This amounts to
the equation $\cB(x)=\cBI(x)$ where $\cBI(x)$ is the (multi)set of inessential boundary
components at $x$, see Definitions \ref{def.iness.bnd} (Case (E)) and \ref{def.ciness.bnd}
 (Case (NE)).)
\end{definition}

Denote by $\Xreg$ (resp. $\XBreg$, resp. $\XBsreg$) the set of the regular
(resp. $\cB$-regular, resp. strongly $\cB$-regular) points of $X$.
These are open subsets of $X$, and dense in $X$ if $X$ is reduced.
Call $\XBsing = X - \XBreg$ the $\cB$-singular locus of $X$.

\medskip
We introduce the following definition for the case of non-reduced schemes.

\begin{definition}\label{def.B-qreg}
(1) Call $x\in X$ quasi-regular, if $X_{red}$ is regular at $x$ and $X$ is normally flat
along $X_{red}$ at $x$.
Call $X$ quasi-regular if it is quasi-regular at all $x\in X$, i.e., if $X_{red}$
is regular and $X$ is normally flat along $X_{red}$. (Compare Definition \ref{def.nf},
but we have reserved the name `permissible' for subschemes not containing any
irreducible component of $X$.)

\smallskip
(2) Call $x\in X$ quasi-$\cB$-regular, if $X_{red}$ is $\cB$-regular at $x$ and $X$ is
normally flat along $X_{red}$ at $x$. Call $X$ quasi-$\cB$-regular, if $X$ is quasi-$\cB$-regular at all $x\in X$,
i.e., if $X_{red}$ is $\cB$-regular and $X$ is normally flat along $X_{red}$. (Similar remark
on comparison with $\cB$-permissibility.)

\smallskip
(3) Call $x\in X$ strongly quasi-$\cB$-regular, if $X_{red}$ is strongly $\cB$-regular at $x$
and $X$ is normally flat along $X_{red}$ at $x$.
\end{definition}

Note that $X$ is regular if and only if $X$ is quasi-regular and reduced. Similarly,
$X$ is $\cB$-regular if and only if $X$ is quasi-$\cB$-regular and reduced. Finally,
$x\in X$ is strongly $\cB$-regular if and only if $x$ is strongly quasi-$\cB$-regular
and $\cO_{X,x}$ is reduced.

\smallskip
Denote by $\Xqreg$, $\XBqreg$ and $\XBsqreg$ the sets of quasi-regular, quasi-$\cB$-regular
and strongly quasi-$\cB$-regular points of $X$, respectively. By Theorem \ref{nf.thm1}
these are dense open subsets of $X$. Moreover, we have inclusions
$$
\begin{matrix}
\Xqreg & \supset & \XBqreg & \supset & \XBsqreg & \supset & \Xqreg\setminus(\Xqreg\cap\cB) \\
\cup && \cup && \cup && \cup \\
\Xreg & \supset & \XBreg & \supset & \XBsreg & \supset & \Xreg\setminus(\Xreg\cap\cB)
\end{matrix}
$$
where the last inclusions of both rows are equalities if no $B\in \cB$ contains any irreducible component of $X$
and the vertical inclusions are equalities if $X$ is reduced.

\begin{lemma}\label{lem.qreg}
Let $X$ be a connected excellent scheme.

\smallskip
(a) For $\nu \in \SigmaXmax$ one has $X(\nu)\cap \Xqreg \neq \emptyset$ if and only if
$X = X(\nu)$. Thus $H_X$ is not constant on $X$ if and only if $X(\nu) \subset X-\Xqreg$
for all $\nu \in \SigmaXmax$.

\smallskip
(b) Let $(\cB,\OB)$ be an admissible boundary with history on $X$. For $\tnu \in \tSigmaXmax$
one has $X^\OB(\tnu) \cap \XBsqreg \neq \emptyset$ if and only if $X = X^\OB(\tnu)$.
Thus $H^\OB_X$ is not constant on $X$ if and only if $X^\OB(\tnu) \subseteq X-\XBsqreg$
for all $\tnu \in \tSigmaXmax$.
\end{lemma}

\textbf{Proof } (a): Let $x\in X(\nu)\cap \Xqreg$, let $Z$ be an irreducible component
of $X$ containing $x$, and let $\eta$ be the generic point of $Z$. Since $\Xqreg$ is open and
dense in $X$, and is quasi-regular, $\eta$ is contained in $\Xqreg$ and $H_X$ is constant
on $\Xqreg$ by Theorem \ref{HSf.perm}. Therefore $\nu = H_X(x) = H_X(\eta)$, i.e., $\eta \in X(\nu)$.
Since $\nu \in \SigmaXmax$, $X(\nu)$ is closed, and we conclude that
$Z = \overline{\{\eta\}} \subset X(\nu)$. By Lemma \ref{lem.constant.HSf} we conclude
that $X = X(\nu)=Z$. This proves the first claim (the other direction is trivial).
The second claim is an obvious consequence.
(b): For the non-trivial direction of the first claim let $\tnu = (\nu,m)$, with $\nu \in \bNN$ and $m \geq 0$.
Then $\nu \in \SigmaXmax$ and $X^\OB(\tnu) \subseteq X(\nu)$. Consequently, if $X^\OB(\tnu) \cap \XBsqreg \neq \emptyset$,
then $X=X(\nu)$ by (a), and $X$ is irreducible. If $\eta$ is the generic point of $X$, then we
conclude as above that $\eta \in X^\OB(\tnu)$, and hence $X=X^\OB(\tnu)$, since the latter
set is closed. Again the second claim follows immediately.

\medskip
Now we study blow-ups. Lemma \ref{lem.Bperm0} implies:

\begin{lemma}\label{lem.Breg}
If $\pi: X' = \Bl_D(X) \rightarrow X$ is the blow-up of $X$ in a
$\cB$-permissible subscheme $D$, and $\cB'$ is the complete transform of
$\cB$, then
$\pi^{-1}(\XBsreg) \subset \XdBdsreg$ and $\pi^{-1}(\XBsqreg) \subset X'_{\cB'\sqreg}$.
\end{lemma}

We first consider the case (NE).

\begin{definition}\label{def.seqBpermbu}
(Case (NE))
A sequence of complete (resp. strict) $\cB$-permissible blowups over $X$ is a diagram
\begin{equation}\label{seqBpermbu1}
\begin{array}{ccccccccccccccc}
 \cB=&\cB_0&& \cB_1&& \cB_2&& \cB_{n-1}&& \cB_{n} & \quad \cdots\\
 & & &  &&  & &  &&  \\
  X=&X_0 & \lmapo{\pi_1} &   X_1 & \lmapo{\pi_2} & X_2 &
 \leftarrow\ldots \leftarrow & X_{n-1} & \lmapo{\pi_n} & X_{n}
&\gets\cdots\\
\end{array}
\end{equation}
where for any $n\geq 0$, $\cB_n$ is a boundary on $X_n$, and
$$
X_{n+1}= \Bl_{D_{n}}(X_{n}) \rmapo{\pi_{n+1}} X_{n}
$$
is the blow-up in a $\cB_n$-permissible center $D_{n}\subset X_{n}$,
and $\cB_{n+1}=\cB_n'$ is the complete transform of $\cB_{n}$ (resp. $\cB_{n+1} = \widetilde{\cB_n}$
is the strict transform of $\cB_n$).

Call a sequence as in \eqref{seqBpermbu1} contracted if none of the
morphisms $\pi_n$ is an isomorphism. For a given sequence of $\cB$-permissible
blowups, define the associated contracted sequence by suppressing all
isomorphisms in the sequence and renumbering as in \eqref{seqBpermbu1} again.

We abbreviate \eqref{seqBpermbu1} as $(X,\cB) = (X_0,\cB_0) \lmapo{\pi_1} (X_1,\cB_1) \leftarrow \ldots\;$, for short.
\end{definition}

We will prove Theorem \ref{thm.intro.1} in the following, more general form.

\begin{theorem}\label{thm.main.1}
Let $(X,\cB)$ be as in (NE), with $X$ dimension at most two.

(a) There is a canonical
finite contracted sequence $S(X,\cB)$ of complete $\cB$-permissible blow-ups over $X$
$$
(X,\cB) = (X_0,\cB_0) \lmapo{\pi_1} (X_1,\cB_1) \lmapo{\pi_2} \ldots \lmapo{\pi_n} (X_n,\cB_n)
$$
such that $\pi_{i+1}$ is an isomorphism over $\XiBsqreg$, $0\leq i <n$, and $X_n$ is quasi-$\cB$-regular.
In particular, the morphism $X_n \rightarrow X$ is an isomorphism over $\XBsqreg$.

\smallskip
Moreover the following functoriality holds:
\begin{itemize}
\item[(F1)]
(equivariance) The action of the automorphism group of $(X,\cB)$ extends to the
sequence in a unique way.
\item[(F2)]
(localization) The sequence is compatible with passing to open subschemes $U \subseteq X$,
arbitrary localizations $U$ of $X$ and \'etale morphisms $U \rightarrow X$ in the following sense:
If $S(X,\cB)\times_XU$ denotes the pullback of $S(X,\cB)$ to $U$, then the associated contracted
sequence $(S(X,\cB)\times_XU)_{contr}$ coincides with $S(U,\cB_U)$.
\end{itemize}
(b) There is also a canonical finite contracted sequence $S_0(X,\cB)$ of strict $\cB$-permissible
blowups with the same properties (except that now each $\cB_{n+1}$ is the strict transform of
$\cB_n$).
\end{theorem}

If $X$ is reduced, then every $X_i$ is reduced, so that $X_n$ is regular and $\cB_n$
is normal crossing with $X_n$; moreover $\XBsqreg = \XBsreg$. In particular,
Theorem \ref{thm.intro.1} can be obtained as the case (b) for $\cB=\emptyset$
(where $\XBsreg = \Xreg$ and $\cB_n= \emptyset$ for all $n$), i.e., as the sequence $S_0(X,\emptyset)$.
If we apply (a) for reduced $X$ and $\cB=\emptyset$, we have $\XBsreg = \Xreg$
as well, but then, for the sequence $S(X,\emptyset)$, $\cB_n$ is not empty for $n>0$,
and we obtain the extra information that the collection
of the strict transforms of all created exceptional divisors is a simple normal crossing
divisor on $X_n$.

\begin{definition}\label{def.str.res.NE}
Let $\cC$ be a category of schemes which is closed under localization. Say that canonical, functorial
resolution with boundaries holds for $\cC$, if the statements in Theorem \ref{thm.main.1} (a) hold
for all schemes in $\cC$ and all boundaries on them. Say that canonical, functorial resolution
holds for $\cC$, if the statements of Theorem \ref{thm.intro.1} (i.e., of Theorem \ref{thm.main.1}
(b) with $\cB=\emptyset$) hold for all schemes in $\cC$.
\end{definition}

Now we will consider the case (E).

\begin{definition}\label{def.seqBpermbuZ}
(Case (E)) A sequence of complete (resp. strict) $\cB$-permissible blowups over $(X,Z)$ is a sequence of blowups:
\begin{equation}\label{seqBpermbu2}
\begin{array}{ccccccccccccccc}
 \cB=&\cB_0&& \cB_1&& \cB_2&& \cB_{n-1}&& \cB_{n} & \dots\\
 \\
  Z=&Z_0 & \lmapo{\pi_1} &   Z_1 & \lmapo{\pi_2} & Z_2 &
 \leftarrow\ldots \leftarrow & Z_{n-1} & \lmapo{\pi_n} & Z_{n} &\ldots\\
  &\cup & & \cup && \cup & & \cup && \cup & \\
 X=&X_0 & \lmapo{\pi_1} &   X_1 & \lmapo{\pi_2} & X_2 &
 \leftarrow\ldots \leftarrow & X_{n-1} & \lmapo{\pi_n} & X_{n} & \ldots \\
 \end{array}
\end{equation}
where for any $i\geq 0$
$$
\begin{matrix}
Z_{i+1}&=& \Bl_{D_{i}}(Z_{i}) &\rmapo{\pi_{i+1}}& Z_{i} \\
&& \cup  && \cup\\
X_{i+1}&=& \Bl_{D_{i}}(X_{i}) &\rmapo{\pi_{i+1}}& X_{i} \\
\end{matrix}
$$
are the blow-ups in a center $D_{i}\subset X_{i}$ which is permissible and n.c. with $\cB_i$,
and where $\cB_{i+1}=\cB_i'$ is the complete transform of $\cB_{i}$
(resp. $\cB_{i+1} = \widetilde{\cB_i}$ is the strict transform of $\cB_i$).

Call a sequence as in \eqref{seqBpermbu2} contracted if none of the
morphisms $\pi_n$ is an isomorphism. For a given sequence of $\cB$-permissible
blowups, define the associated contracted sequence by suppressing all
isomorphisms in the sequence and renumbering as in \eqref{seqBpermbu2} again.

We abbreviate \eqref{seqBpermbu2} as $(X,Z,\cB) = (X_0,Z_0,\cB_0) \lmapo{\pi_1} (X_1,Z_1,\cB_1) \leftarrow \ldots\;$, for short.
\end{definition}

We will prove Theorem \ref{thm.intro.3} in the following form.

\begin{theorem}\label{thm.main.3}
Let $(X,Z,\cB)$ be as in (E), with $X$ of dimension at most two.

(a) There is a canonical finite contracted sequence $S(X,Z,\cB)$ of complete $\cB$-permissible blow-ups over $X$
$$
(X,Z,\cB) = (X_0,Z_0,\cB_0) \lmapo{\pi_1} (X_1,Z_1,\cB_1) \lmapo{\pi_2} \ldots \lmapo{\pi_n} (X_n,Z_n,\cB_n)
$$
such that $\pi_{i+1}$ is an isomorphism over $(Z - X) \cup \XiBsqreg$, $0\leq i <n$ and $X_n$ is
quasi-$\cB_n$-regular.
In particular, the morphism $Z_n \rightarrow Z$ is an isomorphism over $(Z - X) \cup \XBsqreg$.
\smallskip
Moreover the following functoriality holds:
\begin{itemize}
\item[(F1)]
(equivariance) The action of the automorphism group of $(Z,X,\cB)$ (those automorphisms of
$Z$ which respect $\cB$ and $X$) extends to the sequence in a unique way.
\item[(F2)]
(localization) The sequence is compatible with passing to open subschemes $U \subseteq Z$,
arbitrary localizations $U$ of $Z$ and \'etale morphisms $U \rightarrow Z$ in the following sense:
If $S(X,Z,\cB)\times_ZU$ denotes the pullback of $S(X,Z,\cB)$ to $U$, then the associated
contracted sequence $(S(X,Z,\cB)\times_ZU)_{contr}$ coincides with $S(X\times_ZU,U,\cB_U)$.
\end{itemize}
(b) There is also a canonical finite contracted sequence $S_0(X,Z,\cB)$ of strict $\cB$-permissible
blowups over $(X,Z)$ with the same properties (except that now each $\cB_{n+1}$ is the strict transform of
$\cB_n$).
\end{theorem}

Again, for reduced $X$ all $X_i$ are reduced, $\XBsqreg = \XBsreg$, and $X_n$ is regular and
normal crossing with the simple normal crossings divisor $\cB_n$.

\begin{definition}\label{def.str.res.E}
Let $\cC$ be a category of schemes which is closed under localization. Say that canonical, functorial
embedded resolution with boundaries holds for $\cC$, if the statements in Theorem \ref{thm.main.3} (a) hold
for all triples $(X,Z,\cB)$ where $Z$ is a regular excellent scheme, $\cB$ is a simple normal crossing
divisor on $Z$ and $X$ is a closed subscheme of $X$ which is in $\cC$.
\end{definition}

\begin{remark}\label{rem.res.E-NE}
It follows from Lemma \ref{lem.comp.emb-nemb} that Theorem \ref{thm.main.1} implies Theorem \ref{thm.main.3},
in the following way: If $S(X,\cB_X)$ is constructed, one obtains $S(X,Z,\cB)$ by consecutively
blowing up $Z_i$ in the same center as $X_i$, and identifying $X_{i+1}$ with the strict transform
of $X_i$ in $Z_{i+1}$. Conversely, the restriction of $S(X,Z,\cB)$ to $X$ is $S(X,\cB_X)$.
More generally, by the same approach, canonical, functorial embedded resolution with boundaries holds for a category
$\cC$ of schemes as in Definition \ref{def.str.res.E} if canonical, functorial resolution with boundaries
holds for $\cC$.
\end{remark}

\bigbreak\noindent
We set up the strategy of proof for the above theorems in a more general setting.
Let $X$ be an excellent scheme, and recall that we only consider noetherian schemes here.

\begin{definition}\label{def.equising1}
Call an excellent (noetherian) scheme $Y$ equisingular, if $H_Y$ is constant on $Y$.
Call $Y$ locally equisingular if all connected components are equisingular.
\end{definition}

As we will see below, our strategy will be to make $X$ locally equisingular.

\begin{remark}\label{rem.equising1}
(a) For $U\subseteq X$ open and $x\in U$ one has $H_U(x) = H_X(x)$.
Hence $U$ is equisingular if and only if $U \subseteq X(\nu)$
for some $\nu \in \Sigma_X$.

\smallskip
(b) By Lemma \ref{lem.qreg} (a) there are two possibilities for a connected component
$U \subseteq X$: Either $U$ is equisingular (and irreducible),
or $\HSmax U$ is nowhere dense in $U$. In the first case, it follows from
Theorem \ref{HSf.usc} (1) and Theorem \ref{thm.pbu.inv} (1) that,
for any permissible blow-up $\pi: X' \rightarrow X$, $\pi^{-1}(U)$ is equisingular
as well (viz., $\pi^{-1}(U) \subset X'(\nu)$ if $U \subset X(\nu)$),
because by definition, permissible centers are nowhere dense in $X$.

\smallskip
(c) If $X$ is reduced, then a connected component $U \subset X$ is equisingular
if and only if $U$ is regular (cf. Remark \ref{HSregular}). Hence $X$ is locally
equisingular iff it is equisingular iff it is regular.

\smallskip
(d) By way of example, the following situation can occur for non-reduced schemes:
$X$ is the disjoint union of three irreducible components $U_1$, $U_2$ and $U_2$,
where $\Sigma_{U_1} = \{\nu_1\}$, $\Sigma_{U_2} = \{\nu_1,\nu_2\}$ and $\Sigma_{U_3} = \{\nu_3\}$,
such that $\nu_1 < \nu_2 < \nu_3$.
By just blowing up in $\HSmax X$ we cannot make $X$ locally equisingular.
\end{remark}

Motivated by the remarks above, we define:

\begin{definition}\label{def.nu-elim}
Let $X$ be connected and not equisingular. For $\nu\in \SigmaXmax$, a $\nu$-elimination for
$X$ is a morphism $\rho:X'\to X$ that is the composite of a sequence of morphisms:
$$
X=X_0\gets X_1 \gets\cdots \gets X_n=X'
$$
such that
for $0\leq i < n$, $\pi_i:X_{i+1}\to X_i$ is a blowup in a permissible center
$D_i\subseteq \HSa {X_i}$ and $\HSa {X_n}=\emptyset$.
\end{definition}

Let $\nu_1,\dots,\nu_r$ be the elements of $\SigmaXmax$ and assume given
a $\nu_i$-elimination $\rho_i: X_i\to X$ of $X$ for each $i \in \{1,\dots,r\}$.
Noting that $\rho_i$ is an isomorphism over $X-\HS X {\nu_i}$ and that
$\HS X {\nu_i} \cap \HS X {\nu_j} =\emptyset$ if $1\leq i\not= j\leq r$,
we can glue the $\rho_i$ over $X$ -- which is a composition of permissible blow-ups again --
to get a morphism
$\rho: X'\to X$ which is a $\Sigmamax$-elimination where we define:

\begin{definition}\label{def.max-elim1}
Let $X$ be connected and not equisingular.
A morphism $\rho : X' \longrightarrow X$ is called a $\Sigmamax$-elimination for $X$
if the following conditions hold:

\smallskip (ME1)
$\rho$ is the composition of permissible blowups and an isomorphism
over $X-\HSmax X$.

\smallskip (ME2) $\Sigma_{X'} \cap \SigmaXmax = \emptyset$.
\end{definition}

Note that, by Theorem \ref{thm.pbu.inv} (1), (ME1) and (ME2) imply:

\medskip (ME3) For each $\mu \in \Sigma_{X'}$ there exists a $\nu \in \SigmaXmax$
with $\mu < \nu$.

\begin{definition}\label{def.max-elim2}
For any excellent scheme $X$, a morphism $\rho: X' \rightarrow X$ is called a
$\Sigmamax$-elimination, if it is a $\Sigmamax$-elimination after restriction
to each connected component which is not equisingular, and an isomorphism
on the other connected components.
\end{definition}

\begin{theorem}\label{thm.max-elim}
Let $X$ be an excellent (noetherian) scheme, and let $X=X_0\gets X_1\gets \cdots$
be a sequence of morphisms such that $\pi_n:X_{n+1}\to X_n$ is a $\Sigmamax$-elimination
for each $n$. Then there is an $N\in \bN$ such that $X_N$ is locally equisingular.
(So $\pi_n$ is an isomorphism for $n\geq N$.)
\end{theorem}

\textbf{Proof }
(See also Theorem \ref{thm.ordHP} for an alternative, more self-contained proof of Theorem
\ref{thm.max-elim}.)
Suppose there exists an infinite sequence $X = X_0 \leftarrow X_1 \leftarrow \ldots $
of $\Sigmamax$-eliminations such that no $X_n$ is locally equisingular.
For each $n\geq 0$ let $X_n^0 \subseteq X_n$ be the union of those connected
components of $X_n$ which are not equisingular, and let $\Sigma_n= \Sigma_{X_n^0}$ and
$\Sigma^{max}_n = \Sigma^{max}_{X_n^0}$.
For each $n$, choose an element $\nu_n \in \Sigma^{max}_n$. It follows from the finiteness
of $\Sigma_{X_0}$ and Theorem \ref{thm.pbu.inv} that there is an $m\in\bN$ such that all
Hilbert-Samuel functions occurring on the $X_i$ are contained in the set $HF_m$ of all Hilbert
functions $H$ (of standard graded algebras) with $H(1)\leq m$.
Thus it follows from Theorem \ref{thm.ordHF}, i.e., the noetherianess of $HF_m$, that $\nu_i \leq \nu_j$ for some
$i < j$ (see \cite{AP} Proposition 1.3 for several equivalent conditions characterizing
noetherian ordered sets). Choose $x_i \in X_i^0$ and $x_j \in X_j^0$ with $\nu_i = H_{X_i}(x_i)$
and $\nu_j = H_{X_j}(x_j)$, and for $i \leq k \leq j$, let $y_k$ be the image of $x_j$ in $X_k$.
By Remark \ref{rem.equising1} (b) the morphism $X_{n+1} \rightarrow X_n$
maps $X_{n+1}^0$ to $X_n^0$, so we have $y_k \in X_k^0$, and by Theorem \ref{thm.pbu.inv}
and Remark \ref{rem.equising1} we have $\nu_j \leq \mu_k := H_{X_k}(y_k)$.
Since $\nu_i \in \Sigma^{max}_i \subset \Sigma_i$, we conclude that the inequalities
$\nu_i \leq \nu_j \leq \mu_k$ are all equalities. Therefore $\nu_i = \mu_{i+1} \in \Sigma_{i+1}$,
contradicting the assumption that $X_{i+1} \rightarrow X_i$ is a $\Sigma^{max}$-elimination.

\begin{corollary}\label{cor.max-elim} To prove (canonical, functorial) resolution of
singularities for all excellent reduced schemes of dimension $\leq d$, it suffices to prove that
for every connected non-regular excellent reduced scheme $X$ of dimension $\leq d$ there exists a
(canonical functorial) $\Sigmamax$-elimination $X' \rightarrow X$. Equivalently, it suffices to show that for every
such scheme and every $\nu \in \SigmaXmax$, there is a (canonical functorial) $\nu$-elimination for $X$.
Here functoriality means that the analogues of the properties (F1) and (F2) in Theorem \ref{thm.main.1}
hold for the $\Sigmamax$- and $\nu$-eliminations, respectively, where the analogue of property (F2)
for a $\nu$-elimination is the following: Either $U(\nu) = \emptyset$, or $\nu \in \SigmaUmax$
and the pullback of the sequence to $U$ is the canonical $\nu$-elimination
on $U$, after passing to the associated reduced sequence.
\end{corollary}

In fact, under these assumptions one gets a (canonical, functorial) sequence $X \leftarrow X_1 \leftarrow \ldots$
of $\Sigmamax$-eliminations, and by Theorem \ref{thm.max-elim} some $X_n$ is locally equisingular,
which means that $X_n$ is regular (Remark \ref{rem.equising1} (c)).

\medskip
Now we consider the non-reduced case.

\begin{corollary}\label{cor.max-elim2} To prove (canonical, functorial) resolution of
singularities for all excellent schemes of dimension $\leq d$, it suffices to prove that
there exists a (canonical, functorial) $\Sigmamax$-elimination $X' \rightarrow X$
for every connected excellent scheme $X$ of dimension $\leq d$ which is not equi-singular.
Equivalently, it suffices to show that for every
such scheme and every $\nu \in \SigmaXmax$, there is a (canonical, functorial) $\nu$-elimination for $X$.
Here functoriality is defined as in  Corollary \ref{cor.max-elim}.
\end{corollary}

In fact, here we first get a (canonical, functorial) sequence $X \leftarrow X_1 \leftarrow \ldots \leftarrow X_m$
of $\Sigmamax$-eliminations such that $X_m$ is locally equisingular.
By Corollary \ref{cor.max-elim} we get a similar sequence
$(X_m)_{red} \leftarrow X'_{m+1} \leftarrow \ldots \leftarrow X'_n$ such that
$X'_n$ is regular. Blowing up in the same centers we get a sequence of
blow-ups $X_m \leftarrow X_{m+1} \leftarrow \ldots \leftarrow X_n$,
where $X_i'$ is identified with $(X_i)_{red}$, and $X'_{i+1}$ with the strict
transform of $X'_i$ in $X_{i+1}$
(since $D_i \subseteq (X_i)_{red}$, $X_{i+1}'$ is reduced, and homeomorphic to $X_{i+1}$).
For each $i\geq m$, $X_i$ is again equisingular (see Remark \ref{rem.equising1}),
and by Theorem \ref{HSf.perm} the blow-up $X_{i+1} \rightarrow X_i$ is permissible.
It follows that $(X_n)_{red}$ is regular and $X_n$ is normally flat along $(X_n)_{red}$.
Now assume that the first sequence and the sequence $(X_m)_{red} \leftarrow \ldots$
are functorial. Then it is immediate that the sequence $X \leftarrow \ldots \leftarrow X_n$
is functorial for localizations as well. As for automorphisms, it follows
inductively via localization that the automorphisms of $X_i$ $(i\geq m)$
respect the center of the blow-up $X_{i+1} \rightarrow X_i$ and therefore
extend to $X_{i+1}$ in a unique way.

\begin{remark}\label{rem.non-reduced}
We point out the choice of strategy here. It might be tempting to start with desingularizing
$X_{red}$. But if $D \subset X_{red}$ is permissible in $X_{red}$, $D$ will not in general
be permissible in $X$. However, once $X$ is made equi-singular (by the first series of blow-ups),
the arguments above show that we get permissible blow-ups both for $X_{red}$ and $X$, and
can achieve that $X_{red}$ becomes regular and $X$ equi-singular at the same time, so that
$X$ is normally flat along $X_{red}$.
\end{remark}

We now consider a variant of the above for schemes with boundary.
Let $X$ be an excellent scheme, and let $\cB$ be a boundary for $X$, i.e., a boundary
on $X$ (case (NE)) or on $Z$ (case (E)).
In the following we only consider complete transforms for the boundaries, i.e.,
sequences of complete $\cB$-permissible blow-ups, and we will simply speak of
sequences of $\cB$-permissible blow-ups. It is easy to see that the analogous
results also hold for the case of strict transforms, i.e., sequences of strict
$\cB$-permissible blow-ups.

\begin{definition} \label{def.O-equising}
Call $X$ $\OB$-equisingular if $H^\OB_X$ is constant on $X$, and locally $\OB$-equisingular,
if every connected component is $\OB$-equisingular.
\end{definition}

\begin{remark}\label{rem.O-equising}
(a) It follows from Lemma \ref{lem.qreg} (b) that a connected component $U \subseteq X$
is either $\OB$-equisingular, or $U^\OB_{max}$ is nowhere dense in $U$.
In the first case $U \subseteq X^\OB(\tnu)$ for some $\tnu\in\tSigmaX$, and
for every $\cB$-permissible blow-up $\pi: X' \rightarrow X$ one has $\pi^{-1}(U)
\subseteq X'(\tnu)$.

\smallskip
(b) If $X$ is reduced, then $X$ is locally $\OB$-equisingular if and only if $X$ is
$\OB$-regular.

\smallskip
(c) Even for a regular scheme $X$ it can obviously happen that $X$ is the union
of three irreducible components $U_1$, $U_2$ and $U_3$ such that $\Sigma^\OB_{U_1}=\{(\nu_X^{reg},1)\}$,
$\Sigma^\OB_{U_2}=\{(\nu_X^{reg},1), (\nu_X^{reg},2)\}$ and $\Sigma^\OB_{U_3}=\{(\nu_X^{reg},3)\}$,
so that $X$ cannot be made $\OB$-equisingular by blowing up in $X^\OB_{max}$.
\end{remark}

\begin{definition}\label{def.Omax-elim}
Let $\OB$ be a history function for $\cB$ such that $(\cB,\OB)$ is admissible, and let
\begin{equation}\label{eq.Omax-elim}
X=X_0 \lmapo{\pi_1} X_1 \leftarrow \ldots \leftarrow X_{n-1} \lmapo{\pi_n} X_n
\end{equation}
be a sequence of $\cB$-permissible blow-ups (where we have not written the boundaries $\cB_i$,
and neither the regular schemes $Z_i$ in case (E)). For each
$i = 0,\ldots,n-1$ let $(\cB_{i+1},\OB_{i+1})$ be the complete transform of $(\cB_i,\OB_i)$
(where $(\cB_0,\OB_0) = (\cB,\OB)$).
Let $D_i$ be the center of the blow-up $X_i \lmapo{\pi_{i+1}} X_{i+1}$, and $\rho = \pi_n\circ\ldots\circ\pi_1: X_n \rightarrow X$.

\smallskip
(1) If $X$ is connected and not $\OB$-equisingular, and $\tnu \in \tSigmaXmax$, then
\eqref{eq.Omax-elim} or $\rho$ is called a $\tnu$-elimination,
if $D_i \subseteq \tHSa {X_i}$ for $i=1,\ldots,n-1$ and $\tHSa {X_n}=\emptyset$.

\smallskip
(2) If $X$ is connected and not $\OB$-equisingular, then \eqref{eq.Omax-elim} or $\rho$
is called a $\tSigmamax$-elimination for $(X,\cB,\OB)$,
if $D_i \subseteq (X_i)^\OB_{max}$ for $i=0,\ldots,n-1$ and $\tSigmaXn \cap \tSigmaXmax = \emptyset$.

\smallskip
(3) If $X$ is arbitrary, then \eqref{eq.Omax-elim} or $\rho$ is called a $\tSigmamax$-elimination for $(X,\cB,\OB)$,
if it is a $\tSigmamax$-elimination after restriction to each connected component of $X$
which is not $\OB$-equisingular, and an isomorphism after restriction to the
connected components which are $\OB$-singular.

(4) Call the sequence \eqref{eq.Omax-elim} contracted, if none of the morphisms is an isomorphism,
and in general define the associated contracted sequence by omitting the isomorphisms and
renumbering (so the final index $n$ may decrease).
\end{definition}

\begin{remark}\label{rem.tSigmamax.elim}
By glueing, one gets a (canonical, functorial) $\tSigmamax$-elimination for a connected,
not $\OB$-equisingular $X$, if one has (canonical,
functorial) $\tnu$-eliminations for all $\tnu \in \tSigmaXmax$,
and a (canonical, functorial) $\tSigmamax$-elimination for
a non-connected $X$, if one has this for all connected components. Here functoriality is defined
as in Theorem \ref{thm.main.1}.
\end{remark}

In a similar way as in Theorem \ref{thm.max-elim} one proves:

\begin{theorem}\label{thm.Omax-elim}
For any infinite sequence $X=X_0 \leftarrow X_1 \leftarrow X_2 \leftarrow \ldots$ of
$\tSigmamax$-eliminations there is an $n$ such that $(X_n,\cB_n,\OB_n)$ is
$\OB$-equisingular.
\end{theorem}

The following result is now obtained both in the embedded and the non-embedded case.

\begin{corollary}\label{cor.Omax-elim}
Case (NE): To show (canonical, functorial) resolution of singularities
with boundaries for all (noetherian) reduced excellent schemes of dimension $\leq d$ it suffices to show the
existence of (canonical, functorial) $\tSigmamax$-eliminations for all connected reduced excellent schemes $X$ of
dimension $\leq d$ and all admissible boundaries with history $(\cB,\OB)$ for $X$,
for which $X$ is not $\OB$-regular. (Here `functorial' in the last statement means that
the obvious analogues of the conditions (F1) and (F2) in Theorem \ref{thm.main.1} hold for
the sequences considered here.)

Case (E): The obvious analogous statement holds.
\end{corollary}

In fact, if $(X,\cB)$ is given, we start with the history function $\OB(x) = \cB(x)$. Then $X$ is
$\OB$-regular if and only if $X$ is strongly $\cB$-regular at all $x\in X$. If this holds, we are done.
If not, then by assumption there is a (canonical, functorial) $\tSigmamax$-elimination $X_1 \rightarrow X$,
and we let $(\cB_1,\OB_1)$ be
the strict transform of $(\cB,\OB)$ in $X_1$ (which is obtained by successive transforms for
the sequence of $\cB$-permissible blowups whose composition is $X_1 \rightarrow X$). Then
$(\cB_1,\OB_1)$ is admissible by Lemma \ref{lem.cNBncTx}.
If $X_1$ is $\OB_1$-regular, we are done by Lemma \ref{lem.tnureg}.
If not we repeat the process, this time with $(X_1,\cB_1,\OB_1)$, and iterate if necessary. By
Theorem \ref{thm.Omax-elim}, after finitely many steps this process obtains an $X_n$ which is
$\OB_n$-regular and hence achieves the resolution of $X$ by Lemma \ref{lem.tnureg} (case (E))
and Lemma \ref{lem.ctnureg} (case (NE)).

\bigskip
In the non-reduced case we obtain:

\begin{corollary}\label{cor.Omax-elim1}
Case (NE): To show (canonical, functorial) resolution of singularities
with boundaries for all excellent (noetherian) schemes of dimension $\leq d$ it suffices to show the
existence of (canonical, functorial) $\tSigmamax$-eliminations for all connected excellent schemes $X$ of
dimension $\leq d$ and all admissible boundaries with history $(\cB,\OB)$ for $X$,
for which $H^\OB_X$ is not constant. (Here `functorial' in the last statement means that
the obvious analogues of the conditions (F1) and (F2) in Theorem \ref{thm.main.1} hold for
the sequences considered here.)

Case (E): The obvious analogous statement holds.
\end{corollary}

This follows from Corollary \ref{cor.Omax-elim} in a similar way as Corollary \ref{cor.max-elim2}
follows from \ref{cor.max-elim}: First we get a (canonical, functorial) sequence of $\cB$-permissible blow-ups
$X \leftarrow X_1 \leftarrow \ldots \leftarrow X_m$ such that $X_m$ is locally $\OB$-equisingular.
Then we look at the (canonical, functorial) resolution sequence
$(X_m)_{red} \leftarrow X'_{m+1} \leftarrow \ldots \leftarrow X'_n$ from Corollary \ref{cor.Omax-elim}
such that $X'_n$ is $\cB'_n$-regular, where $\cB'_n$ comes from $\cB$ via complete transforms.
By blowing up in the same centers we obtain a sequence of $\cB$-permissible blow-ups
$X_m \leftarrow X_{m+1} \leftarrow \ldots \leftarrow X_n$ such that $(X_n)_{red}$ identifies with
$X'_n$ and thus is $\cB'_n$-regular; moreover, $X_n$ is normally flat along $(X_n)_{red}$,
because $H^\OB_{X_n}$ is constant on all connected components.

\bigskip
We now prove Theorem \ref{thm.main.1}. Then, by Remark \ref{rem.res.E-NE}, Theorem \ref{thm.main.3} follows as well.

\bigskip
By Corollary \ref{cor.Omax-elim1},
it suffices to produce canonical, functorial $\tSigmamax$-eliminations for all connected
excellent (noetherian) schemes $X$ of dimension at most two and all admissible
boundaries with history $(\cB,\OB)$ on $X$ such that $H^\OB_X$ is not constant on $X$.
By the remarks after Definition \ref{def.Omax-elim} it suffices to produce canonical functorial
$\tnu$-eliminations for all $\tnu \in \tSigmaXmax$.
We will slightly modify this procedure and deduce Theorem  \ref{thm.main.1} from the following,
partly weaker and partly more general result.

\begin{theorem}\label{thm.4.3} (Case (NE))
Let $X$ be an excellent connected (noetherian) scheme, let $(\cB,\OB)$ be an admissible
boundary with history on $X$ such that $H^\OB_X$ is not constant on $X$,
and let $\tnu \in \tSigmaXmax$.
Assume the following:
\begin{itemize}
\item[(1)]
$\Char(k(x)) = 0$, or $\Char(k(x)) \geq \dim (X)/2+1$ for any $x\in \tHSa X$,
\item[(2)]
$\dim(\tHSa X)\leq 1$,
\end{itemize}
and there is an integer $e$ with $0\leq e\leq 2$ such that for any
closed point $x\in \tHSa X$,
\begin{itemize}
\item[(3e)]
$\eb_x(X)\leq e$,
\item[(4e)]
either $\NBx \pitchfork \Dirob_x(X)$ or $\eob_x(X) \leq e-1$.
\end{itemize}
Then there exists a canonical reduced $\tnu$-elimination
$S(X,\tnu)$
$$
(X,\cB)=(X_0,\cB_0) \leftarrow (X_1,\cB_1) \leftarrow \ldots \leftarrow (X_n,\cB_n)
$$
for $(X,\cB)$. It satisfies the analogues of properties (F1) and (F2) from Theorem \ref{thm.main.1},
where the analogue of (F2) is the following: Either $U(\tnu) = \emptyset$, or else $\tnu \in \tSigmaUmax$
and the reduced sequence associated to the pullback of the sequence to $U$ is the canonical
$\tnu$-elimination for $(U,\cB_U)$.
\end{theorem}

{\bf Theorem \ref{thm.main.1} follows from this:} If the conditions of Theorem \ref{thm.4.3}
hold for $X$, and $X' \rightarrow X$ is a blowup in a $\cB$-permissible center $D \subseteq \tHSa X$,
the conditions hold for $X'$ as well, with the same $e$. In fact, (1) holds for $X'$ since $\dim(X')=\dim(X)$,
condition (3e) holds for $X'$ by Theorem \ref{thm.pbu.inv}, and condition (4e) holds for $X'$
by Corollary \ref{cor.cNBncDir}. Moreover, condition (2) holds for $X'$ as we will see in the proof
of Theorem \ref{thm.4.3}. Now assume that $X$ is of dimension $d \leq 2$.
Then condition (1) holds, and condition (3e) holds
with $e = d$. If moreover $H^\OB_X$ is not constant on $X$, then condition (2) holds by Lemma \ref{lem.qreg},
and (2) holds for $X'$ as well.
On the other hand, in the presence of condition (1) it suffices to consider admissible boundaries with history
$(\cB,\OB)$ which satisfy condition (4e).
In fact, in the procedure outlined in the proof of Corollary \ref{cor.Omax-elim1}, property $(4e)$ is trivially fulfilled
in the beginning where $\OB(x) = \cB(x)$, i.e., $\NBx = \emptyset$, and as remarked a few lines above,
it is also fulfilled in a sequence of $\cB$-permissible blowups.
Therefore we can keep all assumptions of Theorem \ref{thm.4.3}
in a sequence of $\cB$-permissible blowups, and the arguments for Corollaries \ref{cor.Omax-elim} and \ref{cor.Omax-elim1}
apply in this modified setting.

\begin{remark}\label{rem.strategy}
(1) Before we start with the proof of Theorem \ref{thm.4.3}, we outline the strategy which works under
the conditions (1) to (4) of Theorem \ref{thm.4.3}, but can be stated more generally.
It would be interesting to see if it also works for higher-dimensional schemes.

\smallskip
Let $X$ be an excellent connected noetherian scheme, let $(\cB,\OB)$ be an admissible
boundary with history on $X$ such that $H^\OB_X$ is not constant on $X$,
and let $\tnu \in \tSigmaXmax$. We construct a canonical contracted sequence
$S(X,\tnu)$
\begin{equation}\label{eq.strategy.1}
X=X_0 \leftarrow X_1 \leftarrow X_2\leftarrow X_3 \leftarrow \ldots
\end{equation}
of $\cB$-permissible blowups over $X$ as follows. We introduce some notations. If $X_n$ has
been constructed, then let $Y_n = \tHSa {X_n}$. Give labels (or `years') to the irreducible components
of $Y_n$ in an inductive way as follows. The irreducible components of $Y_0$ all have label $0$.
If an irreducible component of $Y_n$ dominates an irreducible component of $Y_{n-1}$,
it inherits its label. Otherwise it gets the label $n$. Then we can write
\begin{equation}\label{eq.strategy.2}
\tHSa {X_n} = Y_n = Y_n^{(0)} \cup Y_n^{(1)} \cup \ldots \cup Y_n^{({n-1})} \cup Y_n^{(n)}\,,
\end{equation}
where $Y_n^{(i)}$ is the union of irreducible components of $Y_n$ with label $i$.

\smallskip
Now we start the actual definition of \eqref{eq.strategy.1}.
We have $\dim(Y_0)<\dim(X)$. By induction on dimension we have a canonical resolution sequence
$$
Y_0 := Y_{0,0}\leftarrow Y_{0,1} \leftarrow Y_{0,2} \leftarrow \ldots \leftarrow Y_{0,m_0}
$$
for $(Y_0,\cB_{Y_0})$, i.e., a sequence of $\cB$-permissible blowups, so that $Y_{0,m_0}$ is $\cB$-regular.
By successively blowing up $X$ in the same centers -- which are also permissible for the $(X_i,\cB_i)$
by Lemma \ref{lem.comp.incl} (3), we obtain a sequence of $\cB$-permissible blowups
$$
X_0 \leftarrow X_1 \leftarrow X_2 \leftarrow \ldots \leftarrow X_{m_0}
$$
in which $Y_{0,i}$ is the strict transform of $Y_{0,i-1}$ and, moreover, $Y_{0,i} = Y^{(0)}_i$
 so that $Y_{m_0}^{(0)}=Y_{0,m_0}$ is $\cB$-regular.

\smallskip
Then we blow up $Y_{m_0}^{(0)}$ to get $X_{m_0+1}$. We call the obtained sequence (from $X_0$ to $X_{m_0+1}$)
of $\cB$-permissible blow-ups the first (resolution) cycle.
If $Y_{m_0+1}^{(0)}$ is non-empty and not $\cB$-regular,
we proceed as above and have a canonical resolution sequence
$$
Y_{m_0+1}^{(0)}=Y_{0,m_0+1} \leftarrow Y_{0,m_0+2} \leftarrow \ldots \leftarrow Y_{0,m_1}
$$
such that $Y_{0,m_1}$ is $\cB$-regular. Then we get a sequence
$$
X_{m_0+1} \leftarrow X_{m_0+2} \leftarrow \ldots \leftarrow X_{m_1}
$$
of $\cB$-permissible blow-ups by blowing up in the same centers, for which $Y_{0,j}=Y_j^{(0)}$,
and $Y_{m_1}^{(0)}=Y_{0,m_1}$ is $\cB$-regular. Then we blow up $Y_{m_1}^{(0)}$
to get $X_{m_1+1}$. We call the sequence from $X_{m_0+1}$ to $X_{m_1+1}$ the second cycle.
Repeating this process finitely many times, i.e., producing further cycles,
we get $X_{m_\ell}=X_{n_1}$ for which $Y_{n_1}^{(0)}$ is empty (in our situation, which is a non-trivial fact).
Then we proceed with $Y_{n_1}^{(1)}$ as we did before with $Y_0^{(0)}$, in several cycles, until we reach $X_{n_2}$
for which $Y_{n_2}^{(1)}$ is empty, and proceed with $Y_{n_2}^{(2)}$, etc. This procedure ends, i.e.,
there is an $n_r$ such that $Y_{n_r}$ is empty (in our situation, which is a non-trivial fact), so that
we have eliminated the $\tnu$-locus.

\medskip
(2) In the situation of Theorem \ref{thm.4.3}, we will see that $Y_{n}^{(i)}$ is already $\cB$-regular for all $n\geq m_0$ and $i\geq 0$.

\medskip
(3) We will see below in Proposition \ref{prop.funct.strat} that the sequence constructed in (1)
is functorial for automorphisms and morphisms with geometrically regular fibers.

\medskip
(4) By (3), Remark \ref{rem.tSigmamax.elim} and Corollaries \ref{cor.Omax-elim} and \ref{cor.Omax-elim1} we have: If the procedure in (1)
is always finite, i.e., gives $\tnu$-eliminations for $\tSigmamax$-eliminations, then we have a canonical, functorial resolution
with boundaries as formulated in \ref{cor.Omax-elim} and \ref{cor.Omax-elim1}.
\end{remark}

\bigskip
\textbf{Proof of Theorem \ref{thm.4.3}}
Let $X' \rightarrow X$ be a blowup in a $\cB$-permissible center $D \subset \tHSa X$.
As we have seen above, conditions (1), (3e) and (4e) hold for $X'$ and the complete
transform $(\cB',\OB')$ of $(\cB,\OB)$ as well.
Now we will show the same for condition (2). We may assume that $D$ is irreducible.
By (2) $D$ is a closed point or a regular irreducible curve.

\smallskip
\textbf{Step 1} Let $x$ be any closed point in $\tHSa X$ and consider
$\pi:X':=\Bl_x(X) \to X$. Note that $x\hookrightarrow X$ is $\cB$-permissible
for trivial reasons. By Theorem \ref{thm.cBdirectrix} we have
\begin{equation}\label{prop.eq}
D' := \tHSa {X'}\cap\pi^{-1}(x)\subset \bP(\Dirob_x(X))\simeq \bP^t_{k(x)},
\end{equation}
where $t=\eob_x(X)-1\leq e_x(X)-1\leq e-1\leq 1$
(by convention $\bP^t_{k(x)}=\emptyset$ if $t<0$).
Hence condition (2) is also satisfied for $X'$.
Moreover, if $\dim(D')\geq 1$, then $D'=\bP(\Dirob_x(X))$,
so that $D'$ is $\OB'$-permissible,
and condition (3e) implies $e=2$ and $\NBx \pitchfork \Dirob_x(X)$, so that $D'$
is n.c. with $N'$ by Lemma \ref{lem.cNBncDir1}.
Hence $D'$ is a union of closed points or a projective line over $k(x)$, and
in both cases it is $\cB'$-permissible.

\medbreak
\textbf{Step 2} Now let $D\subset \tHSa X$ be regular irreducible of dimension $1$ and n.c.
with $\cB$. By Theorem \ref{thm.cOBHSf.perm}, $D\subset X$ is $\cB$-permissible.
Let $\eta$ be the generic point of $D$.
Consider $\pi: X':=\Bl_D(X) \to X$.
Let $x \in D$ be a closed point. By Theorem \ref{thm.cBdirectrix}, we have
$$
\tHSa {X'}\cap\pi^{-1}(x)
\subset \bP(\Dirob_x(X)/T_x(D))\simeq \bP^s_{k(x)},
$$
where $s=\eob_x(X) - 2 \leq e_x(X)-2\leq e-2\leq 0$ by (3e) for $X$. Hence there is
at most one point in $\tHSa {X'}\cap \pi^{-1}(x)$ so that
$\dim(\tHSa {X'}\cap \pi^{-1}(D))\leq 1$ and
condition (2) is satisfied for $X'$. Similarly we have
$$
\tHSa {X'}\cap \pi^{-1}(\eta) \subset
\bP(\Dirob_\eta(X))\simeq \bP^r_{k(\eta)},
$$
where $r=\eob_\eta(X) -1\leq e_\eta(X)-1\leq e_x(X) -2 \leq 0$ by Theorem \ref{nf.dir}.
Hence, if $\HSa {X'}\cap \pi^{-1}(\eta)$ is not
empty, then it consists of a unique point $\eta'$, and one has
$k(\eta)\simeq k(\eta')$.
This implies that $\pi$ induces an isomorphism $D'\isom D$ where
$D' = \tHSa {X'}\cap \pi^{-1}(D)$. Thus $D'$ is regular, and $\OB'$-permissible
by Theorem \ref{thm.cOBHSf.perm}. Moreover, $e_x(X)=2$ in this case so that
condition (3e) for $X$ implies $e=2$ and $\NBx \pitchfork \Dir_x(X)$.
By Lemma \ref{lem.cNBncDir2}, $D'$ is n.c. with $N'$, hence with $\cB'$. Hence $D'$ is a collection of
closed points or a regular irreducible curve, and in both cases it is $\cB'$-permissible.

\bigskip
\textbf{Step 3} Consider the special case where $\dim(X)=1$. Here $\dim \tHSa X = 0$, so every point
$x \in \tHSa X$ is isolated in $\tHSa X$, and moreover we have $e_x(X) \leq \dim(X)=1$.
The canonical $\tnu$-elimination sequence as defined in Remark \ref{rem.strategy} consists of blowing up all points in $\tHSa X$
and repeating this process as long as $\tHSa X \neq \emptyset$. By Theorem \ref{fu.thm0}
below this process stops after finitely many steps. So Theorem \ref{thm.4.3} holds.
As noticed after Theorem \ref{thm.4.3}, this shows that Theorem \ref{thm.main.1} holds for $\dim(X) = 1$,
i.e., there exists a canonical, functorial resolution sequence for $(X,\cB)$.

\bigskip
\textbf{Step 4} Now we consider the general case and construct a canonical reduced sequence
$S(X,\tnu)$
\begin{equation}\label{eq.can.seq}
X=X_0 \lmapo{\pi_1} X_1 \leftarrow \ldots \leftarrow X_{n-1} \lmapo{\pi_n} X_n \dots
\end{equation}
of $\cB$-permissible blowups over $X$ as follows. Let $Y_0 = \tHSa {X_0}$. If $X_n$ has
been constructed, then let $Y_n = \tHSa {X_n}$. Give labels to the irreducible components
of $Y_n$ in an inductive way as follows. The irreducible components of $Y_0$ all have label $0$.
If an irreducible component of $Y_n$ dominates an irreducible component of $Y_{n-1}$,
it inherits its label. Otherwise it gets the label $n$. Then we can write
\begin{equation}\label{eq.labels}
Y_n = Y_n^{(0)} \cup Y_n^{(1)} \cup \ldots \cup Y_n^{({n-1})} \cup Y_n^{(n)}\,,
\end{equation}
where $Y_n^{(i)}$ is the union of irreducible components of $Y_n$ with label $i$.

\medskip
\textbf{Step 5} By assumption, $\dim(Y_0)\leq 1$. Let $\cB_{Y_o}=\cB\times_X Y_0$ be the pull-back, and let
$$
Y_0 = Y_{0,0} \leftarrow Y_{0,1} \leftarrow \ldots Y_{0,m-1} \leftarrow Y_{0,m}
$$
be the canonical resolution sequence of Theorem \ref{thm.main.1} for $(Y_0,\cB_{Y_0})$
(which exists and is finite by Step 3),
so that $Y_{0,m}$ is regular and normal crossing with $\cB_{0,m}$, where we write $\cB_{0,i}$
for the boundary obtained on $Y_{0,i}$. Let
\begin{equation}\label{eq.can.seq.1}
X=X_0 \lmapo{\pi_1} X_1 \leftarrow \ldots \leftarrow X_{m-1} \lmapo{\pi_m} X_m
\end{equation}
be the sequence of blow-ups obtained inductively by blowing up
$X_i$ in the center $D_i$ of the blowup $Y_{0,i+1} \rightarrow Y_{0,i}$ and identifying
$Y_{0,i+1}$ with the strict transform of $Y_{0,i}$ in $X_{i+1}$.
$D_i$ is a collection of closed points and hence $\cB_i$-permissible, where we
write $\cB_i$ for the boundary obtained on $X_i$. Therefore \eqref{eq.can.seq.1}
is a sequence of $\cB$-permissible blow-ups. (We could also use Lemma \ref{lem.comp.incl} (3)
as in Remark \ref{rem.strategy}.)
By Lemma \ref{lem.comp.incl} (4) we have $\cB_{0,i} = (\cB_i)_{Y_{0,i}}$.
Since each $D_i$ is a nowhere dense subscheme of $Y_{0,i}$, each $Y_{0,j}$ is contained in $\tHSa {X_j} = Y_j$,
and is in fact equal to the label 0 part $Y_j^{(0)}$ of $Y_j$ as defined above.
This is the first stage of \eqref{eq.can.seq}.

\medskip
{\bf Claim 1} For $m$ as above, and all $i\geq 0$, the subschemes $Y_m^{(i)}$ are regular,
of dimension at most 1, and $\cB_m$-permissible.

\medskip
In fact, for $Y_m^{(0)}$ this holds by construction. Moreover, from the statements in Step 1
we conclude that, for $0 < i \leq m$, all schemes $Y_i^{(i)}$ are disjoint unions
of closed points and projective lines and hence regular, moreover they are $\cB_i$-permissible.
Let $[Y_j]_0$ be the union of the 0-dimensional components of $Y_j$.
Since $X_{j+1} \rightarrow X_j$ is a blow-up in closed points of $Y_j^{(0)}$, and
$Y_j^{(0)}\cap [Y_j^{(i)}]_0 = \emptyset$ for $i>0$ (no closed point can dominate the curve $Y_{j-1}^{(0)}$),
the morphism $Y_j^{(i)} \rightarrow Y_i^{(i)}$ is an isomorphism for $j = i,\ldots,m$.
Hence $Y_j^{(i)}$ is regular, and normal crossing with $\cB_j$ (direct check, or application of Lemma \ref{lem.Bperm0}).
Since $Y_j^{(i)}$ is regular and contained in $\tHSa {X_j}$, we conclude
it is $\cB_j$-permissible.

\medskip
\textbf{Step 6} Next we blow up the subscheme $Y_m^{(0)}$, which is regular and $\cB_m$-permissible, and obtain $X_{m+1}$.

\medskip
{\bf Claim 2} For all $i\geq 0$, the subschemes $Y_{m+1}^{(i)}$ are regular, of dimension at most 1, and $\cB_{m+1}$-permissible.
Moreover, the intersection of $Y_{m+1}^{(m+1)}$ with $Y_{m+1}^{(i)}$ is empty for all $i\in \{0,\ldots,m\}$.

\medskip
The first part follows by similar arguments as above. In fact, for $Y_{m+1}^{(0)}$
the arguments are exactly the same as above.
For $Y_{m+1}^{(i)}$ with $i>0$ we have to be careful, since $Y_{m+1}^{(i)}$ consists of irreducible
components of the strict transform of $Y_m^{(i)}$, i.e., the blowup of $Y_m^{(i)}$ in $Y_m^{(i)}\times_{X_m} Y_m^{(0)}$,
which is a zero-dimensional scheme with a possibly non-reduced structure.
But since $Y_m^{(i)}$ is regular of dimension at most $1$ and
$Y_m^{(0)}\cap [Y_m^{(i)}]_0 =\emptyset$ for $i>0$, $Y_{m+1}^{(i)}
\rightarrow Y_m^{(i)}$ is an isomorphism .
As for the second part, $Y_{m+1}^{(m+1)}$ consists of finitely many closed points which,
by definition, are not contained in $Y_{m+1}^{j)}$ for $j\leq m$.

\medskip
\textbf{Step 7} Next we blow up $X_{m+1}$ in $Y_{m+1}^{(0)}$ if this is non-empty, and in $Y_{m+1}^{(1)}$, otherwise, and obtain $X_{m+2}$.
We proceed in this way for $n>m$, blowing up $X_n$ in $Y_n^{(j)}$ where $j \geq 0$ is
the smallest number with $Y_n^{(j)}\neq \emptyset$, to obtain $X_{n+1}$. This is well-defined, because we
always have:

\medskip
{\bf Claim 3} For all $n\geq m$ and $i \geq 0$, the subschemes $Y_n^{(i)}$ are regular, of dimension at most 1,
and $\cB_n$-permissible. For $i \geq m+1$, the intersection of $Y_n^{(i)}$ with $Y_n^{(j)}$ is empty for $j \neq i$.

\medskip
The first part follows as for $n=m+1$. For $i=n$ the second part follows as in claim 2.
For $m+1\leq i < n$, we may assume by induction that
$Y_{n-1}^{(i)}\cap Y_{n-1}^{(j)}=\emptyset$ for all $j=0,\dots, n-1$ with $j\neq i$.
By definition, for all $j=0,\dots, n-1$,
$\pi(Y_n^{(j)}) \subset Y_{n-1}^{(j)}$ where $\pi: X_n\to X_{n-1}$.
This implies $Y_n^{(i)}\cap Y_n^{(j)}=\emptyset$
for all $i\in \{m+1,\dots, n-1\}$ and all $j=0,\dots, n-1$ with $j\neq i$, which proves
the desired assertion.

\medskip
\textbf{Step 8} Thus we have defined the wanted canonical sequence $S(X,\tnu)$, which is reduced
by construction. Now we show the finiteness of this sequence. We have

\begin{lemma}\label{lem.curve.bu}
Let $X=X_0$ be a scheme satisfying the assumptions of Theorem \ref{thm.4.3}. Let $C=C_0$ be an irreducible regular
curve in $\tHSa X$. Let $\pi_1: X_1 = \Bl_C(X_0) \rightarrow X_0$, and let $C_1 = \tHSa {X_1}\cap \pi_0^{-1}(C_0)$.
By Step 2, $\dim C_1\leq 1$ and if $\dim C_1=1$, then $C_1$ is regular,
$C_1\subset X_1$ is $\cB$-permissible and $C_1\simeq C_0$.
In this case we put $X_2=\Bl_{C_1}(X_1)$.
Repeat this procedure to get a sequence
\begin{equation}\label{eq.curve.bu}
\begin{array}{ccccccccccccccc}
X= & X_0 & \lmapo{\pi_1} &   X_1 & \lmapo{\pi_2} & X_2 &
\leftarrow\ldots \leftarrow & X_{m-1} & \lmapo{\pi_m} & X_m & \ldots\\
& \cup & & \cup && \cup & & \cup && \cup & \\
C= & C_0 & \lmapo{\sim} & C_1 & \lmapo{\sim} & C_2 &
\lmapo{\sim} \ldots \lmapo{\sim} & C_{m-1}
& \lmapo{\sim} & C_{m} & \ldots \\
\end{array}
\end{equation}
such that $\pi_i: X_i = \Bl_{C_{i-1}}(X_{i-1})\rightarrow X_{i-1}$ is the blow-up morphism, and $C_i = \tHSa {X_i}\cap \pi_i^{-1}(C_{i-1})$.

Then the process stops after finitely many steps, i.e., there is an $m\geq 0$ with
$C_0 \stackrel{\sim}{\leftarrow} C_1 \stackrel{\sim}{\leftarrow} \ldots \stackrel{\sim}{\leftarrow} C_m$ and
$\dim(C_{m+1}) \leq 0$.
\end{lemma}

\textbf{Proof } Let $\eta$ be the generic point of $C$. As remarked in Step 2, we have $\eob_\eta(X)\leq 1$.
If $\eob_\eta(X) = 0$, then $C_1=\emptyset$ so that $r=1$. If $\eob_\eta(X)=1$, we get a longer sequence,
 which however must be finite by Theorem \ref{fu.thm0}/Corollary \ref{cor.fu.thm0} below, applied to the localization
 $X_\eta = \Spec(\cO_{X,\eta})$ of $X$ at $\eta$, and the point $\eta$ in it, for which $\tHSmax {(X_\eta)} = \{\eta\}$.
 Note that $\eob_\eta(X_\eta) = \eob_\eta(X)$ by definition.

\medskip
By this result, there is an $N\geq 0$ such that $Y_n^{(i)} \cap Y_n^{(j)} = \emptyset$ for all
$i,j \in\{0,\ldots,m\}$ with $i\neq j$, for all $n\geq N$, because $[Y_m^{(i)}]_0 \cap Y_m^{(j)} = \emptyset$
for all $i\neq j$, where $[Y_m^{(i)}]_0$ is the set of zero-dimensional components in $Y_n^{(i)}$.
Note that all schemes $X_n$ satisfy the conditions in Theorem \ref{thm.4.3}.
Therefore we have shown:

\medskip
{\bf Claim 4} There is an $N > 0$ such that $Y_n = \tHSa {X_n}$ is regular for all $n\geq N$.

\medskip
It is clear that the resolution sequence at each step $X_n$ has the following property, because
the centers of the blowups always lie in the subscheme $Y_n$: Let $Y_{n,1},\ldots,Y_{n,s}$ be the
connected components of $Y_n$, and for each $i\in \{1,\ldots,s\}$, let $V_{n,i} \subset X_n$ be an
open subscheme containing $Y_{n,i}$ but not meeting $Y_{n,j}$ for $j\neq i$. Then the resolution
sequence for $X$ is obtained by glueing the resolution sequences for the subsets $V_{n,i}$.
To show finiteness of the resolution sequence we may thus assume that $Y_n$ is
regular and irreducible. Applying Lemma \ref{lem.curve.bu} again, we may assume that $Y_n$ is a
collection of finitely many closed points which are isolated in their $\OB$-Hilbert-Samuel stratum.
Moreover, it is clear that in this case the remaining part of the resolution sequence
$X_n \leftarrow X_{n+1} \leftarrow \ldots$ is just the canonical resolution sequence
$S(X_n,\tnu)$ for $X_n$.

\medskip
\textbf{Step 9} Thus we have reduced to the case of an isolated point $x \in \tHSa X$; in fact, we may
assume that $\tHSa X$ just consists of $x$.
The first step of the canonical sequence then is to form the blowup $X_1 = \Bl_x(X) \rightarrow X$.
If $\eob_x(X) = 0$, then $\tHSa {X_1} = \emptyset$ and we are done. If $\eob_x(X) = 1$, then $\tHSa {X_1}$
is empty or consists of a unique point $x_1$ lying above $x$. In the latter case we have $k(x_1)=k(x)$,
and therefore $e_{x_1}(X_1) \leq e_x(X)$ by Theorem \ref{thm.pbu.inv} (4). If $e_{x_1}(X_1)=1$, then
$\eob_{x_1}(X_1) \leq e_{x_1}(X_1)\leq 1$. Otherwise we must have $e_{x_1}(X_1)=e_x(X)=2$ by assumption $(3e)$.
Then $\eob_{x_1}(X_1) \leq \eob_x(X)=1$ by Theorem \ref{thm.cNBncDir} (1). Thus we obtain a sequence
of blow-ups $X=X_0 \leftarrow X_1 \leftarrow \ldots$ in points $x_i\in X_i(\tnu)$ such that either
$e_{x_n}(X_n)=0$ for some $n$ so that $X_{n+1}(\tnu)=\emptyset$ and the sequence stops, or we have
a sequence where $e_{x_i}(X_i)=1$ for all $i$. But this sequence must be finite by
Theorem \ref{fu.thm0}/Corollary \ref{cor.fu.thm0} below. It remains the case where
$\eob_x(X) = e_x(X) = \eb_x(X) = 2$. This follows from Theorem \ref{Thm2} below.

\medskip
\textbf{Step 10} As a last step we show the functoriality, i.e., the properties (F1) and (F2) in
Theorem \ref{thm.4.3}. Property (F1) is of course an easier special case of (F2), and was only
written for reference. We show more generally:

\begin{proposition}\label{prop.funct.strat}
The canonical sequence defined in Remark
\ref{rem.strategy} is functorial for arbitrary flat morphisms with geometrically regular fibers $\alpha: X^\ast \rightarrow X$,
in the following sense: Let $\tnu \in \tSigmaXmax$, and let $\mathcal X = (X_n,\cB_n,\OB_n,Y_n,D_n)$ be the sequence $S(X,\tnu)$
as defined in Remark \ref{rem.strategy}, so that
$$
Y_n = X_n(\tnu) \qaq X_{n+1} =  \Bl_{D_n}(X_n)\,.
$$
If $X^\ast(\tnu)\neq \emptyset$, then $\tnu \in \Sigma_{X^\ast}^{\OB^\ast,max}$, where $\OB^\ast$ is the natural history function
induced on $X^\ast$, see Lemma \ref{lem.Operm.nemb}, and we let ${\mathcal X}^\ast = (X^\ast_n,\cB^\ast_n,\OB^\ast_n,Y^\ast_n,D^\ast_n)$ be the analogous sequence $S(X^\ast,\tnu)$ for $X^\ast$ and $\tnu$. If $X^\ast(\tnu) = \emptyset$, let
$\mathcal X^\ast = (X^\ast,\cB^\ast,\OB^\ast,\emptyset,\emptyset)$ (sequence consisting of one term), where $\cB^\ast$ and $\OB^\ast$ are obtained by pull-back of $\OB$ and $\cB$ to $X^\ast$ (notations as in Lemma \ref{lem.Operm.nemb}).
For any closed subscheme $Z \subseteq X_n$, let
$$
Z_\ast = X^\ast \times _X Z
$$
be its base change with $X^\ast \rightarrow X$. (Note that all $X_n$ are $X$-schemes in a canonical way.)
In particular let
$$
{\mathcal X}_\ast = (X_{n,\ast},\cB_{n,\ast},\OB_{n,\ast},Y_{n,\ast},D_{n,\ast})
$$
be the base change of $\mathcal X$ with $X^\ast \rightarrow X$, so that
$X_{n,\ast} = X^\ast\times_X X_n, Y_{n,\ast}= X^\ast\times_XY_n$ etc..\\
Then $\cX^\ast$ is canonically isomorphic to the associated contracted
sequence $(\cX_\ast)_{contr}$ (obtained by omitting the isomorphisms).
More precisely we claim the following.

(i) There is a canonical isomorphism $X^\ast(\tnu) \rmapo{\sim} X(\tnu)_\ast$.
Define the function $\varphi: \bN \rightarrow \bN$
inductively by $\varphi(0)=0$ and the following property for all $n\geq 1$:
\begin{equation}\label{def.varphi}
\varphi(n) =\left.\left\{\gathered
\varphi(n-1) \\
\varphi(n-1)+1 \\
\endgathered\right.\quad
\begin{aligned}
&\text{if $D_{n,\ast} = \emptyset$}\\
&\text{otherwise}.
\end{aligned}\right.
\end{equation}
Then the following holds:
\begin{itemize}
\item[(a)]
There are canonical isomorphisms $X_{n+1,\ast} \cong \Bl_{D_{n,\ast}}(X_{n,\ast})$.
\item[(b)]
Let $\beta_0: X^\ast_0 = X^\ast \rightarrow X^\ast\times_XX = X_{0,\ast}$ be the canonical
isomorphism induced by $\alpha: X^\ast \rightarrow X$. For all $n\geq 1$ there are unique morphisms
$\beta_n: (X^\ast)_{\varphi(n)} \rightarrow X_{n,\ast}$ such that all diagrams
\begin{equation}\label{eq.alpha}
\begin{CD}
X^\ast_{\varphi(n)} @>{\beta_n}>> X_{n,\ast} \\
@VVV @VVV \\
X^\ast_{\varphi(n-1)} @>{\beta_{n-1}}>> X_{n-1,\ast}
\end{CD}
\end{equation}
are commutative. Moreover all $\beta_n$ are isomorphisms. Here the left hand morphism
is the identity for $\varphi(n) = \varphi(n-1)$, and the morphism occurring in $\cX^\ast$
for $\varphi(n) = \varphi(n-1)+1$. The right hand morphism is the base change of $X_{n} \rightarrow X_{n-1}$
with $X^\ast \rightarrow X$.
\item[(c)]
The morphism $\beta_n$ induces isomorphisms
\begin{equation}\label{eq.Y.n}
Y^\ast_{\varphi(n)} \rmapo{\sim} Y_{n,\ast}
\end{equation}
and
\begin{equation}\label{eq.Y.n.j}
(Y^\ast_{\varphi(n)})^{(\varphi(j))} \rmapo{\sim} X^\ast\times_X Y_n^{(j)}=(Y_n^{(j)})_\ast
\end{equation}
if $(Y_n^{(j)})_\ast \neq \emptyset$ and
\begin{equation}\label{eq.D.n}
D^\ast_{\varphi(n)} \rmapo{\sim} D_{n,\ast}
\end{equation}
if $D_{n,\ast} \neq \emptyset$.
\end{itemize}
(ii) The analogous statement holds, if the $\tnu$-elimination sequence of Remark \ref{rem.strategy}
is replaced by the $\Sigma^{\OB,max}$-elimination sequence deduced from these $\tnu$-eliminations by gluing,
and for the canonical resolution sequences outlined in the proof of Corollary \ref{cor.Omax-elim}
or Corollary \ref{cor.Omax-elim1} using the mentioned max-eliminations.
\end{proposition}

\bigskip
\textbf{Proof }
We prove both claims by induction on dimension. It is quite clear how to deduce
(ii) from (i) once the dimension is fixed, and both claims are trivial for dimension zero.
Therefore it suffices to prove (i) when (ii) is assumed for smaller dimension.
Since $X^\ast \rightarrow X = X_0$ is flat with regular fibers, the morphism
\begin{equation}\label{eq.X}
X^\ast(\tnu) \rightarrow X^\ast\times_X X(\tnu) = Y_{0,\ast}
\end{equation}
is an isomorphism by Lemma \ref{lem.Operm.nemb} (2), and property (a) follows since $X_{n+1} = \Bl_{D_n}(X_n)$,
and $X^\ast \rightarrow X$ is flat (compare the proof of Lemma \ref{lem.Operm.nemb} (4)).
For the other claims we use induction on $n$.

\smallskip
(b) is empty for $n=0$, and we show (c) for $n=0$. The isomorphism \eqref{eq.Y.n} follows immediately from \eqref{eq.X},
since $Y_0^\ast = X^\ast(\tnu)$ by definition.
Moreover, $Y_0 = Y_0^{(0)}$ and $Y^\ast_0 = (Y^\ast_0)^{(0)}$ so that \eqref{eq.Y.n.j} follows trivially for $n=0$.
Now we show \eqref{eq.D.n} for $n=0$.
By the strategy of Remark \ref{rem.strategy}, the following holds. If $Y_0$ is $\cB$-regular, then the first
resolution cycle is already completed, and we have to blow up in $Y_0$, so that $D_0= Y_0$.
Since $Y^\ast_0 \rightarrow Y_0$ is flat with geometrically regular fibers,
we conclude that $Y^\ast_0$ is regular as well, and $D^\ast_0=Y^\ast_0$, so that
\eqref{eq.D.n} follows from \eqref{eq.Y.n}. If $Y_0$ is not $\cB$-regular, then $D_0$ is determined as the first
center of the canonical resolution sequence for $(Y_0,\cB_{Y_0})$, so that $D_0 = \cHSmax Y$. Similarly,
$D^\ast_0 = \cHSmax {(Y^\ast)}$.
Now we note again that $Y^\ast \rightarrow Y$ is flat with regular fibers and that
$D_{0,\ast}=X^\ast\times_X D_0 = Y^\ast\times_Y D_0$. Hence, if $D_{0,\ast} \neq \emptyset$, then
$\cHSmax {(Y^\ast)} = Y^\ast\times_Y \cHSmax Y$ by Lemma \ref{lem.Operm.nemb} (2), which shows
\eqref{eq.D.n} for $n=0$.

\smallskip
For $n\geq 1$ suppose we have already constructed the morphisms $\beta_m$ for $m\leq n-1$ and proved (b) and (c) for them.
For $m\leq n-1$ define the morphisms
$$
\alpha_m: X^\ast_{\varphi(m)} \rmapo{\beta_m} X_{m,\ast} = X^\ast\times_X X_m \rmapo{pr_2} X_m
$$
as the composition of the isomorphism $\beta_m$ with the canonical projection.
Then $\alpha_m$ is flat with regular fibers, because this holds for $pr_2$.

If $D_{n-1,\ast}=\emptyset$, then $X_{n,\ast} \rmapo{\sim} X_{n-1,\ast}$ is an isomorphism
by (a), we have $\varphi(n) = \varphi(n-1)$, and there is a unique morphism $\beta_n$,
necessarily an isomorphism, making the diagram in (c) commutative, and trivially cartesian.
Define $\alpha_n$ as the composition of $\beta_n$ with the projection to $X_n$.

If $D_{n-1,\ast}\neq \emptyset$, then $\varphi(n)=\varphi(n-1)+1$, and $\beta_{n-1}$ induces an isomorphism
$$
D^\ast_{\varphi(n-1)} \rmapo{\sim} D_{n-1,\ast}=X^\ast\times_X D_{n-1} \cong X^\ast_{\varphi(n-1)}\times_{X_{n-1}} D_{n-1}\,,
$$
by induction. As $\alpha_{n-1}$ is flat, there is a unique morphism $\alpha_n$ making the diagram
\begin{equation}\label{eq.beta}
\begin{CD}
\Bl_{D^\ast_{n-1}}(X^\ast_{\varphi(n-1)}) = X^\ast_{\varphi(n)} @>{\alpha_n}>> X_n = \Bl_{D_{n-1}}(X_{n-1}) \\
@VVV @VVV \\
X^\ast_{\varphi(n-1)} @>{\alpha_{n-1}}>> X_{n-1}
\end{CD}
\end{equation}
commutative, and the diagram is cartesian (see Lemma \ref{lem.Operm.nemb} (4)). By base change
of the right column with $X^\ast$ over $X$ we obtain a diagram as wanted in (b).

By flatness of $\alpha_n$ and regularity of its fibers, and by Lemma \ref{lem.Operm.nemb} (2), we have the isomorphism
\begin{equation}\label{eq.Y.n+1}
Y^\ast_{\varphi(n)} = X^\ast(\tnu) \mathop{\longrightarrow}\limits^{\sim} X^\ast_{n-1} \times_{X_{n-1}} Y_n \cong X^\ast\times_X Y_n = Y_{n,\ast}
\end{equation}
in \eqref{eq.Y.n} in both cases ($D_{n-1,\ast}$ empty or non-empty). We note that we have a cartesian diagram
\begin{equation}\label{eq.base.change.Y}
\begin{CD}
Y^\ast_{\varphi(n)} @>{\alpha_n}>> Y_n \\
@V{\pi^\ast}VV @VV{\pi}V \\
Y^\ast_{\varphi(n-1)} @>{\alpha_{n-1}}>> Y_{n-1}\,.
\end{CD}
\end{equation}
In fact, we have
$$
Y^\ast_{\varphi(n-1)} \times_{Y_{n-1}} Y_n \cong (X^\ast\times_X Y_{n-1})\times_{Y_{n-1}} Y_n \cong X^\ast\times_X Y_n \cong
Y^\ast_{\varphi(n)}\,.
$$
Now we show \eqref{eq.Y.n.j} for $n$, i.e., that the isomorphism \eqref{eq.Y.n+1} induces isomorphisms
\begin{equation}\label{eq.Y.n+1.j}
(Y^\ast_{\varphi(n)})^{(\varphi(j))} \mathop{\longrightarrow}\limits^{\sim} X^\ast_{\varphi(n-1)}\times_{X_{n-1}} Y_n^{(j)} \cong X^\ast\times_X Y_n^{(j)}
\end{equation}
for all $j\leq n$. By the factorization
$$
\alpha_{n-1}: Y^\ast_{\varphi(n-1)} \mathop{\longrightarrow}\limits^{\sim} X^\ast\times_X Y_{n-1} \mathop{\longrightarrow}\limits^{pr_2} Y_{n-1}\,,
$$
condition \eqref{eq.Y.n.j} for $n-1$ means that $(\alpha_{n-1})^{-1}(Y_{n-1}^{(j)}) = (Y^\ast_{\varphi(n-1)})^{(\varphi(j))}$ for all $j\leq n-1$,
or, equivalently, that all generic points of $(Y^\ast_{\varphi(n-1)})^{(\varphi(j))}$ map to generic points of $Y_{n-1}^{(j)}$,
and we have to show the corresponding property for $n$ in place of $n-1$, where we can assume that $D_{n-1,\ast}\neq \emptyset$.
Since $\alpha_{n-1}$ and $\alpha_n$ are flat, the generic points of $Y^\ast_{\varphi(n-1)}$ are those points
which lie over generic points $\eta$ of $Y_{n-1}$ and are generic points in the fiber over $\eta$,
and the analogous statement holds for $n$ in place of $n-1$.

Let $\varphi(j)\leq \varphi(n)=\varphi(n-1)+1$, and let $\eta^\ast_{n}$ be a generic point of $(Y^\ast_{\varphi(n)})^{(\varphi(j))}$.
If $\varphi(j) \leq \varphi(n-1)$, then $\eta^\ast_{n}$ maps to a generic point $\eta^\ast_{n-1}$ of $(Y^\ast_{\varphi(n-1)})^{(\varphi(j))}$
which in turn maps to a generic point $\eta_{n-1}$ of $Y_{n-1}^{(j)}$. Let $\eta_n$ be the image of $\eta^\ast_n$ in $Y_n$. Then $\eta_n$
is a generic point of $Y_n$ and lies in $Y_n^{(j)}$ since it maps to $\eta_{n-1}$. Note that in \eqref{eq.Y.n+1} the fibers of $\pi^\ast$
are obtained from those of $\pi$ by base change with a field extension. If $\varphi(j)=\varphi(n-1)+1 (=\varphi(n))$, then $\eta_n$
is still a generic point of $Y_n$, but $\eta^\ast_{n-1}$ is not a generic point of $Y^\ast_{n-1}$. By \eqref{eq.base.change.Y}, the fiber $F_n$ of $\alpha_n$ over $\eta_n$ is flat over the fiber $F_{n-1}$ of $\alpha_{n-1}$ over $\eta_{n-1}$, because it is obtained by base change with
the residue field extension $k(\eta_n)/k(\eta_{n-1})$. Therefore $\eta^\ast_n$ is a generic point of the fiber over $\eta_n$. This implies
that $\eta_{n-1}$ is not a generic point of $Y_{n-1}$, and hence that $\eta_n$ lies in $Y_n^{(n)}$.

\medskip
Finally we show that $\alpha_n$ induces an isomorphism
$$
D^\ast_{\varphi(n)} \rmapo{\sim} D_{n,\ast}=X^\ast\times_X D_n \cong X^\ast_{\varphi(n)}\times_{X_n} D_{n}\,,
$$
provided that $D_{n,\ast}\neq \emptyset$. Suppose that
\begin{equation}\label{eq.D.Y}
D_n \subset Y_n^{(j)}\,.
\end{equation}
By the prescription
in Remark \ref{rem.strategy} this means that either $j=0$ or that $j>0$ and $Y_n^{(j-1)}=\emptyset$.
In the first case we have $\varphi(0)=0$, and in the second case we have $(Y^\ast_{\varphi(n)})^{(\varphi(j-1))}=\emptyset$
by \eqref{eq.Y.n+1.j}. As $\varphi(j)-1 \leq \varphi(j-1)$, we then also
have $(Y^\ast_{\varphi(n)})^{(\varphi(j)-1)}=\emptyset$. By \eqref{eq.D.Y} we get $Y_n^{(j)}\neq\emptyset$ and
hence also $(Y^\ast_{\varphi(n)})^{(\varphi(j))}\neq \emptyset$ by \eqref{eq.Y.n+1.j}. Therefore in both cases we have
$D^\ast_n \subset (Y^\ast_n)^{(\varphi(j))}$ by our strategy.

\bigskip
If $Y_n^{(j)}$ is regular, then $D_n=Y_n^{(j)}$ by our strategy, and $(Y^\ast_n)^{(\varphi(j))}$ is regular,
too, because $(Y^\ast_n)^{(\varphi(j))}\rightarrow Y_n^{(j)}$ is flat with regular fibers by \eqref{eq.Y.n.j}.
Then also $D^\ast_n = (Y^\ast_n)^{(\varphi(j))}$, and the claim follows from \eqref{eq.Y.n+1.j}.


\bigskip
If $Y_n^{(j)}$ is not regular, then $D_n$ lies in the singular locus of $Y_n^{(j)}$ by our strategy,
and we are in one of the resolution cycles for $Y_?^{(j)}$. Suppose this cycle has started with $Y_r^{(j)}$
($r \leq n$). Then we obtain the centers from the canonical resolution sequence for $Y_r^{(j)}$.

\bigskip
First assume that $r=0$. Then $j=0$, and we start a resolution cycle with $Y^\ast_0=(Y^\ast_0)^{(0)}$ as well.
By induction on dimension we can apply (ii) to $Y=Y_0$, and we have a canonical commutative diagram with cartesian squares
$$
\begin{array}{cccccccccc}
Y_n  & \rightarrow &    Y_{n-1}   & \rightarrow &  \ldots  & \rightarrow &   Y_1    & \rightarrow & Y_0 & = Y\\
\uparrow &             & \uparrow &             &  \ldots  &             & \uparrow &             & \uparrow & \\
Y^\ast_{\psi(n)} & \rightarrow & Y^\ast_{\psi(n-1)} & \rightarrow & \ldots & \rightarrow & Y^\ast_{\psi(1)} & \rightarrow & Y^\ast_0 & = Y^\ast
\end{array}
$$
where the upper line is part of the canonical resolution sequence for $Y$, the $Y^\ast_m$ in the lower line
come from the canonical resolution sequence of $Y^\ast_0$, and the function $\psi$ is defined in an analogous way
as $\varphi$ in \eqref{def.varphi}, but using $Y^\ast \rightarrow Y$ and the $D_{i-1} \subseteq Y_{i-1}$ for all $i\leq n$
instead of $X^\ast\rightarrow X$ and the $D_{i-1} \subseteq X_{i-1}$: Note that, by definition of the resolution cycle, the center of the
blow-up $Y_i \rightarrow Y_{i-1}$ is the same center $D_{i-1}$ as for the blow-up $X_i \rightarrow X_{i-1}$.
Again by induction the vertical morphisms induce isomorphisms $D_{\psi(i)} \rmapo{\sim} Y^\ast\times_Y D_{i}$
for all $i\leq n$ if the last scheme is non-empty.
We claim that $\psi(i) = \varphi(i)$, and that $D^\ast_{\varphi(i)} \rightarrow D_i$ coincides with the
morphism induced by $\alpha_i: X^\ast_{\varphi(i)} \rightarrow X_i$, for all $i\leq n$, which then
proves the claim, since $Y^\ast\times_Y D_i = X^\ast\times_X D_i$.

The claim is true for $i=0$, and we assume it is true up to some $i\leq n-1$. But then we have
$Y^\ast\times_Y D_i=X^\ast\times_X D_i = D_{i,\ast}$, and hence the first scheme is empty if and
only if the last one is, which shows $\psi(i+1)=\varphi(i+1)$ since $\psi(i) = \varphi(i)$.
Moreover the existence and uniqueness of the morphisms $Y^\ast_{\varphi(i+1)}\rightarrow Y_{\ast,i+1}$
induced by the above diagram and by $\alpha_i$ implies the remaining claim.

\bigskip
If $r>0$, then $Y_{r-1}^{(j)}$ is regular, $D_{r-1}=Y_{r-1}^{(j)}$ and $X_r = \Bl_{D_{r-1}}(X_{r-1})$,
and by induction and \eqref{eq.Y.n.j},
$(Y^\ast_{\varphi(r-1)})^{(\varphi(j))}$ is regular as well (note that $r \leq n$ and
that $(Y^\ast_{\varphi(r-1)})^{(\varphi(j))}\rightarrow Y_{r-1}^{(j)}$ is flat with regular fibers).
Hence, by our strategy, with $(Y^\ast_{\varphi(r)})^{(\varphi(j))}$ we start a resolution cycle for $X^\ast$ as well.
Now the proof is the same as for the cycle starting with $Y_0$, up to some renumbering.
$\square$

\begin{corollary}\label{cor.strategy}
If in the situation of Proposition \ref{prop.funct.strat} the $\tnu$-elimination sequence
(resp. the $\Sigma^{\OB,max}$-elimination sequence, resp. the resolution sequence) is finite
for $X$, then it is also finite for $X^\ast$.
\end{corollary}

\bigskip\noindent
We now turn to the two key theorems used in the proof above. Let $X$ be an excellent scheme, and
let $(\cB,\OB)$ be a boundary with history on $X$. Let $x \in X$ and assume that
\begin{itemize}
\item[$(F1)$]
$\Char(k(x)) = 0\text{  or  } \Char(k(x)) \geq \dim (X)/2+1.$
\item[$(F2)$]
$\NBx\pitchfork \Dirob_x(X).$
\item[$(F3)$]
$\eob_x(X)\leq 1$ or $e_x(X)=\eb_x(X)$.
\end{itemize}

\medbreak
Consider
$$
\pi_1 : X_1=\Bl_x(X)\to X,
$$
$$
C_1:=\bP(\Dirob_x(X))\subset \bP(\Dir_x(X)) \subset X_1.
$$
Let $\eta_1$ be the generic point of $C_1$. We note
$C_1\simeq \bP_k^{t-1}$, where $t=\eob_x(X)$.
By Theorem \ref{thm.cBdirectrix}, any point of $X_1$
which is $\OB$-near to $x$, lies in $C_1$.

\begin{lemma}\label{lem0.fupb}
\begin{itemize}
\item[(1)]
If $\eta_1$ is $\OB$-near to $x$, then so is any point of $C_1$.
\item[(2)]
If $\eta_1$ is very $\OB$-near to $x$, then so is any point of $C_1$.
\end{itemize}
\end{lemma}
\textbf{Proof }
Take any point $y\in C_1$. By Theorems \ref{thm.cBHSf.usc} and
\ref{thm.cBpbu.inv}, we have
$$
\tHSfXd {\eta_1} \leq \tHSfXd y \leq \tHSfX x.
$$
(1) follows from this.
By Theorems \ref{nf.dir} and \ref{thm.pbu.inv}, we have
$$
e_{\eta_1}(X')\leq e_y(X') -\dim(\cO_{C_1,y})
\leq e_x(X)-\delta_{y/x} -\dim(\cO_{C_1,y})=e_x(X)-\delta_{\eta_1/x},
$$
where we used $(F3)$ for the second inequality. In fact, if $\eob_x(X)\leq 1$ then $k(y)=k(x)$,
so in both cases of $(F3)$ we can apply \ref{thm.pbu.inv} (4). Hence the assumption of (2)
implies that $y$ is very near to $x$. Then,
by Theorems \ref{nf.cBdir} and \ref{thm.cBpbu.inv}, we get
$$
\eob_{\eta_1}(X')\leq \eob_y(X') -\dim(\cO_{C_1,y})
\leq \eob_x(X)-\delta_{y/x} -\dim(\cO_{C_1,y})=\eob_x(X)-\delta_{\eta_1/x},
$$
which implies the conclusion of (2).
$\square$

\bigskip
We now assume that $t\geq 1$ and $\eta_1$ is very $\OB$-near to $x$.
This implies
\begin{equation}\label{eq1.esbpbu}
\eob_{\eta_1}(X_1)=\eob_x(X)-\dim(C_1)=1.
\end{equation}
By $(F2)$, Lemma \ref{lem.cNBncDir1} and Remark \ref{rem.cOperm}, $C_1$ is $\cB$-permissible with
respect to the complete transform $(\cB_1,\OB_1)$ of $(\cB,\OB)$ for $X_1$. Consider the blow-up
$$
\pi_2: X_2=\Bl_{C_1}(X_1)\to X_1
$$
and the complete transform $(\cB_2,\OB_2)$ of $(\cB_1,\OB_1)$ for $X_2$.
By \eqref{eq1.esbpbu} and Theorem \ref{thm.cBdirectrix},
there is at most one point $\eta_2\in X_2$ which is $\OB$-near to $\eta_1$.
If it exists, let $C_2$ be the closure of $\eta_2$ in $X_2$. Then $k(\eta_1)=k(\eta_2)$,
$\pi_2$ induces an isomorphism $C_2\cong C_1$, and $C_2$ is $\OB_2$-permissible.
\def\NBy{N(y)}
By Lemma \ref{lem0.fupb} and Theorem \ref{thm.cNBncDir}, $(F2)$ implies
$\NBy\pitchfork \Dir_y(X_1)$ for any point $y\in C_1$.
By Lemma \ref{lem.cNBncDir2}, $C_2$ is n.c. with $N\cB_2(y') = \cB_2(y')-\OB_2(y')$
at the unique point $y' \in C_2$ above $y$, so that $C_2$ is $\cB_2$-permissible.
Consider
$$
\pi_3 : X_3=\Bl_{C_2}(X_2)\to X_2
$$
and proceed in the same way as before.
This construction (which occurred in the proof of Theorem \ref{thm.4.3} for
$t = 1,2$) leads us to the following:

\begin{definition}\label{Def.fupb}
Assume $\eob_x(X) \geq 1$ and let $m$ be a non-negative integer or $\infty$.
The fundamental sequence of $\cB$-permissible blowups over $x$
of length $m$
is the canonical (possibly infinite) sequence of permissible blowups:
\begin{equation}\label{esbpbu.eq}
\begin{array}{ccccccccccccccc}
\cB=&\cB_0&& \cB_1&& \cB_2&& \cB_{n-1}&& \cB_{n} \\
\\
X= & X_0 & \lmapo{\pi_1} &   X_1 & \lmapo{\pi_2} & X_2 &
\leftarrow\ldots \leftarrow & X_{n-1} & \lmapo{\pi_n} & X_n & \leftarrow \ldots\\
& \uparrow & & \cup && \cup & & \cup && \cup \\
& x & \leftarrow & C_1 & \stackrel{\sim}{\leftarrow} & C_2 &
\stackrel{\sim}{\leftarrow} \ldots \stackrel{\sim}{\leftarrow} & C_{n-1}
&\stackrel{\sim}{\leftarrow} & C_{n} & \leftarrow \ldots\\
\end{array}
\end{equation}
which satisfies the following conditions:
\begin{itemize}
\item[$(i)$]
$X_1=\Bl_x(X)$ and
$$
C_1=\bP(\Dirob_x(X)) \cong\bP^{t-1}_{k(x)}\quad (t=\eob_x(X))\,.
$$
\item[$(ii)$]
For $1\leq q < m$,
$$
C_q=\{\xi\in \phi^{-1}_q(x)|\; \tHSf{X_q}(\xi)=\tHSf{X}(x)\}\qwith
\phi_q:X_q\rightarrow X \,.
$$
Let $\eta_q$ be the generic point of $C_q$.
\item[$(iii)$]
For $2\leq q < m$, $X_q=\Bl_{C_{q-1}}(X_{q-1})$
and $\pi_q: C_q \isom C_{q-1}$ is an isomorphism.
\item[$(iv)$]
If $m=1$, then the generic point of $C_1$ is not $\OB$-near to $x$.
If $1<m<\infty$, then $X_m=\Bl_{C_{m-1}}(X_{m-1})$ and there is no point in $X_m$ which is near to $\eta_{m-1}$.
If $m=\infty$, then the sequence is infinite.
\end{itemize}
Here $(\cB_q,\OB_q)$ is the complete transform of $(\cB_{q-1},\OB_{q-1})$.
\end{definition}

\bigskip

The proof of the following first key theorem will be given in \S\ref{sec:fund.seq}.

\begin{theorem}\label{fu.thm0}
Assume that there is no regular closed subscheme
$$
D \subseteq \{\xi\in \Spec(\cO_{X,x}) \; |\; \tHSfX {\xi}\geq \tHSfX {x}\}
$$
of dimension $\eob_x(X)$.
Then, for the sequence \eqref{esbpbu.eq}, we have $m<\infty$, i.e., it stops in finitely many steps.
\end{theorem}

\begin{remark}\label{rem.fu.thm0}
We note that the assumption of the theorem holds in particular, if
\begin{equation}\label{fu.thm0.eq1}
\dim\big(\{\xi\in \Spec(\cO_{X,x}) \; |\; \tHSfX {\xi}\geq
\tHSfX {x}\}\big)< \eob_x(X).
\end{equation}
\end{remark}

Thus a special case of Theorem \ref{fu.thm0} is the following.

\begin{corollary}\label{cor.fu.thm0}
If $x$ is isolated in the $\OB$-Hilbert-Samuel locus of $X$ and $\eob_x(X) = 1$, then the
fundamental sequence \eqref{esbpbu.eq} consists of a sequence of blowups in closed points
and is finite.
\end{corollary}

In particular, if $\dim(X)=2$, we obtain the canonical sequences constructed in steps 8 and 9 of
the proof of Theorem \ref{thm.4.3}.
Hence we obtain their finiteness as needed in that proof.

\bigskip
Now we consider the fundamental sequence of $\cB$-permissible blowups over $x$
as in Definition \ref{Def.fupb} for the second case needed in the proof
of Theorem \ref{thm.4.3}, namely where $x$ is isolated in $X(\tnu)$ and
$$
\eob_x(X)=e_x(X)=\eb_x(X)=2.
$$
Here we have $C_q \cong \bP^1_{k(x)}$ for $1\leq q < m$.
Again by Theorem \ref{fu.thm0} (and Remark \ref{rem.fu.thm0} (a))
we deduce that the fundamental sequence \eqref{esbpbu.eq}
is finite, i.e., there exists an $m<\infty$ such that there is no
point in $X_m$ which is $\OB$-near to $\eta_{m-1}$.
Let $\tnu = \tHSfX {x}$. If $X_m(\tnu)$ is non-empty, all its points
lie above $X_{m-1}(\tnu)$ (by maximality of $\tnu$), and the image in
$X_{m-1}(\tnu)$ consists of finitely many closed points (since it does
not contain the generic point of in $X_{m-1}(\tnu)$). By the argument used
in Step 2 of the proof of Theorem \ref{thm.4.3}, each of these points has
at most one point in $X_m(\tnu)$ above it, so that the latter set consist
of finitely many closed points as well. Pick one such point $y$. Since we
already treated the case where $y$ is isolated with $\eob_y(X_m)=1$, and
we always have $\eob_y(X_m)\leq 2$, we are led to the
following definition.

\begin{definition}\label{Def.fupb2}
A sequence of $\cB$-permissible blowups $(\cX,\cB)$:
$$
\begin{array}{ccccccccccccccc}
 \cB=&\cB_0&& \cB_1&& \cB_2&& \cB_{m-1}&& \cB_{m} \\
 \\
X= & X_0 & \lmapo{\pi_1} &   X_1 & \lmapo{\pi_2} & X_2 &
\leftarrow\ldots \leftarrow & X_{m-1} & \lmapo{\pi_m} & X_m\\
& \uparrow & & \cup && \cup & & \cup && \uparrow \\
& x & \leftarrow & C_1 & \stackrel{\sim}{\leftarrow} & C_2 &
\stackrel{\sim}{\leftarrow} \ldots \stackrel{\sim}{\leftarrow} & C_{m-1}
&\gets & x_m \\
\end{array}
$$
is called a fundamental unit of $\cB$-permissible blowups of length $m$
if the following conditions are satisfied:
\begin{itemize}
\item[$(i)$]
$x$ is a closed point of $X$ such that $\eob_x(X)=e_x(X)=\eb_x(X)=2$.
\item[$(ii)$]
$X_1=\Bl_x(X)$ and
$$
C_1=\bP(\Dirob_x(X)) \cong\bP^1_{k(x)}\,.
$$
\item[$(iii)$]
For $2\leq q\leq m$, $C_q = \phi^{-1}_q(x)\cap X_q(\tnu)$ and $X_q = \Bl_{C_{q-1}}(X_{q-1})$,
where  $\tnu = \tHSf{X}(x)$ and $\phi_q:X_q\rightarrow X$ is the natural morphism\,.
\item[$(iv)$]
For $2\leq q < m$, $\pi_q$ induces an isomorphism $C_q \cong C_{q-1}$.
\item[$(v)$]
$\pi_m: C_m \rightarrow C_{m-1}$ is not surjective.
\item[$(vi)$]
$x_m$ is a closed point of $X_m$ above $x$ such that
$$
\tHSf{X_m}(x_m)=\tHSf{X}(x) \qaq
\eob_{x_m}(X_m)=e_{x_m}(X_m)=\eb_{x_m}(X_m)=2.
$$
\end{itemize}
Here $\tHSf{X_q}$ is considered for the successive complete transform $(\cB_q,\OB_q)$ of $(\cB,\OB)$ for $X_q$.
\medbreak

By convention, a fundamental unit of $\cB$-permissible blowups of length 1 is
a sequence of $\cB$-permissible blowups such as
$$
\begin{array}{ccccccc}
 \cB&=&\cB_0&& \cB_1&& \\
 \\
X&=&X_0&\stackrel{\pi_1}{\leftarrow} & X_1&=&\Bl_x(X)\\
\uparrow && \uparrow & & \uparrow\\
x& =& x_0 & \leftarrow & x_1
\end{array}
$$
where $x\in X$ is as in $(i)$ and $x_1$ is as in $(vi)$ with $m=1$.\\
We call $(x,X,\cB)$ (resp. $(x_m,X_m,\cB_m)$) the initial
(resp. terminal) part of $(\cX,\cB)$.
\end{definition}
\medbreak

We remark that, in this definition, we have not assumed that $x$ (resp. $x_m$) is isolated
in $\tHSmax X$ (resp. $\tHSmaxXm$).

\begin{definition}\label{Def.fspb}
A chain of fundamental units of $\cB$-permissible blowups is
a sequence of $\cB$-permissible blowups:
$$
\cX_1\gets\cX_2\gets\cX_3\gets \ldots
$$
where $\cX_i=(\cX_i,\cB_i)$ is a fundamental unit of $\cB$-permissible blowups such that the terminal part of $\cX_i$ coincides with initial part
of $\cX_{i+1}$ for $\forall \, i\geq 1$.
\end{definition}

The finiteness of the canonical $\tnu$-elimination $S(X,\tnu)$ for the case where
$\tHSa X = \{x\}$ and $\eob_x(X)=e_x(X)=\eb_x(X)=2$, as needed in the proof
of Theorem \ref{thm.4.3}, is now a consequence of the following second key theorem
whose proof will be given in \S\ref{sec:e2II} and \S\ref{sec:e2III}.

\begin{theorem}\label{Thm2} Let
$\cX_1\gets\cX_2\gets\cX_3\gets \ldots $
be a chain of fundamental units of $\cB$-permissible blowups. Let
$(x^{(i)}, X^{(i)}, B^{(i)})$ be the initial part of $(\cX_i)$ for
$i\geq 0$. Assume that, for each $i$, there is no regular closed subscheme $C\subseteq \tHSmaxXi$
of dimension 1 with $x^{(i)} \in C$ (which holds if $x^{(i)}$ is isolated in $\tHSmaxXi$).
Then the chain must stop after finitely many steps.
\end{theorem}

In fact, to show the finiteness of $S(X,\tnu)$ in the considered case, we have to show
that there is no infinite sequence of closed points $x_n\in X_n(\tnu)$ such that $x_0=x$
and $x_{n+1}$ lies above $x_n$. This can only happen if
$\eob_{x_n}(X_n)=2\; (=e_{x_n}(X_n) = \eb_{x_n}(X_n))$ for all $n$, and by construction,
the canonical sequence $S(X,\tnu)$ would then give rise to an infinite chain of fundamental
units.

\bigbreak
We remark that the claims on the fundamental sequences, fundamental units and
chains of fundamental units depend only on the localization $X_x = \Spec(\cO_{X,x})$
of $X$ at $x$.
Moreover, by the results in Lemma \ref{completion.nu.2} and Lemma \ref{completion.HS}
we may assume that $X = \Spec(\cO)$ for a complete local ring. Thus we assume that there
is an embedding $X \hookrightarrow Z$ into a regular excellent scheme $Z$, and moreover,
by Lemma \ref{lem.comp.emb-nemb}, that there is a simple normal crossings
boundary $\cB_Z$ on $Z$ whose pull-back to $X$ is $\cB$. Thus we may consider an embedded
version of the constructions above, where each blowup $X_{n+1} = \Bl_{C_n}(X_n) \rightarrow X_n$
(where $C_0 = \{x\}$) can be embedded into a diagram
$$
\begin{matrix}
Z_{m+1}&=& \Bl_{C_m}(Z_m) &\rmapo{\pi_{m+1}}& Z_m \\
&& \cup  && \cup\\
X_{m+1}&=& \Bl_{C_m}(X_m) &\rmapo{\pi_{m+1}}& X_m \,.
\end{matrix}
$$
In the proofs of Theorems \ref{fu.thm0} and \ref{Thm2}, this situation will be assumed.

\newpage
\section{$(u)$-standard bases}\label{sec:u-st.base}

\bigskip

Let $R$ be a regular noetherian local ring with maximal ideal $\fm$ and
residue field $k=R/\fm$, and let $J \subseteq \fm$ be an ideal.
It turns out that the directrix $\Dir(R/J)$ is an important
invariant of the singularity of $X=\Spec(R/J)$, and that it is
useful to consider a system $(y_1, \dots, y_r, u_1, \dots, u_e)$
of regular parameters for $R$ such that:
\begin{equation}\label{eq.adim}
\text{
$\IDir(R/J)=\langle Y_1, \dots, Y_r \rangle \;\subset\grm(R)$, where
$Y_i := \inm(y_i)\in \grm^1(R)$.
}
\end{equation}
Then $(U_1, \dots, U_e)$ with $U_j := \inm(u_j)$ form coordinates of
the affine space $\Dir(R/J) \cong \bA^e_k$.
Consequently, it will be useful to distinguish the $Y$- and $U$-coordinates
in $\grm(R)=k[Y,U]$. This observation leads us to the following:

\begin{definition}\label{def2.11}
\begin{itemize}
\item[(1)]
 A system $(y,u)=(y_1,\ldots , y_r, u_1,\ldots , u_e)$ of regular
parameters for $R$ is called strictly admissible for $J$
if it satisfies the above condition \eqref{eq.adim}.
\item[(2)]
A sequence $u=(u_1,\ldots , u_e)$ of elements in $\fm\subseteq R$ is called
admissible (for $J$), if it can be extended to a strictly admissible system
$(y,u)$ for $J$.
\item[(3)]
 Let $(y,u)$ be strictly admissible for $J$, and let $f=(f_1,\dots, f_m)$ be a system of elements in $J$.
 Then $(f,y,u)$ is called admissible if $in_\fm(f_i)\in k[Y]$ for all $i=1,\dots,m$.
\end{itemize}
\end{definition}

Let $T_1,\dots,T_e$ be a tuple of new variables over $k$.
Note that $(u)=(u_1,\ldots , u_e)$ is admissible if and only if
we have an isomorphism of $k$-algebras
$$
k[T_1,\dots,T_e] \isom \grm(R)/\IDir(R/J)
\;;\;T_i\to in_{\fm}(u_i)\mod \IDir(R/J).
$$
The map induces the following isomorphism which we will use later.
\begin{equation}\label{eq.u-coord}
\psi_{(u)}:
\bP(\Dir(R/J)) \isom \Proj(k[T_1,\dots,T_e])=\bP^{e-1}_k
\end{equation}
\medbreak

The admissibility will play an essential role in the next section.
For the moment we shall work in the following general setup:

\vspace{0,5cm} \textbf{Setup A:} \quad
Let $J\subset\fm\subset R$ be as above.
Let $(u) = (u_1, \dots, u_e)$ be a system of elements in $\fm$ such that
$(u)$ can be extended to a system of regular parameters $(y, u)$ for some
$y = (y_1, \dots, y_r)$.
In what follows we fix $u$ and work with various choices of $y$ as above.
Such a choice induces an identification
$$
\gr_\fm(R) = k[Y,U]=k[Y_1\dots, Y_r, U_1,\dots, U_e].\quad
(Y_i=in_\fm(y_i),\; U_j=in_\fm(u_j)).
$$
Let
$\tilde{R} = R/\langle u \rangle$ and $\tilde{\fm} = {\fm}/\langle u \rangle$
where $\langle u \rangle=\langle u_1, \dots, u_e\rangle \; \subset R$.
For $f\in R-\{0\}$ put
$$
n_{(u)}(f) = v_{{\tilde{\fm}}} (\tilde{f})
\qwith \tilde{f} = f \mod \langle u \rangle\;\in \tilde{R}.
$$
Note $n_{(u)}(f)\geq v_{\fm}(f)$ and $n_{(u)}(f)=\infty$ if and only if
$f\in \langle u \rangle$.
Let $f\in R-\{0\}$ and write an expansion in the $\fm$-adic completion $\hat{R}$ of $R$:
\begin{equation}\label{eq1.2.1}
f = \sum\limits_{(A, B)} C_{A, B} \; y^B u^A \quad \mbox{with}
\quad C_{A, B} \in R^{\times} \cup \{ 0 \}
\end{equation}
where for $A= (a_1, \dots, a_e) \in {\mathbb Z}_{\geq 0}^e$ and
$B = (b_1, \dots,b_r) \in {\mathbb Z}_{\geq 0}^r$,
$$
y^B = y_1^{b_1} \dots y_r^{b_r}\quad \text{ and }\quad
u^A = u_1^{a_1} \dots u_e^{a_e}.
$$
If $n_{(u)}(f)<\infty$, then we define the $0$-initial form of $f$ by:
\begin{equation}\label{eq.0intial}
in_0(f) = in_0(f)_{(y,u)} =
 \underset{ \overset{B}{|B|=n_{(u)}(f)} }{\sum}
\overline{C}_{0, B}\; Y^B \in k[Y],
\end{equation}
where $\overline{C_{A, B}}=C_{A, B}\mod\fm \;\in k =R/\fm$.
If $n_{(u)}(f) =\infty$, we define $in_0(f)_{(y,u)} = 0$.
It is easy to see that $in_0(f)$ depends only on $(y,u)$, not on
the presentation \eqref{eq1.2.1}.
\medbreak

We will need to make the expansion \eqref{eq1.2.1} uniquely determined by $f$.
By \cite{EGAIV} Ch. 0 Th. (19.8.8), we can choose a ring $S$
of coefficients of $\hat{R}$: $S$ is a subring of $\hat{R}$
which is a complete local ring with the maximal ideal
$pS$ where $p=\Char(k)$ such that $\fm\cap S=pS$ and $S/pS=R/\fm$.
We choose a set $\Gamma\subset S$ of representatives of $k$.
Note that $S\simeq k$ and $\Gamma=k$ if
$\Char(k)=\Char(K)$ where $K$ is the quotient field of $R$.
Back in the general situation each $f\in R$ is expanded in $\hat{R}$ in a unique way as:
\begin{equation}\label{eq1.2.1.c}
f = \sum\limits_{(A, B)} C_{A, B} \; y^B u^A \quad \mbox{with}
\quad C_{A, B} \in \Gamma.
\end{equation}
We will use the following map of sets
\begin{equation}\label{eq1.2.1.d}
\omega=\omega_{(y,u,\Gamma)}:
k[[Y,U]] \to \hat{R}\;;\;
\underset{(A, B)}{\sum} c_{A, B} \; Y^B U^A \to
\underset{(A, B)}{\sum} C_{A, B} \; y^B u^A ,
\end{equation}
where $C_{A,B}\in \Gamma$ is the representative of $c_{A,B}\in k$.
For $F,G\in k[[Y,U]]$ we have
\begin{equation}\label{eq1.2.1.e}
\omega(F+G)-\omega(F)-\omega(G)\in p\hat{R} \qaq
\omega(F\cdot G)-\omega(F)\cdot \omega(G)\in p\hat{R}.
\end{equation}
\bigskip

We now introduce the notion of a $(u)$-standard base
(see Definition \ref{def1.4}),
which generalizes that of a standard base (cf. Definition \ref{def0.5}).
The following facts are crucial: Under a permissible blowup a standard base
is not necessarily transformed into a standard base but into
a $(u)$-standard base (see Theorem \ref{prop2.1}), on the other hand
there is a standard procedure to transform a $(u)$-standard base into
a standard base (see Theorem \ref{cor.wellprepared}).
\medbreak

A linear form $L : {\mathbb R}^e \longrightarrow {\mathbb R}$ given by
$$
L(A) = \sum\limits_{i = 1}^e c_i a_i
\text{  with } c_i \in \bR \quad (A = (a_i) \in {\mathbb R}^e)
$$
is called positive (resp. semi-positive) if $c_i>0$ (resp. $c_i\geq 0$)
for all $1\leq i\leq e$.

\begin{definition}\label{def1.3}
Let $(y,u)$ be as in Setup A and let $L$ be a non-zero semi-positive linear
form $L$ on $\mathbb R^e$.
\begin{itemize}
\item[(1)]
For $f\in \hat{R}-\{0\}$ define the $L$-valuation of $f$
with respect to $(y,u)$ as:
$$
v_L (f) = v_L(f)_{(y, u)} :=
\min \{ |B| + L(A) \; \vert \; C_{A, B} \ne 0 \},
$$
where the $C_{A, B}$ come from a presentation \eqref{eq1.2.1.c} and
$|B|=b_1+\cdots b_r$ for $B=(b_1,\dots,b_r)$.
We set $v_L(f)=\infty$ if $f=0$.
\item[(2)]
Fix a representative $\Gamma$ of $k$ in $\hat{R}$ as in \eqref{eq1.2.1.c}.
The initial form of $f\in \hat{R}-\{0\}$ with respect to
$L$, $(y,u)$ and $\Gamma$ is defined as:
$$
in_L (f) = in_L (f)_{(y, u,\Gamma)} :=
\sum\limits_{A, B} \overline{C_{A, B}} \; Y^B U^A $$
where
$A,B$ range over $\bZ_{\geq 0}^e \times \bZ_{\geq 0}^r$ satisfying
$|B|+ L (A) = v_L (f)$. We set $in_L(f)=0$ if $f=0$.
\item[(3)]
For an ideal $J\subset R$, we define
$$
\InhLJ = \InhLJ_{(y, u,\Gamma)} = \langle in_L(f)|\; f\in J \rangle
\subset k[[U]][Y].
$$
In case $L$ is positive we define
$$
\InLJ = \InLJ_{(y, u,\Gamma)} = \langle in_L(f)|\; f\in J \rangle
\subset k[U,Y]=\grm(R).
$$
\end{itemize}
\end{definition}
\medskip
It is easy to see that this is well-defined, i.e., one has the following:
\begin{itemize}
\item[$(i)$]
$in_L (f)$ is an element of $k[[U]][Y]$,
the polynomial ring of $Y$ with coefficients in the formal power series ring
$k[[U]]$,
\item[$(ii)$]
If $L$ is positive, $in_L (f)\in k[U,Y]=\grm(R)$, and is independent
of the choice of $\Gamma$.
\end{itemize}

\begin{remark}\label{rem1.1}
Note that $v_{\fm} (f) = v_{L_0} (f)$ and
$in_{\fm} (f) = in_{L_0} (f)$, where
$$
L_0 (A) = |A| = a_1 + \dots + a_e \quad\text{ for }
A = (a_1, \dots,a_e).
$$
\end{remark}
\medbreak

The proofs of the following lemmas \ref{lem1.1} and \ref{L0initial} are easy
and left to the readers.

\begin{lemma}\label{lem1.1}
Let the assumptions be as in Definition \ref{def1.3}.
\begin{itemize}
\item[(1)]
$v_L (f)$ is independent of the choice of $\Gamma$. We have
$$v_L (fg)  = v_L (f) + v_L (g) \quad\text{ and }\quad
v_L (f+g)  \geq \min \{ v_L (f), v_L (g)\}.$$
\item[(2)]
Assume $v_L(\Char(k))>0$ (which is automatic if $L$ is positive).
If $f = \sum\limits_{i = 1}^m f_i$ and $v_L (f_i) \geq
v_L (f)$ for all $i = 1, \dots, m$, then
$in_L (f) = \underset{1\leq i\leq m}{\sum} in_L (f_i)$, where
the sum ranges over such $i$ that $v_L (f_i) = v_L (f)$.
\item[(3)]
Let $z=(z_1,\dots,z_r)\subset R$ be another system of parameters such that
$(z,u)$ is a system of regular parameters of $R$.
Assume $v_L(z_i-y_i)_{(z,u)}\geq 1$ for all $i=1,\dots,r$.
Then, for any $f\in R$, we have
$v_L(f)_{(z,u)}\geq v_L(f)_{(y,u)}$.
\end{itemize}
\end{lemma}
\medbreak

\begin{definition}\label{def.effective}
Let the assumption be as Definition \ref{def1.3}. Let $f=(f_1,\dots.f_m)$
be a system of elements in $R-\{0\}$.
A non-zero semi-positive linear form $L$ on $\mathbb R^e$ is called
effective for $(f,y,u)$ if $in_L(f_i)\in k[Y]$ for all $i=1,\dots, m$.
\end{definition}

\begin{lemma}\label{L0initial}
Let the assumption be as in Definition \ref{def.effective}.
\begin{itemize}
\item[(1)]
The following conditions are equivalent:
\begin{itemize}
\item[(i)]
$L$ is effective for $(f,y,u)$.
\item[(ii)]
$v_L(f_i)=n_{(u)}(f_i)<\infty$ and $in_L(f_i)=in_0(f_i)$
for all $i=1,\dots, m$.
\item[(iii)]
For $1\leq i\leq m$ write as \eqref{eq1.2.1}
\begin{equation*}
f_i = \sum\limits_{(A, B)} C_{i,A, B} \; y^B u^A \quad \mbox{with}
\quad C_{i,A, B} \in R^{\times} \cup \{ 0 \}\;.
\end{equation*}
Then, for $1\leq i\leq m$ and $A\in \bZ^e_{\geq 0}$ and $B\in \bZ^r_{\geq 0}$,
we have
$$
|B|+ L(A)> n_{(u)}(f_i)
\text{   if $|B|< n_{(u)}(f_i)$ and $C_{i,A,B}\not=0$.}
$$
\end{itemize}
\item[(2)]
There exist a positive linear form $L$ on $\bR^e$ effective for $(f,y,u)$
if and only if $f_i$ is not contained in $\langle u \rangle\subset R$ for any
$1\leq i\leq m$.
\item[(3)]
If $L$ is effective for $(f,y,u)$ and $\Lambda$ is a linear form such that
$\Lambda \geq L$, then $\Lambda$ is effective for $(f,y,u)$. More precisely
one has the following for $1 \leq i \leq m$:
$$ v_{\Lambda} (f_i) = v_L (f_i)=n_{(u)}(f_i) \quad\text{ and }
 \quad in_{\Lambda} (f_i) = in_L (f_i) =in_0(f_i)
$$
\end{itemize}
\end{lemma}

\begin{definition}\label{def1.4}
Let $u$ be as in Setup A.
Let $f=(f_1, \dots, f_m)\subset J$ be a system of elements in $R-\{0\}$.
\begin{itemize}
\item[(1)]
$f$ is called a $(u)$-effective base of $J$, if there is a tuple
$y =(y_1, \dots, y_r)$ as in Setup A and a positive form $L$ on $\mathbb R^e$
such that $L$ is effective for $(f,y,u)$ and
\begin{equation*}\label{eq1.6.2}
In_L (J) = \langle in_0 (f_1),\dots, in_0 (f_m)\rangle \; \subset \grmR.
\end{equation*}

\item[(2)]
$f$ is called a $(u)$-standard base, if in addition
$(in_0 (f_1), \dots, in_0 (f_m))$ is a standard base of $In_L(J)$.
\medbreak

In both cases (1) and (2), $(y, L)$ is called a reference datum for the
$(u)$-effective (or $(u)$-standard) base $(f_1, \dots, f_m)$.
\end{itemize}
\end{definition}

\begin{lemma}\label{rem1.2}
Let $u$ be as in Setup A.
\begin{itemize}
\item[(1)]
Let $f=(f_1,\dots, f_m)$ be a standard base of $J$ such that
$\inm(f_i)\in k[Y]$ for $i=1,\dots,m$.
Then $f$ is a $(u)$-standard base of $J$
with reference datum $(y, L_0)$, where $L_0$ is as in Remark \ref{rem1.1}.
\item[(2)]
Assume that $(u)$ be admissible for $J$ (cf. Definition \ref{def2.11}).
A standard base $f=(f_1,\dots,f_m)$ of $J$ is a $(u)$-standard base of $J$.
\end{itemize}
\end{lemma}
\medbreak\noindent
\textbf{Proof }
(1) is an immediate consequence of Definition \ref{def2.11}(3) and
Remark \ref{rem1.1}. We show (2). First we note that
$$
\nu^*(J)=(n_1,\dots,n_m,\infty,\dots)\qwith n_i=v_\fm(f_i).
$$
and that $(\inm(f_1),\dots,\inm(f_m))$ is a standard base of $\InmJ$.
Choose $y=(y_1,\dots,y_r)$ such that $(y,u)$ is strictly admissible for
$J\subset R$ and identify $\grmR=k[Y,U]$. Then there exist
$\psi_1,\dots,\psi_m\in k[Y]$ which form a standard base of $\InmJ$.
Note that $\psi_i$ is homogeneous of degree $n_i$ for $i=1,\dots,m$. We have
$$ \langle \inm(f_1),\dots,\inm(f_m)\rangle =\InmJ=
\langle \psi_1,\dots,\psi_m\rangle.$$
Writing
$$
\inm(f_i)=\phi_i +\underset{A\in \bZ_{\geq 0}^e-\{0\}}{\sum} U^A P_{i,A},
\quad \text{with } \phi_i,\; P_{i,A}\in k[Y]\;\; (i=1,\dots,m),
$$
this implies that $(\phi_1,\dots,\phi_m)$ is a standard base of $\InmJ$.
Hence it suffices to show that there exists a positive linear form
$L:\bR^e\to \bR$ such that $in_L(f_i)=\phi_i$ for all $i=1,\dots,m$.
We may write
$$
P_{i,A}=\underset{B\in \bZ_{\geq 0}^r}{\sum} c_{i,A,B} Y^B,\quad
(c_{i,A,B}\in k),
$$
where the sum ranges over $B\in \bZ_{\geq 0}^r$ such that
$|B|+|A|=n_i:=v_\fm(f_i)$. It is easy to see that there exists
a positive linear form $L$ satisfying the following for all $A\in \bZ_{\geq 0}^e-\{0\}$:
$$
L(A)>|A|\qaq L(\frac{A}{n_i-|B|})>1\text{  if } c_{i,A,B}\not=0.
$$
Then $in_L(f_i)=\phi_i$ and the proof of Lemma \ref{rem1.2} is complete.
$\square$
\bigskip

A crucial fact on $(u)$-standard bases is the following:

\begin{theorem}\label{refdatum0}
Let $f=(f_1,\dots, f_m)$ be a $(u)$-effective (resp. standard) base of $J$.
Then, for any $y=(y_1\dots,y_r)$ as in Setup A and for any positive linear
form $L$ on $\mathbb R^e$ effective for $(f,y,u)$, $(y,L)$ is a reference
datum for $f$.
\end{theorem}

Before going to the proof of Theorem \ref{refdatum0},
we deduce the following:

\begin{corollary}\label{cor.refdatum0}
Let $f=(f_1,\dots, f_m)$ be a $(u)$-effective (resp. standard) base of $J$.
Let $g=(g_1,\dots, g_m)\subset J$ be such that
$in_0(g_i)=in_0(f_i)$ for all $i=1,\dots,m$.
Then $g$ is a $(u)$-effective (resp. standard) base of $J$.
\end{corollary}
\textbf{Proof}\quad
The assumption implies that no $f_i$ or $g_i$ is contained in
$\langle u \rangle\subset R$. By Lemma \ref{L0initial} (2) and (3) there
exists a positive
linear form $L$ on $\bR^e$ effective for both $(f,y,u)$ and $(g,y,u)$.
By the assumption on $f$, Theorem \ref{refdatum0} implies that
$in_L(f_1),\dots,in_L(f_m)$ generate (resp. form a standard base of) $in_L(J)$.
By Lemma \ref{L0initial} (1)$(ii)$, we get
$$
in_L(g_i)=in_0(g_i)=in_0(f_i)=in_L(f_i)\qfor i=1,\dots,m,
$$
which implies Corollary \ref{cor.refdatum0}.
$\square$

\bigskip
\textbf{Proof of Theorem \ref{refdatum0}}\quad
Let $(z,\Lambda)$ be a reference datum for $f$ which exists by the assumption.
By definition $(z,u)$ is a system of regular parameters of $R$ and $\Lambda$
is a positive linear form on $\mathbb R^e$
such that $in_\Lambda(f_1)_{(z,u)},\dots,in_\Lambda(f_m)_{(z,u)}$
generate (resp. form a standard base of) $in_\Lambda(J)_{(z,u)}$.
First we assume $y=z$. Then the theorem follows from Proposition \ref{prop1.1}
below in view of Lemma \ref{L0initial} (Note that condition (4) of
Proposition \ref{prop1.1} is always satisfied if $L$ and $\Lambda$ are
both positive). We consider the general case.
Since $(y,u)$ and $(z,u)$ are both systems of regular parameters of $R$,
there exists $M=(\alpha_{ij})\in GL_r(R)$ such that
$$
y_i=l_i(z) + d_i,\quad\text{ where }
l_i(z) = \sum_{j=1}^r \alpha_{ij} z_j \text{ and } d_i\in \langle u \rangle.
$$
Take any positive linear form $L'$ on $\mathbb R^e$ such that
$L'(A)> |A|$ for $\forall A\in \mathbb R^e_{\geq 0}-\{0\}$.
An easy computation shows that for $A\in \mathbb Z^e_{\geq 0}$ and
$B=(b_1,\dots,b_r) \in \mathbb Z^r_{\geq 0}$ we have
$$
y^B u^A = u^A \cdot l_1(z)^{b_1}\cdots l_1(z)^{b_1} + w \quad\text{ with}
\quad v_{L'}(w)>v_{L'}(y^B u^A ) = |B|+ L'(A).
$$
This implies $in_{L'}(g)_{(z,u)} =\phi(in_{L'}(g)_{(y,u)})$ for
$g\in R-\{0\}$, where
$$
\phi:k[Y]\simeq k[Z] \;;\; Y_i \to \sum_{j=1}^r \overline{\alpha}_{ij} Z_j.
\quad
(Z_j=in_\fm(z_j),\; \overline{\alpha}_{ij}=\alpha_{ij}\mod\fm \in k)
$$
In view of Lemma \ref{L0initial}, the proof of Theorem \ref{refdatum0} is
now reduced to the case $y=z$.
$\square$

\begin{proposition}\label{prop1.1}
Let $(y,u)$ be as in Setup A, and let $f=(f_1,\dots, f_m)\subset J$.
Let $\Lambda$ and $L$ be semi-positive linear forms on
${\mathbb R}^e$. Assume $v_\Lambda(\Char(k))>0$ (cf. Lemma \ref{lem1.1})
Assume further the following conditions:
\begin{itemize}
\item[(1)]
$v_{\Lambda}(f_i)=v_L(f_i)=n_{(u)}(f_i)<\infty$ for $i=1,\dots,m$.
\item[(2)]
$in_{\Lambda} (f_i) = in_0 (f_i):=F_i(Y)\in k[Y]$ for $i=1,\dots,m$.
\item[(3)]
$\Inh_{\Lambda} (J) = \langle F_1(Y), \dots, F_m(Y)\rangle \subset k[[U]][Y]$.
\end{itemize}
Then for any $g \in J$ and any $M \geq 0$, there exist
$\lambda_1, \dots,\lambda_m \in \hat{R}$ such that
\begin{equation*}\label{eq1.9.1}
v_L (\lambda_i f_i) \geq v_L (g),\quad
v_\Lambda(\lambda_i f_i) \geq v_\Lambda(g),\quad
v_\Lambda \left(g - \sum\limits_{i = 1}^m \lambda_i f_i \right) > M.
\end{equation*}
If $\Lambda$ is positive, one can take $\lambda_i\in R$ for $i=1,\dots,m$.
Assume further:
\begin{itemize}
\item[(4)]
there exist $c>0$ such that $L \geq c \Lambda$.
\end{itemize}
Then we have
$$\Inh_L(J) = \langle in_L(f_1), \dots, in_L(f_m)\rangle \subset k[[U]][Y].$$
If $L$ is positive, then
$$In_L(J) = \langle in_L(f_1), \dots, in_L(f_m)\rangle \subset k[U,Y].$$
\end{proposition}
\medbreak
\textbf{Proof  }
Let $g \in J$ and expand as in \eqref{eq1.2.1.c}:
$$
g = \sum\limits_{A, B} C_{A, B}y^B u^A \quad \mbox{in} \; \hat{R},
\quad C_{A, B} \in \Gamma .
$$
Then we have
$$
in_{\Lambda} (g) =
\underset{\stackrel{A, B}{|B| + \Lambda (A)= v_{\Lambda} (g)}}{\sum}
 \overline{C_{A, B}} \quad Y^B U^A.
$$
Put
$$
B_{\max} = B_{\max} (g, \Lambda)
= \max \{ B \; \vert \; \overline{C_{A,B}} \ne 0,
\; |B| + \Lambda (A) = v_{\Lambda} (g) \text{ for some } A\in \bZ^e_{\geq0}\}
$$
where the maximum is taken with respect to the lexicographic order.

\begin{lemma}\label{lem1.2}
Under the assumption (1) and (2) of Proposition \ref{prop1.1},
there exist $\lambda_1, \dots, \lambda_m \in \hat{R}$ such that:
\begin{itemize}
\item[(1)] $v_{\Lambda} (\lambda_i f_i) = v_{\Lambda} (g)$ if
$\lambda_i \ne 0$, and $v_L (\lambda_i f_i) \geq v_L (g)$.
\item[(2)]
 $v_L (g_1) \geq v_L (g)$ where $g_1 =g-\sum\limits_{i = 1}^m \lambda_i f_i$.
\item[(3)]
 $v_{\Lambda} (g_1) \geq v_{\Lambda} (g)$ and $B_{\max} (g_1, \Lambda) <
B_{\max} (g, \Lambda)$ if $v_{\Lambda} (g_1) = v_{\Lambda} (g)$.
\end{itemize}
If $\Lambda$ is positive, one can take $\lambda_i\in R$ for $i=1,\dots,m$.
\end{lemma}
\medbreak \textbf{Proof  }
By the assumptions of Proposition \ref{prop1.1}, we can write
$$
in_{\Lambda} (g) =
\underset{\stackrel{A, B}{|B| + \Lambda (A)= v_{\Lambda} (g)}}{\sum}
\overline{C_{A, B}} \quad Y^B U^A
= \underset{1\leq i\leq m}{\sum} H_i F_i(Y)
$$
for some $H_1,\dots, H_m \in k[[U]][Y]$, where $F_i(Y)\in k[Y]$ is homogeneous
of degree $n_i: = v_{\Lambda}(f_i)=v_L(f_i)$ for $i=1,\dots,m$.
Writing
$H_i=\underset{A\in \bZ^e_{\geq 0}}{\sum} h_{i,A}(Y) U^A$ with
$h_{i,A}(Y)\in k[Y]$, this implies
\begin{equation}\label{lem1.2.eq1}
\underset{\stackrel{B}{|B| + \Lambda (A)= v_{\Lambda} (g)}}{\sum}
\overline{C_{A, B}} \quad Y^B
= \underset{1\leq i\leq m}{\sum} h_{i,A}(Y) F_i(Y)
\quad\text{ for each } A\in \bZ^e_{\geq 0}.
\end{equation}
Take any $A_0$ with
$|B_{\max}|+\Lambda(A_0)=v_\Lambda(g)$ and $\overline{C_{A_0, B_{\max}}}\ne 0$.
Looking at the homogeneous part of degree $|B_{\max}|$ in \eqref{lem1.2.eq1}
with $A=A_0$, we get
$$
\underset{\overset{B}{|B|= |B_{\max}|}}{\sum} \overline{C_{A_0, B}} \; Y^B
= \underset{1\leq i\leq m}{\sum} S_i (Y) F_i (Y)
$$
where $S_i (Y)\in k[Y]$ is the homogenous part of degree
$|B_{\max}|-n_i$ of $h_{i,A_0}(Y)$. Therefore
\begin{equation}\label{eq1.10.1}
Y^{B_{\max}} - \sum\limits_{i = 1}^m P_i (Y) F_i(Y)
=\sum\limits_{\stackrel{\scriptstyle |B|=|B_{\max}|}{B\neq
B_{\max}}}a_B Y^B,
\end{equation}
where
$P_i (Y) = \left( \overline{C_{A_0, B_{\max}}} \right)^{-1}S_i (Y)\in k[Y]$
and $a_B=-(\overline{C_{A_0,B_{\max}}})^{-1}\overline{C_{A_0B}}\in k$.
Now put
$$
g_1 = g - \sum\limits_{i = 1}^m \lambda_i f_i\qwith
\lambda_i = \left\{
\begin{aligned}
\tilde{P}_i(y)\tilde{Q}(u) & \quad\text{ if  } P_i (Y) \ne 0, \\
0  & \quad\text{ if  }P_i (Y) = 0, \\
\end{aligned} \right.
$$
where $\tilde{P}_i(Y)\in R[Y]$ is a lift of $P_i(Y)\in k[Y]$ and
$$
\tilde{Q}(u) =
\underset{\overset{A}{|B_{\max}|+ \Lambda (A)=v_{\Lambda}(g)}}{\sum}
C_{A, B_{\max}} u^A\;\in \hat{R}.
$$
Note that if $\Lambda$ is positive, then the sum is finite and
$\tilde{Q}(u) \in R$ and $\lambda_i\in R$.
For $1\leq i\leq m$ with $P_i (Y) \ne 0$, we have
$$v_{\Lambda} (\lambda_i) =
(|B_{\max}| - n_i) + (v_{\Lambda} (g) -|B_{\max}|)=v_{\Lambda}(g)-n_i,
$$
$$
\begin{aligned}
v_L (\lambda_i)  = & (|B_{\max}| - n_i) + v_L (\tilde{Q}(u)) \\
 \geq & (|B_{\max}|-n_i) + v_L (g) - |B_{\max}|= v_L(g) - n_i\; , \\
\end{aligned}
$$
which shows Lemma \ref{lem1.2} (1) in view of Proposition \ref{prop1.1} (1).
Here the last inequality holds because
\begin{align*}
v_L (g) =& \min \{ |B|+L(A) \; \vert \;
\overline{C_{A, B}} \ne 0 \},\\
v_L (\tilde{Q}(u))+|B_{\max}| =& \min \{|B|+ L(A) \; \vert
\; \overline{C_{A, B}} \ne 0, \; B=B_{\max}, \;
|B_{\max}|+\Lambda (A) = v_{\Lambda} (g) \}.
\end{align*}
Therefore, by Lemma \ref{lem1.1} (1) we get
$$
v_{\Lambda} (g_1) \geq
\min \{ v_{\Lambda} (g), \; v_{\Lambda}(\lambda_if_i) \; (1\leq i\leq m)\}
\geq v_{\Lambda} (g).
$$
If $v_{\Lambda} (g_1) = v_{\Lambda} (g)$, then Lemma \ref{lem1.1} (2) implies
\begin{align*}
in_{\Lambda} (g_1) & =  in_{\Lambda} (g) + in_{\Lambda} \left(
\sum\limits_{i = 1}^m \lambda_i f_i \right) \\
& =  in_{\Lambda} (g) +
\underset{\overset{1\leq i\leq m}{\lambda_i \neq 0}}{\sum}
in_{\Lambda} (\lambda_i) in_{\Lambda} (f_i) \\
& =  in_{\Lambda} (g) - Q(U) \sum\limits_{i = 1}^m P_i(Y) F_i (Y),
\end{align*}
where
$$
Q(U)= \underset{\overset{A}{|B_{\max}|+ \Lambda (A)=v_{\Lambda}(g)}}{\sum}
 \overline{C_{A, B_{\max}}} U^A\;\in k[[U]].
$$
Hence, by \eqref{eq1.10.1}, we get
$B_{\max} (g_1, \Lambda) < B_{\max} (g, \Lambda)$, which proves
Lemma \ref{lem1.2} (3). Finally,
$$
v_L (g_1) \geq \min \{ v_L (g), v_L (\lambda_i f_i) \;(1\leq i\leq m)\}
\geq v_L (g),
$$
because $v_L (\lambda_i f_i) \geq v_L (g)$. This proves Lemma \ref{lem1.2} (2)
and the proof of the lemma is complete.
$\square$

\bigskip
We now proceed with the proof of Proposition \ref{prop1.1}.
From $g = g_0$ we construct $g_1$ as in Lemma
\ref{lem1.2}, and by applying Lemma \ref{lem1.2} repeatedly, we
get a sequence $g_0, g_1, g_2, \dots, g_{\ell}, \dots$ in $J$ such
that for all $\ell \geq 1$ we have
$$
g_{\ell} = g_{\ell -1} - \sum\limits_{i = 1}^m \lambda_{\ell,i}
f_i \quad \mbox{with} \quad \lambda_{\ell,i} \in R,
$$
$$
v_{\Lambda} (\lambda_{\ell,i} f_i)  \geq  v_{\Lambda} (g_{\ell-1})
\geq v_{\Lambda} (g) \; \mbox{and} \; v_L (\lambda_{\ell,i} f_i)
\geq v_L (g_{\ell-1}) \geq v_L (g),
$$
$$
v_L (g_{\ell}) \geq v_L (g_{\ell-1}) \;\mbox{and} \; v_{\Lambda}
(g_{\ell})  \geq  v_{\Lambda} (g_{\ell-1}),
$$
$$
B_{\max}(g_\ell,\Lambda)< B_{\max}(g_{\ell-1},\Lambda)\text{  if  }
v_{\Lambda} (g_{\ell})  = v_{\Lambda} (g_{\ell-1}).
$$
Then we have
$$
g_{\ell} = g - \sum\limits_{i = 1}^m \mu_{\ell,i} f_i\quad\mbox{
with } \mu_{\ell,i} = \sum\limits_{q = 1}^{\ell} \lambda_{q ,i},
$$
\begin{equation}\label{eq1.10.2}
\begin{aligned}
v_{\Lambda} (\mu_{\ell,i} f_i)  \geq & \min\limits_{1 \leq q \leq
\ell} \{v_{\Lambda} (\lambda_{q i} f_i)\} \geq v_{\Lambda} (g),\\
v_L (\mu_{\ell,i} f_i)  \geq & \min\limits_{1 \leq q \leq \ell}
\{v_L (\lambda_{q i} f_i)\} \geq v_L (g).\\
\end{aligned}
\end{equation}

Note that $B_{\max} (g_{\ell},\Lambda)$ cannot drop forever in the
lexicographic order so that we must have
$v_{\Lambda} (g_{\ell})  \not= v_{\Lambda} (g_{\ell-1})$ for infinitely
many $\ell$. Noting $v_\Lambda(R)$ is a discrete subset of $\bR$,
this implies that for given $M \geq 0$, taking $\ell$ sufficiently large,
$$
v_\Lambda(g_{\ell}) = v_\Lambda(g-\sum\limits_{i=1}^m \mu_{\ell,i} f_i)>M.
$$
This shows the first assertion of Proposition \ref{prop1.1}
in view of \eqref{eq1.10.2}.
It implies by (4) that for any $g \in J$, there exist
$\lambda_1, \dots,\lambda_m \in \hat{R}$ such that
\begin{equation*}
v_L (\lambda_i f_i) \geq v_L (g),\quad
v_L(g - \sum\limits_{i = 1}^m \lambda_i f_i ) > v_L(g),
\end{equation*}
which implies, by Lemma \ref{lem1.1} (2), that
$in_L (g) = \sum\limits_{i = 1}^m in_{L} (\lambda_i) in_{L} (f_i),$
where the sum ranges over all $i$ for which $v_L (\lambda_i f_i) =
v_L (g)$. This shows
$\Inh_L(J)=\langle in_L(f_1),\dots,in_L(f_m)\rangle$,
which also implies the last assertion of Proposition \ref{prop1.1}
by the faithful flatness of $k[[U]][Y]$ over $k[U,Y]$.
This completes the proof of Proposition \ref{prop1.1}
$\square$
\medbreak

\begin{lemma}\label{def1.4cor}
Let $u$ be as in Setup A and assume that $(u)$ is admissible for $J$
(cf. Definition \ref{def2.11}). Choose $y=(y_1,\dots,y_r)$ such that $(y,u)$ is strictly admissible for $J$.
\begin{itemize}
\item[(1)]
If $f=(f_1,\dots,f_m)$ is a $(u)$-effective base of $J$, then
$$
\InmJ=\langle in_0(f_1),\dots,in_0(f_m)\rangle \qaq J=(f_1,\dots,f_m).
$$
Here the equality makes sense via the isomorphism $gr_m(R) \simeq k[Y,U]$ induced by the above choice of $(y,u)$.
\item[(2)]
If $f=(f_1,\dots,f_m)$ is a $(u)$-standard base, then
$(in_0(f_1),\dots,in_0(f_m))$ is a standard base of $\InmJ$ and
$\nu^*(J)=(n_{(u)}(f_1),\dots, n_{(u)}(f_m),\infty,\infty,\dots).$
\end{itemize}
\end{lemma}
\medbreak\noindent
\textbf{Proof }
(2) follows at once from (1).
The second assertion of (1) follows from the first in view of
\cite{H3}, (2.21.d).
We now show the first assertion of (1).
Take a positive linear form $\Lambda$ which is effective for $(f,y,u)$.
By Theorem \ref{refdatum0}, $(y,\Lambda)$ is a reference datum for $f$.
By Definition \ref{def1.4} and Lemma \ref{L0initial} (1) this implies
$v_\Lambda(f_i)=n_{(u)}(f_i)<\infty$ for $i=1,\dots,m$ and
$$
In_\Lambda(J)=\langle F_1(Y),\dots,F_m(Y)\rangle
\qwith
F_i(Y)=in_\Lambda(f_i)=in_0(f_i)\in k[Y].
$$
It suffices to show $In_\Lambda(J)=In_\fm(J)$.
By the strict admissibility of $(y,u)$, Lemma \ref{stbTJ} implies that
there exists a standard base $g=(g_1,\dots,g_s)$
which is admissible for $(y,u)$ (cf. Definition \ref{def2.11}).
By Lemma \ref{rem1.2} (2) this implies
$v_{\fm}(g_i)=v_{L_0}(g_i)=n_{(u)}(g_i)<\infty$ for $i=1,\dots,s$ and
$$
In_\fm(J)=In_{L_0}(J)=\langle G_1(Y),\dots,G_s(Y)\rangle
\qwith
G_i(Y)=in_{L_0}(g_i)\in k[Y].
$$
Take a positive linear form $L$ such that $L\geq\Lambda$ and $L\geq L_0$.
Then Proposition \ref{prop1.1} and Lemma \ref{L0initial} (3) imply
$In_\Lambda(J)=In_L(J)=In_{L_0}(J)=In_\fm(J).$
This completes the proof.
$\square$
\medbreak

\newpage
\section{Characteristic polyhedra of $J\subset R$}
\def\ut{u_{\leq t}}
\def\Jsa{J_{s+1}}
\def\vsa{v_{s+1}}
\def\usa{u_{\leq s+1}}
\def\tlam{\widetilde{\lambda}}

\bigskip

In this section we are always in Setup A (beginning of \S\ref{sec:u-st.base}).
We introduce a polyhedron $\Delta(J,u)$ which plays a crucial role
in this paper. It will provide us with useful invariants of singularities
of $\Spec(R/J)$ (see \S\ref{sec:ad.inv}). It also give us a natural way
to transform a $(u)$-standard base of $J$ into a standard base of $J$
(see Corollary \ref{cor.wellprepared}).

\begin{definition}\label{def1.1}
\begin{itemize}
\item[(1)]
An $F$-subset $\Delta \subseteq {\bR}_{\geq 0}^e$ is a
closed convex subset of ${\bR}_{\geq 0}^e$ such that $v \in
\Delta$ implies $v + {\bR}_{\geq 0}^e \subseteq \Delta$.
The essential boundary $\partial\Delta$ of an $F$-subset $\Delta$ is
the subset of $\Delta$ consisting of those $v\in \Delta$ such that $v\not\in
v'+ {\bR}^e_{\geq 0}$ with $v'\in \Delta$ unless $v=v'$.
We write
$\Delta^+=\Delta - \partial\Delta$.
\item[(2)]
For a semi-positive positive linear form $L:\bR^e\to \bR$, put
$$
\delta_L(\Delta)=\min\{L(v)\;|\; v\in \Delta\}.
$$
Then
$E_L=\Delta\cap \{v\in \bR^e|\; L(v)=\delta_L(\Delta)\}$
is called a face of $\Delta$ with slope $L$.
One easily sees that $E_L$ is bounded if and only if $L$ is positive.
If $E_L$ consists of a unique point $v$, we call $v$ a vertex of $\Delta$.
\item[(3)]
When $L=L_0$ as in Remark \ref{rem1.1}, we call
$\delta(\Delta) = \delta_{L_0}(\Delta)=\min\{a_1+\ldots +a_e\;|$ $\; (a_1,\ldots,a_e)\in \Delta\}$
the $\delta$-invariant of $\Delta$ and $E_{L_0}$ the $\delta$-face of $\Delta$.
\end{itemize}
\end{definition}

\begin{definition}\label{def.intial}
Let $(y,u)$ be as in Setup A in \S\ref{sec:u-st.base}.
Let $f \in \fm$ be not contained in $\langle u \rangle \; \subset R$.
Write as in \eqref{eq1.2.1}:
\begin{equation*}
f = \sum\limits_{(A, B)} C_{A, B} \; y^B u^A \quad \mbox{with}
\quad C_{A, B} \in R^{\times} \cup \{ 0 \}.
\end{equation*}

\begin{itemize}
\item[(1)]
The polyhedron
$$\Delta (f, y, u) \subseteq {\bR}_{\geq 0}^e$$
is defined as the smallest $F$-subset containing all points of
$$
\left\{ v = \frac{A}{n_{(u)}(f) - \vert B \vert} \; \big\vert \; C_{A, B}
\ne 0, \vert B \vert < n_{(u)}(f) \right\}.
$$
This is in fact a polyhedron in ${\bR}_{\geq 0}^e$, which depends only on
$f,y,u$, and does not depend on the presentation \eqref{eq1.2.1}.
\item[(2)]
For $v\in \bR^e-\Delta(f,y,u)^+$, the $v$-initial of $f$ is defined as
$$
in_v(f) = in_v(f)_{(y,u)} = in_0(f) + in_v(f)^+ \; \in k[Y,U],
$$
where writing as \eqref{eq1.2.1},
$$
in_v(f)^+ =in_v(f)^+_{(y,u)} =
\sum\limits_{(A, B)} \overline{C}_{A, B} \; Y^B U^A \in k[Y,U]
$$
where the sum ranges over such $(A,B)$ that
$\vert B \vert < n_{(u)}(f)$ and
$\displaystyle{\frac{A}{n_{(u)}(f) - \vert B \vert} =v}$.
\item[(3)]
For a semi-positive linear form $L:\bR^e\to \bR$, we write
$\delta_L(f,y,u)=\delta_L(\Delta(f,y,u))$. By definition
$$
\delta_L(f,y,u) =
\min \left\{ \frac{L(A)}{n_{(u)}(f) - |B|} \;
\big\vert \; C_{A,B} \ne 0, \; \vert B \vert < n_{(u)}(f) \right\}.
$$
and $E_L=\Delta(f,y,u)\cap \{A\in \bR^e|\; L(A)=\delta_L(f,y,u)\}$
is a face of $\Delta(f,y,u)$ of slope $L$.
When $E_L$ is the $\delta$-face of $\Delta(f,y,u)$ (namely
$L=L_0$ as in Definition \ref{def1.1} (3)), we write simply
$\delta(f,y,u)=\delta_L(f,y,u)$.
\item[(4)]
Let $E_L$ be as in (3). We define the $E_L$-initial of $f$ by
$$
in_{E_L}(f) = in_{E_L}(f)_{(y,u)} = in_0(f) +
\sum\limits_{(A, B)} \overline{C}_{A, B} \; Y^B U^A \in k[[U]][Y]
$$
where the sum ranges over such $(A,B)$ that
$$
|B| < n_{(u)}(f) \qaq L(A)=\delta_L(f,y,u)(n_{(u)}(f) - |B|).
$$
We note that $in_{E_L}(f)$ is different from $in_L(f)$ in
Definition \ref{def1.3} (2).
When $E_L$ is the $\delta$-face of $\Delta(f,y,u)$, we write
$in_\delta(f)$ for $in_{E_L}(f)$.
\end{itemize}
\end{definition}

One easily sees the following:

\begin{lemma}\label{initial.lem}
Let the notation be as in Definition \ref{def.intial}.
\begin{itemize}
\item[(1)]
$in_v(f)_{(y,u)}$ is independent of the presentation \eqref{eq1.2.1}.
\item[(2)]
If $E_L$ is bounded, $in_{E_L}(f)\in k[U,Y]$ and it is independent of
the presentation \eqref{eq1.2.1}.
Otherwise it may depend on \eqref{eq1.2.1} (so there is an abuse of
notation).
\item[(3)]
If $E_L$ is bounded,
$$
in_{E_L}(f)=in_0(f)+\underset{v\in E_L}{\sum} in_v(f)^+.
$$
\item[(4)]
$in_v(f)=in_0(f)$ if $v\not\in \Delta(f,y,u)$ and
$in_v(f)\not=in_0(f)$ if $v$ is a vertex of $\Delta(f,y,u)$.
\end{itemize}
\end{lemma}
\medbreak

\begin{lemma}\label{lem2.2}
\begin{itemize}
\item[(1)]
$\delta(f,y,u) \geq 1$ if and only if $n_{(u)}(f)=v_\fm(f)$.
\item[(2)]
$\delta(f,y,u) = 1$ if and only if $\inm(f)=in_\delta(f)$.
\item[(3)]
$\delta(f,y,u) > 1$ if and only if $\inm(f) =in_0(f)\in k[Y]$.
\end{itemize}
\end{lemma}
\textbf{Proof}
\quad By definition $\delta(f,y,u)\geq 1$ is equivalent to the condition:
$$
C_{A,B}\not=0 \text{ and } |B| < n_{(u)}(f) \Rightarrow
\vert A \vert + \vert B \vert \geq n_{(u)}(f),
$$
which is equivalent to
$$
\vert A \vert + \vert B \vert <  n_{(u)}(f) \Rightarrow
C_{A,B}=0 \text{ or } |B| = n_{(u)}(f).
$$
Lemma \ref{lem2.2} (1) follows easily from this.
(2) and (3) follow by a similar argument and
the details are omitted.
$\square$
\bigskip

\begin{definition}\label{def.initialform}
Let $f=(f_1,\dots, f_m) \subset \fm $ be a system of elements such that
$f_i\not\in \langle u \rangle $.
\begin{itemize}
\item[(1)]
Define the polyhedron
$$\Delta((f_1, \dots, f_m), y, u) = \Delta (f, y, u)
\subseteq {\bR}_{\geq 0}^e$$ as the smallest $F$-subset
containing
$\underset{1\leq i\leq m}{\cup}\Delta (f_i, y, u)$.
\item[(2)]
For $v\in \bR^e- \Delta (f,y,u)^+$, the $v$-initial of $f$ is defined as
$$in_v(f)=(in_v(f_1),\ldots , in_v(f_m))$$
by noting
\begin{equation*}\label{poly.rem}
\Delta (f, y, u)^+\supset \underset{1\leq i\leq m}{\cup} \Delta (f_i, y, u)^+,
\end{equation*}
The $E_L$-initial $in_{E_L}(f)$ of $f$ for a face $E_L$ of $\Delta(f,y,u)$ is
defined similarly.
\item[(3)]
For a semi-positive linear form $L:\bR^e\to \bR$, we put
$$
\delta_L(f,y,u) = \min\{\delta_L(f_i,y,u)|\; 1\leq i\leq m\}.
$$
\item[(4)]
We let $V(f,y,u)$ denote the set of vertices of $\Delta(f,y,u)$. We put
$$
\EP(f,y,u)=\{v\in \bR^e-\Delta(f,y,u)^+|\; in_v(f_i)\not=in_0(f_i)
\text{ for some } 1\leq i\leq m \}.
$$
We call it the set of the essential points of $\Delta(f,y,u)$.
By definition $\Delta(f,y,u)$ is the smallest $F$-subset of $\bR^e$
which contains $\EP(f,y,u)$. By Lemma \ref{initial.lem} we have
$$
V(f,y,u) \subset \EP(f,y,u)\subset \partial\Delta(f,y,u).
$$

\end{itemize}
\end{definition}
\medbreak

The following fact is easily seen:

\begin{lemma}\label{Vlattice}
We have
$$
\EP(f,y,u)\subset \frac{1}{d!} {\bZ}_{\geq 0}^e\subseteq \bR^e
\qwith  d =\max\{n_{(u)}(f_i)|\; 1\leq i\leq m\}.
$$
In particular $\EP(f,y,u)$ is a finite set.
\end{lemma}
\medbreak

\begin{theorem}\label{refdatum.cor}
Let the assumption be as in Definition \ref{def.initialform}.
The following conditions are equivalent:
\begin{itemize}
\item[(1)]
$f$ is a $(u)$-standard base of $J$ and $\delta(f,y,u)>1$.
\item[(2)]
$f$ is a standard base of $J$ and $in_\fm(f_i)\in k[Y]$ for $\forall i$.
\end{itemize}
If $(u)$ is admissible for $J$, the conditions imply that $(y,u)$ is
strictly admissible for $J$.
\end{theorem}
\textbf{Proof  }
The implication (2)$\Rightarrow$(1) follows from Lemma \ref{lem2.2}
in view of Remark \ref{rem1.1}. We show (1)$\Rightarrow$(2).
When $\delta(f,y,u)>1$, $L=L_0$ is effective for $(f,y,u)$ by
Lemma \ref{lem2.2}, where $L_0(A)=|A|$ (cf. Remark \ref{rem1.1}).
Thus the desired assertion follows from Theorem \ref{refdatum0}.
Assume that $(u)$ is admissible for $J$. By Lemma \ref{def1.4cor},
the conditions imply that $\InmJ$ is generated by polynomials in $k[Y]$.
Thus we must have $\IDir(R/J)=\langle Y_1,\dots, Y_r\rangle$ by the
assumption that $(u)$ is admissible for $J$. This proves the last assertion.
$\square$
\medskip

\begin{definition}\label{def1.5}
Let the assumption be as in Setup A in \S\ref{sec:u-st.base}.
The polyhedron $\Delta(J,u)$ is the intersection of all
$\Delta (f,y,u)$ where $f=(f_1,\dots, f_m)$ is a $(u)$-standard basis
with reference datum $(y, L)$ for some $y$ and $L$.
\end{definition}
\medskip

\begin{remark}
\begin{itemize}
\item[(1)]
This is not the original definition given in \cite{H3}, (1.12)
(which is formulated more intrinsically), but it follows from \cite{H3}, (4.8)
that the definitions are equivalent.
\item[(2)]
As a polyhedron, $\Delta(f,y,u)$ and $\Delta(J,u)$ are defined by equations
$$ L_1 (A) \geq d_1, \dots, L_t (A) \geq d_t $$
for different non-zero semi-positive linear forms
$L_i$ on ${\bR}^e$.
\end{itemize}
\end{remark}
\medskip

Another important result of Hironaka provides a certain condition
under which we have $\Delta(J,u)=\Delta(f,y,u)$ (see Theorem
\ref{thm.wellprepared2}).
First we introduce the notion of normalizedness.

\begin{definition}\label{Esubset}
Let $S=k[X_1,\ldots , X_n]$ be a polynomial ring over a field $k$
and $I\subset S$ be a homogeneous ideal. We define:
$$
E(I)=\{LE(\varphi)\in \bZ^n_{\geq 0} \mid \varphi\in I\mbox{ homogeneous}\}\,,
$$
where for a homogeneous polynomial $\varphi\in S$, $LE(\varphi)$ is
its leading exponent, i.e., the biggest exponent
(in the lexicographic order on $\bZ^n_{\geq 0}$) occurring in $\varphi$:
For $\varphi=\sum c_A X^A$ we have $LE(\varphi)=\max\{A\mid c_A\neq 0\}$.
If $I$ is generated by homogeneous elements $\varphi_1,\ldots ,
\varphi_m$, we also write $E(I)=E(\varphi_1,\ldots , \varphi_m)$.
We note $E(I) +\bZ^n_{\geq 0} \subset E(I)$.
\end{definition}
\medbreak

\begin{definition}\label{def.normalized0}
Assume given $G_1,\dots,G_m\in k[[U]][Y]=k[[U_1,\dots,U_e]][Y]$:
$$
G_i= F_i(Y) + \underset{|B|<n_i}{\sum} Y^B P_{i,B}(U),\quad
(P_{i,B}(U)\in k[U])
$$
where $F_i(Y)\in k[Y]$ is homogeneous of degree $n_i$ and
$P_{i,B}(U)\not\in k-\{0\}$.
\begin{itemize}
\item[(1)] $(F_1,\dots,F_m)$ is normalized if writing
$$
F_i(Y)=\underset{B}{\sum} C_{i,B} Y^B \qwith  C_{i,B}\in k,
$$
$C_{i,B}=0$ if $B\in E(F_1,\dots,F_{i-1})$ for $i=1,\dots,m$ .
\item[(2)]
$(G_1,\dots,G_m)$ is normalized if $(F_1,\dots,F_m)$ is normalized and
$P_{i,B}(U)\equiv 0$ if $B\in E(F_1,\dots,F_{i-1})$ for $i=1,\dots,m$ .
\end{itemize}
\end{definition}
\medbreak

It is easy to see that if $(F_1,\dots,F_m)$ is normalized, then
it is weakly normalized in the sense of Definition \ref{def0.2}.
There is a way to transform a weakly normalized standard base of
a homogeneous ideal $I\subset k[Y]$ into a normalized standard base of $I$
(cf. \cite{H3}, Lemma 3.14 and Theorem \ref{thm.normalized} below).

\begin{definition}\label{def.normalized}
Let the assumption be as in Definition \ref{def.initialform}.
\begin{itemize}
\item[(1)]
$(f,y,u)$ is weakly normalized if $(in_0(f_1),\ldots ,in_0(f_m))$
is weakly normalized.
\item[(2)]
$(f,y,u)$ is $0$-normalized if $(in_0(f_1),\ldots ,in_0(f_m))$
is normalized in the sense of Definition \ref{def.normalized0} (1).
\item[(3)]
$(f,y,u)$ is normalized at $v\in \bR^e - \Delta(f,y,u)^+$ if so is
$(in_v(f_1),\ldots ,in_v(f_m))$ in the sense of Definition
\ref{def.normalized0} (2).
\item[(4)]
$(f,y,u)$ is normalized along a face $E_L$ of $\Delta(f,y,u)$ if so is
$(in_{E_L}(f_1),\ldots ,in_{E_L}(f_m))$ in the sense of Definition
\ref{def.normalized0} (2).
\end{itemize}
\end{definition}

Now we introduce the notions of (non-) solvability and preparedness.

\begin{definition}\label{def.solvable}
Let the assumption be as Definition \ref{def.initialform}.
For $v\in V(f,y,u)$, $(f, y, u)$ is called solvable at $v$
if there are $\lambda_1,\ldots , \lambda_r\in k[U]$ such that
$$
in_v(f_i)_{(y,u)}=F_i(Y+\lambda) \quad\text{ with  } F_i(Y)=in_0(f_i)_{(y,u)}
\quad (i=1,\ldots , m)\, ,
$$
where $Y+\lambda=(Y_1+\lambda_1,\ldots , Y_r+\lambda_r)$.
In this case the tuple
$\lambda=(\lambda_1,\ldots , \lambda_r)$ is called a solution for
$(f,y,u)$ at $v$.
\end{definition}

\begin{remark}
For $v \in V(f,y,u)$, it is not possible that $in_v(f_i)\in k[Y]$
for all $i=1,\ldots , m$ (cf. Definition \ref{def.initialform} (4)); hence
$\lambda\neq 0$ if $v$ is solvable.
\end{remark}

\begin{definition}\label{def.wellprepared}
Let the assumption be as Definition \ref{def.initialform}.
\begin{itemize}
\item[(1)]
Call $(f,y,u)$ prepared at $v\in V(f,y,u)$ if $(f,y,u)$ is
normalized at $v$ and not solvable at $v$.
\item[(2)]
Call $(f,y,u)$ prepared along a face $E_L$ of $\Delta(f,y,u)$ if $(f,y,u)$
is normalized along $E_L$ and not solvable at any $v\in V(f,y,u)\cap E_L$.
\item[(3)]
Call $(f,y,u)$ $\delta$-prepared if it is prepared along the $\delta$-face
of $\Delta(f,y,u)$.
\item[(4)]
Call $(f,y,u)$ well prepared
if it is prepared at any $v\in V(f,y,u)$.
\item[(5)]
Call $(f,y,u)$ totally prepared if it is well prepared and normalized
along all bounded faces of $\Delta(f,y,u)$.
\end{itemize}
\end{definition}

We can now state Hironaka's crucial result (cf. \cite{H3}, (4.8)).

\begin{theorem}\label{thm.wellprepared2}
Let the assumption be as in Definition \ref{def.initialform}.
Assume that $f$ is a $(u)$-standard base of $J$ and the following condition
holds, where $\tR=R/\langle u \rangle $, $\tfm=\fm/\langle u \rangle$,
$\tJ=J\tR$:
\begin{enumerate}
\item[(*)]
There is no proper $k$-subspace $T\subsetneq gr_{\tfm}^1(\tR)$ such that
$$
(In_{\tfm}(\tJ)\cap k[T])\cdot gr_{\tfm}(\tR)= In_{\tfm}(\tJ).
$$
\end{enumerate}
Let $v$ be a vertex of $\Delta(f,y,u)$ such that $(f,y,u)$ is prepared at
$v$. Then $v$ is a vertex of $\Delta(J,u)$.
In particular, if $(f,y,u)$ is well-prepared, then
$\Delta(J,u)=\Delta(f,y,u)$.
\end{theorem}

We note that condition $(*)$ is satisfied if $(u)$ is admissible
(Definition \ref{def2.11}).
\medbreak

\begin{corollary}\label{cor.wellprepared2}
Let the assumption be as Theorem \ref{thm.wellprepared2}.
Assume further that $(u)$ is admissible for $J$.
Then the following conditions are equivalent:
\begin{itemize}
\item[(1)]
$(f,y,u)$ is prepared at any $v\in V(f,y,u)$ lying in
$\{A\in \bR^e|\; |A|\leq 1\}$.
\item[(2)]
$\delta(f,y,u)>1$.
\item[(3)]
$(y,u)$ is strictly admissible for $J$ and $f$ is
a standard base of $J$ admissible for $(y,u)$.
\end{itemize}
The above conditions hold if $(f,y,u)$ is $\delta$-prepared.
\end{corollary}
\textbf{Proof}
Clearly (2) implies (1). The equivalence of (2) and (3) follows from
Theorem \ref{refdatum.cor}. It remains to show that (1) implies (2).
By Lemma \ref{stbTJ} we can find a strictly admissible $(z,u)$ and
a standard base $g$ of $J$ admissible for $(z,u)$.
By Lemma \ref{lem2.2} we have $\delta(g,z,u)>1$ and hence
$$
\Delta(J,u)\subset \Delta(g,z,u) \subset \{A\in \bR^e|\; |A| >1\}.
$$
By Theorem \ref{thm.wellprepared2} this implies $\delta(f,y,u)>1$ since
$(f,y,u)$ is well-prepared at any vertex $v$ with $|v|\leq 1$.
Finally, if $(f,y,u)$ is $\delta$-prepared, Theorem \ref{thm.wellprepared2}
implies $\delta(f,y,u)=\delta(\Delta(J,u))$ and the same argument as
above shows $\delta(f,y,u)>1$. This completes the proof.
$\square$
\bigskip

We have the following refinement of Theorem \ref{nf.thm1} (2)$(iv)$.

\begin{theorem}\label{C.thm}
Let the assumption be as in Definition \ref{def.initialform}.
Assume that $(u)$ is admissible for $J$ and that
$(f)$ is a $(u)$-standard base of $J$.
Let $X=\Spec(R/J)$ and $D=\Spec(R/\fp)$ for
$\fp=(y,u_1,\dots,u_s)\subset \fm=(y,u_1,\dots,u_e)$.
Assume that $D\subset X$ is permissible and that
there exists a vertex $v$ on the face $E_L$ such that $(f,y,u)$ is
prepared at $v$, where
$$
L:\bR^e\to\bR;\; (a_1,\dots,a_e)\to \underset{1\leq i\leq s}\sum a_i,
$$
Then
$v_\fp(f_i)=v_\fm(f_i)=n_{(u)}(f_i)$ for $i=1,\dots,m$.
In particular we have $\delta(f,y,u)\geq 1$.
\end{theorem}
\textbf{Proof of Theorem \ref{C.thm} }
By Theorem \ref{thm.wellprepared2}, the last assumption implies
\begin{equation}\label{compl.eq-1}
\delta_L(f,y,u)=\delta_L(\Delta(J,u)).
\end{equation}
Let $n_i=n_{(u)}(f_i)=\nu^{(i)}(J)$ for $i=1,\dots,m$.
By Lemma \ref{def1.4cor} (2) and \eqref{eq.ordersc} we have
$$v_\fp(f_i)\leq v_\fm(f_i)\leq n_{(u)}(f_i)=n_i
\qfor i=1,\dots,m.
$$
Thus it suffices to show that \eqref{compl.eq-1} implies $v_\fp(f_i)\geq n_i$.
Let $g=(g_1,\dots,g_m)$ be a $(u)$-standard base of $J$. As usual write
\begin{equation*}
g_i = \sum\limits_{i,A, B} C_{i,A, B} \; y^B u^A \quad \mbox{with}
\quad C_{i,A, B} \in R^{\times} \cup \{ 0 \},\;
A\in {\bZ}_{\geq 0}^e,\; B \in {\bZ}_{\geq 0}^r.
\end{equation*}
Note that $n_{(u)}(g_i)=n_{(u)}(f_i)=n_i$ by
Lemma \ref{def1.4cor} (2). We have
\begin{equation}\label{compl.eq0}
\begin{aligned}
v_{\fp} (g_i)\geq n_i  \Leftrightarrow &
|B|+L(A)\geq n_i \text{  if $C_{i,A, B}\not=0$ and $|B|<n_i$} \\
 \Leftrightarrow &
L(\frac{A}{n_i-|B|})\geq 1 \text{  if $C_{i,A, B}\not=0$ and $|B|<n_i$} \,. \\
\end{aligned}
\end{equation}
Hence we get the following equivalences for a $(u)$-standard base $g$ of $J$:
\begin{equation}\label{compl.eq1}
\begin{aligned}
\delta_L(g,y,u)\geq 1
& \Leftrightarrow v_{\fp} (g_i)\geq \nu^{(i)}(J) \qfor i=1,\dots,m \\
& \Leftrightarrow v_{\fp} (g_i) =v_{\fm}(g_i)=n_{(u)}(g_i) = \nu^{(i)}(J)
\qfor i=1,\dots,m \\
\end{aligned}
\end{equation}
Noting that any standard base is a $(u)$-standard base by
Lemma \ref{rem1.2} (2) (here we used the assumption (1)),
Theorem \ref{nf.thm1} (2)$(iv)$ implies that there exists a $(u)$-standard
base $g$ which satisfies the conditions of \eqref{compl.eq1}.
Since $\Delta(g,y,u)\supset\Delta(J,u)$, this implies
$\delta_L(f,y,u)\geq \delta_L(g,y,u)\geq 1$, which
implies the desired assertion by \eqref{compl.eq1}.
Finally the last assertion follows from Lemma \ref{lem2.2}.
$\square$

\bigskip

We can not expect to get the desirable situation of Theorem \ref{thm.wellprepared2}
right away. So we need procedures to attain this situation.
This is given by the following results (cf. \cite{H3}, (3.10), (3.14) and (3.15)).
First we discuss {\bf normalizations}.

\begin{theorem}\label{thm.normalized}
Let the assumption be as in Definition \ref{def.initialform}.
Assume that $(f,y,u)$ is weakly normalized. For $v\in V(f,y,u)$, there exist
$x_{ij}\in \langle u \rangle \;\subset R$ ($1\leq j <i\leq m$)
such that the following hold for
$$
h=(h_1,\ldots , h_m),\quad\text{ where  }
h_i=f_i-\sum\limits^{i-1}_{j=1}x_{ij}f_j.
$$
\begin{itemize}
\item[$(i)$]
$\Delta (h,y,u)\subseteq \Delta(f,y,u)$.
\item[$(ii)$]
If $v\in \Delta(h,y,u)$, then $v\in V(h,y,u)$ and
$(h,u,y)$ is normalized at $v$.
\item[$(iii)$]
$V(f,y,u)-\{v\}\subset V(h,y,u)$.
\item[(iv)]
For $v'\in V(f,y,u)-\{v\}$, we have
$in_{v'}(f)_{(y,u)}=in_{v'}(h)_{(y,u)}$.
\end{itemize}
\end{theorem}

\begin{remark}\label{remnormal}
In the above situation, the passage from $(f,y,u)$ to $(h,y,u)$ is
called a normalization at $v$. It is easy to see that
$ in_0(h_i) = in_0(f_i)$ for all $i$, $1\leq i \leq m$.
\end{remark}
\medbreak

We will need the following slight generalization of
Theorem \ref{thm.normalized}:

\begin{theorem}\label{thm.normalizededge}
Let the assumption be as in Definition \ref{def.initialform}.
Let $E$ be a bounded face of $\Delta(f,y,u)$.
Assume that $(f,y,u)$ is normalized at any $v\in E\cap V(f,y,u)$.
Then there exist
$x_{ij}\in (u)R$ ($1\leq j <i\leq m$)
such that putting
$$
h_i=f_i-\sum\limits^{i-1}_{j=1}x_{ij}f_j\qfor 1\leq i\leq m,
$$
$\Delta (h,y,u)=\Delta(f,y,u)$ and $(h,y,u)$ is normalized along $E$.
\end{theorem}
\medbreak\noindent
\textbf{Proof  }
Write $E=\Delta(f,y,u)\cap \{A\in \bR^e|\; L(A)=1\}$
for a positive linear form $L:\bR^e\to \bR$.
For $i=1,\dots,m$, we can write
$$
in_E(f_i)=F_i(Y)+ \underset{|B|<n_i}{\sum} Y^B P_{i,B}(U)\in k[Y,U],\quad
P_{i,B}(U)= \underset{|B|+L(A)=n_i}{\sum} c_{i,A,B} U^A\in k[U],
$$
where $n_i=n_{(u)}(f_i)=v_L(f_i)$ and $F_i(Y)=in_0(f_i)\in k[Y]$ which is
homogeneous of degree $n_i$. For each $i\geq 1$ put
$$
\Sigma(f_1,\dots,f_i)=\{B|\; P_{i,B}(U)\not\equiv 0,\; B\in
E(F_1,\dots,F_{i-1})\}.
$$
If $\Sigma(f_1,\dots,f_i)=\emptyset$ for all $i\geq 1$, there is nothing to be
done. Assume the contrary and let
$j=\min\{i|\; \Sigma(f_1,\dots,f_i)\not=\emptyset\}$ and $\Bm$ be the
maximal element of $\Sigma(f_1,\dots,f_j)$ with respect to the lexicographic
order. Note that
\begin{equation}\label{thm.normalizededge.eq1}
\{\frac{A}{n_j-|\Bm|}\;|\; c_{j,A,\Bm}\not=0\}\subset E_L\backslash V(f,y,u)
\end{equation}
by the assumption that $(f,y,u)$ is normalized at any $v\in E_L\cap V(f,y,u)$.
By the construction there exist $G_i(Y)\in k[Y]$, homogeneous of degree
$|\Bm|-n_i$, for $1\leq i\leq j-1$ such that
$$
H(Y):= Y^{\Bm}-\underset{1\leq i\leq j-1}{\sum} G_i(Y) F_i(Y)
$$
has exponents smaller than $\Bm$. Note that $H(Y)$ is homogeneous of degree
$|\Bm|$. For $i=1,\dots,j-1$, take $g_i\in R$ such that $\inm(g_i)=G_i(Y)$
and take
$$
\tilde{P}_{j,\Bm} =\underset{|\Bm|+L(A)=n_i}{\sum}
\tilde{c}_{j,A,\Bm} u^A\;\in R
\qwith c_{j,A,\Bm}= \tilde{c}_{j,A,\Bm}\mod \fm.
$$
Put
$$
h_j=f_j- \tilde{P}_{j,\Bm} \underset{1\leq i\leq j-1}{\sum} g_i f_i\qaq
h_i=f_i \qfor 1\leq i\not=j\leq m.
$$
By \eqref{thm.normalizededge.eq1}, we have
\begin{align*}
in_E(h_j)&=in_E(f_j) - P_{j,\Bm}(U)(Y^{\Bm}-H(Y)),\\
in_v(h_j)&=in_v(f_j) \qfor \forall v\in V(f,y,u).
\end{align*}
Hence $\Sigma(h_1,\dots,h_i)=\Sigma(f_1,\dots,f_i)$ for all $i=1,\dots,j-1$
and all elements of $\Sigma(h_1,\dots,h_j)$ are smaller than $\Bm$ in the
lexicographical order. This proves Theorem \ref{thm.normalizededge} by
induction.
$\square$
\medbreak

Now we discuss {\bf dissolutions}.

\begin{theorem}\label{thm.solvable}
Let the assumption be as in Definition \ref{def.initialform} and let $v\in V(f,y,u)$.
\begin{itemize}

\item[$(a)$] Any solution for $(f,y,u)$ at $v$ is of the form
$\lambda$ with $\lambda_i=c_iU^v$, where $c_i\in k^\times$. In particular,
if it exists, a solution is always non-trivial and $v\in \bZ^e$.

\item[$(b)$] Let $d=(d_1,\ldots , d_r)\subset R$ with
$d_i\in \langle u^v \rangle $ be
such that the image of $d_i$ in $\grm^{|v|}(R)$ is
$\lambda_i$. Let $z=y-d=(y_1-d_1,\dots, y_r-d_r)$. Then
\begin{itemize}
\item[$(i)$] $\Delta(f,z,u)\subseteq\Delta (f,y,u)$.
\item[$(ii)$]
$v\not\in \Delta(f,z,u)$ and $V(f,y,u)-\{v\}\subset V(f,z,u)$.
\item[$(iii)$]
For $v'\in V(f,y,u)-\{v\}$, we have
$$ in_{v'}(f)_{(z,u)}={in_{v'}(f)_{(y,u)}}_{|Y=Z}\; \in k[Z,U]. \quad
(Z= \inm(z) \in gr^1_{\mathfrak m}(R))$$
\end{itemize}
\end{itemize}
\end{theorem}

\begin{remark}\label{remdiss}
In the above situation, the passage from $(f,y,u)$ to $(f,z,u)$ is
called the dissolution at $v$. It is easy to see
$ in_0(f_i)_{(z,u)} = {in_0(f_i)_{(y,u)}}_{| Y=Z}$
for $1\leq \forall i\leq m$.
\end{remark}
\bigskip

We now come to the \textbf{preparation} of $(f,y,u)$:\quad
Let the assumption be as in Definition \ref{def.initialform}.
We apply to $(f,y,u)$ alternately and repeatedly
normalizations and dissolutions at vertices of polyhedra.
To be precise we endow $\bR^e_{\geq 0}$ with an order defined by
$$
v>w \Leftrightarrow |v|>|w|\text{ or } |v|=|w|\text{ and }
v>w \text{ in the lexicographical order}.
$$
Let $v\in V(f,y,u)$ be the smallest point and apply the normalization at $v$ from
Theorem \ref{thm.normalized} and then the dissolution at $v$ from Theorem \ref{thm.solvable}
if $v$ is solvable, to get $(g,z,u)$. Then $(g,z,u)$ is prepared at $v$.
Repeating the process, we arrive at the following conclusion
(cf. \cite{H3}, (3.17)).

\begin{theorem}\label{thm.wellprepared}
Let the assumption be as in Definition \ref{def.initialform}.
Assume that $(f,y,u)$ is weakly normalized. For any integer $M>0$,
there exist
$$
x_{ij} \in \langle u \rangle \quad (1\leq j < i\leq m),\quad
d_\nu  \in \langle u \rangle  \quad (\nu=1,\ldots, r),
$$
such that putting
\begin{align*}
&z=(z_1,\ldots , z_r)\text{ with  } z_\nu=y_\nu-d_\nu, \\
&g=(g_1,\ldots , g_m)\text{ with  } g_i=f_i-\sum^{i-1}_{j=1}x_{ij}f_j\:,
\end{align*}
we have $\Delta(g,z,u)\subseteq \Delta (f,y,u)$, and $(g,z,u)$ is prepared along
all bounded faces contained in $\{A\in \bR^e|\; |A|\leq M\}$.
If $R$ is complete, we can obtain the stronger conclusion that
$(g,z,u)$ is well-prepared.
\end{theorem}

\begin{remark}\label{remark.normalizededges}
By Theorem \ref{thm.normalizededge} we can make $(g,z,u)$ in Theorem
\ref{thm.wellprepared} satisfy the additional condition that it is
normalized along all bounded faces of the polyhedron contained in
$\{A\in \bR^e|\; |A|\leq N\}$. If $R$ is complete, we can make
$(g,z,u)$ totally prepared.
\end{remark}

\begin{corollary}\label{cor.wellprepared}
Let the assumption and notation be as in Theorem \ref{thm.wellprepared}.
\begin{itemize}
\item[(1)]
If $f$ is a $(u)$-standard base of $J$, then so is $g$.
\item[(2)]
If $f$ is a $(u)$-standard base of $J$ and $(u)$ is admissible
(cf. \ref{def2.11}), then $\delta(g,z,u)>1$ and $(z,u)$ is strictly
admissible and $g$ is a standard base admissible for $(z,u)$.
\end{itemize}
\end{corollary}
\textbf{Proof} \quad
(1) follows from Corollary \ref{cor.refdatum0}.
(2) follows from (1) and Theorem \ref{cor.wellprepared2}.
$\square$
\bigskip

\def\us{u_{\leq s}}
\def\usa{u_{\leq {s+1}}}
\def\Us{U_{\leq s}}
\def\Ps{P^{\leq s}}
\def\Pbs{\overline{P}^{\leq s}}
\def\Pb{\overline{P}}

At the end of this section we prepare a key result which relates certain localizations
of our ring to certain projections for the polyhedra.
Let $(f,y,u)$ be as in Definition \ref{def.initialform}.
For $s=1,\dots,e$, we let
$$
\fp_s= \langle y,\us\rangle=\langle y_1,\dots,y_r,u_1,\dots,u_s\rangle.
$$
Let $R_s$ be the localization of $R$ at $\fp_s$, let $J_s=JR_s$, and
$\fm_s=\fp_s R_s$ (the maximal ideal of $R_s$).
We want to relate $\Delta(f,y,u)\subset \bR^e$ to $\Delta(f,y,\us)\subset \bR^s$,
the characteristic polyhedron for $J_s\subset R_s$.
Assume given a presentation as in \eqref{eq1.2.1}:
\begin{equation}\label{eq1.proj}
f_i = \sum\limits_{(A, B)} P_{i,A, B} \; y^B u^A \quad \mbox{with}
\quad P_{i,A, B} \in R^\times\cup \{0\} \quad (i=1,\dots,m).
\end{equation}
\eqref{eq1.proj} can be rewritten as:
\begin{equation}\label{eq2.proj}
f_i = \sum\limits_{(C, B)} \Ps_{i,C, B} \; y^B \us^C,
\end{equation}
where for $C=(a_1,\dots,a_s)\in \bZ^s_{\geq 0}$,
$\us^C=u_1^{a_1}\cdots u_s^{a_s}$ and
\begin{equation}\label{eq2.1.proj}
\Ps_{i,C, B}=\underset{A=(a_1,\dots,a_s,a_{s+1},\dots,a_e)}{\sum}
P_{i,A, B}
u_{s+1}^{a_{s+1}}\cdots u_e^{a_e}\;\in \hat{R}.
\end{equation}
We now introduce some conditions which are naturally verified in the
case where $\Spec(R/\langle u_1\rangle)$ is the exceptional divisor of a blowup
at a closed point (see Lemma \ref{ekeylemma2} below). Assume
\begin{enumerate}
\item[$(P0)$]
$J\subset \fp_1$ and there is a subfield $k_0\subset R/\fp_1$ such that
$P_{i,A, B} \mod \fp_1\;\in k_0$.
\end{enumerate}
By $(P0)$ we get the following for $C=(a_1,\dots,a_s)\in \bZ^s_{\geq 0}$:
$$
\Pbs_{i,C,B}:=\Ps_{i,C,B}\mod \fp_s =
\underset{A=(a_1,\dots,a_s,a_{s+1},\dots,a_e)}{\sum}  \Pb_{i,A, B}
\ub_{s+1}^{a_{s+1}}\cdots \ub_e^{a_e}\;\in
k_0[[\ub_{s+1},\dots,\ub_e]]\subset \hat{R}/\fp_s,
$$
$$
\Pb_{i,C, B}:=P_{i,A, B} \mod\fp_s \in k_0\hookrightarrow R/\fp_s,
\quad
\ub_j=u_j\mod\fp_s\in R/\fp_s.
$$
Hence we have the following equivalences for $C\in \bZ^s_{\geq 0}$
\begin{enumerate}
\item[$(P1)$]
$\begin{aligned}
\Ps_{i,C,B}\in \fp_s
&\Leftrightarrow P_{i,A, B}=0
\quad\text{ for all $A\in \bZ^e_{\geq 0}$ such that $\pi_s(A)=C$}\\
&\Leftrightarrow \Ps_{i,C,B}=0 \\
\end{aligned}$
\end{enumerate}
where
$\pi_s:\bR^e\to \bR^s\;;\; (a_1,\dots,a_e)\to (a_1,\dots,a_s).$
We further assume
\begin{enumerate}
\item[$(P2)$]
For fixed $B$ and $a\in \bZ_{\geq0}$, there are only finitely many $A$
such that $P_{i,A,B}\not=0$ and $\pi_1(A)=a$.
\end{enumerate}
This condition implies $\Ps_{i,C, B}\in R$.

\begin{theorem}\label{lem.projection}
Let $L:\bR^s\to \bR$ be a semi-positive linear form and $L_s=L\circ \pi_s$.
\begin{itemize}
\item[(1)]
If $(P0)$ holds, then
$$\Delta(f,y,\us)=\pi_s(\Delta(f,y,u))\qaq
\delta_L(f,y,\us)=\delta_{L_s}(f,y,u).
$$
\item[(2)]
If $(P0)$ and $(P2)$ hold and $L(1,0,\dots,0)\not=0$,
then the initial form along the face $E_L$ of $\Delta(f,y,\us)$ (with respect to the presentation \eqref{eq2.proj},
cf. Definition \ref{def.intial}) lies in the polynomial ring $R/\fp_s[Y,\Us]$
and we have
$$
in_{E_L}(f)_{(y,\us)}= {in_{E_{L_s}}(f)_{(y,u)}}_{|U_i=\ub_i\;
(s+1\leq i\leq e)}\;,
$$
considered as an equation in
$k_0[\ub_{s+1},\dots,\ub_e][Y,\Us]\subset R/\fp_s[Y,\Us]$ (cf. $(P2)$).
\item[(3)]
Assume $(P0)$ and $(P2)$ and $L(1,0,\dots,0)\not=0$. Assume further:
\begin{itemize}
\item[$(i)$]
There is no proper $k$-subspace
$T\subset \underset{1\leq i\leq r}{\bigoplus} k\cdot Y_i$
($k=R/\fm$) such that

$F_i(Y):=in_0(f_i)\in k[T]\subset k[Y]$ for all $j=1,\dots,m$.
\item[$(ii)$]
$(f,y,u)$ is prepared along the face $E_{L_s}$ of $\Delta(f,y,u)$.
\end{itemize}
Then $(f,y,\us)$ is prepared along the face $E_{L}$ of $\Delta(f,y,\us)$.
\end{itemize}
\end{theorem}
\textbf{Proof }
By $(P1)$ we get
$$
\begin{aligned}
\delta_L(f,y,\us)
&=\min\{\frac{L(C)}{n_i-|B|}\;|\; |B|<n_i,\; \Ps_{i,C,B}\not=0\}\\
&=\min\{\frac{L_s(A)}{n_i-|B|}\;|\; |B|<n_i,\; P_{i,A,B}\not=0\} =\delta_{L_s}(f,y,u)
\end{aligned}
$$
and (1) follows from this. To show (2), we use \eqref{eq2.proj} to compute
\begin{equation}\label{eq3.5.proj}
\begin{aligned}
in_{E_L}(f)_{(y,\us)}
&=F_i(Y)+\underset{B,C}{\sum}\Pbs_{i,C,B} Y^B \Us^C\\
&=F_i(Y)+\underset{B,A}{\sum}\Pb_{i,A,B} \ub_{s+1}^{a_{s+1}}\cdots \ub_e^{a_e} \cdot Y^B \Us^{\pi_s(A)} \;,\\
\end{aligned}
\end{equation}
where the first (resp. second) sum ranges over those $B,C$ (resp. $B,A$) for which
$|B|<n_i$ and $L(C)=\delta_L(f,y,\us)(n_i-|B|)$ (resp. $|B|<n_i$ and
$L_s(A)=\delta_{L_s}(f,y,u)(n_i-|B|)$). (2) follows easily from this.

We now show (3). By (2) the assumption $(ii)$ implies that
$in_{E_L}(f)_{(y,\us)}$ is normalized. It suffices to show the following:

\begin{claim}\label{claim0.proj}
Let $v$ be a vertex on $E_L$. Then $(f,y,\us)$ is not solvable at $v$.
\end{claim}
\medbreak

We prove the claim by descending induction on $s$ (the case $s=e$ is obvious).
From (1) we easily see that there exists a vertex $w\in \Delta(f,y,u)$
such that $\pi_s(w)=v$ and $v_t:=\pi_t(w)$ is a vertex of $\Delta(f,y,,\ut)$
for all $s\leq t\leq e$. By induction hypothesis, we may assume
\begin{enumerate}
\item[$(*)$]
$(f,y,\usa)$ is not solvable at $\vsa$.
\end{enumerate}
Assume $(f,y,\us)$ solvable at $v$. Then $v\in \bZ_{\geq 0}^s$ and there
exist elements $\lambda_1,\dots,\lambda_r$ of the fraction field $K$ of
$R/\fp_s$ such that
\begin{equation}\label{eq5.proj}
in_{v}(f)_{(y,\us)}=
F_i(Y_1+\lambda_1 \Us^v,\dots,Y_r+\lambda_r \Us^v)
\quad\text{ for all }i=1,\dots,m.
\end{equation}

\begin{claim}\label{claim1.proj}
$\lambda_j$ lies in the localization $S$ of $R/\fp_s$ at
$\fp_{s+1}/\fp_s$ for all $j=1,\dots,r$.
\end{claim}
\medbreak

Admit the claim for the moment.
By the claim we can lift $\lambda_j\in S$ to $\tlam_j\in R_{s+1}$,
the localization of $R$ at $\fp_{s+1}$.
Set $z_j=y_j + \tlam_j \us^v \in R_{s+1}\subset R_s$.
Take a positive linear form $L_v:\bR^s\to \bR$ such that
$E_{L_v}=\{v\}$ and hence $L_v(v)=\delta_{L_v}(f,y,\us)$.
By Theorem \ref{thm.solvable},
$\Delta(f,z,\us)\subset\Delta(f,y,\us)-\{v\}$ so that
\begin{equation}\label{eq1.claim1.proj}
\delta_{L_v}(f,z,\us) > \delta_{L_v}(f,y,\us).
\end{equation}
Now we apply (1) to $J_{s+1}\subset R_{s+1}$ and $(f,z,\usa)$ instead of
$J\subset R$ and $(f,y,u)$. Note that in the proof of (1) we have used only $(P0)$
which carries over to the replacement. We get
$$
\Delta(f,z,\us)=\pi\big(\Delta(f,z,\usa)\big)\qaq
\delta_{L_v\circ \pi}(f,z,\usa)= \delta_{L_v}(f,z,\us),
$$
where $\pi:\bR^{s+1}\to \bR^s\;;\; (a_1,\dots,a_{s+1})\to (a_1,\dots,a_s)$.
By the assumption, $v=\pi(\vsa)$ and
$$
L_v\circ \pi(\vsa)=L_v(v)=\delta_{L_v}(f,y,\us).
$$
By Theorem \ref{thm.wellprepared2}, $(*)$ implies
\begin{equation}\label{eq3.claim1.proj}
\vsa\in \Delta(\Jsa,\usa)\subset\Delta(f,z,\usa)).
\end{equation}
Thus we get
$$
\delta_{L_v}(f,z,\us)=\delta_{L_v\circ \pi}(f,z,\usa)\leq
\delta_{L_v\circ \pi}(\Delta(\Jsa,\usa))\leq L_v\circ \pi(\vsa) =\delta_{L_v}(f,y,\us),
$$
where the inequalities follow from \eqref{eq3.claim1.proj}.
This contradicts \eqref{eq1.claim1.proj} and the proof of (3) is complete.
\medbreak

Now we show Claim \ref{claim1.proj}.
Note that $S$ is a discrete valuation ring with a prime element
$\pi:=u_{s+1}\mod\fp_s$.
Thus it suffices to show $v_\pi(\lambda_j)\geq 0$
for $j=1,\dots,r$. Assume the contrary. We may assume
$$
v_\pi(\lambda_1)=-\epsilon <0,\quad v_\pi(\lambda_1)\leq v_\pi(\lambda_j)
\qfor j=1,\dots,r.
$$
Set
$$
Z_j=Y_j + \mu_j V\quad\text{where $V=\Us^v$ and
$\mu_j=\lambda_j \pi^\epsilon\in S$},
$$
and recall that $F_i(Y)=in_0(f_i)\in k_0[Y]$ and $k_0\subset S$ by $(P0)$ .
Consider
$$
F_i(\overline{Z})=F_i(Y_1+\overline{\mu}_1 V,\dots,Y_r+\overline{\mu}_r V)
\in \k[Y,V],\quad \k:=S/\langle \pi\rangle=\k(\fp_{s+1}),
$$
where $\overline{\mu}_i=\mu_j\mod \pi\;\in \k$.
We claim that there is some $i$ for which $F_i(\overline{Z})\not \in \k[Y]$.
Indeed, by the structure theorem of complete local rings, $(P0)$ implies
$\hat{R}/\fp_{s+1}\simeq k[[u_{s+2},\dots,u_e]]$ so that
$\k$ is contained in $k((u_{s+2},\dots,u_e))$ which is a separable
extension of $k$. By (3)$(i)$ and Lemma \ref{dimDir} (2),
$$
T=\underset{1\leq j\leq r}{\bigoplus} \k\cdot (Y_1+\overline{\mu}_j V)
\;\subset \;
\underset{1\leq j\leq r}{\bigoplus} \k\cdot Y_j \oplus \k\cdot V
$$
is the smallest $\k$-subspace such that
$F_i(\overline{Z})\in \k[T]$ for all $i=1,\dots,m$.
Thus the claim follows from the fact that $\overline{\mu}_1\not=0$.
For the above $i$, we expand
$$
F_i(Y_1+\lambda_1 V,\dots,Y_r+\lambda_r V)=
\underset{B}{\sum} \gamma_B {Y}^B V^{n_i-|B|}\quad (\gamma_B\in K),
$$
and we get
$$
F_i(Z)=F_i(Y_1+\mu_i V,\dots,Y_r+\mu_i V)
=\underset{B}{\sum} \gamma_B \pi^{\epsilon(n_i-|B|)}\cdot {Y}^B V^{n_i-|B|}
\;\in S[Y,V].
$$
Since $F_i(\overline{Z})\not\in \k[Y]$, there is some $B$ such that
$|B|<n_i$ and $\gamma_B \pi^{\epsilon(n_i-|B|)}$ is a unit of $S$.
Noting $\epsilon>0$, this implies $\gamma_B\not\in S$.
On the other hand, \eqref{eq3.5.proj} and \eqref{eq5.proj} imply
$\gamma_B\in S$, which is absurd.
This completes the proof of Claim \ref{claim1.proj}.
$\square$

\newpage
\section{Transformation of standard bases under blow-ups}
\label{sec:pr.usb}

\bigskip
In this section we will study the transformation of a standard base
under permissible blow-ups, in particular with respect to near points
in the blow-up. We begin by setting up a local description of the situation
in Theorem \ref{thmdirectrix} in \S2.

\bigskip\noindent
\textbf{Setup B }
\medbreak

Let $Z$ be an excellent regular scheme and $X\subset Z$ be a closed subscheme
and take a closed point $x\in X$. Put $R=\cO_{Z,x}$ with the maximal ideal
$\fm$ and put $k=R/\fm=k(x)$.
Write $X\times_Z \Spec(R)=\Spec(R/J)$ for an ideal $J\subset \fm$.
Define the integers $n_1\leq n_2\leq \cdots\leq n_m$ by
$$
\nu^*_x(X,Z)=\nu^*(J)=(n_1,\dots,n_m,\infty,\infty,\dots).
$$
Let $D\subset X$ be a closed subscheme permissible at $x\in D$
and let $J\subset \fp$ be the prime ideal defining $D\subset Z$.
By Theorem \ref{nf.thm1} (2) we have $T_x(D)\subset \Dir_x(X)$ so that
we can find a system of regular parameters for $R$,
$$(y, u, v ) = (y_1, \dots, y_r, u_1, \dots, u_s, v_1,\dots, v_t)$$
such that $\fp=(y,u)$ and
$(y,(u,v))$ is strictly admissible for $J$ (cf. Definition \ref{def2.11}).
This gives us an identification
$$
\grmR = k [Y, U, V]=k[Y_1,\dots,Y_r,U_1,\dots,U_s,V_1,\dots,V_t]
\quad (k = R/{\fm})
$$
where
$Y_i = \inm(y_i), \; U_i = \inm(u_i),\; V_i =\inm(v_i)\in \grmR$.
Consider the diagram:
$$
\begin{matrix}
\Bl_D(X)=& X'&\subset& Z'
& = & B\ell_D (Z)  \hookleftarrow  \pi_Z^{-1} (x)  = E_x \\
&\downarrow\rlap{$\pi_X$}&&\downarrow\rlap{$\pi_Z$} \\
&X&\subset& Z  \\
\end{matrix}
$$
and note that
$$
E_x:= \pi_Z^{-1} (x) = {\bP}(T_x(Z)/T_x(D)) =
\Proj(k [Y, U]) \cong {\bP}_k^{r+s-1}.
$$
We fix a point
$$
x'\in \Proj(k [U]) \subset \Proj(k [Y,U])= E_x.
$$
(By Theorem \ref{thmdirectrix}, if $\Char(k(x)) = 0$ or
$\Char(k(x)) \geq \dim (X)/2+1$, any point of $X'$ near to $x$ lies in
$\Proj(k [U])$. Moreover, if $x'$ is near to $x$, Theorem \ref{thm.pbu.inv} implies
$$
\nu^*_{x'}(X',Z') = \nu^*_x(X,Z).
$$
\medbreak

Without loss of generality we assume further that $x'$ lies in the chart
$\{U_1 \ne 0 \}\subset E_x$.
Let $R'=\cO_{Z',x'}$ with the maximal ideal $\fm'$ and let $J'\subset \fm'$
be the ideal defining $X'\subset Z'$ at $x'$. Put $k'=k(x')=R'/\fm'$.
Then $\fm R' = (u_1,v)=(u_1,v_1,\dots,v_t)$, and
$$
(y', u_1, v),\qwith
y' =(y'_1, \dots, y'_r), \; y'_i = y_i/u_1, \; v=(v_1, \dots, v_t)
$$
is a part of a system of regular parameters for $R'$.
Choose any $\phi_2, \dots, \phi_{s'} \in \fm'$ such that $(y',u_1, \phi,v)$,
with $\phi =(\phi_2, \dots, \phi_{s'})$, is a system of regular parameters
for $R'$ (note $s-s'=\trdeg_{k}(k(x'))$). Then
$$
\gr_{\fm'} (R') = k'[Y', U_1, \Phi, V] =
k'[Y'_1,\dots,Y'_r, U_1, \Phi_1,\dots,\Phi_{s'}, V_1,\dots,V_t],
$$
where
$Y'_i =\inmd(y'_i),\; \Phi_i = \inmd(\phi_i),\; V_i = \inmd(v_i)\in \grmd(R')$.
Assume now given
$$
\text{a standard base $f=(f_1,\dots,f_m)$ of $J$ which is admissible for
$(y,(u,v))$.}
$$
By definition
\begin{equation}\label{eqB1}
F_i(Y):=in_{\fm}(f_i)\in k[Y]\;(i=1,\dots, m)\qaq
In_{\fm}(J)=\langle F_1(Y),\dots, F_m(Y) \rangle  .
\end{equation}
By Lemma \ref{def1.4cor} (2) we have
\begin{equation}\label{eqB1.5}
v_\fm(f_i)= n_{(u,v)}(f_i) = n_i \qfor  i=1,\dots, m.
\end{equation}
Finally we assume
\begin{equation}\label{eqB3}
v_\fm(f_i)=v_\fp(f_i) \qfor i=1,\dots, m.
\end{equation}
This assumption is satisfied if $D=\{x\}$ (for trivial reasons) or under
the conditions of Theorem \ref{C.thm} (for example if
$(f,y,(u,v))$ is well-prepared). The assumptions imply
\begin{equation}\label{eqB1.6}
f_i = \sum\limits_{A, B} C_{i, A, B} \; y^B u^{A_u} v^{A_v},
\quad C_{A, B} \in \Gamma, \quad A=(A_u,A_v),\;
A_u\in \bZ_{\geq 0}^s,\; A_v\in \bZ_{\geq 0}^t,
\end{equation}
where the sum ranges over $A_u,A_v,B$ such that
\begin{equation}\label{eqB1.8}
|B|+|A_u|\geq n_i.
\end{equation}

By \cite{H1} Ch. III.2 p. 216 Lemma 6 we have $J'=\langle f'_1,\dots,f'_m \rangle  $ with
\begin{equation}\label{eqB1.7}
f'_i = f_i/u_1^{n_i}=
\sum\limits_{A, B} \big(C_{A, B} \; u^{\prime A_u}\big) y^{\prime^B}
u_1^{|A_u|+ |B|-n_i} v^{A_v},
\end{equation}
$$
u^{\prime A_u} = u_2^{\prime a_2} \dots u_s^{\prime a_s}
\qfor A_u = (a_1, \dots, a_s)\quad(u'_i = u_i/u_1).
$$
This implies
\begin{equation}\label{eqB2}
f'_i= \tilde{F}'_i \; \mod \;  \langle u_1,v \rangle  =
\langle u_1,v_{1},...,v_t \rangle
\quad \text{ for $i=1,\dots,m$, where}
\end{equation}
$$
\tilde{F}'_i = \sum\limits_{|B|=n_i} C_{i,0, B} y^{\prime^B}
\qaq
n_i= v_{\fm'}(\tilde{F}'_i)= n_{(u_1,\phi,v)}(f'_i)
$$
so that
\begin{equation}\label{eqB2.1}
in_0(f'_i)_{(y',(u_1,\phi,v))}=F_i(Y')
\quad \text{ for $i=1,\dots,m$}
\end{equation}

For later use, we choose $S$ (resp. $S'$), a ring of coefficients of
$\hat{R}$ (resp. $\hat{R'}$) (cf. \eqref{eq1.2.1.c}).
We also choose a set $\Gamma\subset S$ of representatives of $k$
(resp. a set $\Gamma'\subset S'$ of representatives of $k'$).
We note that the choices for $R$ and $R'$ are independent:
We do not demand $S\subset S'$ nor $\Gamma\subset \Gamma'$.

\medskip
(end of Setup B).

\bigbreak
We want to compare the properties of $(f, y, (u, v))$ (downstairs)
and $(f', y', (u_1, \phi, v))$ (upstairs), especially some
properties of the polyhedra and initial forms.

\medbreak
Let $e=s+t$ and $e'=s'+t$. For a semi-positive linear form $L$ on ${\bR}^e$
(downstairs) (resp. on $\bR^{e'}$ (upstairs)), $v_L$ and $in_L(*)$ denote the
$L$-valuation of $R$ (resp. $R'$) and the corresponding initial form of
$*\in R$ (resp. $*\in R'$) with respect to $(y,(u,v),\Gamma)$
(resp. $(y',(u_1,\phi,v),\Gamma')$).

\medbreak

\begin{theorem}\label{prop2.1}
In Setup B, $f'=(f'_1, \dots, f'_m)$ is a $(u_1,\phi,v)$-effective basis
of $J'$.
If $(f_1, \dots, f_m)$ is a standard base of $J$, then $(f'_1, \dots, f'_m)$
is a $(u_1,\phi,v)$-standard basis of $J'$.
More precisely there exists a positive linear form $L'$ on ${\bR}^e$
(upstairs) such that:
$$
in_{L'} (f'_i) = F_i (Y') \quad (1 \leq i \leq m)
\quad\text{ and }\quad
in_{L'} (J') = \langle F_1 (Y'), \dots, F_m (Y') \rangle  .
$$
\end{theorem}
\medbreak

First we need to show the following:

\begin{lemma}\label{lem2.1}
Let the assumption be as in Setup B.
Choose $d > 1$ and consider the linear forms on ${\bR}^e$ and ${\bR}^{e'}$:
$$
\begin{aligned}
L(A) & = & \frac{1}{d}
\big(\sum_{1\leq i\leq s} a_i + \sum_{1\leq j\leq t}^t a'_j\big)
& \qquad \mbox{(downstairs)} \\
\\
\Lambda' (A) & = & \frac{1}{d-1}
\big(a_1 + \sum_{1\leq j\leq t} a'_j\big) & \qquad \mbox{(upstairs)} \\
\end{aligned}
$$
respectively, where
$A=(a_1,\dots, a_*,a'_1,\dots,a'_t)$ with
$*=s$ (downstairs) and $*=s'$ (upstairs).
Then the following holds for $g \in R$.
\begin{itemize}
\item[(1)]
$v_{\Lambda'} (g) = \frac{d}{d-1} v_L (g)$.
\item[(2)]
Assuming that $g' := g/u_1^n \in R'$ with $v_L (g) = n$, we have
$v_{\Lambda'}(g') =n$. Assuming further $in_L(g) = G(Y) \in k[Y]$, we have
$in_{\Lambda'} (g') = G(Y')$.
\end{itemize}
\end{lemma}
\textbf{Proof} \quad
Let $g\in \hat{R}$ and write
$$
g = \sum\limits_{A, B} C_{A, B} \; y^B u^{A_u} v^{A_v}
$$
as in \eqref{eqB1.6}. Then, in $\widehat{R'}$ we have
\begin{equation}\label{eq2.2.1}
g = \sum\limits_{A, B} \big(C_{A, B} \; u^{\prime A_u}\big) y^{\prime^B}
u_1^{\vert A_u \vert + \vert B \vert} v^{A_v},
\end{equation}
where the notation is as in \eqref{eqB1.7}. Note that
$$
C_{A, B} \; u^{\prime A_u} = 0 \quad \mbox{in} \quad
\widehat{R'}/(y', u_1, v) \Longleftrightarrow C_{A, B} = 0
$$
because $C_{A, B} \ne 0$ implies $C_{A, B} \in R^{\times}$ so that
$C_{A, B} \in (R')^{\times}$. Hence
\begin{align*}
v_{\Lambda'} (g)  = & \min \left\{ \vert B \vert + \frac{\vert B
\vert + \vert A
\vert}{d-1} \; \big\vert \; C_{A, B} \ne 0 \right\} \\
 = & \frac{d}{d-1} \; \min \left\{ \vert B \vert + \frac{\vert A
\vert}{d}
\; \big\vert \; C_{A, B} \ne 0 \right\} = \frac{d}{d-1} \; v_L (g)\;,\\
\end{align*}
which proves (1) of Lemma \ref{lem2.1}.
Next assume $v_L (g) = n$ and $g' = g/u_1^n \in R'$. Then
$$
v_{\Lambda'} (g') = v_{\Lambda'} (g) - v_{\Lambda'} (u_1^n) =
\frac{d}{d-1} \; v_L (g) - \frac{n}{d-1} = n.
$$
Note that
$$
g' = \sum\limits_{A, B} \big(C_{A, B} \; u^{\prime A_u}\big)
y^{\prime^B}  u_1^{|A_u| + |B| -n }v^{A_v},
$$
$$
\vert B \vert + \frac{\vert A \vert}{d} = n \Longleftrightarrow \vert B
\vert + \frac{\vert A \vert + \vert B \vert - n}{d-1} = n.
$$
Therefore, with \eqref{eq2.2.1} we see
\begin{align*}
in_L (g) \in k [Y]  \Longleftrightarrow & \; C_{A, B} = 0 \quad
\mbox{if} \quad \vert B \vert + \frac{\vert A \vert}{d} = n \quad
\mbox{and} \quad \vert B
\vert < n \\
 \Longleftrightarrow & \; u^{\prime A_u} \; C_{A, B} = 0 \quad
\mbox{if} \quad \vert B \vert + \frac{\vert A \vert + \vert B
\vert - n}{d-1} = n \quad
\mbox{and} \quad \vert B \vert < n \\
 \Longrightarrow & \; in_{\Lambda'} (g') \in k' [Y']. \\
\end{align*}
(Note that the last implication is independent of the choice of
a representative $\Gamma'$ of $k'$.)
Moreover, if these conditions hold, we have
$$
in_L (g) = \sum\limits_{\stackrel{B}{\vert B \vert = v_L(g)}}
\overline{C_{0,B}} \; Y^B
\qaq
in_{\Lambda'} (g') = \sum\limits_{\stackrel{B}{\vert B \vert = v_L(g)}}
\overline{C_{0, B}} \; {Y'}^B.
$$
This completes the proof of Lemma \ref{lem2.1}.
$\square$

\bigskip
\textbf{Proof of Theorem \ref{prop2.1}}\quad
By \eqref{eqB1} and Lemma \ref{lem2.2} we have $\delta :=\delta(f,y,(u,v))>1$.
Choose $d$ with $1 < d < \delta$, and consider the linear forms $L$ on $\bR^e$
(downstairs) and $\Lambda'$ on $\bR^{e'}$ (upstairs) on ${\bR}$ as
in Lemma \ref{lem2.1}. As for the desired positive linear form in
Theorem \ref{prop2.1}, we take
$$
L'(A) = \frac{\sum a_i + \sum a'_j}{d-1} \quad \mbox{(upstairs)}
\quad (A=((a_i)_{1\leq i\leq s'},(a'_j)_{1\leq i\leq t}).
$$
By Lemma \ref{L0initial}(3) and Proposition \ref{prop1.1}, we have
$v_L (f_i) =v_\fm(f_i)= n_i$ and
\begin{equation}\label{eq2.3.2}
in_L (f_i) = in_{\fm} (f_i) = F_i (Y),\quad
In_L(J)=\langle F_1(Y),\dots,F_m(Y) \rangle  .
\end{equation}
Clearly $L' \geq \Lambda'$, so that by Proposition \ref{prop1.1}
it suffices to show
\begin{equation}\label{eq2.3.2.5}
in_{\Lambda'} (f'_i)  =  F_i (Y') \quad \text{ and }\quad
\Inh_{\Lambda'} (J)  =  \langle F_1 (Y'), \dots, F_m (Y') \rangle  .
\end{equation}
($\Lambda'$ satisfies the condition $v_{\Lambda'}(\Char(k'))>0$ in the
proposition since $\Char(k')=\Char(k) \in \fm R'=\langle u_1,v \rangle $).
The first part follows from \eqref{eq2.3.2} and Lemma \ref{lem2.1} (2)
To show the second part, choose any $g' \in J'$.
Take an integer $N > 0$ such that $g := u_1^N g' \in J$.
By \eqref{eq2.3.2}, Proposition \ref{prop1.1} implies that
there exist $\lambda_1, \dots,\lambda_m \in R$ such that,
\begin{equation}\label{eq2.3.3}
v_{\fp} (\lambda_i f_i) \geq v_{\fp} (g)
\end{equation}
\begin{equation}\label{eq2.3.4}
v_L (\lambda_i f_i) \geq v_L (g),
\end{equation}
\begin{equation}\label{eq2.3.5}
v_L (g - \sum\limits_i \lambda_i f_i) \gg N.
\end{equation}
where $\fp=(y,u)=(y_1,\dots,y_r,u_1,\dots,u_s)$. Here we used the fact that
$v_{\fp}=v_{L_\fp}$ with
$$
L_{\fp}(A) = \underset{1\leq i\leq s}{\sum} a_i  \quad \mbox{(downstairs)}
\quad (A=((a_i)_{1\leq i\leq s},(a'_j)_{1\leq i\leq t}).
$$
Let $v_{u_1}$ be the discrete valuation of $R'$ with respect to the ideal
$\langle u_1 \rangle \subset R'$. Because $\fp R' = \langle u_1 \rangle$,
\eqref{eq2.3.3} implies
$$
v_{u_1} (\lambda_i) = v_{\fp} (\lambda_i) \geq v_{\fp} (g) - n_i
= v_{u_1} (g) - n_i \geq N-n_i,
$$
where $n_i = v_{\fp} (f_i) = v_{\fm} (f_i)$ (cf. \eqref{eqB3}). Therefore
$$\lambda'_i := \lambda_i /u_1^{N-n_i} \in R'.$$
We calculate
\begin{equation}\label{eq2.3.6}
\begin{aligned}
v_{\Lambda'} (\lambda'_i f'_i)
 = & \; v_{\Lambda'} (\lambda_i f_i/u_1^N) \\
 = & \; \frac{d}{d-1} v_L (\lambda_i f_i) - \frac{N}{d-1} &
\mbox{(by Lemma \ref{lem2.1}(1))} \\
 \geq &\; \frac{d}{d-1} v_L (g) - \frac{N}{d-1}&\text{(by \eqref{eq2.3.4})} \\
 = & \; v_{\Lambda'} (g/u_1^N) = v_{\Lambda'} (g') &
\mbox{(by Lemma \ref{lem2.1}(1))} . \\
\end{aligned}
\end{equation}
\eqref{eq2.3.5} implies
\begin{align*}
v_{\Lambda'} \left( g' - \sum\limits_{1\leq i\leq m} \lambda'_i f'_i \right)
= & v_{\Lambda'} \left( \left(
g - \sum\limits_{1\leq i\leq m} \lambda_i f_i\right)/u_1^N\right) \\
 = & \frac{d}{d-1} v_L \left( g - \sum\limits_{1\leq i\leq m}
\lambda_i f_i\right) -\frac{N}{d-1} \gg 0\;. \\
\end{align*}
Therefore we may assume
\begin{equation}\label{eq2.3.7}
v_{\Lambda'} \left( g' - \sum\limits_i \lambda'_i f'_i \right) >
v_{\Lambda'} (g').
\end{equation}
By Lemma \ref{lem1.1} (2), \eqref{eq2.3.6} and \eqref{eq2.3.7} imply
$$
in_{\Lambda'} (g') = in_{\Lambda'} \left( \sum\limits_{1\leq i\leq m}
\lambda'_i f'_i\right) =
\sum\limits_i in_{\Lambda'} (\lambda'_i) in_{\Lambda'}(f'_i),
$$
where the last sum ranges over all $i$ such that
$v_{\Lambda'}(\lambda'_i f'_i)=v_{\Lambda'}(g')$.
This proves the second part of \eqref{eq2.3.2.5} and the proof of Theorem \ref{prop2.1} is complete.
\bigskip

We keep the assumptions and notations of Setup B and assume that $(f_1, \dots, f_m)$ is a standard base
of $J$. So by Theorem \ref{prop2.1} $f'=(f'_1,\dots,f'_m)$ is a
$(u_1,\phi,v)$-standard base of $J'=JR'$. Then Corollary \ref{cor.wellprepared}
assures that we can form a standard base of $J'$ from $f'$ by preparation if
$(u_1,\phi,v)$ is admissible for $J'$. Hence the following result is important.

\begin{theorem}\label{thm3.1} Let $k'=k(x')$.
If $x'$ is very near to $x$ (cf. Definition \ref{def.nearpoint}),
there exist linear forms
$L_1(U_1,V),\dots,L_r(U_1,V)\in k'[U_1,V]$ such that
$$
\IDir(R'/J')= \langle   Y'_1+L_1(U_1,V),\dots,Y'_r+L_r(U_1,V)\rangle \;
\subset \grmdRd=k'[Y,U_1,\Phi,V].
$$
In particular $(u_1,\phi,v)$ is admissible for $J'$.
\end{theorem}
\medbreak

Let $K/k'$ be a field extension. Consider the following map
$$
\psi: K[Y] \isom K[Y']  \hookrightarrow K[Y',\Phi] =
\grmdRd_K/\langle U_1,V\rangle ,
$$
where the first isomorphism maps $Y_i$ to $Y'_i$ for $i=1,\dots,r$.
Recall that
$$
\IDir(R/J)=\langle Y_1,\dots,Y_r\rangle \;\subset\grmR=k[Y,U,V],
$$
so that
$$
\IDir(R/J)_{/K}^{(1)}:=\IDir(R/J)_{/K}\cap \grm^1(R)_K\subset
\underset{1\leq i\leq r}{\bigoplus} K\cdot Y_i\; \subset K[Y].
$$
Theorem \ref{thm3.1} is an immediate consequence of the following more general
result.

\begin{theorem}\label{thm3.1a}
Assume that $x'$ is near to $x$. Then we have
$$
\psi(\IDir(R/J)_{/K}^{(1)}) \subset \IDir(R'/J')_{/K}
\mod \langle U_1,V\rangle  \text{  in }
\grmdRd_K/\langle U_1,V\rangle .$$
If $e(R/J)=e(R/J)_K$, we have
$$
\IDir(R'/J')_{/K} \supset \;
\langle Y'_1+L_1(U_1,V), \dots, Y'_r + L_r(U_1,V)\rangle
$$
for some linear forms
$L_1(U_1,V),\dots,L_r(U_1,V)\in K[U_1,V]$.
\end{theorem}
\medbreak

For the proof, we need the following.

\begin{proposition}\label{prop3.1}
If $x'$ is near to $x$, there exist
$h_{ij}\in  R'$ for $1\leq j<i\leq m$ such that setting
$$
g_i=f'_i-\sum\limits^{i-1}_{j=1}h_{ij}f'_j,
$$
we have the following for $i=1,\dots, m$:
\begin{equation}\label{eq.prop3.1}
v_{\fm'}(g_i) =n_i, \qaq in_{\fm'}(g_i) \equiv F_i(Y')
\mod \langle U_1,V\rangle
\text{  in } \gr_{{\frak m}'}(R').
\end{equation}
In particular, $(g_1,\dots,g_m)$ is a standard base of $J'$.
\end{proposition}

\bigskip
\textbf{Proof}
First we note that since $x'$ is near to $x$, we have
\begin{equation}\label{eq3.3.2}
\nu^*(J') = \nu^*(J) = (n_1,n_2,\ldots , n_m,\infty,\infty,\ldots )
\quad (n_1\leq n_2 \leq \ldots \leq n_m)
\end{equation}
The last assertion of the proposition follows from
\eqref{eq.prop3.1} and \eqref{eq3.3.2} in view of \cite{H1}, Ch. III Lemma 3.
Put
$$
\tFidg:= \omega_{(y',u_1,\phi,v,\Gamma')}(F_i(Y'))\,,
\quad\text{where  }\quad
 \omega_{(y',u_1,\phi,v,\Gamma')}:k'[[Y',U_1,\Phi,V]] \to \hat{R'}
$$
is the map \eqref{eq1.2.1.d} for $(y',u_1,\phi,v)$ and $\Gamma'$
(notations as in Setup B). \eqref{eqB2} implies
\begin{equation}\label{eqB2d}
f'_i \equiv \tFidg + \lambda_i\mod (\fm')^{n_i+1} \qwith \lambda_i\in
\langle u_1,v \rangle \qfor i=1,\dots,m.
\end{equation}
To prove \eqref{eq.prop3.1}, it suffices to show that there are
$h_{ij}\in \hat{R'}$ for $1\leq j<i,\leq m$ such that
letting $g_i$ be as in the proposition, we have
\begin{equation}\label{eq2.prop3.1}
v_{\fm'}(g_i) =n_i, \qaq  g_i- \tFidg \in\;\langle u_1,v\rangle+(\fm')^{n_i+1}
\subset \hat{R'}.
\end{equation}
Indeed \eqref{eq.prop3.1} follows from \eqref{eq2.prop3.1}
by replacing the $h_{ij}$ with elements of $R'$ sufficiently close to them.
\eqref{eqB2d} implies
$$
v_{\fm'}(f_i')\leq n_i=v_{\fm'}(\tFidg)\qfor i=1,\dots,m.
$$
\eqref{eq3.3.2} and Corollary \ref{cordef0.1} imply $v_{\fm'}(f_1')= n_1$ so that
one can take $g_1=f_1'$ for \eqref{eq2.prop3.1}.
Let $\ell$ be the maximal $t\in \{1,\ldots , m\}$ for which the following
holds.
\begin{itemize}
\item[($*_{t}$)] There exist $h_{ij}\in \hat{R'}$ for
$1\leq j<i\leq t$, such that \eqref{eq2.prop3.1} holds for
\par\noindent
$g_i= f_i'-\sum\limits^{i-1}_{j=1} h_{ij}f'_j$ with $i=1,\dots,t$.
\end{itemize}
We want to show $\ell=m$. Suppose $\ell<m$.
Then $v_{\fm'}(f'_{\ell+1})< n_{\ell+1}$ since otherwise
\eqref{eq2.prop3.1} holds for $g_{\ell+1}=f'_{\ell+1}$,
which contradicts the maximality of $\ell$. This implies
\begin{equation}\label{eqq0}
in_{\fm'}(f'_{\ell+1})=in_{\fm'}(\lambda_{\ell+1})\; \in \;
\langle U_1,V\rangle \subset k'[Y',U_1,\Phi,V].
\end{equation}
By the assumption we have
$$
G_i:=\inmd(g_i) \equiv F_i(Y')\mod \langle U_1,V\rangle \qfor i=1,\dots, \ell.
$$
Since $(F_1,\dots, F_m)\subset k[Y]$ is normalized
(cf. Definition \ref{def0.2}), $(G_1,\ldots, G_\ell)$ is normalized
in $\grmdRd$. Therefore \eqref{eq3.3.2} and Corollary \ref{cordef0.1}
imply that there exist $H_i\in \grmdRd$ homogeneous of degree
$v_{\fm'}(f'_{\ell+1})-n_i$ for $i=1,\dots,\ell$, such that
\begin{equation}\label{eqq1}
in_{\fm'}(f'_{\ell+1})= \sum_{1\leq i\leq \ell} H_i G_i.
\end{equation}
Setting
$$
g^{(1)}_{\ell+1}=f'_{\ell+1}-\sum_{1\leq i\leq \ell} \tilde{H}_i g_i
\qwith \tilde{H}_i=\omega_{(y',u_1,\phi,v,\Gamma')}(H_i),
$$
we have
$v_{\fm'}(g^{(1)}_{\ell+1})>v_{\fm'}(f'_{\ell+1}).$
We claim
\begin{equation}\label{eq.claim}
g^{(1)}_{\ell+1}- \tilde{F}'_{n_{\ell+1},\Gamma'}
\in \langle u_1,v\rangle +(\fm')^{n_{\ell+1}+1},
\end{equation}
which completes the proof. Indeed \eqref{eq.claim} implies
$$
v_{\fm'}(g^{(1)}_{\ell+1})\leq v_{\fm'}(\tilde{F}'_{n_{\ell+1},\Gamma'})
=n_{\ell+1}.
$$
If $v_{\fm'}(g^{(1)}_{\ell+1})=n_{\ell+1}$,
\eqref{eq2.prop3.1} holds for $g^{(1)}_{\ell+1}$,
which contradicts the maximality of $\ell$.
If $v_{\fm'}(g^{(1)}_{\ell+1})<n_{\ell+1}$, we apply the same argument
to $g^{(1)}_{\ell+1}$ instead of $f'_{\ell+1}$ and find
$h_1,\dots,h_\ell\in R'$ such that setting
$$
g^{(2)}_{\ell+1}=g^{(1)}_{\ell+1}-\sum_{1\leq i\leq \ell} h_i g_i,
\quad\text{  we have}\quad
$$
$$
v_{\fm'}(g^{(1)}_{\ell+1})<v_{\fm'}(g^{(2)}_{\ell+1})\leq n_{\ell+1}
\qaq
g^{(2)}_{\ell+1}- \tilde{F}'_{n_{\ell+1},\Gamma'}
\in \langle u_1,v\rangle +(\fm')^{n_{\ell+1}+1}.
$$
Repeating the process, we get $g_{\ell+1}$ for which \eqref{eq2.prop3.1}
holds, which contradicts the maximality of $\ell$.
We now show claim \eqref{eq.claim}.
Noting
$$
f'_{\ell+1} -\tilde{F}'_{n_{\ell+1},\Gamma'},\;\;
g_i - \tFidg \in \;\langle u_1,v\rangle +(\fm')^{n_i+1}\quad (i=1,\dots,\ell),
$$
it suffices to show
$$
\sum_{1\leq i\leq \ell} \tilde{H}_i \tFidg \in \langle u_1,v\rangle .
$$
For $i=1,\dots, \ell$, write
$H_i=H_i^- + H_i^+$, where $H_i^-\in k'[Y',\Phi]$ and $H_i^+\in
\langle U_1,V\rangle $
and both are homogeneous of degree $v_{\fm'}(f'_{\ell+1})-n_i$.
Similarly we write
$G_i= F_i(Y') + G_i^+$ with $G_i^+\in \langle U_1,V\rangle $.
Then \eqref{eqq0} and \eqref{eqq1} imply
$$
\sum_{1\leq i\leq \ell} H_i^- F_i(Y')=0.
$$
Let
$\tilde{H}_i^\pm = \omega_{(y',u_1,\phi,v,\Gamma')}(H_i^\pm)$.
Noting $\Char(k)\in \fm R'=\langle u_1,v\rangle $, property \eqref{eq1.2.1.e} implies
$$
\sum_{1\leq i\leq \ell} \tilde{H}_i^- \tFidg,\;
\tilde{H}_i^+, \; \tilde{H}_i- (\tilde{H}_i^+ + \tilde{H}_i^-)
\;\in \;\langle u_1,v\rangle ,
$$
which shows the desired assertion.
$\square$

\bigskip{\bf Proof of Theorem \ref{thm3.1a}}
The second assertion follows at once from the first one, and we show the first.
By Proposition \ref{prop3.1} there exist homogeneous
$G_1,\dots,G_m \in \grmdRd$ such that
\begin{equation}\label{Th3.1.eq1}
\begin{aligned}
&G_j\equiv F'_j:=F_j(Y')\mod \langle U_1,V\rangle
\quad\text{for all }j=1,\dots,m,\\
&\InmJd_K=\langle G_1,\dots,G_m\rangle \; \subset \grmdRd_K = K[Y',U_1,\Phi,V].
\end{aligned}
\end{equation}
Let $W'=\IDir(R'/J')_{/K}\cap \grmd^1(R')$, then there exist
$H_1,\dots,H_n\in K[W']\cap \InmJd_K$ such that
$$
G_i=\underset{1\leq j\leq n}{\sum} h_{i,j} H_j
\quad\text{ for some homogeneous  }h_{i,j}\in \grmdRd.
$$
On the other hand, \eqref{Th3.1.eq1} implies
$$
H_i=\underset{1\leq j\leq m}{\sum} g_{i,j} G_j
\quad\text{ for some homogeneous  }g_{i,j}\in \grmdRd.
$$
Regarding everything $\mod \langle U_1,V\rangle $, we get
$$
F'_i=\underset{1\leq j\leq n}{\sum} \ol{h}_{i,j} \ol{H}_j
\qaq
\ol{H}_i=\underset{1\leq j\leq m}{\sum} \ol{g}_{i,j} F'_j
\;\text{ in  } K[Y',\Phi]=\grmdRd/\langle U_1,V\rangle .
$$
The second equality implies
$$
\ol{H}_i\in K[\ol{W'}]\;\cap \langle F'_1,\dots,F'_m\rangle \;
\subset K[Y',\Phi]\quad\text{ where } \ol{W'}=W'\mod \langle U_1,V\rangle,
$$
so that the first equality implies
\begin{equation}\label{thm8.2.eq1}
\big(K[\ol{W'}]\;\cap \langle F'_1,\dots,F'_m\rangle \big) \cdot
K[Y',\Phi]=\langle F'_1,\dots,F'_m\rangle .\end{equation}
Since $W=\IDir(R/J)_{/K}\cap \grm^1(R)$ is the minimal subspace of
$\underset{1\leq i\leq r}{\bigoplus} K\cdot Y_i$ such that
$$
\big(K[W]\;\cap \langle F_1(Y),\dots,F_m(Y)\rangle \big) \cdot K[Y]=
\langle F_1(Y),\dots,F_m(Y)\rangle ,
$$
\eqref{thm8.2.eq1} implies $\psi(W)\subset \ol{W'}$,
which is the desired assertion.
$\square$

\medbreak\medbreak
We conclude this section with the following useful criteria for the nearness
and the very nearness of $x'$ and $x$.
We keep the assumptions and notations of Setup B.

\begin{theorem}\label{deltanearcriteria}
\begin{itemize}
\item[(1)]
If $\delta(f',y',(u_1,\phi,v))\geq 1$, $x'$ is near to $x$.
The converse holds if $(f',y',(u_1,\phi,v))$ is prepared at any vertex lying in
$\{A\in \bR^{s'+t}|\; |A|\leq 1\}$.
\item[(2)]
Assume $e_x(X)=e_x(X)_{k(x')}$.
If $\delta(f',y',(u_1,\phi,v))>1$, $x'$ is very near to $x$.
The converse holds under the same assumption as in (1).
\end{itemize}
\end{theorem}
\textbf{Proof }
Write $k'=k(x')$.
Assume $\delta(f',y',(u_1,\phi,v))\geq 1$.
By Lemma \ref{lem2.2} (1) and \eqref{eqB2.1}, $v_{\fm'}(f'_i)=n_i$ and
we can write
$$
in_{\fm'}(f_i')= F(Y')+ \underset{|B|<n_i}{\sum} {Y'}^B P_B(U_1,\Phi,V),
\quad P_B(U_1,\Phi,V)\in k'[U_1,\Phi,V].
$$
Put $I=\langle in_{\fm'}(f_1'),\dots, in_{\fm'}(f_m') \rangle\subset
\grmdRd$. We have
$$
\nu^*(J) \geq \nu^*(J')=\nu^*(\InmJd) \geq \nu^*(I),
$$
where the first inequality follows from Theorem \ref{thm.pbu.inv}
and the last from \cite{H1} Lemma Ch.II Lemma 3.
By the assumption $(F_1(Y),\dots,F_m(Y))$ is weakly normalized
(cf. Definition \ref{def0.2}), which implies that
$(in_{\fm'}(f_1'),\dots, in_{\fm'}(f_m'))$ is weakly normalized so that
$\nu^*(I)=\nu^*(J)$. Therefore  we get $\nu^*(J) = \nu^*(J')$ and
$x'$ is near to $x$. We also get $\nu^*(\InmJd) =\nu^*(I)$,
which implies $\InmJd = I$ by loc.cit..
Now assume $x'$ is near to $x$. Let $g=(g_1,\dots,g_m)$ be as in Proposition
\ref{prop3.1}. By \eqref{eq.prop3.1} we have
$\delta(g,y',(u_1,\phi,v))\geq 1$.
We claim that $g$ is a $(u_1,\phi,v)$-standard basis of $J'$.
Indeed $f'$ is a $(u_1,\phi,v)$-standard basis of $J'$
by Theorem \ref{prop2.1}, so the claim follows from the fact that $in_0(f')=in_0(g)$
by using Corollary \ref{cor.refdatum0}.
By the claim $\Delta(J',(u_1,\phi,v))\subset \Delta(g,y',(u_1,\phi,v))$
so that there exists no vertex $w$ of $\Delta(J',(u_1,\phi,v))$ such that
$|w|<1$. If $\Delta(f',y',(u_1,\phi,v))$ is prepared at any vertex $w$ with $|w|\leq 1$,
we obtain $\delta(f',y',(u_1,\phi,v))\geq 1$ from Theorem \ref{thm.wellprepared2}.

\medbreak
Assume $\delta(f',y',(u_1,\phi,v))> 1$.
By Lemma \ref{lem2.2} (1) and \eqref{eqB2.1},
$in_{\fm'}(f_i')= F(Y')$ so that
\begin{equation}\label{deltanearcriteria.eq1}
\inmd(J') =I=\langle F_1(Y'),\dots,F_m(Y')\rangle.
\end{equation}
By the assumed equality $e_x(X)=e_x(X)_{k'}$, we have
$I\Dir_{k'}(R/J)=\langle Y_1,\dots, Y_r\rangle$ (cf. Remark \ref{rem.dir}).
Since $\InmJ =\langle F_1(Y),\dots,F_m(Y)\rangle$, this implies
$I\Dir(R'/J')=\langle Y'_1,\dots, Y'_r\rangle$
by \eqref{deltanearcriteria.eq1} so that $x'$ is very near to $x$.
Finally assume $x'$ very near to $x$.
By Theorem \ref{prop2.1}, $f'$ is $(u_1,\phi,v)$-standard basis of $J'$.
By Theorem \ref{thm3.1}, $(u_1,\phi,v)$ is admissible for $J'$.
Thus Corollary \ref{cor.wellprepared2} implies
$\delta(f',y',(u_1,\phi,v)) >1$. This completes the proof of Theorem
\ref{deltanearcriteria}.
$\square$
\medskip

\begin{remark}\label{deltanearcriteria.rem}
By the above theorem it is important to compute $\delta(f',y',(u_1,\phi,v))$
from $\Delta(f,y,(u,v))$. It is also important to
see if the well-preparedness of $\Delta(f,y,(u,v))$ implies that of
$\Delta(f',y',(u_1,\phi,v))$. These issues are discussed later in this paper
in various situations (e.g., see Lemma \ref{ekeylemma0}).
\end{remark}

\newpage
\section{Termination of the fundamental sequences of $B$-permissible blowups, and the case $e_x(X)=1$}\label{sec:fund.seq}

\def\tL{\widetilde{L}}
\medbreak
\def\Retaq{R_{\eta_q}}
\def\Jetaq{J_{\eta_q}}
\def\fmetaq{\fm_{\eta_q}}
\bigskip

In this section we prove the Key Theorem \ref{fu.thm0} in \S\ref{sec:strategy}, by deducing it
from a stronger result, Theorem \ref{fu.thm} below. Moreover we will give an explicit bound on
the length of the fundamental sequence, by the $\delta$-invariant of the polyhedron at the beginning.
First we introduce a basic setup.
\medbreak

\textbf{Setup C: }
Let $Z$ be an excellent regular scheme, let $X\subset Z$ be a closed subscheme
and take a point $x\in X$. Let $R=\cO_{Z,x}$ with maximal ideal
$\fm$ and residue field $k=R/\fm=k(x)$, and
write $X\times_Z \Spec(R)=\Spec(R/J)$ with an ideal $J\subset \fm$.
Define the integers $n_1\leq n_2\leq \cdots\leq n_m$ by
$$
\nu^*_x(X,Z)=\nu^*(J)=(n_1,\dots,n_m,\infty,\infty,\dots).
$$
We also assume given a simple normal crossing boundary $\cB$ on $Z$
and a history function
$\OB: X \to \{\text{subsets of $\cB$}\}$
for $\cB$ on $X$ (Definition \ref{def.histfcnZ}).
Note that $\cB$ may be empty.
\bigskip

We introduce some notations.

\begin{definition}\label{def.label}
\begin{itemize}
\item[(1)]
A prelabel of $(X,Z)$ at $x$ is
$$
(f,y,u)=(f_1\dots, f_m,\; y_1,\dots,y_r, u_1,\dots,u_e),
$$
where $(y,u)$ is a system of regular parameters of $R$ such that
$(u)$ is admissible for $J$ (cf. Definition \ref{def2.11}) and $f$ is a
$(u)$-standard base of $J$. By Lemma \ref{def1.4cor} we have
$n_{(u)}(f_i)=n_i$ for $i=1,\dots, N$.
\item[(2)]
A prelabel $(f,y,u)$ is a label of $(X,Z)$ at $x$ if
$\delta(f,y,u)>1$. By Corollary \ref{cor.wellprepared2}, this means that
$(y,u)$ is strictly admissible for $J$ and $f$ is a standard base of $J$
such that $(f,y,u)$ is admissible. By Lemma \ref{def1.4cor} we have
\begin{equation*}\label{ec.eq2}
\begin{aligned}
&v_{\fm}(f_i)=n_i \qaq
\inm(f_i) =in_0(f_i) \qfor i=1,\dots, N,\\
& \langle Y_1,\ldots,Y_r \rangle =\IDir(R/J)\; \subset \grmR,\\
& In_{\fm}(J)=\langle F_1(Y),\dots,F_N(Y) \rangle \qwith F_i(Y)=in_0(f_i) \in k[Y],\\
\end{aligned}
\end{equation*}
where $k[Y]= k[Y_1\dots,Y_r] \subset \grmR$ with $Y_i=\inm(y_i)\in \grmR$.
\item[(3)]
A label $(f,(y,u))$ is well-prepared (resp.
totally prepared) if so is $(f,y,u)$ in the sense of Definition
\ref{def.wellprepared} (resp. Remark \ref{remark.normalizededges}).
By Theorem \ref{thm.wellprepared2}, for a well-prepared label $(f,y,u)$,
we have $\Delta(J,u)=\Delta(f,y,u)$.
\end{itemize}
\end{definition}
\medbreak

For each $B\in \cB$ choose an element $l_B\in R$ such that
$B\times_Z\Spec(R)=\Spec(R/\langle l_B \rangle)$.
For a positive linear form $L:\bR^e\to \bR$ let
$\delta_L(l_B,y,u)$ be as in Definition \ref{def.intial} (3).
Writing
\begin{equation}\label{fu.eq9}
l_B=\underset{1\leq i\leq r}{\sum} a_i y_i +
\underset{A\in \bZ^e_{\geq 0}}{\sum} b_A u^A + g
\quad (a_i,\; b_A\in R^\times\cup \{0\}),
\end{equation}
where $g\in {\langle y_1,\dots,y_r \rangle}^2$, we have
$$
\delta_L(l_B,y,u)=\min\{L(A)\;|\; b_A\not\in \fm\}.
$$
It is easy to see that $\delta_L(l_B,y,u)$ depends only on $B$ and not on
the choice of $l_B$. We define
$$
\deltaob_L(y,u)=\min\{\delta_L(l_B,y,u)\;|\;
\text{$B$ irreducible component of $\OBx$}\}.
$$

\begin{theorem}\label{fu.thm}
Assume $\Char(k(x)) = 0$ or $\Char(k(x)) \geq \dim (X)/2+1$, and assume there is
a fundamental sequence of length $\geq m$ starting with $(X,\cB_X,\OB)$ and $x$
as in \eqref{esbpbu.eq} of Definition \ref{Def.fupb}, where $m\geq 1$.
Let $(f,y,u)$ be a label of $(X,Z)$ at $x$.
In case $e_x(X)>\eob_x(X)$, assume that $u=(u_1,\dots,u_e)$ ($e=e_x(X)$)
satisfies the following condition:
\begin{itemize}
\item[]
There exist $B_j\in \OBx$ for
$2\leq j \leq s:=e_x(X)-\eob_x(X)+1$ such that
$$
B_j\times_Z \Spec(R)=\Spec(R/\langle u_{j}\rangle)
\qaq
\Dirob_x(X)=\Dir_x(X)\cap \underset{2\leq j\leq s}{\cap} T_x(B_j).
$$
\end{itemize}

Assume further $(f,y,u)$ is prepared along the faces $E_{L_q}$
for all $q=0,1,\dots, m-1$, where
$$
L_q:\bR^e\to \bR\;;\; A=(a_1,\dots,a_e)\to |A|+q\cdot \sum_{j=2}^s a_j.
$$
Then we have
\begin{equation}\label{fu.thm.eq}
 \delta_{L_q}(f,y,u) \geq q+1\qaq \deltaob_{L_q}(y,u)\geq q+1
\qfor q=1,\dots, m-1.
\end{equation}
\end{theorem}
\bigskip

First we show how to deduce Theorem \ref{fu.thm0} from Theorem \ref{fu.thm}.
It suffices to show that under the assumption of \ref{fu.thm0},
there is no infinite fundamental sequence of $\cB_X$-permissible blowups over $x$.
Assume the contrary.
As in \eqref{eq1.2.1} write
\begin{equation}\label{fu.eq.-1}
f_i = \sum\limits_{(A, B)} C_{i, A, B} \; y^B u^A \quad \mbox{with}
\quad C_{i, A, B} \in R^{\times} \cup \{ 0 \}\,.
\end{equation}
Write $|A|_1=\sum_{j=2}^s a_j$ for $A=(a_1,\dots,a_e)$.
By Definition \ref{def.intial} (3), we have
$$
\begin{aligned}
\delta_{L_q}(f,y,u) >q+1
&\Leftrightarrow |A|+ q |A|_1 \geq (q+1)(n_i-|B|)
\text{  if  }|B|<n_i,\; C_{i, A, B}\not=0\\
&\Leftrightarrow |B|+ |A|_1 \geq n_i- \frac{|A|-|A|_1}{q+1}
\text{  if  } |B|<n_i,\; C_{i, A, B}\not=0 \,.\\
\end{aligned}
$$
Hence by \eqref{fu.thm.eq} and the assumption we have the last statement for all $q\in \bN$,
which implies
$|B|+ |A|_1 \geq n_i$ for all $A$ and $B$ such that $C_{i, A, B}\not=0$
in \eqref{fu.eq.-1}.
Setting $\fq=\langle y_1,\dots,y_r,u_2,\dots,u_s \rangle \in \Spec(R)$,
this implies $J\subset \fq$ and $v_\fq(f_i)=n_i=v_\fm(f_i)$ for all $i$.
By Theorems \ref{nf.thm1} $(iv)$ and \ref{HSf.perm} we conclude
$\HSfX \omega= \HSfX x$, where $\omega\in X$ is the image of $\fq\in \Spec(R)$.
By a similar argument, the second part of \eqref{fu.thm.eq} implies
$\omega\in B$ for all $B\in \OBx$ so that $\tHSfX \omega= \tHSfX x$.
This contradicts the assumption of Theorem \ref{fu.thm0} since $\dim(R/\fq)=e-s+1=\eob_x(X)$.
This completes the proof of Theorem \ref{fu.thm0}.
$\square$
\bigskip

Now we prepare for the proof of Theorem \ref{fu.thm}. Consider
$$
\pi: Z'=\Bl_x(Z)\to Z \qaq \pi: X'=\Bl_x(X)\to X.
$$
By Theorem \ref{thmdirectrix}, any point of $X'$ near to $x$ is contained in
$\bP(\Dir_x(X)) \subset X'$. We now take a label $(f,y,u)$ of $(X,Z)$ at
$x$ as in Definition \ref{def.label} (2).
By \eqref{eq.u-coord} we have the identification determined by $(u)$:
\begin{equation}\label{ec.eq1.5}
\psi_{(u)} : \bP(\Dir_x(X)) = \Proj(k[T]) =\bP^{e-1}_k\quad
(k[T]=k[T_1,\dots,T_e])
\end{equation}

Let $x'=(1:0:\cdots:0) \in \bP(\Dir_x(X))$. If $x'$ is near to $x$,
Theorems \ref{thm.pbu.inv} and \ref{def.nearpoint} imply
$$
\nu^*_{x'}(X',Z') = \nu^*_x(X,Z) \qaq e_{x'}(X')\leq e.
$$
Let $R'=\cO_{Z',x'}$ with maximal ideal $\fm'$, and let $J'\subset R'$
be the ideal defining $X'\subset Z'$. Note that
$$
\bP(\Dir_x(X))\times_{Z'} \Spec(R')=\Spec(R'/\langle u_1\rangle ),
$$
and denote
$$
y'=(y'_1,\dots,y'_r), \;  u'=(u'_2,\dots,u'_e), \; f'=(f'_1,\dots,f'_N) \;,
$$
where
$y'_i=y_i/u_1,\; u'_i=u_i/u_1,\;  f_i'=f_i/u_1^{n_i}$. As is seen in
Setup B in \S8, $(y', u_1,u')$ is a system of regular parameters of $R'$
and $J'=\langle f'_1,\dots,f'_N\rangle$ and $R'$ is the localization of
$R[y'_1,\dots,y'_r,u'_2,\dots,u'_e]$ at $(y', u_1,u')$.
We will use the usual identifications
$$
\grmdRd=k[Y',U_1,U'_2,\dots,U'_e],\quad
Y'_i=\inmd(y'_i),\; U'_j=\inmd(u'_j).
$$
Moreover we will consider the maps:
$$\Psi : \bR^e\rightarrow\bR^e\;;\; (a_1,a_2,\dots,a_e)\mapsto
(\sum_{i=1}^e a_i-1,a_2,\dots,a_e),$$
$$\Phi : \bR^e\rightarrow\bR^e\;;\; (a_1,a_2,\dots,a_e)\mapsto
(\sum_{i=1}^e a_i,a_2,\dots,a_e).$$
A semi-positive linear from $L:\bR^e\to \bR$ is called monic if
$L(1,0,\dots,0)=1$.

\begin{lemma}\label{ekeylemma0}
\begin{itemize}
\item[(1)]
$\Delta(f', y', (u_1, u'))$ is the minimal $F$-subset containing
$\Psi(\Delta(f,y,u))$.
\item[(2)]
For any monic semi-positive linear form
$L:\bR^e\to \bR$ and $\tL=L\circ\Phi$,
we have
$$
\delta_L(f',y',(u_1,u'))=\delta_{L\circ\Phi}(f,y,u)-1,
$$
$$
\begin{aligned}
&U_1^{n_i}\cdot in_{E_L}(f'_i)_{(y',(u_1,u'))}=
{in_{E_{\tL}}(f_i)_{(y,u)}}_{|Y=U_1Y',U_i=U_1U'_i\; (2\leq i\leq e)}\;,\\
&in_{E_{\tL}}(f_i)_{(y,u)}=U_1^{n_i}\cdot
{in_{E_L}(f_i)_{(y',(u_1,u'))}}_{|Y'=Y/U_1,U'_i=U_i/U_1\; (2\leq i\leq e)}\;.\\
\end{aligned}
$$
If $(f,y,u)$ is prepared along $E_{\tL}$, then
$(f', y', (u_1, u'))$ is prepared along $E_{L}$.
\item[(3)]
Assume $(f,y,u)$ prepared along $E_{\tL_0}$ with $\tL_0=L_0\circ \Phi$.
Then $(f', y', (u_1, u'))$ is $\delta$-prepared and $x'$
is near (resp. very  near) if and only if
$\delta_{\tL_0}(f,y,u)\geq 2$ (resp. $\delta_{\tL_0}(f,y,u) >2$).
If $(f,y,u)$ is totally prepared, so is $(f', y', (u_1, u'))$.
\item[(4)]
Assume $x'$ very near to $x$. Then $(f',y',(u_1,u'))$ is a prelabel of $(X',Z')$ at $x'$. Assume further $(f,y,u)$ prepared along
$E_{\tL_0}$ (resp. totally prepared), then $(f',y',(u_1,u'))$ is
a $\delta$-prepared (resp. totally prepared) label of $(X',Z')$ at $x'$.
\end{itemize}
\end{lemma}
\medbreak
\textbf{Proof  }
From \eqref{fu.eq.-1} we compute
\begin{equation}\label{fu.eq0}
f'_i=f_i/u_1^{n_i} = \sum\limits_{(A, B)} C_{i, A, B} \;
(y')^B {u'}^A u_1^{|A| + (|B|-n_i)} ,
\end{equation}
$$
{u'}^A={u'_2}^{a_2}\cdots {u'_e}^{a_e} \qfor A=(a_1,a_2,\dots,a_e).
$$
(1) follows at once from this.
For a semi-positive linear form $L:\bR^e\to \bR$, we have
$$
\begin{aligned}
\delta_L(f',y',(u_1,u'))
&=\min\{\frac{L(|A|+|B|-n_i,\; a_2,\dots,a_e)}{n_i-|B|}\;|\; |B|<n_i,\; C_{i,A,B}\not=0 \}\\
&=\min\{\frac{L\circ\Phi(A)}{n_i-|B|}-1\;|\; |B|<n_i,\; C_{i,A,B}\not=0\}\\ &=\delta_{L\circ\Phi}(f,y,u)-1.
\end{aligned}
$$
From \eqref{fu.eq0}, we compute
$$
in_{E_L}(f')_{(y',(u_1,u'))}
=F_i(Y)+\underset{B,A}{\sum}\Cb_{i,A,B} {Y'}^B {U'}^A U_1^{|A|+|B|-n_i},
$$
where the sum ranges over such $B, A$ that $|B|<n_i$ and
$L\circ\Phi(A)=\delta_{L\circ\Phi}(n_i-|B|)$. (2) follows easily from this.
(3) follows from (2) and Theorem \ref{deltanearcriteria}.
The first assertion of (4) follows from Theorem \ref{prop2.1}
and Theorem \ref{thm3.1}. The other assertion of (4) follows from (3) and
Corollary \ref{cor.wellprepared2}.
$\square$

\bigskip

We now consider
$$
C':=\bP(\Dirob_x(X))\subset \bP(\Dir_x(X))=\Proj(k[T_1,\dots,T_e]) \subseteq X'.
$$
Let $\eta'$ be the generic point of $C'$, and let
$$
t:=\eob_x(X)\qaq s=e-t+1.
$$
We assume $t\geq 1$ so that $1\leq s\leq e$.
By making a suitable choice of the coordinate $(u)=(u_1,\dots,u_e)$,
we may assume:
\begin{equation}\label{ekeylemma2.cond}
\text{there exists $B_j\in \OBx$ for $2\leq j\leq s$ such that}
\end{equation}
$$
B_j\times_Z \Spec(R)=\Spec(R/\langle u_{j}\rangle)
\qaq
\Dirob_x(X)=\Dir_x(X)\cap \underset{2\leq j\leq s}{\cap} T_x(B_j).
$$
Then
$$
C'\times_{Z'} \Spec(R')=\Spec(R'/\fp')\qwith
\fp'=\langle y',u_1,u'_2,\dots,u'_s \rangle.
$$
Note $\delta_{k(\eta')/k}:=\trdeg_k(k(\eta'))=t-1$.
By Theorems \ref{def.nearpoint} and \ref{nf.dir}, if $\eta'$ is near $x$,
we have $e_{\eta'}(X')\leq e-\delta_{k(\eta')/k}=s$. It also implies
that $C' \subset X'$ is permissible by Theorems \ref{HSf.usc} and
\ref{HSf.perm}. By Theorem \ref{def.Bnearpoint}, if $\eta'$ is very near $x$,
we have $\eob_{\eta'}(X')\leq \eob_x(X)-\delta_{k(\eta')/k}=1$.
\medbreak

Write $\Rdetad=\cO_{Z',\eta'}$ and let $\fmdetad$ be its maximal ideal and
$J'_{\eta'}=J'\Rdetad$.
Note that $\Rdetad$ is the localization of $R'$ at $\fp'$ and
$(y',u_1,u'_2,\dots,u'_s)$ is a system of regular parameters of $\Rdetad$.

\begin{lemma}\label{ekeylemma2}
Let $\Delta(f',y',(u_1,u'_2,\dots,u'_s))$ be the characteristic polyhedron for
$J'_{\eta'}\subset \Rdetad$.
\begin{itemize}
\item[(1)]
$\Delta(f', y', (u_1,u'_2,\dots,u'_s))$ is the minimal $F$-subset containing
$\pi\cdot\Psi(\Delta(f,y,u))$, where $\Psi$ is as in Lemma \ref{ekeylemma0} and
$$
\pi : \bR^e\rightarrow\bR^s\;;\;
A=(a_1,a_2,\dots,a_e)\mapsto (a_1,a_2,\dots,a_s).
$$
For any monic semi-positive linear form $L:\bR^s\to \bR$, we have
$$
\delta_L(f',y',(u_1,u'_2,\dots,u'_s))=\delta_{L\circ\pi\circ\Phi}(f,y,u)-1.
$$
In particular,
$$
\delta(f', y', (u_1,u'_2,\dots,u'_s))=\delta_{L_1}(f,y,u)-1,
$$
$$
L_1 : \bR^e\rightarrow\bR\;;\; A\mapsto
|A|+ \underset{2\leq i\leq s}{\sum} a_i.
$$
\end{itemize}
\begin{itemize}
\item[(2)]
If $(f,y,u)$ is prepared along the face $E_{L_1}$ of $\Delta(f,y,u)$, then
$(f', y', (u_1,u'_2,\dots,u'_s))$ is $\delta$-prepared.
If $(f,y,u)$ is totally prepared, so is $(f', y', (u_1,u'_2,\dots,u'_s))$.
\item[(3)]
Assume $(f,y,u)$ is prepared along the face $E_{L_1}$ of $\Delta(f,y,u)$.
Then $\eta$ is near $x$ if and only if $\delta_{L_1}(f,y,u)\geq 2$.
if this holds, we have
$v_{\fmdeta}(f'_i)=v_{\fm}(f_i)=n_i$ for $i=1,\dots,N$.
\item[(4)]
Assume $(f,y,u)$ is prepared along the face $E_{L_1}$ of $\Delta(f,y,u)$.
Then $\eta$ is very near $x$ if and only if $\delta_{L_1}(f,y,u)>2$.
if this holds, then $(f',(y',(u_1,u'_2,\dots,u'_s)))$ is
a $\delta$-prepared label of $(X',Z')$ at $\eta$.
\end{itemize}
\end{lemma}
\medbreak
\textbf{Proof  }
(1) and (2) follows from Lemma \ref{ekeylemma0} and Theorem
\ref{lem.projection} applied to $(f', y', (u_1,u'))$ and
$\fp'=(u_1,u'_2,\dots,u'_s)$ and \eqref{fu.eq0} in place of
$(f,y,u)$ and $\fp_s=(u_1,u_2\dots,u_s)$ and \eqref{eq1.proj}.
We need to check the conditions $(P0)$ and $(P2)$ as well as
Theorem \ref{lem.projection} (3)$(i)$ and $(ii)$ for the replacement.
$(P0)$ holds since $k=R/\fm\hookrightarrow R'/\langle u_1\rangle$.
In view of the presentation \eqref{fu.eq0}, $(P2)$ holds by the fact that
for fixed $B$ and $a\in \bZ$, there are only finitely many
$A\in \bZ^e_{\geq 0}$ such that $|A|+|B|-n_i=a$.
Theorem \ref{lem.projection} (3)$(i)$ is a consequence of the assumption
that $(f,y,u)$ is a label of $(X,Z)$ at $x$
(cf. Definition \ref{def.label} (2)).
Finally Theorem \ref{lem.projection} (3)$(ii)$ follows from Lemma
\eqref{ekeylemma0}(2).
(3) and (4) are consequences of (2) and Theorems \ref{deltanearcriteria}
and \ref{prop2.1} and \ref{thm3.1} and Lemma \ref{lem2.2}.
This completes the proof of Lemma \ref{ekeylemma2}.
$\square$
\bigskip

\textbf{Proof of Theorem \ref{fu.thm}.}
Write $Z_1=Z'$, $X_1=X'$, $C_1=C'$ and assume $\pi:Z'\to Z$ extends to
a sequence \eqref{esbpbu.eq}.
Let $(f,y,u)$ be a label of $(X,Z)$ at $x$.
For $1\leq q\leq m$, write $\Retaq=\cO_{Z_q,\eta_q}$ with the maximal ideal
$\fmetaq$, and let $\Jetaq\subset\Retaq$ be the ideal defining $X_q\subset Z_q$
at $\eta_q$. Write
$$
f_i^{(q)}=f_i/u_1^{q n_i},\; y_i^{(q)}=y_i/u_1^{q},\;
u_i^{(q)}=u_i/u_1^{q}\; (2\leq i \leq s),\; u_i'=u_i/u_1\; (s+1\leq i\leq e).
$$

\begin{claim}\label{fu.claim1}
Let $(f,y,u)$ be a label of $(X,Z)$ at $x$ prepared along the faces $E_{L_q}$
for all $q=1,\dots, m-1$. Then, for $q=1,\dots, m$, $\Retaq$ is the localization of
$$
R[y_1^{(q)},\dots,y_r^{(q)},u_2^{(q)},\dots,u_s^{(q)},
u'_{s+1},\dots,u'_e]
\quad \text{ at  }
\langle y_1^{(q)},\dots,y_r^{(q)},u_1,u_2^{(q)},\dots,u_s^{(q)}\rangle,
$$
and $\Jetaq =\langle f_1^{(q)},\dots,f_N^{(q)}\rangle$,
and $(f^{(q)},y^{(q)},(u_1,u_2^{(q)},\dots,u_s^{(q)}))$
is $\delta$-prepared and
$$
\delta(f^{(q)},y^{(q)},(u_1,u_2^{(q)},\dots,u_s^{(q)}))=
\delta_{L_q}(f,y,u)-q.
$$
For $q\leq m-1$, $(f^{(q)},y^{(q)},(u_1,u_2^{(q)},\dots,u_s^{(q)}))$
is a $\delta$-prepared label of $(X_q,Z_q)$ at $\eta_q$.
\end{claim}

\textbf{Proof }
For $q=1$ the claim follows from Lemma \ref{ekeylemma2}.
For $q>1$, by induction it follows from loc.cit. applied to
$$
\Spec(\cO_{X_{q-1},\eta_{q-1}})\gets X_{q}\times_{X_{q-1}}\Spec(\cO_{X_{q-1},\eta_{q-1}})
$$
and $(f^{(q-1)},y^{(q-1)},(u_1,u_2^{(q-1)},\dots,u_s^{(q-1)}))$
in place of $X\gets X_1$ and $(f,y,u)$ (note that the condition
\eqref{ekeylemma2.cond} is satisfied for
$\Spec(\cO_{X_{q-1},\eta_{q-1}})$ by Lemma \ref{ekeylemma2} (5)).
$\square$

\bigbreak

Recalling now that $\eta_q$ is near to $\eta_{q-1}$ for $1\leq q\leq m-1$
($\eta_0=x$ by convention), we get
$\delta_{L_q}(f,y,u)-q \geq 1$ for these $q$ by Theorem \ref{deltanearcriteria} and
Claim \ref{fu.claim1}. It remains to show $\deltaob_{L_q}(y,u)\geq q+1$.
For this we rewrite \eqref{fu.eq9} in $\Retaq$ as follows:
\begin{equation}\label{fu.eq10}
l_B= u_1^q \underset{1\leq i\leq r}{\sum} a_i y^{(q)}_i +
\underset{C=(c,a_2,\dots,a_s) \in \bZ^s_{\geq 0}}{\sum} P_C \cdot
u_1^c (u_2^{(q)})^{a_2} \cdots (u_s^{(q)})^{a_s}\;\; + u_1^{2q} g',
\end{equation}
where $g'=g/u_1^{2q} \in {\langle y^{(q)}_1,\dots,y^{(q)}_r \rangle}^2$
and
$$
P_C= \underset{\Omega(A)=C}{\sum} b_A\cdot
{u'}_{s+1}^{a_{s+1}} \cdots {u'_e}^{a_e},
$$
for the map
$\Omega:\bR^e\to \bR^s\;;\; A=(a_1,\dots,a_e)\to (L_{q-1}(A),a_2,\dots,a_s)$.
We easily see
$$
P_C\in \fmetaq =\langle y^{(q)},u_1,u_2^{(q)},\dots,u_s^{(q)}\rangle
\;\Leftrightarrow\; b_A\in \fm\text{ for all $A$ such that $\Omega(A)=C$}\,,
$$
by noting that $R\to \Retaq/\fmetaq$ factors through $R\to k=R/\fm$ and
$k[u'_{s+1},\dots,u'_e]\hookrightarrow k(\eta_q)$.
The strict transform $\tB_q$ of $B$ in $\Spec(\Retaq)$ is defined by
$$
l'_B=l_B/u_1^{\gamma}\qwith \gamma = v_{\langle u_1\rangle}(l_B),
$$
where $v_{\langle u_1\rangle}$ is the valuation of $\Retaq$ defined by
the ideal $\langle u_1\rangle\subset \Retaq$.
From \eqref{fu.eq10}, we see
$$
\begin{aligned}
l'_B\not\in \fmetaq &\Leftrightarrow
\text{$P_C\not\in \fmetaq$ and $c\leq q$ for some
$C=(c,0,\dots,0)\in \bZ^s_{\geq 0}$}\\
&\Leftrightarrow
\text{$b_A\not\in \fm$ and $|A|\leq q$ for some
$A=(c',0,\dots,0)\in \bZ^s_{\geq 0}$.}\\
\end{aligned}
$$
For $A\in \bZ^e_{\geq 0}-\{0\}$, we have
$$
L_q(A)=|A|+\underset{2\leq j\leq s}{\sum} a_j\geq q+1  \Leftrightarrow
(a_2,\dots,a_s)\not=(0,\dots,0) \text{ or } |A|\geq q+1.
$$
Hence we get
$$
l'_B\not\in \fmetaq \Leftrightarrow
 \text{$L(A)< q+1$ for some $A\in \bZ^s_{\geq 0}-\{0\}$.}
 \Leftrightarrow \delta_{L_q}(l_B,y,u)<q+1.
$$
It implies
$\deltaob_{L_q}(y,u)=\min\{\delta(l_B,y,u)\;|\; B\in \OBx\}\geq q+1$
since $\eta_q\in \tB_q$ for $q\leq m-1$ by the assumption.
This shows the desired assertion and the proof of Theorem \ref{fu.thm}
is complete.

\begin{corollary}\label{cor.fu}
Let $(f,y,u)$ be a $\delta$-prepared label of $(X,Z)$ at $x$.
Assume $e_x(X)=\eob_x(X)$ (for example $\cB=\emptyset$), and that the assumptions of
Theorem \ref{fu.thm0} hold. Then for the length $m$ of the fundamental unit
(Definition \ref{Def.fupb}) we have
$$
m = \llcorner\delta(f,y,u)\lrcorner \quad(\;:= \mbox{ greatest integer }\leq \delta(f,y,u)\;)\,.
$$
\end{corollary}
\textbf{Proof }
By Claim \ref{fu.claim1} in case $s=1$ and $L_q=L_0$
together with Theorem \ref{deltanearcriteria} (1) and by the assumption that $\eta_{m-1}$ is near to $\eta_{m-2}$ and $\eta_m$ is
not near to $\eta_{m-1}$, we have
$$
 \delta(f^{(m-1)},y^{(m-1)},u_1)= \delta(f,y,u)-(m-1) \geq 1, \quad
 \delta(f^{(m)},y^{(m)},u_1)= \delta(f,y,u)-m < 1\,.
$$

\newpage
\section{Additional invariants in the case $e_x(X)=2$}\label{sec:ad.inv}
\bigskip

In order to show key Theorem \ref{Thm2} in \S5, we recall further invariants
for singularities, which were defined by Hironaka (\cite{H6}).
The definition works for any dimension, as long as the directrix is
$2$-dimensional.

\begin{definition}\label{def.invpol}
For a polyhedron $\Delta\subset \bR^2_{\geq 0}$ we define
\begin{align*}
\alpha(\Delta)  :=& \inf\left\{ v_1\mid (v_1,v_2)\in \Delta \right\} \\
\beta(\Delta)  :=& \inf\left\{ v_2 \mid (\alpha(\Delta),v_2)\in \Delta
\right\} \\
\delta(\Delta)  :=& \inf\left\{ v_1+v_2 \mid (v_1,v_2)\in \Delta \right\} \\
\gamma^+(\Delta)  :=& \sup\left\{v_2 \mid
(\delta(\Delta)-v_2,v_2)\in \Delta \right\} \\
\gamma^-(\Delta)  :=& \inf\left\{ v_2 \mid
(\delta(\Delta)-v_2,v_2)\in \Delta \right\} \\
\epsilon(\Delta)  :=& \inf\left\{ v_2\mid (v_1,v_2)\in \Delta \right\} \\
\zeta(\Delta) :=& \inf\left\{ v_1 \mid (v_1,\epsilon(\Delta))\in \Delta\right\}
\end{align*}
\end{definition}

\bigskip
The picture is as follows:

\vspace{0.5cm}
\begin{figure}[h]
\begin{center}
\scalebox{1} {\setlength{\unitlength}{0.00087489in}
\begingroup\makeatletter\ifx\SetFigFontNFSS\undefined%
\gdef\SetFigFontNFSS#1#2#3#4#5{%
  \reset@font\fontsize{#1}{#2pt}%
  \fontfamily{#3}\fontseries{#4}\fontshape{#5}%
  \selectfont}%
\fi\endgroup%
{\renewcommand{\dashlinestretch}{30}
\begin{picture}(5424,3639)(0,-10)
\path(462,12)(462,3612)
\path(492.000,3492.000)(462.000,3612.000)(432.000,3492.000)
\path(12,462)(5412,462)
\path(5292.000,432.000)(5412.000,462.000)(5292.000,492.000)
\dashline{60.000}(462,2262)(912,2262)(912,462)
\dashline{60.000}(462,912)(3162,912)(3162,462)
\dashline{60.000}(462,1362)(1812,1362)(1812,462)
\dashline{60.000}(1137,1812)(462,1812)
\path(1542,1587)(1812,1677)
\path(3747,912)(3882,1137)
\dashline{60.000}(462,2487)(2892,462)
\thicklines
\path(912,2712)(912,2262)(1272,1812)
	(1812,1362)(3162,912)(4512,912)
\put(237,867){\makebox(0,0)[lb]{\smash{{\SetFigFontNFSS{12}{14.4}{\rmdefault}{\mddefault}{\updefault}$\epsilon$}}}}
\put(1000,2262){\makebox(0,0)[lb]{\smash{{\SetFigFontNFSS{12}{14.4}{\rmdefault}{\mddefault}{\updefault}$\mathbf{v}$}}}}
\put(1250,1857){\makebox(0,0)[lb]{\smash{{\SetFigFontNFSS{12}{14.4}{\rmdefault}{\mddefault}{\updefault}$\mathbf{w^+}$}}}}
\put(1860,1700){\makebox(0,0)[lb]{\smash{{\SetFigFontNFSS{12}{14.4}{\rmdefault}{\mddefault}{\updefault}slope $-1$}}}}
\put(1800,1432){\makebox(0,0)[lb]{\smash{{\SetFigFontNFSS{12}{14.4}{\rmdefault}{\mddefault}{\updefault}$\mathbf{w^-}$}}}}
\put(3837,1227){\makebox(0,0)[lb]{\smash{{\SetFigFontNFSS{12}{14.4}{\rmdefault}{\mddefault}{\updefault}$\partial\Delta$}}}}
\put(5097,147){\makebox(0,0)[lb]{\smash{{\SetFigFontNFSS{12}{14.4}{\rmdefault}{\mddefault}{\updefault}$a_1$}}}}
\put(102,3432){\makebox(0,0)[lb]{\smash{{\SetFigFontNFSS{12}{14.4}{\rmdefault}{\mddefault}{\updefault}$a_2$}}}}
\put(237,2487){\makebox(0,0)[lb]{\smash{{\SetFigFontNFSS{12}{14.4}{\rmdefault}{\mddefault}{\updefault}$\delta$}}}}
\put(237,2172){\makebox(0,0)[lb]{\smash{{\SetFigFontNFSS{12}{14.4}{\rmdefault}{\mddefault}{\updefault}$\beta$}}}}
\put(192,1317){\makebox(0,0)[lb]{\smash{{\SetFigFontNFSS{12}{14.4}{\rmdefault}{\mddefault}{\updefault}$\gamma^-$}}}}
\put(192,1812){\makebox(0,0)[lb]{\smash{{\SetFigFontNFSS{12}{14.4}{\rmdefault}{\mddefault}{\updefault}$\gamma^+$}}}}
\put(867,192){\makebox(0,0)[lb]{\smash{{\SetFigFontNFSS{12}{14.4}{\rmdefault}{\mddefault}{\updefault}$\alpha$}}}}
\put(1585,192){\makebox(0,0)[lb]{\smash{{\SetFigFontNFSS{12}{14.4}{\rmdefault}{\mddefault}{\updefault}$\delta-\gamma^-$}}}}
\put(2847,192){\makebox(0,0)[lb]{\smash{{\SetFigFontNFSS{12}{14.4}{\rmdefault}{\mddefault}{\updefault}$\delta$}}}}
\put(3117,192){\makebox(0,0)[lb]{\smash{{\SetFigFontNFSS{12}{14.4}{\rmdefault}{\mddefault}{\updefault}$\zeta$}}}}
\put(2667,2350){\makebox(0,0)[lb]{\smash{{\SetFigFontNFSS{12}{14.4}{\rmdefault}{\mddefault}{\updefault}$\Delta=\Delta (f,y,u)$}}}}
\end{picture}
}}
\end{center}
\end{figure}

\medskip
There are three vertices of $\Delta(f,y,u)$ which play crucial roles:
\begin{align*}
\bv:=\bv(\Delta) :=& (\alpha(\Delta),\beta(\Delta)),\\
\bw^+:=\bw^+(\Delta) :=&(\delta(\Delta)-\gamma^+(\Delta) , \gamma^+(\Delta)),\\
\bw^-:=\bw^-(\Delta) :=&(\delta(\Delta)-\gamma^-(\Delta) , \gamma^-(\Delta)).
\end{align*}
We have
\begin{equation}\label{bcineq}
\beta(\Delta) \; \geq \; \gamma^+(\Delta)\; \geq \; \gamma^-(\Delta)\, .
\end{equation}

\bigskip

Now let us consider the situation of Setup C in \S\ref{sec:fund.seq}
and assume $e_x(X)=2$.

Let $(f,y,u)=(f_1,\dots,f_N,y_1,\dots,y_r,u_1,u_2)$ be a prelabel of
$(X,Z)$ at $x$. Recall $n_{(u)}(f_i)=n_i$ for $i=1,\dots, N$.
Write as \eqref{eq1.2.1}:
$$
f_i=\sum\limits_{A,B}C_{i,A,B}\; y^B u^A\;,\qwith
A=(a_1,a_2), \; B=(b_1,\ldots , b_r),\; C_{i,A,B}\in R^\times\cup \{0\}.
$$
\medbreak

For $*=\alpha,\;\beta,\; \delta,\; \gamma^{\pm},\; \bv,\; \bw^{\pm}$,
we write $*(f,y,u)$ for $*(\Delta(f,y,u))$. Then we see
$$
\begin{aligned}
\alpha(f,y,u) =&
\inf\left\{\frac{a_1}{n_i-|B|}\mid 1\leq i\leq m, \; C_{i,A,B}\neq 0\right\}\\
\beta(f,y,u) =&
\inf\left\{\frac{a_2}{n_i-|B|}\mid 1\leq i\leq m,\;
\frac{a_1}{n_i-|B|}=\alpha(f,y,u), \; C_{i,A,B}\neq 0\right\} \\
\delta(f,y,u) =&
\inf\left\{\frac{|A|}{n_i-|B|}\mid 1\leq i\leq m,\; C_{i,A,B}\neq 0\right\} \\
\gamma^+(f,y,u) =&
\sup\left\{\frac{a_2}{n_i-|B|}\mid 1\leq i\leq m,\;
\frac{|A|}{n_i-|B|}=\delta,\; C_{i,A,B}\neq 0\right\} \\
\gamma^-(f,y,u) =&
\inf\left\{\frac{a_2}{n_i-|B|}\mid 1\leq i\leq m,\;
\frac{|A|}{n_i-|B|}=\delta, \; C_{i,A,B}\neq 0\right\}\\
\epsilon(f,y,u) =&
\inf\left\{\frac{a_2}{n_i-|B|}\mid 1\leq i\leq m, \; C_{i,A,B}\neq 0\right\}\\
\end{aligned}
$$

\begin{equation*}\label{bvbw}
\begin{aligned}
\bv:=\bv(f,y,u) =& (\alpha(f,y,u),\beta(f,y,u)),\\
\bw^+:=\bw^+(f,y,u) = & (\delta(f,y,u)-\gamma^+(f,y,u) , \gamma^+(f,y,u)),\\
\bw^-:=\bw^-(f,y,u) = & (\delta(f,y,u)-\gamma^-(f,y,u),\gamma^-(f,y,u))\, .
\end{aligned}
\end{equation*}

\begin{definition}\label{def4.1}
\begin{itemize}
\item[(1)]
A prelabel $(f,y,u)$ (cf. Definition \ref{def.label}) is $\bv$-prepared
if $(f,y,u)$ is prepared at $\bv(f,y,u)$.
\item[(2)]
We say that $(X,Z)$ is $\bv$-admissible at $x$ if there exists a $\bv$-prepared
prelabel $(f,y,u)$ of $(X,Z)$ at $x$. By Theorem \ref{thm.wellprepared},
$(X,Z)$ is $\bv$-admissible at $x$ if $R=\cO_{Z,x}$ is complete.
\end{itemize}
\end{definition}

We now extend the above definition to the situation where
the old components of $\cB$ at $x$ are taken into account.
We assume
\begin{equation}\label{logu.eq}
\Spec(R/\langle u\rangle )\not\subset B\quad\text{ for any } B\in \OB\cBx.
\end{equation}

\begin{definition}\label{loginv}
Let $(y,u)$ be a system of regular parameters of $R$ such that
$(u)$ is admissible for $J$.
\begin{itemize}
\item[(1)]
For each $B\in \OBx$, choose $l_B\in R$ such that
$B=\Spec(R/\langle l_B\rangle)\subset Z= \Spec(R)$. \eqref{logu.eq} implies
$l_B\not\in \langle u\rangle$ so that $\Delta(l_B,y,u)$ is well-defined
(cf. Definition \ref{def.intial}). We define
$$
\Deltaob(y,u)=\text{ the minimal $F$-subset containg }
\underset{B\in \OBx}{\cup} \Delta(l_B,y,u).
$$
It is easy to see that $\Deltaob(y,u)$ is independent of the choice of $l_B$.
For a prelabel $(f,y,u)$ of $(X,Z)$ at $x$, let
$$
\Deltaob(f,y,u)=\text{ the minimal $F$-subset containg }
\Delta(f,y,u) \cup \Deltaob(y,u),
$$
$$
*^{\OB}(f,y,u)=*(\Deltaob(f,y,u))
\qfor *=\bv,\; \alpha,\; \beta,\; \gamma^{\pm},\dots.
$$
Note that
\begin{equation}\label{loginv.eq}
\Deltaob(f,y,u)=\Delta(\fob,y,u)\quad\text{ where  }
\fob=(f,\; l_B\; (B\in \OBx)).
\end{equation}
\item[(2)]
Assume that $(X,Z)$ is $\bv$-admissible at $x$. Then we define
$$
\betaob_x(X,Z)=\betaob(J) :=\underset{(f,y,u)}{\min}\{\betaob(f,y,u)|\;
\text{$(f,y,u)$ is $\bv$-prepared}\}.
$$
\item[(3)]
A prelabel $(f,y,u)$ of $(X,Z)$ at $x$ is called $\OB$-admissible if
$(f,y,u)$ is $\bv$-prepared and $\betaob_x(X,Z)=\betaob(f,y,u)$.
Such a prelabel exists if and only if $(X,Z)$ is $\bv$-admissible.
\end{itemize}
\end{definition}

By Lemma \ref{Vlattice} and Lemma \ref{def1.4cor} (2),
for $\bv$-prepared $(f,y,u)$ we have
\begin{equation}\label{Vlattice2}
\betaob(f,y,u)\in \frac{1}{n_N !} \bZ_{\geq 0} \subseteq \bR\;,
\end{equation}

\begin{lemma}\label{loginv.lem}
\begin{itemize}
\item[(1)]
Let $(f,y,u)$ be a $\bv$-prepared prelable of $(X,Z)$ at $x$.
Then, for a preparation $(f,y,u)\to (g,z,u)$ at a vertex $v\in \Delta(f,y,u)$,
we have $\betaob(f,y,u)=\betaob(g,z,u)$.
\item[(2)]
Assume that $(X,Z)$ is $\bv$-admissible at $x$.
For any integer $m\geq 1$, there exists an $\delta$-prepared and
$\OB$-admissible label $(f,y,u)$ of $(X,Z)$ at $x$.
Moreover, if $R$ is complete, one can make $(f,y,u)$ totally prepared.
\end{itemize}
\end{lemma}
\textbf{Proof  }
(2) is a consequence of (1) in view of Corollary \ref{cor.wellprepared2}.
We prove (1). Setting
$$
\bvob(y,u)=\bv(\Deltaob(y,u))=(\alphaob(y,u),\betaob(y,u)),
$$
we have
\begin{equation}\label{loginv.lem.eq}
\vob(f,y,u)=(\alphaob(f,y,u),\betaob(f,y,u))=
\left.\left\{\gathered
 \bvob(y,u) \\
 \bv(f,y,u) \\
\endgathered\right.\quad
\begin{aligned}
&\text{if $\alphaob(y,u) \leq \alpha(f,y,u)$}\\
&\text{if $\alphaob(y,u) > \alpha(f,y,u)$}.
\end{aligned}\right.
\end{equation}
By the $\bv$-preparedness of $(f,y,u)$, any vertex $v\in \Delta(f,y,u)$
which is not prepared lies in the range
$\{(a_1,a_2)\in \bR^2|\; a_1>\alpha(f,y,u)\}$.
Theorem \ref{thm.normalized} implies $\vob(f,y,u)=\vob(g,y,u)$ for
the normalization $(f,y,u)\to (g,y,u)$ at such $v$. Thus it suffices to show
$\vob(f,y,u)=\vob(f,z,u)$ for the dissolution $(f,y,u)\to (f,z,u)$ at $v$.
Write $v=(a,b)$. By the above remark, we have $a>\alpha(f,y,u)$.
The dissolution is given by a coordinate transformation:
$$
y=(y_1,\dots,,y_r)\to z=(z_1,\dots,z_r)\qwith z_i=y_i+ \lambda_i u_1^a u_2^b
\quad(\lambda_i\in R)
$$
Write $\aol=\alphaob(y,u)$ for simplicity.
For each $B\in \OBx$ choose $l_B\in R$ such that
$B=\Spec(R/\langle l_B\rangle)$. We may write
$$
l_B=\Lambda_B(y) + u_1^{\aol}\phi_B\qwith
\Lambda_B(y)=\underset{1\leq i\leq r}{\sum}c_{B,i}\cdot y_i\quad
(c_{B,i},\; \phi_B\in R)
$$
and $\phi_B\not=0$ for some $B\in \OBx$. Then we get
\begin{equation}\label{loginv.lem.eq2}
l_B=\Lambda_B(z) + u_1^{\aol}\phi_B + u_1^{a}\psi_B \quad
\text{ for some $\psi_B\in R$}
\end{equation}
In case $\aol=\alphaob(y,u)\leq \alpha(f,y,u)$, \eqref{loginv.lem.eq} implies
$$\bvob(f,y,u)=\bvob(y,u)=\bvob(z,u)=\bvob(f,z,u),$$
where the second equality follows from \eqref{loginv.lem.eq2} because
$a>\alpha(f,y,u)\geq \aol$, and the third follows from $\bv(f,y,u)=\bv(f,z,u)$
by the $\bv$-preparedness of $(f,y,u)$.
In case $\aol=\alphaob(y,u)> \alpha(f,y,u)$, \eqref{loginv.lem.eq} implies
$$
\bvob(f,y,u)=\bv(f,y,u)=\bv(f,z,u)=\bvob(f,z,u),
$$
where the second equality follows from the $\bv$-preparedness of $(f,y,u)$.
The third equality holds since by \eqref{loginv.lem.eq2}, we have
$\alphaob(z,u)\geq a>\alpha(f,y,u)$ if $\psi_B\not=0$ for some $B\in \OBx$, and
$\alphaob(z,u)=\aol>\alpha(f,y,u)$ if $\psi_B=0$ for all $B\in \OBx$.
This completes the proof of Lemma \ref{loginv.lem}.
$\square$
\medbreak

\begin{lemma}\label{case1.lem3}
Let $(f,y,u)$ be a prelabel of $(X,Z)$ at $x$.
Assume that there is no regular closed subscheme
$D \subseteq \{\xi\in X|\; \tHSf{X}(\xi)\geq \tHSf{X}(x)\}$ of
dimension 1 with $x\in D$. (In particular this holds if $x$ is isolated in
$\{\xi\in X|\; \tHSf{X}(\xi)\geq \tHSf{X}(x)\}$.) Then
$\alphaob(f,y,u)<1$ and $\epsilonob(f,y,u)<1$.
\end{lemma}
\medbreak
\textbf{Proof  }
By corollary \ref{cor.wellprepared2}, we prepare $(f,y,u)$ at the vertices
in $\{A\in \bR^2|\; |A|\leq 1\}$ to get a label $(g,z,u)$ of $(X,Z)$ at $x$.
Then $\alphaob(g,z,u)\geq \alphaob(f,y,u)$ since
$\Deltaob(g,z,u)\subset \Deltaob(f,y,u)$. Thus we may replace $(f,y,u)$
with $(g,z,u)$ to assume that $f$ is a standard base of $J$.
\medbreak
Assume $\alphaob(f,y,u)\geq 1$. Then, letting
$\fp=(y_1,\dots,y_r,u_1)\subset R$, we have $v_\fp(f_j)\geq  n_j$ for
$j=1,\dots,N$ (cf. \eqref{compl.eq0}).
Since $n_j=n_{(u)}(f_j)\geq v_\fm(f_j)\geq v_\fp(f_j)$,
this implies $v_\fp(f_j)=v_\fm(f_j)$ for $j=1,\dots,N$.
This implies by Theorems \ref{nf.thm1} $(iv)$ and \ref{HSf.perm} that
$\eta\in X$ and $\HSfX \eta= \HSfX x$, where $\eta\in Z$ is the point
corresponding to $\fp$. By the same argument we prove
$v_\fp(l_B)= 1= v_{\fm}(l_B)$ for $B\in \OBx$ so that
$\tHSfX \eta= \tHSfX x$. Thus $\overline{\{\eta\}} = \Spec(R/\fp)$ is
$\OB$-permissible, which contradicts the assumption of the lemma.
The assertion $\epsilonob(f,y,u)<1$ is shown in the same way.
$\square$
\medskip

\begin{lemma}\label{case1.lem4}
Let $(y,u)$ be a system of regular parameters of $R$ which is strictly admissible for
$J$. Assume
$$
\eob_x(X)=\dim_k\big(\Dir_x(X)\cap \underset{B\in \OBx}{\cap}T_x(B)\big)
=2.\leqno(*)
$$
Then we have $\deltaob(y,u)>1$. Assume in addition that $\delta(f,y,u)>1$
(so that $(f,y,u)$ is a label). Then we have $\deltaob(f,y,u)>1$.
\end{lemma}
\medbreak
\textbf{Proof  }
By the assumption on $(y,u)$, we have
$I\Dir_x(X)=\langle \inm(y_1),\dots,\inm(y_r)\rangle$.
Hence $(*)$ implies
$$
l_B\in \langle y_1,\dots, y_r\rangle +\fm^2\qfor B\in \OBx.
$$
We then easily deduce the first assertion of the lemma.
The second assertion is an obvious consequence of the first.
$\square$

\newpage
\section{Proof in the case $e_x(X)=\eb_x(X)=2$, I: some key lemmas}\label{sec:e2I}
\def\OBd{O'}
\bigskip

In this section we prepare some key lemmas for the proof of Theorem \ref{Thm2}.
\medbreak

Let the assumption be those of Setup C in \S\ref{sec:fund.seq}.
We assume
\begin{itemize}
\item
$\Char(k(x)) = 0\text{  or  } \Char(k(x)) \geq \dim (X)/2+1.$
\item
$\eob_x(X)=e_x(X)=\eb_x(X)=2.$
\end{itemize}

We fix a label $(f,y,u)$ of $(X,Z)$ at $x$ and adopt the notations
of Definition \ref{def.label} (1) and (2). We recall
$$
F_i(Y)=\inm(f_i)\in k[Y]=k[Y_1,\dots,Y_r] \subset \grmR .
\quad(Y_j=in_{\fm}(y_j))
$$
By Lemma \ref{case1.lem4}, the assumption $\eob_x(X)=e_x(X)=2$ implies:
\begin{equation}\label{eq.eOB2}
\deltaob(f,y,u)>1.
\end{equation}
It also implies that \eqref{logu.eq} is always satisfied
so that $\Deltaob(f,y,u)$ is well-defined.

\medbreak
For each $B\in \OBx$, we choose $l_B\in R$ such that $B\times_Z\Spec(R) =\Spec(R/\langle l_B \rangle)$.
We study two cases.

\medskip\medskip\noindent
\textbf{Case 1 (point blowup):}
Consider
$$
\pi: Z'=\Bl_x(Z)\to Z \qaq \pi: X'=\Bl_x(X)\to X.
$$
Note that $x\hookrightarrow X$ is $\cB$-permissible for trivial reasons.
Let $(\cB',\OBd)$ be the complete transform of $(\cB,\OB)$ in $Z'$
(cf. Definition \ref{def.compl.transf.cBO}).
In this case we have
$\cB'=\tilde{\cB}\cup \{\pi^{-1}(x)\}$, where $\tilde{\cB}=\{\tB|\; B\in \cB\}$
with $\tB$, the strict transform of $B$ in $Z'$.
\medbreak

By Theorem \ref{thmdirectrix}, any point of $X'$ near to $x$ is contained in
\begin{equation*}
C' :=\bP(\Dir_x(X)) \subset \bP(T_x(Z))=\pi_{Z}^{-1}(x) \subset Z'\, .
\end{equation*}
Let $T_1,T_2$ be a pair of new variables over $k$ and let
\begin{equation}\label{e2c.eq1.5}
\psi_{(u)} : C' = \bP(\Dir(R/J)) \isom \Proj(k[T_1, T_2]) =\bP^1_k
\end{equation}
be the isomorphism \eqref{eq.u-coord} which is determined by $(u)$.
\medbreak

Take a closed point $x'\in C'$ near to $x$.
By Theorems \ref{thm.pbu.inv} and \ref{def.nearpoint}, we have
$$
\nu^*_{x'}(X',Z') = \nu^*_x(X,Z)\qaq
e_{x'}(X')\leq 2.
$$
Put $R'=\cO_{Z',x'}$ with the maximal ideal $\fm'$.
Let $J'$ and $\fp'$ be the ideals of $R'$ such that
\begin{equation}\label{eq.Jdpd}
X'\times_{Z'}\Spec(R')=\Spec(R'/J') \qaq
C' \times_{Z'} \Spec(R')=\Spec(R'/\fp').
\end{equation}

\medbreak

\begin{lemma}\label{e2keylemma0}
Assume $\psi_{(u)}(x')=(1:0)\in \Proj(k[T_1,T_2])$. Let
$$
y'=(y'_1,\dots,y'_r)\text{ with }
y'_i=y_i/u_1, \;  u'_2=u_2/u_1,\;
f'=(f'_1,\dots,f'_N) \text{ with }f_i'=f_i/u_1^{n_i}.
$$
\begin{itemize}
\item[(1)]
$(y', u_1,u'_2)$ is a system of regular parameters of $R'$ such that
$\fp'=(y',u_1)$, and $J'=\langle f'_1,\dots,f'_N\rangle$.
\item[(2)]
If $x'$ is very near to $x$, then $(f',y',(u_1,u_2'))$ is
a prelabel of $(X',Z')$ at $x'$.
\item[(3)]
$\Delta(f', y', (u_1, u_2'))$ is the minimal $F$-subset containg
$\Psi(\Delta(f,y,u))$, where
\par\noindent
$$\Psi : \bR^2\rightarrow\bR^2, (a_1,a_2)\mapsto (a_1+a_2-1,a_2)$$

\vspace{0,5cm}
\begin{center}
\myfig{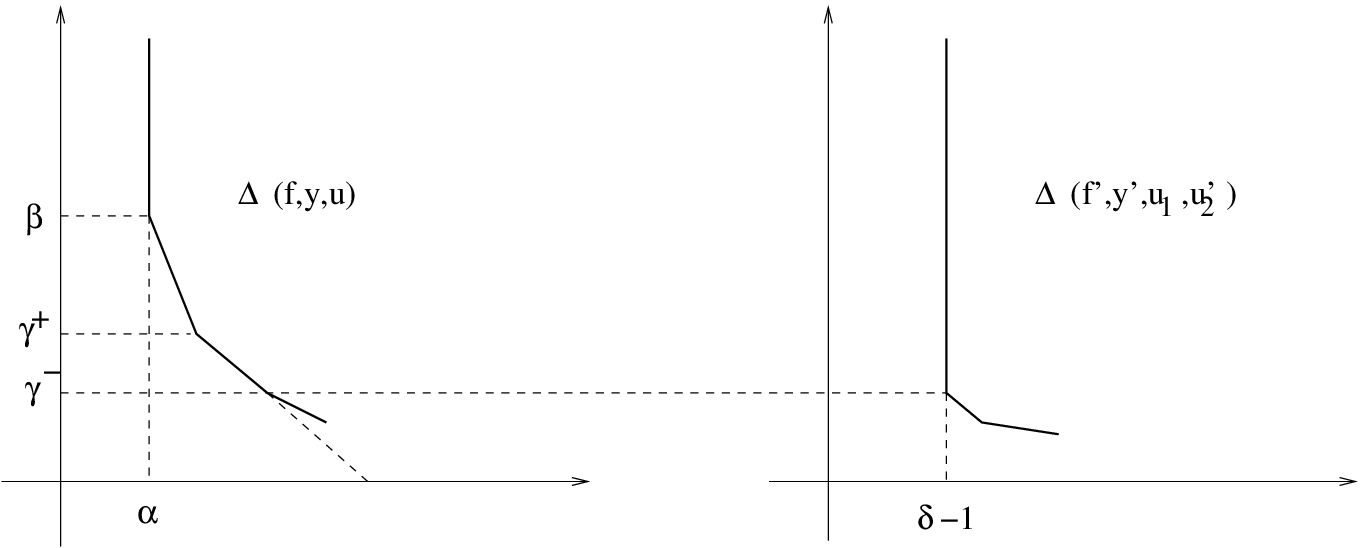}{0,9}
\end{center}
for which all vertices move horizontally. We have
$$
\beta(f',y',(u_1,u_2'))  =  \gamma^-(f,y,u)\leq \beta(f,y,u).
$$
$$
\alpha(f',y', (u_1,u_2'))  =  \delta(f,y,u)-1.
$$
\item[(4)]
If $(f,y,u)$ is prepared at $\bwm(f,y,u)$, then
$(f',y',(u_1,u_2'))$ is $\bv$-prepared.
If $(f,y,u)$ is prepared along the face $E_L$, then
$(f',y',(u_1,u_2'))$ is $\delta$-prepared, where
$$
L:\bR^2\to \bR\;;\; (a_1,a_2)\to a_1+2a_2.
$$
If $(f,y,u)$ is totally prepared, so is $(f',y',(u_1,u_2'))$.
\item[(5)]
Assume that $x'$ is $\cB$-near to $x$.
Putting $l'_B=l_B/u_1\in R'$ for $B\in \OBx$, we have
$$
\Deltaob(f',y',(u_1,u'_2))=\Delta(\fdob,y,(u_1,u_2'))\qwith
\fdob=(f',\; l'_B\; (B\in \OBx)).
$$
The same assertions as (3) hold
replacing $\Delta$ by $\Deltaob$ and
$*$ by $*^{OB}$ for $*=\alpha,\; \beta,\; \gamma^-,\; \delta$.
\end{itemize}
\end{lemma}
\medbreak
\textbf{Proof  }
By Definition \ref{def.label}, $(y,u)$ is strictly admissible for $J$ and
$f$ is a standard base of $J$ which is admissible for $(y,u)$. Hence
(1) has been seen in Setup B in \S\ref{sec:pr.usb}. (2) follows from Theorem
\ref{prop2.1} and Theorem \ref{thm3.1}. (3) and (4) follow from Lemma
\ref{ekeylemma0}. As for (5), the assumption implies
$\OBdxd=\{B'\;|\; B\in \OBx\}$ with $B'$, the strict transform of $B$ in $Z'$
and we have
$$
B'\times_{Z'}\Spec(R')=\Spec(R'/\langle l'_B \rangle) \subset \Spec(R').
$$
This implies the first assertion of (5) by \eqref{loginv.eq}.
The second assertion of (5) then follows from the first.
$\square$
\medbreak

The following lemma is shown in the same way as the previous lemma
except that the last assertion of (5) follows from Lemma \ref{case1.lem3}.

\begin{lemma}\label{e2keylemma0.1}
Assume $\psi_{(u)}(x')=(0:1)\in \Proj(k[T_1,T_2])$ and put
$$
z'=(z'_1,\dots,z'_r)\text{ with }
z'_i=y_i/u_2, \;  u'_1=u_1/u_2,\;
g'=(g'_1,\dots,g'_N) \text{ with }f_i'=f_i/u_2^{n_i}.
$$
\begin{itemize}
\item[(1)]
$(z', u'_1,u_2)$ is a system of regular parameters of $R'$ such that
$\fp=(z',u_2)$, and $J'=\langle g'_1,\dots, g'_N\rangle$.
\item[(2)]
If $x'$ is very near to $x$, then $(g',z',(u'_1,u_2))$ is
a prelabel of $(X',Z')$ at $x'$.
\item[(3)]
$\Delta(g', z', (u'_1, u_2))$ is the minimal $F$-subset containg
$\Psi(\Delta(f,y,u))$, where
\par\noindent
$$\Psi : \bR^2\rightarrow\bR^2, (a_1,a_2)\mapsto (a_1,a_1+a_2-1)$$
We have
$$
\alpha(g',z',(u'_1,u_2))  =  \alpha(f,y,u).
$$
$$\beta(g', z',(u'_1, u_2)) \leq \beta(f,y,u)+\alpha(f,y,u)-1.$$
\item[(4)]
If $(f,y,u)$ is prepared at $\bwp(f,y,u)$, then
$(g',z',(u'_1,u_2))$ is $\bv$-prepared.
If $(f,y,u)$ is prepared along the face $E_{L'}$, then
$(f',y',(u_1,u_2'))$ is $\delta$-prepared, where
$$
L':\bR^2\to \bR\;;\; (a_1,a_2)\to 2 a_1+a_2.
$$
If $(f,y,u)$ is totally prepared, so is $(g',z',(u'_1,u_2))$.
\item[(5)]
Assume that $x'$ is $\cB$-near to $x$. Then
$$\betaob(g', z',(u'_1, u_2)) \leq \betaob(f,y,u)+\alphaob(f,y,u)-1.$$
If there is no regular closed subscheme
$D \subseteq \{\xi\in X|\; \tHSf{X}(\xi)\geq \tHSf{X}(x)\}$ of
dimension 1 with $x\in D$, then $\alphaob(f,y,u)<1$ so that
$$\betaob(g', z',(u_1, u'_2)) < \betaob(f,y,u).$$
\end{itemize}
\end{lemma}
\medskip

Now let $\eta'$ be the generic point of $C'$.
By Theorem \ref{def.nearpoint} and Theorem \ref{nf.dir}, if $\eta'$ is near
to $x$, we have $e_{\eta'}(X')\leq 1$.
Write $\Rdetad=\cO_{X',\eta'}$ with the maximal ideal $\fmdetad$.
Note $(y',u_1)=(y'_1,\dots,y'_r,u_1)$ is a system of regular parameters of
$\Rdetad$.

\begin{lemma}\label{e2keylemma2}
Let $\Delta(f', y', u_1)$ be the characteristic polyhedron for
$(\Rdetad, J\Rdetad)$.
\begin{itemize}
\item[(1)]
$\Delta(f', y', u_1)=[\delta(f,y,u)-1,\infty)\subset \bR_{\geq 0}$.
\end{itemize}
Assume $(f,y,u)$ is $\delta$-prepared.
\begin{itemize}
\item[(2)]
$\Delta(f', y', u_1)$ is well prepared.
\item[(3)]
$\eta'$ is near to $x$ if and only if $\delta(f,y,u)\geq 2$.
If $\eta'$ is near to $x$, we have
$$
v_{\fmdetad}(f'_i)=v_{\fm}(f_i)=n_i \qfor i=1,\dots, N.
$$
\item[(4)]
$\eta'$ is very near to $x$ if and only if $\delta(f,y,u)> 2$.
If $\eta'$ is very near to $x$, then $(f',(y',u_1))$ is a well-prepared label of
$(X',Z')$ at $\eta'$.
\end{itemize}
\end{lemma}
\medbreak
\textbf{Proof  }
The lemma is a special case of Lemma \ref{ekeylemma2}.
$\square$

\medskip\medskip\noindent
\textbf{Case 2 (curve blowup):}
Let the assumption be as in the beginning of this section.
Assume given a regular curve $C\subset X$ containing $x$ which is
$\cB$-permissible (cf. Definition \ref{def.Bperm}) and such that
$$
C\times_Z \Spec(R)= \Spec(R/\fp)\qwith \fp=(y,u_1)=(y_1,\dots,y_r, u_1).
$$
By Theorem \ref{nf.dir}, the assumption $e_x(X)=2$ implies
$$
e_{\eta}(X)\leq 1\quad\text{ where $\eta$ is the generic point of $C$.}
$$
Consider
$$
\pi: Z'=\Bl_C(Z)\to Z \qaq \pi: X'=\Bl_C(X)\to X.
$$
Let $(\cB',\OBd)$ be the complete transform of $(\cB,OI)$ in $Z'$
(cf. Definition \ref{def.compl.transf.cBO}). In this case we have
$\cB'=\tilde{\cB}\cup \{\pi^{-1}(C)\}$, where $\tilde{\cB}=\{\tB|\; B\in \cB\}$
with $\tB$, the strict transform of $B$ in $Z'$.
By Theorem \ref{thmdirectrix}, there is the unique point $x'\in \pi^{-1}(x)$
possibly near to $x$, given by
\begin{equation*}
x':=\bP(\Dir_x(X)/T_x(C)) \subset \Proj(T_x(Z)/T_x(C))) =
\pi_{Z}^{-1}(x) \subset Z'\, .
\end{equation*}
In what follows we assume $x'$ near to $x$.
By Theorems \ref{thm.pbu.inv} and \ref{def.nearpoint}, we have
$$
\nu^*_{x'}(X',Z') = \nu^*_x(X,Z)
\qaq e_{x'}(X')\leq 2.
$$
Let $R'=\cO_{Z',x'}$ with the maximal ideal $\fm'$, and let $J' \subset R'$
be the ideal such that
$$
X'\times_{Z'} \Spec(R')=\Spec(R'/J').
$$
As is seen in Setup B in \S\ref{sec:pr.usb},
$$
(y',u)=(y'_1,\dots,y'_r, u_1,u_2)=(y_1/u_1,\dots,y_r/u_1, u_1,u_2)
$$
is a system of regular parameters of $R'$.

\begin{lemma}\label{e2keylemma0b}
\begin{itemize}
\item[(1)]
If $(f,y,u)$ is $\bv$-prepared, then
$v_{\fp}(f_i)=v_{\fm}(f_i)=n_{(u)}(f_i)$ and
$$
J'=\langle f'_1,\dots,f'_N\rangle \qwith f'_i:=f_i/u_1^{n_i}\in R'.
$$
\item[(2)]
If $x'$ is very near to $x$, then $(f',y',u)$ is a prelabel of $(X',Z')$ at $x'$.
\item[(3)]
If $(f,y,u)$ is $\bv$-prepared (resp. $\delta$-prepared, resp.
totally prepared), so is $(f',y',u)$.
\item[(4)]
We have
$$
\Delta(f', y', u)=\Psi(\Delta(f,y,u))
\qwith
\Psi : \bR^2\rightarrow\bR^2, (a_1,a_2)\mapsto (a_1-1,a_2),
$$
$$
\beta(f',y',u)  =  \beta(f,y,u)\qaq \alpha(f',y', u)  =  \alpha(f,y,u)-1.
$$
If $x'$ is $\cB$-near to $x$, the same assertions hold replacing $\Delta$
by $\Deltaob$ and $*$ by $*^{OB}$ for $*=\alpha,\; \beta$.
\end{itemize}
\end{lemma}
\medbreak
\textbf{Proof  }
The first assertion of (1) follows from Theorem \ref{C.thm} and the second
from the first (cf. Setup B in \S\ref{sec:pr.usb}). (2) follows from Theorems
\ref{prop2.1} and \ref{thm3.1}. To show (3) and (4), write, as in \eqref{eq1.2.1}:
\begin{equation}\label{eqf.Thm2-1}
f_i = \sum\limits_{(A, B)} C_{i, A, B} \; y^B u^A \quad \mbox{with}
\quad C_{i, A, B} \in R^{\times} \cup \{ 0 \}
\end{equation}
We compute
\begin{equation}\label{case1.eq0b}
f'_i=f_i/u_1^{n_i} = \sum\limits_{(A, B)} C_{i, A, B} \;
(y')^B {u_2}^{a_2} u_1^{a_1 + (|B|-n_i)} \qwith A=(a_1,a_2).
\end{equation}
This immediately implies the first assertion of (4).
Then (3) is shown in the same way as Lemma \ref{ekeylemma0} (2).
Finally the last assertion of (4) is shown in the same way as
Lemma \ref{e2keylemma0} (5).
This completes the proof of Lemma \ref{e2keylemma0b}.
$\square$
\medbreak

\begin{lemma}\label{e2keylemma2b}
Assume that $(f,(y,u_1))$ is a well-prepared label of $(X,Z)$ at $\eta$
(note that this implies $e_{\eta}(X)=1$).
Let $C'=\Spec(R'/\fp')$ with $\fp'=(y',u_1)\subset R'$ and
let $\eta'$ be the generic point of $C'$.
Then $C'\subset X'$ and $\eta'$ is the unique point of $X'$
possibly near to $\eta$. If $\eta'$ is very near to $\eta$, then
$(f',(y',u_1))$ is a well-prepared label of $(X',Z')$ at $\eta'$.
\end{lemma}
\medbreak
\textbf{Proof  }
The first assertion is a direct consequence of Theorem \ref{thmdirectrix}.
The second assertion follows from Lemma \ref{ekeylemma0} applied to
the base change via $\eta\to C$ of the diagram
$$
\begin{matrix}
C' &\hookrightarrow& X' &\hookrightarrow& Z' \\
\downarrow && \downarrow\rlap{$\pi$} && \downarrow\rlap{$\pi$} \\
C &\hookrightarrow& X &\hookrightarrow& Z \\
\end{matrix}
$$
$\square$

\newpage
\section{Proof in the case $e_x(X)=\eb_x(X)=2$, II: separable residue extensions}\label{sec:e2II}
\bigskip

In this section we prove Theorem \ref{Thm2-1.Step4.thm} below,
which implies Key Theorem \ref{Thm2} under the assumption that
the residue fields of the initial points of $\cX_n$ are separably algebraic
over that of $\cX_1$. The proof is divided into two steps.

\bigskip\noindent
\textbf{Step 1 (one fundamental unit):}
Let the assumptions and notations be as in the beginning of
the previous section. Assume given a fundamental unit of $\cB$-permissible
blowups as in Definition \ref{Def.fupb}:
\begin{equation}\label{fupb.eq2-1}
\begin{array}{ccccccccccccccc}
 \cB=&\cB_0&& \cB_1&& \cB_2&& \cB_{m-1}&& \cB_{m} \\
 \\
Z= & Z_0 & \lmapo{\pi_1} &   Z_1 & \lmapo{\pi_2} & Z_2 &
\leftarrow\ldots \leftarrow & Z_{m-1} & \lmapo{\pi_m} & Z_m\\
& \cup & & \cup && \cup & & \cup && \cup \\
X= & X_0 & \lmapo{\pi_1} &   X_1 & \lmapo{\pi_2} & X_2 &
\leftarrow\ldots \leftarrow & X_{m-1} & \lmapo{\pi_m} & X_m\\
& \uparrow & & \cup && \cup & & \cup && \uparrow  \\
x= & x_0 & \leftarrow & C_1 & \stackrel{\sim}{\leftarrow} & C_2 &
\stackrel{\sim}{\leftarrow} \ldots \stackrel{\sim}{\leftarrow} & C_{m-1}
& \leftarrow & x_m\\
\end{array}
\end{equation}
We denoted it by $(\cX,\cB)$.
For $2\leq q\leq m-1$, let $\eta_q$ be the generic point
of $C_q$ and let $x_q\in C_q$ be the image of $x_m$.
By definition the following conditions hold:
\begin{itemize}
\item
For $1\leq q \leq m$, $x_q$ is near to $x_{q-1}$
and $k(x_{q-1})\simeq k(x_q)$.
\item
For $1\leq q \leq m-1$,
$C_q=\{\xi\in \phi_q^{-1}(x)|\; \tHSf{X_q}(\xi)=\tHSf{X}(x)\}$
with $\phi_q:X_q\rightarrow X$.
\item
For $1\leq q\leq m$,
$\tHSf {X_q}(x_q)=\tHSf{X}(x)$ and $\eob_{x_q}(X_q)=e_{x_q}(X_q)=2$.
\item
For $1\leq q\leq m-1$, $\tHSf {X_q}(\eta_q)=\tHSf{X}(x)$.
\item
For $1\leq q \leq m-2$, $\eob_{\eta_q}(X_q)=e_{\eta_q}(X_q)=1$.
\end{itemize}
Let $R_q=\cO_{Z_q,x_q}$ with the maximal ideal $\fm_q$, and
let $J_q$ and $\fp_q$ be the ideals of $R_q$ such that
\begin{equation}\label{XqZqCq}
X_q\times_{Z_q} \Spec(R_q)=\Spec(R_q/J_q) \qaq
C_q\times_{Z_q} \Spec(R_q)=\Spec(R_q/\fp_q).
\end{equation}
Let $R_{\fp_q}$ be the localization of $R_q$ at $\fp_q$
and $J_{\fp_q}=J_q R_{\fp_q}$.
Let $T_1,T_2$ be a pair of new variables over $k$ and consider
the isomorphism \eqref{e2c.eq1.5}:
\begin{equation}\label{e2c.eq1.5.5}
\psi_{(u)} : C_1=\bP(\Dir_x(X)) \isom \Proj(k[T_1, T_2]) =\bP^1_k
\end{equation}

\begin{definition}\label{def.label2}
\begin{itemize}
\item[(1)]
A prelabel (resp. label) $\Lambda$ of $(\cX,\cB)$ is a prelabel (resp. label)
$\Lambda=(f,y,u)$ of $(X,Z)$ at $x$.
When $x_1$ is a $k$-rational point of $C_1$,
the homogeneous coordinates of $\psi_{(u)}(x_1)\in \Proj(k[T_1,T_2])$
are called the coordinates of $(\cX,\Lambda)$.
\item[(2)]
We say that $(\cX,\cB)$ is $\bv$-admissible if $(X,Z)$ is $\bv$-admissible
in the sense of Definition \ref{def4.1} (2).
A prelabel $(f,y,u)$ of $(\cX,\cB)$ is $\OB$-admissible if it is
$\OB$-admissible as a prelabel of $(X,Z)$ at $x$ (cf. Definition \ref{loginv} (3)).
\end{itemize}
\end{definition}

We remark that the coordinates of $(\cX,\cB;\Lambda)$ depend only on $(u)$,
not on $(f,y)$.
\medbreak

\begin{lemma}\label{e2keylemma3}
Let $\Lambda=(f,y,u)$ be a label of $(\cX,\cB)$ which is
$\bv$-prepared and $\delta$-prepared and prepared along the face $E_L$, where
$L:\bR^2\to \bR\;;\; (a_1,a_2)\to a_1+2a_2.$
Assume the coordinates of $(\cX,\cB;\Lambda)$ are $(1:0)$. Let
$$
u_2'=u_2/u_1,\;
y^{(q)}=(y_1^{(q)},\dots,y_r^{(q)})\; (y_i^{(q)}=y_i/u_1^q),\;
f^{(q)}=(f_1^{(q)}, \dots, f_p^{(q)})\; (f^{(q)}_i=f_i/u_1^{q n_i}).
$$
For $1\leq q\leq m$, the following holds:
\begin{itemize}
\item[(1)]
$(f^{(q)},(y^{(q)},(u_1, u_2')))$ is a $\bv$-prepared and $\delta$-prepared
label of $(X_q,Z_q)$ at $x_q$.
\item[(2)]
$\Deltaob(f^{(q)}, y^{(q)}, (u_1, u_2'))$ is the minimal $F$-subset containg
$\Psi_q(\Deltaob(f,y,u))$, where
\par\noindent
$$\Psi_q : \bR^2\rightarrow\bR^2\;, \quad(a_1,a_2)\mapsto (a_1+a_2-q,a_2)$$
We have
$$
\betaob(f^{(q)},y^{(q)},(u_1,u_2'))  =  \gammamob(f,y,u)\leq \betaob(f,y,u).
$$
$$\alphaob(f^{(q)},y^{(q)}, (u_1,u_2'))  =  \deltaob(f,y,u)-q.$$
\item[(3)]
For $q \leq m-1$,
$\fp_q=(y^{(q)},u_1)=(y^{(q)}_1,\dots,y^{(q)}_r,u_1)$
and $v_{\fp_q}(f^{(q)}_i)=n_i$ for $i=1,\dots,N$.
For $q \leq m-2$, $(f^{(q)},y^{(q)},u_1)$ is a well prepared label of
$(X_q,Z_q)$ at $\eta_q$.
\end{itemize}
\end{lemma}
\textbf{Proof }
By Lemma \ref{e2keylemma0}, $(f^{(1)},(y^{(1)},(u_1, u_2')))$ is a label of $(X_1,Z_1)$ at $x_1$ which is $\bv$-prepared and $\delta$-prepared.
By Lemma \ref{e2keylemma2},  $(f^{(1)},(y^{(1)},(u_1, u_2')))$ is a
well prepared label of $(X_1,Z_1)$ at $\eta_1$ if $\eta_1$ is very near to
$x_1$. Then the lemma follows from Lemmas \ref{e2keylemma0b} and
\ref{e2keylemma2b}, applied to $X_q\gets X_{q+1}$ in place of $X\gets X'$.
$\square$
\medbreak

\begin{definition}\label{def.quasi.isolated}
Call $x$ or $(\cX,\cB)$ quasi-isolated if there is no regular closed subscheme
$$
D\subseteq \{\xi\in X|\; \tHSf{X}(\xi)\geq \tHSf{X}(x)\}
$$
of dimension 1 with $x\in D$.
\end{definition}

We note that this holds in particular if $x$ is isolated in $\{\xi\in X|\; \tHSf{X}(\xi)\geq \tHSf{X}(x)\}$,
and especially if $x$ lies in the $\OB$-Hilbert-Samuel locus and is isolated in it.

\begin{lemma}\label{e2keylemma4}
Let $\Lambda=(f,(y,u))$ be a label of $(\cX,\cB)$ which is
$\bv$-prepared and $\delta$-prepared and prepared along the face $E_{L'}$, where
$L':\bR^2\to \bR\;;\; (a_1,a_2)\to 2a_1+a_2.$
Assume the coordinates of $(\cX,\cB;\Lambda)$ are $(0:1)$. Set
$$
u_1'=u_1/u_2,\;
z^{(q)}=(z_1^{(q)},\dots,z_r^{(q)})\; (z_i^{(q)}=y_i/u_2^q),\;
g^{(q)}=(g_1^{(q)}, \dots, g_N^{(q)})\; (g^{(q)}_i=f_i/u_2^{q n_i}).
$$
For $1\leq q\leq m$, the following hold:
\begin{itemize}
\item[(1)]
$(g^{(q)},(z^{(q)},(u'_1, u_2)))$ is a $\bv$-prepared and $\delta$-prepared
label of $(X_q,Z_q)$ at $x_q$.
\item[(2)]
$\Deltaob(g^{(q)}, z^{(q)}, (u'_1, u_2))$ is the
minimal $F$-subset containg $\Phi_q(\Deltaob(f,y,u))$, where
\par\noindent
$$\Phi_q : \bR^2\rightarrow\bR^2, (a_1,a_2)\mapsto (a_1,a_1+a_2-q)$$
and we have
$$\alphaob(g^{(q)}, z^{(q)},(u'_1, u_2))= \alphaob(f,y,u),$$
$$\betaob(g^{(q)}, z^{(q)},(u'_1, u_2))
\leq \betaob(f,y,u)+\alphaob(f,y,u)-q.$$
\item[(3)] If $x$ is quasi-isolated, then $\alphaob(f,y,u)<1$ so that
$$\betaob(g^{(q)}, z^{(q)},(u'_1, u_2)) < \betaob(f,y,u).$$
\end{itemize}
\end{lemma}
\textbf{Proof }
This is shown in the same way as Lemma \ref{e2keylemma3} using
Lemma \ref{e2keylemma0.1} instead of \ref{e2keylemma0}.
$\square$
\medbreak

\begin{proposition}\label{e2keylemma5.cor}
Assume that $(\cX,\cB)$ is $\bv$-admissible and quasi-isolated, and that $k(x)=k(x_1)$.
Then, for all $q=1,\dots, m$, $(X_q,Z_q)$ is $\bv$-admissible at $x_q$
and
$$
\betaob_{x_q}(X_q,Z_q) \leq \betaob_x(X,Z).
$$
\end{proposition}
\textbf{Proof }
By the assumption we can take an $\OB$-admissible prelabel $\Lambda=(f,y,u)$
of $(X,Z)$, and the coordinates of $(\cX,\Lambda)$ are either $(1:0)$, or
$(0:1)$, or $(1:\lambda)$ for some $\lambda\in k$.

Assume we are in the first case.
By Lemma \ref{loginv.lem}, after preparation we may assume that
$(f,y,u)$ is prepared along the $\delta$-face and the face $E_L$
in Lemma \ref{e2keylemma3}.
By Lemmas \ref{e2keylemma3} and \ref{e2keylemma4}, applied to
$X_q\gets X_{q+1}$ for $q=1,\dots,m$ in place of $X\gets X'$, we get
a label $\Lambda_q:=(f^{(q)},(y^{(q)},(u_1, u_2')))$ of $(X_q,Z_q)$ at
$x_q$ which is $\bv$-prepared and $\delta$-prepared and deduce
$$
\begin{aligned}
\betaob_{x_q}(X_q,Z_q)\leq\betaob(f^{(q)},y^{(q)},(u_1,u_2'))
& =\betaob(f^{(1)},y^{(1)},(u_1,u_2')) \\
& = \gammamob(f,y,u)\leq \betaob(f,y,u)=\betaob_x(X,Z).
\end{aligned}
$$
where the first inequality comes from $\bv$-preparedness of $\Lambda_q$.
The case that the coordinates of $(\cX,\Lambda)$ are $(0:1)$ is shown in
the same way by using Lemma \ref{e2keylemma4} instead of Lemma
\ref{e2keylemma3}.
\medbreak

Now assume that the coordinates of $(\cX,\Lambda)$ are $(1:\lambda)$.
Let $\tu_2:=u_2 - \phi u_1$ for some $\phi\in R=\cO_{Z,x}$ such that
$\phi\mod \fm=\lambda$. Then $(f,y,(u_1, \tu_2))$ is a prelabel of $(X,Z)$ at
$x$ and $\psi_{(u_1,\tu_2)}(x_1)=(1:0)\in \Proj(k[T_1,T_2])$ so that the
coordinate of $(f,y,(u_1, \tu_2))$ is $(1:0)$. Hence the proof is reduced to
the first case in view of the following lemma.

\begin{lemma}\label{claim-case4}
Let $\tu_2=u_2 - \phi u_1$ for $\phi\in R=\cO_{Z,x}$.
Then $\bv(f,y,u)=\bv(f,y,(u_1,\tu_2))$.
For $v=\bv(f,y,u)$ we have
$$
in_v(f)_{(y,(u_1,\tu_2))}={in_v(f)_{(y,u)}}_{|U_2=\tU_2}
\in k[Y,U_1,\tU_2]=\grmR,
$$
where $\tU_2=in_{\fm}(\tu_2)\in \grmR$. Hence, if $(f,y,u)$ is $\bv$-prepared,
then so is $(f,y,(u_1,\tu_2))$.

We also have $\bvob(f,y,u)=\bvob(f,y,(u_1,\tu_2))$ and hence if
$(f,y,u)$ is $\OB$-admissible, then so is $(f,y,(u_1,\tu_2))$. .
\end{lemma}
\medbreak\noindent
\textbf{Proof }
We compute
$$
\begin{array}{rcl}
u_1^{a_1}u_2^{a_2} & = & u_1^{a_1}(\tu_2+ u_1\phi)^{a_2}\\
&=& u_1^{a_1}\sum\limits^{a_2}_{m=0}\binom{a_2}{m} u_1^m \tu_2^{a_2-m}\phi^m\\
& = & u_1^{a_1}\tu_2^{a_2}+\sum\limits^{a_2}_{m=1} \binom{a_2}{m}
u_1^{a_1+m} \tu_2^{a_2-m}\phi^m.
\end{array}
$$
This implies that the vertices on the line
$\{(a_1,a_2)|\; a_1=\alpha(f,y,u)\}$ together with
the initial forms at it, are not affected by the transformation
$(f,y,u)\to (f,y,(u_1,\tu_2))$.
Thus the first assertion follows.
The last assertion is shown by the same argument applied to
$\fob$ instead of $f$ (cf. \eqref{loginv.eq}) .
$\square$

\bigskip\noindent
\textbf{Step 2 (a chain of fundamental units):}
In this step we consider the following situation: Assume given a chain of
fundamental units of $\cB$-permissible blowups (cf. Definition \ref{Def.fspb}):
\begin{equation}\label{fspb.Thm2-1}
 (\cX_0,\cB_0) \leftarrow (\cX_1,\cB_1) \leftarrow (\cX_2,\cB_2) \gets \ldots
\end{equation}
where each $(\cX_q,\cB_q)$ is as \eqref{fupb.eq2-1}.
For $q\geq 0$, let $(x^{(q)}, X^{(q)}, Z^{(q)},\cB^{(q)})$ be the initial
part of $(\cX_q,\cB_q)$ and $m_q$ be the length of $\cX_q$. Let
$R_q=\cO_{Z^{(q)},x^{(q)}}$ with the maximal ideal $\fm_q$ and
$$
X^{(q)}\times_{Z^{(q)}} \Spec(R_q)=\Spec(R_q/J_q)
\quad\text{for an ideal $J_q\subset R_q$}.
$$

\begin{theorem}\label{Thm2-1.Step4.thm}
Assume that for all $q\geq 0$, $k(x^{(q)})$ is separably algebraic over
$k(x^{(0)})$ and $x_q$ is quasi-isolated (see Definition \ref{def.quasi.isolated}).
Then the sequence \eqref{fspb.Thm2-1} stops after finitely many steps.
\end{theorem}
\textbf{Proof  }
By Lemmas \ref{completion.nu.2} and \ref{completion.HS},
it suffices to show the claim after replacing each $\cX_q$ by its base changes via
$\Spec(\hat{R}^{ur}) \to Z^{(0)}$, where $\hat{R}^{ur}$ is
the maximal unramified extension of the completion of $R=R_0$
(here we use that $Z$ is excellent).
Thus we may assume that $k(x^{(q)})=k(x)$ and that $R_0$ is complete and hence that
$(\cX_0,\cB_0)$ is $\bv$-admissible.
Proposition \ref{e2keylemma5.cor} implies that for all $q\geq 0$,
$(\cX_q,\cB_q)$ is $\bv$-admissible so that
$\beta_q:=\betaob_{x^{(q)}}(X^{(q)},Z^{(q)})$ is well-defined and
\begin{equation}\label{Thm2-1.Step4.thm.eq}
\beta_{q+1}\leq\beta_q\quad \text{ for all }q\geq 0.
\end{equation}
Note
$\beta_q\in 1/n_N! \cdot \bZ^2\subset\bR^2$,
where
$\nu^*(J_0)=\nu^*(J_q)=(n_1,\dots,n_N,\infty,\dots)$
(cf. Lemma \ref{Vlattice}, Lemma \ref{def1.4cor} and
Theorem \ref{thm.pbu.inv} (3)).
Hence the strict inequality may occur in \eqref{Thm2-1.Step4.thm.eq} only for
finitely many $q$. Hence we may assume $\beta_q=\beta_0$ for all $q\geq 0$.
\medbreak

In what follows, for a prelabel $(g,(z,u))$ of $\cX_0$ with
$g=(g_1,\dots,g_N)$ and $z=(z_1,\dots,z_r)$) and for $q\geq 0$, we write
$$
g^{(q)}=(g_1^{(q)},\dots,g_N^{(q)}) \qwith
g_i^{(q)} =
\left\{ \begin{array}{ccc}
g_i/u_1^{n_i(m_1+\cdots + m_q)} & , & q \geq 1, \\
g_i & , & q=0.\\
\end{array} \right.
$$
$$
z^{(q)}=(z_1^{(q)},\dots, z_r^{(q)}) \qwith z_i^{(q)}=
\left\{ \begin{array}{ccc}
z_i/u_1^{m_0+\cdots + m_{q-1}} & , & q \geq 1, \\
z_i & , & q=0.\\
\end{array} \right.
$$
By the completeness of $R_0$, we can choose a totally prepared and
$\OB$-admissible label $\Lambda_0=(f,y,u))$ of $(\cX_0,\cB_0)$. By the
assumption (a), $\beta_0=\beta_1$ implies by Lemma \ref{e2keylemma4} (2)
that the coordinate of $(\cX_0,\Lambda_0)$ must be $(1:-\lambda_0)$ for
some $\lambda_0\in k$.
Put $v_1=u_2 + \phi_0 u_1$ for a lift $\phi_0\in R$ of $\lambda_0$ and
prepare $(f,y,(u_1,v_1))$ to get a totally prepared label
$\Lambda'_0=(g,z,(u_1,v_1))$ of $(\cX_0,\cB_0)$.
By Lemma \ref{claim-case4}, $\Lambda'_0$ is $\OB$-admissible and
the coordinate of $(\cX_0,\Lambda'_0)$ is $(1:0)$. Lemma \ref{e2keylemma3} (1)
implies that $\Lambda_1=(g^{(1)},z^{(1)},(u_1,v_1/u_1))$ is a totally
prepared label of $\cX_1$. By the assumption (a),
$\beta_0=\beta_1=\beta_2$ implies by Lemma
\ref{e2keylemma4} (2) that the coordinate of $(\cX_1,\Lambda_1)$
is $(1:-\lambda_1)$ for some $\lambda_1\in k$.
Put $v_2=v_1 + \phi_1 u_1^2=u_2 + \phi_0 u_1+ \phi_1 u_1^2$ for
a lift $\phi_1\in R$ of $\lambda_1$. Prepare $(f,y,(u_1,v_2))$ to get
a totally prepared label $\Lambda''_0=(h,w,(u_1,v_2))$ of
$(\cX_0,\Lambda'_0)$. Then $\Lambda''_0$ is $\OB$-admissible by Lemma
\ref{claim-case4} and $\Lambda'_1=(h^{(1)},(w^{(1)},(u_1,v_2/u_1)))$ is
a totally prepared label of $\cX_1$ by Lemma \ref{e2keylemma3} (1).
Moreover, the coordinate of $(\cX_0,\Lambda''_0)$ and that of
$(\cX_1,\Lambda'_1)$ are both $(1:0)$. Lemma \ref{e2keylemma3} (1) then
implies that $\Lambda_2=(h^{(2)},w^{(2)},(u_1,v_2/u_1^2))$ is a
totally prepared label of $\cX_2$. The same argument repeats itself to imply
the following:

\begin{claim}\label{step4.claim1}
Assume that the sequence \eqref{fspb.Thm2-1} proceeds in infinitely many steps
and that $\beta_q=\beta_0$ for all $q\geq 0$.
Then there exists a sequence $\phi_0,\phi_1,\phi_2,\dots$ of elements in $R$
for which the following holds: Recalling that $R$ is complete, set
$$
v= \lim_{q\to \infty}
\big(u_2 + \phi_0 u_1+ \phi_1 u_1^2 + \ldots + \phi_{q-1} u_1^q\big) \;
\in \; R\;.
$$
Prepare $(f,y,(u_1,v)))$ to get a totally prepared label
$(\hf,\hy,(u_1,v))$ of $(\cX_0,\cB_0)$. Then $(\hf,\hy,(u_1,v))$ is
$\OB$-admissible and $\hLam_q=(\hf^{(q)},\hy^{(q)},(u_1,v^{(q)}))$ is
a totally prepared label of $(\cX_q,\cB_q)$ and the coordinate of
$(\cX_q,\hLam_q)$ is $(1:0)$ for all $q\geq 0$.
Here $v^{(q)}=v/u_1^q$.
\end{claim}
\bigskip

We now write $(f,y,(u_1,u_2))$ for $(\hf,\hy,(u_1,v))$.
Lemma \ref{e2keylemma3} implies that for all $q\geq 0$,
$(f^{(q)},y^{(q)},(u_1,u_2^{(q)}))$ is a totally prepared label of
$(\cX_q,\cB_q)$.
Moreover $\Deltaob(\ff q, \yy q, (u_1, u_2^{(q)}))$ is
the minimal $F$-subset of $\bR^2$ containing
$$
T_q(\Deltaob(\ff {q-1}, \yy {q-1},(u_1, u_2^{(q-1))}))),
\quad
T_q:\bR^2\rightarrow \bR^2\;;\; (a_1,a_2)\to (a_1+a_2-m_{q-1}, a_2),
$$
where $m_{q-1}$ is the length of $\cX_{q-1}$.
This implies that $\epsilon^\OB(\ff {q},\yy {q}, (u_1, u_2^{(q)}))=\epsilon^\OB(f,y,u)=: \epsilon^\OB$
for all $q\geq 0$ (see Definition \ref{def.invpol}), and that
\begin{eqnarray*}
\zeta^\OB(f^{(q)},y^{(q)},(u_1,u_2^{(q)}))& = & \zeta^\OB(f^{(q-1)},y^{(q-1)},(u_1,u_2^{(q-1)}))
+ \epsilon^\OB -m_{q-1} \\
& < & \zeta^\OB(f^{(q-1)},y^{(q-1)},(u_1,u_2^{(q-1)}))\,,
\end{eqnarray*}
for all $q\geq 1$, because $m_{q-1}\geq 1$, and $\epsilon^\OB < 1$ by Lemma \ref{case1.lem3}
and the assumption (a) of Theorem \ref{Thm2-1.Step4.thm}.
This implies that the sequence must stop after finitely many steps as claimed.

\newpage
\section{Proof in the case $e_x(X)=\overline e_x(X)=2$, III: inseparable residue extensions}\label{sec:e2III}

\bigskip

In this section we complete the proof of key Theorem \ref{Thm2}
(see Theorem \ref{Thm.nonrational} below).
\bigskip

Let the assumptions and notations be as in the beginning of Case 1 of
\S\ref{sec:e2I}.
Let $\Phi(T_1,T_2)\in k[T_1,T_2]$ be an irreducible homogeneous
polynomial corresponding to $x'\in C=\Proj(k[T_1, T_2])$
(cf. \eqref{e2c.eq1.5}). We set
$$
\Phi=\Phi(U_1,U_2)\in k[U_1,U_2]\subset k[Y,U_1,U_2]=\grmR.
$$
We assume $x'\not=(0:1)$ so that $U_1$ does not divide $\Phi$. Let
$$
d= \deg\Phi=[k(x'):k(x)].
$$
Choosing a lift $\tilde{\Phi}(U_1,U_2)\in R[U_1,U_2]$ of
$\Phi\in k[U_1,U_2]$, set
$$
\phi=\tilde{\Phi}(u_1,u_2)\in R \qaq \phi'=\phi/u_1^{deg\Phi}\in R'\;.
$$
The following two lemmas are shown in the same way as
Lemma \ref{e2keylemma0} (1), (2) and (5).

\begin{lemma}\label{lem.13pre1}
Let
$$
y'=(y'_1,\dots,y'_r)\text{ with } y'_i=y_i/u_1, \;
f'=(f'_1,\dots,f'_N) \text{ with }f_i'=f_i/u_1^{n_i}.
$$
\begin{itemize}
\item[(1)]
$(y', (u_1,\phi'))$ is a system of regular parameters of $R'$ such that
$\fp=(y',u_1)$ (cf. \eqref{eq.Jdpd}), and $J'=\langle f'_1,\dots,f'_N\rangle$.
\item[(2)]
If $x'$ is very near to $x$, then $(f',y',(u_1,\phi'))$ is a prelabel of
$(X',Z')$ at $x'$.
\end{itemize}
\end{lemma}

\begin{lemma}\label{lem.13pre2}
Assume $x'$ is $\OB$-near to $x$. Setting $l'_B=l_B/u_1\in R'$ for $B\in \OBx$,
we have
$$
\Deltaob(f',y',(u_1,\phi'))=\Delta(\fdob,y,(u_1,\phi'))\qwith
\fdob=(f',\; l'_B\; (B\in \OBx)).
$$
\end{lemma}
\medbreak

Now we state the main result of this section.

\begin{proposition}\label{pr.13}
Let $\delta=\delta(f,y,u)$. Assume the following conditions:
\begin{itemize}
\item[(a)]
$d \geq 2$ and there is no regular closed subscheme $D \subseteq \tHSmax X$
of dimension 1 with $x\in D$,
\item[(b)]
$(f,y,u)$ is $\bv$-prepared and prepared at $\bw^+(f,y,u)$,
\item[(c)]
$in_\delta(f)_{(y,u)}$ is normalized.
\end{itemize}
Then there exists a part of a system of regular parameters
$z'=(z'_1,\dots,z'_r)\subset \fm'$ such that the following holds:
\begin{itemize}
\item[(1)]
$(f',z',(u_1,\phi'))$ is a $\bv$-prepared prelabel of $(X',Z')$ at $x'$,
\item[(2)]
$\betaob(f',z',(u_1,\phi'))<\betaob(f,y,u)$,
$\alpha(f',z',(u_1,\phi'))=\alpha(f',y',(u_1,\phi'))=\delta-1$,
\item[(3)]
$z'=y'$ unless $\delta\in \bZ$ and $\delta\geq 2$.
In the latter case we have
$z'_i- y'_i\in \langle u_1^{\delta-1}\rangle$ for $i=1,\dots,r$.
In particular $\langle y',u_1\rangle=\langle z',u_1\rangle $.
\end{itemize}
\end{proposition}
\bigskip

Before proving Proposition \ref{pr.13}, we now complete the proof of
Theorem \ref{Thm2}.
In view of Proposition \ref{e2keylemma5.cor} and
Theorem \ref{Thm2-1.Step4.thm} (and its proof),
Theorem \ref{Thm2} is an obvious consequence of the following.

\begin{theorem}\label{Thm.nonrational}
Assume given a fundamental unit of $\cB$-permissible blowups
\eqref{fupb.eq2-1}. Assume that $k(x_1)\not=k(x)$ and that there is no
regular closed subscheme $D \subseteq \tHSmax X$ of dimension 1 with $x\in D$.
(E.g., this holds if $x$ is isolated in $\tHSmax X$).
Assume that $(X,Z)$ is $\bv$-admissible at $x$ (cf. Definition \ref{loginv} (3)).
Then, for $q=1,\dots, m$, $(X_q,Z_q)$ is $\bv$-admissible at $x_q$ and
$$
\betaob_{x_q}(X_q,Z_q) < \betaob_x(X,Z).
$$
\end{theorem}
\textbf{Proof  }
By the assumption we can take an $\OB$-admissible prelabel $\Lambda=(f,y,u)$.
By Lemma \ref{loginv.lem}, after preparation we may assume that
$(f,y,u)$ is $\bv$-prepared and $\delta$-prepared.
Let $(f',y',(u_1,\phi'))$ and $(f',z',(u_1,\phi'))$ be as in
Lemma \ref{lem.13pre1} and Proposition \ref{pr.13}, applied to $X\gets X_1$
in place of $X\gets X'$. By Claim \ref{fu.claim1}, we have (cf. \eqref{XqZqCq})
$$
\fp_q=\langle y'_1/u_1^{q-1},\dots,y'_r/u_1^{q-1},u_1\rangle
\qfor q=1,\dots, m-1.
$$
By Proposition \ref{pr.13} (3) and since $\delta(f,y,u)\geq m$ by
Corollary \ref{cor.fu}, this implies
\begin{equation}\label{pr.13.eq0}
\fp_q=\langle z'_1/u_1^{q-1},\dots,z'_r/u_1^{q-1},u_1\rangle
\qfor q=1,\dots, m-1.
\end{equation}
We prepare $(f',z',(u_1,\phi'))$ at all vertices and the faces lying in
$$
\{A\in \bR^2\;|\; |A|\leq |\bv(f',z',(u_1,\phi'))|\}
$$
to get a $\bv$-prepared and $\delta$-prepared label $(g,w,(u_1,\phi'))$ of $(X_1,Z_1)$ at $x_1$. Note that
\begin{equation}\label{pr.13.eq1}
w_i-z'_i\in \langle u_1^\gamma\rangle \; (1\leq i\leq r)
\qfor \gamma\in \bZ_{\geq 0},\; \gamma>\alpha(f',z',(u_1,\phi'))=\delta-1.
\end{equation}
By Lemma \ref{loginv.lem}, we have
$$
\betaob(g,w,(u_1,\phi'))=\betaob(f',z',(u_1,\phi')).
$$
For $q=1,\dots,m$, let
$$
g^{(q)}=(g_1^{(q)},\dots,g_N^{(q)})\; (g_i^{(q)}=g/u_1^{n_i(q-1)}),\quad
w^{(q)}=(w_1^{(q)},\dots,w_N^{(q)})\; (w_i^{(q)}=w/u_1^{(q-1)}).
$$
Then \eqref{pr.13.eq0} and \eqref{pr.13.eq1} imply
$\fp_q=\langle w_1^{(q)},\dots,w_r^{(q)},u_1\rangle$ for $q=1,\dots, m-1$.
For $q\geq 1$, Lemma \ref{e2keylemma0b}, applied to $X_q\gets X_{q+1}$ in place of
$X\gets X'$, implies that $\Lambda_q:=(g^{(q)},w^{(q)},(u_1, \phi'))$ is a label
of $(X_q,Z_q)$ at $x_q$ which is $\bv$-prepared and $\delta$-prepared.
Then we get
$$
\begin{aligned}
\betaob_{x_q}(X_q,Z_q)\leq \betaob(g^{(q)},w^{(q)},(u_1,\phi'))
&=\betaob(g,w,(u_1,\phi'))=\betaob(f',z',(u_1,\phi'))\\
&<\betaob(f,y,u)=\betaob_x(X,Z),\\
\end{aligned}
$$
where the first inequality (resp. equality) comes from $\bv$-preparedness of
$\Lambda_q$ (resp. Lemma \ref{e2keylemma0b} (4)).
This completes the proof of the theorem.
$\square$
\bigskip

Now we start the proof of Proposition \ref{pr.13}.
We may write
\begin{equation}\label{eq.13-4}
in_{\delta}(f_i)_{(y,u)}
= F_i(Y)+\underset{|B|<n_i}{\sum} P_{i,B}(U) \cdot Y^B \;\in k[Y,U_1,U_2],\\
\end{equation}
where $P_{i,B}(U)\in k[U_1,U_2]$, for $B\in \bZ^r_{\geq 0}$ with $|B|<n_i$,
is either $0$ or homogeneous of degree $\delta(n_i-|B|)$.
Write
\begin{equation}\label{eq.13-3.5}
P_{i,B}(U)=\Phi^{e_i(B)}\cdot Q_{i,B}(U),
\end{equation}
where $e_i(B)\in \bZ_{\geq 0}$ and $Q_{i,B}(U) \in k[U_1,U_2]$ is either $0$, or
homogeneous of degree $(\delta-d\cdot e_i(B))(n_i-|B|)$ and not divisible by
$\Phi$. Then we get
\begin{equation}\label{eq.13-4.1}
in_{\delta}(f_i)_{(y,u)}= F_i(Y)+{\sum} Q_{i,B}(U) Y^B \Phi^{e_i(B)}
 \;\in k[Y,U_1,U_2],
\end{equation}
From this we compute
\begin{equation}\label{eq.gammap}
\begin{aligned}
\gamma^+ (f,y,u)
&=\sup \left\{ \frac{\deg_{U_2} P_{i,B}(1,U_2)}{n_i-|B|}\left| \;
1\leq i\leq N,\; P_{i,B}(U)\not\equiv 0
\right.\right\}\,.\\
&=\sup \left\{ \frac{d\cdot e_i(B)+\deg_{U_2} Q_{i,B}(1,U_2)}{n_i-|B|}\left|\;
1\leq i\leq N, \; Q_{i,B}(U)\not\equiv 0
\right.\right\}\,.\\
&\geq d \cdot
\sup \left\{ \frac{e_i(B)}{n_i-|B|}\left|\;
1\leq i\leq N,\; Q_{i,B}(U)\not\equiv 0
\right.\right\}\,.
\end{aligned}
\end{equation}
where $\deg_{U_2}$ denotes the degree of a polynomial in $k[U_2]$. Set
$$
q_{i,B}=\tilde{Q}_{i,B}(u_1,u_2)\in R \qaq
q'_{i,B}=q_{i,B}/u_1^{\deg Q_{i,B}}\in R'\;,
$$
where $\tilde{Q}_{i,B}(U_1,U_2)\in R[U_1,U_2]$ is a lift of
$Q_{i,B}(U)\in k[U_1,U_2]$.
Letting $\tilde{F}(Y)\in R[Y]$ be a lift of $F(Y)\in k[Y]$,
\eqref{eq.13-4} and \eqref{eq.13-3.5} imply
\begin{equation}\label{eq.3stardownstairs}
f_i=\tilde F_i(y)+\sum\limits_{|B|<n_i}
y^B u_1^{a_1}(\phi)^{e_i(B)} q_{i,B} +g
\qwith v_L(g)_{(y,u)}> v_L(f_i)_{(y,u)} = n_i\,,
\end{equation}
By Lemma \ref{lem2.1}(2) this implies
\begin{equation}\label{eq.3star}
f_i'=f_i/u_1^{n_i}=\tilde F_i(y')+\sum\limits_{|B|<n_i}
{y'}^B u_1^{(\delta-1)(n_i-|B|)}(\phi')^{e_i(B)} q'_{i,B} +g'
\qwith v_\Lambda(g')_{(y',(u_1,\phi'))} >n_i\,,
\end{equation}
where
$\displaystyle{\Lambda:\bR^2\to \bR\;;\; (a_1,a_2) \to \frac{a_1}{\delta-1}}$
and $g'=g/u_1^{n_i}$. Noting that
\begin{equation}\label{eq.Qq}
q'_{i,B}\in (R')^\times \Leftrightarrow Q_{i,B}(U)\not\equiv 0\in k[U_1,U_2],
\end{equation}
we get
\begin{equation}\label{eq.alphabeta}
\begin{aligned}
\alpha' := & \alpha( f',y',u_1,\phi')=\delta-1
\quad \text{where }\delta:=\delta(f,y,u),\\
\beta':= & \beta ( f',y',u_1,\phi')=
\inf\left\{\frac{e_i(B)}{n_i-|B|}\left|\;
1\leq i\leq N,\; Q_{i,B}(U)\not\equiv 0 \right.\right\}\,.
\end{aligned}
\end{equation}
The same argument, applied to $l_B$ with $B\in \OBx$, shows
\begin{equation}\label{eq.alphaob}
\alphaob( f',y',u_1,\phi')=\deltaob(f,y,u)-1.
\end{equation}
Let $v'=\bv(f',y',(u_1,\phi'))$. Under the identification
$$
\grmdRd=k'[Y',U_1,\Phi']
\quad\big(Y'_i=\inmd(y'_i),\; U_1=\inmd(u_1),\;\Phi'=\inmd(\phi')\big),
$$
we get
\begin{equation}\label{eq.13-5}
in_{v'}(f_i')=F_i(Y')+\underset{e_i(B)=\beta'(n_i-|B|)}{\sum} \qbdiB \cdot
{Y'}^B (U_1^{\alpha'} {\Phi'}^{\beta'})^{n_i-|B|}\;
\in k'[Y',U_1,\Phi']
\end{equation}
where
$\qbdiB \in k':=k(x')$ is the residue class of $\qdiB\in R'$.

\begin{lemma}\label{lem.claim1}
We have
$$
\begin{aligned}
\beta(f,y,u)\geq \gamma^+(f,y,u)&\geq d \cdot \beta(f',y',(u_1,\phi')),\\
\betaob(f,y,u)\geq \gammapob(f,y,u)&\geq d \cdot \betaob(f',y',(u_1,\phi')).\\
\end{aligned}
$$
\end{lemma}
\textbf{Proof }
The first inequality holds in general (cf. the picture below Definition \ref{def.invpol})
and the second follows from \eqref{eq.alphabeta} and \eqref{eq.gammap}.
This proves the first assertion. The second assertion follows
by applying the same argument to $l_B$ for $B\in \OBx$ in view of
Lemma \ref{lem.13pre2}. This completes the proof.
$\square$

\begin{corollary}\label{cor.cor2}
If $(f',y',(u_1,\phi'))$ is not solvable at $v'$, Proposition \ref{pr.13}
holds.
\end{corollary}
\textbf{Proof  }
Indeed, it suffices to take $z'=y'$ in this case. Proposition \ref{pr.13} (3)
follows from Lemma \ref{e2keylemma0} (1).
As for (1), Proposition \ref{pr.13} (c) implies, in view of
\eqref{eq.13-4.1}, \eqref{eq.13-5}, and \eqref{eq.Qq}, that $(f',y',(u_1,\phi'))$ is normalized at $v'$.
Hence the assumption implies $(f',y',(u_1,\phi'))$ is prepared at $v'$.
It remains to show (2). By Lemma \ref{lem.claim1} it suffices to show
$\betaob(f,y,u)\not=0$. By \eqref{eq.eOB2} we have
$$
\betaob(f,y,u)+\alphaob(f,y,u) \geq \deltaob(f,y,u)>1.
$$
By the assumption (a) and Lemma \ref{case1.lem3},
$\alphaob(f,y,u)<1$ and hence $\betaob(f,y,u)>0$.
$\square$

\bigskip\medskip

So for the proof of Proposition \ref{pr.13}, it remains to treat the case
where $(f',y',(u_1,\phi'))$ is solvable at $v'$.
Assume that we are now in this case. This implies
\begin{equation}\label{betadeltaZ}
\beta':=\beta(f',y',(u_1,\phi'))\in \bZ_{\geq 0},\quad
\delta:=\delta(f,y,u)\in \bZ,\quad \delta\geq 2.
\end{equation}
Indeed $\delta\in \bZ$ since $\alpha'=\delta-1\in \bZ$ by
\eqref{eq.alphabeta} and $\delta>1$ by the assumption (d).
It also implies that there exist $\lambda_1,\ldots , \lambda_r\in k'$ such that
\begin{equation}\label{eq.dissolve}
in_{v'}(f'_i)=F_i(Y'+\lambda\cdot \Phi^{\beta'}U_1^{\delta-1})
\qfor i=1,\dots,N,
\end{equation}
where $\lambda=(\lambda_1,\ldots , \lambda_r)$.
For $1\leq j \leq r$, let $A_j(U)=A_j(U_1,U_2)\in k[U_1, U_2]$ be
homogeneous polynomials such that:
\begin{itemize}
\item[(C1)]
$\alpha_j:= \deg A_j<d=\deg \Phi$,
\item[(C2)]
$A_j$ is not divisible by $U_1$
(which implies $\alpha_j=\deg_{U_2} A_j(1, U_2))$),
\item[(C3)]
$\lambda_j\equiv A_j(1,U_2)\mod\Phi(1, U_2)$ in
$k'=k(x') \cong k[U_2]/\langle \Phi(1, U_2) \rangle$.
\end{itemize}
\medskip
Choose a lift of $\tilde A_j(U_1,U_2)\in R[U_1,U_2]$ of
$A_j(U)\in k[U_1,U_2]$ and set
$$
a_j=\tilde A_j(u_1,u_2)\in R \qaq a'_j=a_j/u_1^{\alpha_j}\in R'.
$$
Then $\lambda_j \equiv a'_j\mod \fm'$. Define
\begin{equation}\label{star.on.6}
z_j:= y_j+\phi^{\beta'}\cdot u_1^{\delta-(d\beta'+\alpha_j)}
\cdot a_j \in R[1/u_1],
\end{equation}
\begin{equation}\label{star.on.6up}
z'_j:=z_j/u_1=y'_j+{\phi'}^{\beta'}\cdot u_1^{\delta-1}\cdot a'_j \in R'.
\end{equation}
Note that $z'_i-y'_i\in \langle u_1\rangle$ for $i=1,\dots, r$ since
$\delta\geq 2$ as noted in \eqref{betadeltaZ}.
Consider the following condition
\begin{equation}\label{fis}
\gamma^ +:=\gamma^+(f,y,u) \geq d(\beta'+1)
\quad (\beta':=\beta(f',y',(u_1,\phi'))).
\end{equation}

\begin{lemma}\label{lem.claim4}
Assume \eqref{fis} holds. Then the following is true.
\begin{itemize}
\item[(1)]
We have $z_j\in R$ for $j=1,\dots,r$ and  $(z,u)$ is a system of regular
parameters of $R$ which is strictly admissible for $J$.
The conditions (b) and (c) of Proposition \ref{pr.13} are satisfied for $(z,u)$
in place of $(y,u)$.
\item[(2)]
$\bw^+(f,z,u)=\bw^+(f,y,u)$, $\bv(f,z,u)=\bv(f,y,u)$,

$\bvob(f,z,u)=\bvob(f,y,u)$, $\delta(f,z,u)=\delta(f,y,u)$.
\item[(3)]
$\alpha(f',y',(u_1,\phi')) = \alpha(f',z',(u_1,\phi')),$

$\beta(f',y',u_1,\phi')  <  \beta(f',z',u_1,\phi') \leq
\gamma^+(f,y,u)=\gamma^+(f,z,u).$
\end{itemize}
\end{lemma}

By Lemma \ref{lem.claim4}, if \eqref{fis} holds for $(f,y,u)$,
we may replace $(f,y,u)$ by $(f,z,u)$ to show Proposition \ref{pr.13}.
If $(f,z,u)$ is not solvable at $\bv(f,z,u)$, then we are done by
Corollary \ref{cor.cor2}. If $(f,z,u)$ is solvable at $\bv(f,z,u)$ and
\eqref{fis} holds for $(f,z,u)$, then we apply the same procedure to $(z,u)$
to get a new system of regular parameters of $R$. This process must stop
after finitely many steps, by the last inequality in Lemma \ref{lem.claim4} (3).
Thus Proposition \ref{pr.13} follows from Corollary \ref{cor.cor2}, Lemma \ref{lem.claim4}
and the following Lemma \ref{lem.claim5}.

\begin{lemma}\label{lem.claim5}
Assume that \eqref{fis} does not hold for $(f,y,u)$ and that
$(f',y',(u_1,\phi'))$ is solvable at $v'=\bv(f',y',(u_1,\phi'))$.
Dissolve $v'$ as in \eqref{eq.dissolve} and let $z'_j$ be as in
\eqref{star.on.6up}. Then $(f',z',(u_1,\phi'))$ is $\bv$-prepared and we have
$$
\betaob(f',z',(u_1,\phi'))<\betaob(f,y,u),\quad
\alpha(f',z',(u_1,\phi'))=\alpha(f',y',(u_1,\phi'))=\delta-1.
$$
\end{lemma}

\bigskip
\textbf{Proof of Lemma \ref{lem.claim4}}.
Conditions \eqref{fis} and (C1) imply
\begin{equation}\label{eq.13-6}
\delta-(d \beta'+\alpha_j) > \delta-d(\beta'+1) \geq
\delta -\gamma^+\geq 0 \qfor j=1,\dots, r,
\end{equation}
where the last inequality holds in general (cf. the picture below Definition \ref{def.invpol}).
This implies $z_j\in R$, and the second assertion of (1) is obvious from
\eqref{star.on.6}. Let $v_{u_2}$ be the valuation on $R$ with respect to
$\langle u_2\rangle $. By \eqref{fis} and (C1) we have
$$
v_{u_2}(\phi^{\beta'}a_j) = d \beta'+\alpha_j <
d \beta'+ d  \leq  \gamma^+ \,.
$$
This together with \eqref{eq.13-6} implies that the coordinate transformation
\eqref{star.on.6} affects only those vertices of $\Delta(f,y,u)$ lying in
$\{(a_1,a_2)\in \bR^2|\; a_1> \delta-\gamma^+,\; a_2< \gamma^+\}$.
This shows (2), and that condition (b) of Proposition \ref{pr.13} holds for $(f,z,u)$.
\medbreak

The first assertion of (3) follows from \eqref{eq.alphabeta} and
the equality $\delta(f,y,u)=\delta(f,z,u)$ implied by (2).
The first inequality in the second assertion of (3) is a consequence of
the first assertion since
$\bv(f',y',(u_1,\phi'))\not\in \Delta(f',z',(u_1,\phi'))$.
Lemma \ref{lem.claim1} implies $\gamma^+(f,z,u)\geq \beta(f',z',(u_1,\phi'))$
and (2) implies $\gamma^+(f,z,u)=\gamma^+(f,y,u)$, which completes the proof
of (3).
\medbreak

It remains to show that condition (c) of Proposition \ref{pr.13} holds for $(f,z,u)$.
Introduce $\rho=(\rho_1,\dots,\rho_r)$, a tuple of independent
variables over $k$. For $i=1,\dots, N$, write
$$
in_{\delta}(f_i)_{(y,u)}=\underset{B}{\sum} Y^B T_{i,B}(U)
\qwith T_{i,B}(U)\in k[U]
$$
and substitute $Y_i=Z_i-\rho_i$ in $in_{\delta}(f_i)_{(y,u)}\in k[Y,U]$ to get
$$
G_i(Z,U,\rho)=\underset{B}{\sum}(Z-\rho)^B T_{i,B}(U)\,
\in k[Z,U,\rho].
$$
By \eqref{eq.13-4} and \eqref{star.on.6}
we have $in_{\delta}(f_i)_{(z,u)}= G_i(Z,U,\rho)_{|\rho=s}$, where
$$
s=(s_1,\dots,s_r) \qwith
s_i= U_1^{\delta-(d \beta'+\alpha_j)}
\Phi^{\beta'}A_j(U_1, U_2)\in k[U_1,U_2]\,.
$$
Write
$$
G_i(Z,U,\rho) = \sum\limits_C Z^C S_{i,C}(U,\rho)
\quad(S_{i,C}(U,\rho)\in k[U,\rho])\,.
$$
By condition \ref{pr.13} (c) for $(f,y,u)$, i.e.,  the normalizedness of $in_{\delta}(f_i)_{(y,u)}$,
we have $T_{i,B}(U)\equiv 0$ if $B\in E^r(F_1(Y),\ldots , F_{i-1}(Y))$.
By Lemma \ref{lem.N} below, this implies
$S_{i,C}(U,\rho)\equiv 0$ for $C\in E^r(F_1(Z),\ldots , F_{i-1}(Z))$,
which is the normalizedness of $in_{\delta}(f)_{(z,u)}$, i.e., \ref{pr.13} (c)
for $(f,z,u)$.
$\square$

\begin{lemma} \label{lem.N}
Assume given
$$
G(Y)=\sum\limits_AC_A\cdot Y^A \in k[Y]=k[Y_1,\ldots , Y_r],
$$
and a subset $E\subset \bZ^r_{\geq 0}$ such that
$E+\bZ^r_{\geq 0}\subset E$ and that $C_A=0$ if $A\in E$.
Let $\rho = (\rho_1,\ldots , \rho_r)$ be a tuple of independent variables
over $k$ and write
$$
G(Y+\rho)=\sum\limits_K S_K\cdot Y^K\quad\mbox{with}\quad S_K\in k[\rho]\,.
$$
Then $S_K\equiv 0$ if $K\in E$.
\end{lemma}
\textbf{Proof} We have
$$
\begin{array}{rcl}
G(Y+\rho) & = & \sum\limits_A C_A
\left[\prod\limits^r_{i=1}(Y_i+\rho_i)^{A_i}\right]\\
& = & \sum\limits_A C_A \prod\limits^r_{i=1}
\prod\limits^{A_i}_{k_i}
\begin{pmatrix} A_i\\ k_i\end{pmatrix}
\rho_i^{A_i-k_i}\cdot Y^{k_i}\\
& = & \sum\limits_{K=(k_1,\ldots , k_r)}
\left[\sum\limits_{A\in K+\mathbb Z_{\geq 0}}
\left(C_A\prod\limits^m_{i=1}\begin{pmatrix}
A_i\\k_i\end{pmatrix} \right)\rho^{A-K}\right]Y^K\,.
\end{array}
$$
Thus
$$
S_K=\sum\limits_{A\in K+\bZ_{\leq 0}}
\left(C_A\prod\limits^m_{i=1}
\begin{pmatrix}A_i\\k_i\end{pmatrix} \right) \cdot\rho^{A-K}\,.
$$
Now, if $S_K\not\equiv 0$ for $K\in E$, then there is an
$A\in K+\bZ^r_{\geq 0}\subset E$ with $C_A\neq 0$.
This completes the proof of the lemma.
$\square$

\bigskip
\textbf{Proof of Lemma \ref{lem.claim5}}.
Let $\tilde{F}_i(Y)\in R[Y]$ be a lift of $F_i(Y)\in k[Y]$.
For a tuple of independent variables $\rho=(\rho_1,\dots,\rho_r)$ over $R$,
write
\begin{equation}\label{eq.2star.on.9-1}
\tilde{F_i}(Y + \rho)=\tilde{F_i}(Y)+
\sum\limits_{\overset{|B|<n_i}{|B+D|=n_i}}
K_{i,B,D}\cdot {Y}^B\rho^D \qwith K_{i,B,D}\in R,
\end{equation}
By \eqref{star.on.6up} we have $z'=y'+\mu$, where
\begin{equation}\label{eq.13-7}
\mu = (\mu_1,\ldots,\mu_r) \qwith
\mu_i = u_1^{\delta-1}(\phi')^{\beta'}\cdot a_i'.
\end{equation}
Hence \eqref{eq.2star.on.9-1} implies:
\begin{equation*}\label{eq.2star.on.9-1.1}
\tilde{F_i}(z')=\tilde{F_i}(y')+
\sum\limits_{\overset{|B|<n_i}{|B+D|=n_i}}
K_{i,B,D}\cdot {y'}^B\mu^D \qwith K_{i,B,D}\in R,
\end{equation*}
By \eqref{eq.3star} this implies
\begin{equation}\label{eq.2star.on.9-1.5}
f'_i=\tilde{F_i}(z')+\sum\limits_{|B|< n_i}(z'-\mu)^B \theta_{i,B}+g'
\qwith v_\Lambda(g')_{(y',(u_1,\phi'))} > n_i,
\end{equation}
where
$$
\begin{array}{rcl}
\theta_{i,B} & = & (\phi')^{e_i(B)}q'_{i,B}u_1^{(\delta-1)(n_i-|B|)}-
\sum\limits_{|D|=n_i-|B|}K_{i,B,D}\cdot \mu^D \\
& = & ({\phi'}^{\beta'} u_1^{\delta-1})^{n_i-|B|} \omega_{i,B}\;.
\end{array}
$$
Here we set
\begin{equation}\label{eq.omega}
\omega_{i,B}=(\phi')^{b_i(B)}q'_{i,B}-
\sum\limits_{|D|=n_i-|B|}K_{i,B,D}\cdot {a'}^D
\quad (a'=(a'_1,\ldots, a'_r)),
\end{equation}
$$
b_i(B)=e_i(B)-(n_i-|B|)\beta'.
$$
By \eqref{eq.alphabeta} we have $b_i(B)\geq 0$. For each $B$
write (in $k[U_2]$)
\begin{equation} \label{eq.3star.on.9-2}
\Phi(1,U_2)^{b_i(B)}Q_{i,B}(1,U_2)-\sum\limits_{|D|=n_i-|B|}
\overline{K}_{i,B,D}A(1,U_2)^D=\Phi(1,U_2)^{c_i(B)}\cdot R_{i,B}(1,U_2),
\end{equation}
$$
\mathrm{with}\quad A(1,U_2)=\big(A_1(1,U_2),\dots,A_r(1,U_2)\big),\quad
\overline{K}_{i,B,D}=K_{i,B,D}\mod \fm\;\in k ,
$$
where $c_i(B)\in \bZ_{\geq 0}$ and $R_{i,B}(U_1,U_2)\in k[U_1,U_2]$ is
either $0$ or homogeneous and not divisible by $\Phi$ nor by $U_1$.
Choose a lift $\tilde{R}_{i,B}(U_1,U_2)\in R[U_1,U_2]$ of
$R_{i,B}(U_1,U_2))\in k[U_1,U_2]$ and set
$$
r_{i,B}=\tilde{R}_{i,B}(u_1,u_2)\in R, \quad
r'_{i,B}=r_{i,B}/u_1^{\deg R_{i,B}}\in R'.
$$
Then \eqref{eq.omega} implies
$$
\omega_{i,B} - (\phi')^{c_i(B)}r_{i,B}'\in \fm R'=\langle u_1\rangle\subset R'
$$
so that \eqref{eq.2star.on.9-1.5} gives
\begin{equation}\label{eq.2star.on.9-3}
f'_i = \tilde{F}_i(z')+\sum\limits_{|B| < n_i}
(z'-\mu)^B(\phi')^{(n_i-|B|)\beta'+c_i(B)}
u_1^{(n_i-|B|)(\delta-1)}\cdot r_{i,B}'+ h',
\end{equation}
where $v_\Lambda(h')_{(y',(u_1,\phi'))} > n_i$.
By Lemma \ref{lem1.1}(3) this implies
$v_\Lambda(h')_{(z',(u_1,\phi'))} > n_i$ by noting
$z'_i-y'_i\in \langle u_1^{\delta-1}\rangle$ by \eqref{star.on.6up}.
Now we need the following lemma:

\begin{lemma}\label{lem.claim5-2}
There exist $i \in \{1,\ldots, N\}$ and $B$ with $|B| < n_i$ such that
$R_{i,B}(U)\not\equiv 0$ in $k[U_1,U_2]$ (which is equivalent to
$r'_{i,B}\not\in \fm'$).
\end{lemma}
\medbreak
The proof of Lemma \ref{lem.claim5-2} will be given later.
Using the lemma, we see from \eqref{eq.2star.on.9-3} that there is a vertex of
$\Delta (f',z',(u_1,\phi'))$ on the line
$\{(a_1,a_2)\in \bR^2|\; a_1=\delta-1\}$.
Since $\Delta(f',z',(u_1,\phi'))\subset \Delta(f',y',(u_1,\phi'))$ and
$\bv(f',y',(u_1,\phi'))\notin \Delta(f',z',(u_1,\phi'))$, this implies
\begin{equation}\label{eq.claim5-3}
\begin{aligned}
\alpha(f',z',(u_1,\phi'))&=\alpha(f',y',(u_1,\phi'))=\delta-1,\\
\beta':=\beta(f',y',(u_1,\phi'))&<\beta(f',z',(u_1,\phi')).
\end{aligned}
\end{equation}
Moreover there exist $1\leq i\leq N$ and $B$ such that
\begin{equation}\label{eq.2star.on.9-5.5}
\beta(f',z',(u_1,\phi')) = \beta'+ \frac{c_i(B)}{n_i-|B|},
\end{equation}
where $c_i(B)$ is defined by equation \eqref{eq.3star.on.9-2}.

\begin{lemma}\label{lem.claim5-1}
If condition \eqref{fis} does not hold, $c_i(B)<n_i-|B|$ for all $i=1,\dots,N$.
\end{lemma}
\textbf{Proof }
By the assumption we have $\gamma^+(f,y,u)<d(\beta' +1)$.
Then we claim that for all $(i,B)$ such that $Q_{i,B}\not\equiv 0$,
\begin{equation*}\label{eq.2star.on.9-6}
d\cdot b_i(B) + \deg_{U_2} Q_{i,B}(1,U_2)<d\cdot (n_i-|B|)\, .
\end{equation*}
Indeed, recalling $b_i(B)=e_i(B)-(n_i-|B|)\beta'$, the assertion is
equivalent to
$$
\frac{d\cdot e_i(B)+\deg_{U_2} Q_{i,B}(1,U_2)}{n_i-|B|}< d(\beta'+1)
$$
and this follows from \eqref{eq.gammap}.
On the other hand, in \eqref{eq.3star.on.9-2} we have
$$
\deg_{U_2} A(1,U_2)^D \leq |D|\max\{\deg_{U_2} A_j(1,U_2)\;|\;
1\leq j\leq r\} < d\cdot (n_i-|B|)\;,
$$
because $|D|=n_i-|B|$ and $\deg_{U_2} A_j(1,U_2)= \alpha_j < d$ by (C1).
By \eqref{eq.3star.on.9-2}, this implies
$$
\deg_{U_2}\big(\Phi(1,U_2)^{c_i(B)} R_{i,B}(1,U_2)\big)< d\cdot (n_i-|B|),
$$
which implies $c_i(B)< n_i- |B|$ since $d=\deg_{U_2} \Phi(1,U_2)$.
$\square$
\medbreak

By Lemma \ref{lem.claim5-1}, \eqref{eq.2star.on.9-5.5} and \eqref{eq.claim5-3}
imply
\begin{equation}\label{eq.2star.on.9-4}
\beta'<\beta(f',z',(u_1,\phi'))<\beta'+1\,.
\end{equation}
By \eqref{betadeltaZ} we have $\beta'\in \bZ$ so that
$\beta(f',z',(u_1,\phi'))\notin\bZ$. Hence $(f',z',(u_1,\phi'))$
is not solvable at $\bv(f',z',(u_1,\phi'))$.
\medbreak

We now show $(f',z',(u_1,\phi'))$ is normalized at
$\bv(f',z',(u_1,\phi'))$. Let
$$
v''=\bv(f',z',(u_1,\phi'))\qaq \beta''=\beta(f',z',(u_1,\phi')).
$$
Setting $Z'=(Z'_1,\dots,Z'_r)$ with $Z'_i=in_{\fm'}(z'_i)$,
\eqref{eq.2star.on.9-3} implies
\begin{equation}\label{infz}
\begin{aligned}
in_{v''}(f'_i)_{(z',(u_1,\phi'))}
&= F_i(Z') + \underset{B}{\sum} (Z'-\overline{\mu})^B
\big({\Phi'}^{\beta''}U_1^{\delta-1}\big)^{n_i-|B|}
\cdot\overline{r}'_{i,B},\\
&= F_i(Z') + \underset{C}{\sum} {Z'}^C S_{i,C}\qwith
S_{i,C}\in k'[U_1,\Phi'],\\
\end{aligned}
\end{equation}
where the first sum ranges over such $B$ that
$|B|<n_i$ and $\displaystyle{\beta'+\frac{c_i(B)}{n_i-|B|}=\beta''}$
and
$$
\begin{aligned}
&\overline{\mu}=(\overline{\mu}_1,\dots,\overline{\mu}_r)\qwith
\overline{\mu}_i= U_1^{\delta-1}{\Phi}^{\beta'}\cdot \overline{a}'_i,\\
&\overline{r}'_{i,B}=r'_{i,B} \mod\fm' = R_{i,B}(1,U_2)\mod \Phi(1,U_2)\;\in
k'\simeq k[U_2]/\langle \Phi(1,U_2)\rangle.\\
\end{aligned}
$$
Let $B\in E^r(F_1,\ldots , F_{i-1})$.
Proposition \ref{pr.13} (c) implies $Q_{i,B}(U)\equiv 0$
(cf. \eqref{eq.13-4} and \eqref{eq.13-3.5}).
This implies that $(F_1,\ldots , F_N)$ is normalized in the sense of
Definition \ref{def.normalized0} (1), and Lemma \ref{lem.N} implies that
$K_{i,B,D}\in \fm$ for all $D$ in \eqref{eq.2star.on.9-1}.
Hence $R_{i,B}(1,U_2)\equiv 0$ in \eqref{eq.3star.on.9-2} so that
$\overline{r}'_{i,B}=0$ in \eqref{infz}.
By Lemma \ref{lem.N} this implies that $S_{i,C}\equiv 0$ if
$C\in E(F_1,\ldots , F_{i-1})$, which proves the desired assertion.
\medbreak

Finally it remains to show
\begin{equation}\label{eq.betabeta}
\betaob(f',z',(u_1,\phi'))<\betaob(f,y,u).
\end{equation}
The proof is divided into the following two cases:
\medbreak

\textbf{Case (1)}  $\deltaob(f,y,u)=\delta(f,y,u)$,

\medbreak
\textbf{Case (2)}  $\deltaob(f,y,u)< \delta(f,y,u)$.
\medbreak

Note that we always have $\deltaob(f,y,u)\leq \delta(f,y,u)$ since
$\Delta(f,y,u) \subset \Deltaob(f,y,u)$.
\medbreak

\textbf{Assume we are in Case (1).} We have
\begin{equation}\label{eq.case1-1}
\gamma^+(f,y,u)\leq \gammapob(f,y,u)\leq \betaob(f,y,u),
\end{equation}
where the first inequality holds by the assumption and the second holds
in general. By \eqref{eq.alphaob} and \eqref{eq.alphabeta}, the assumption
implies
\begin{equation*}\label{eq.case1-2}
\alphaob(f',y',(u_1,\phi')) = \alpha(f',y',(u_1,\phi'))=\delta-1.
\end{equation*}
Since the coordinate transformation $y'\to z'$ in \eqref{star.on.6up} affects
only those vertices lying in $\{(a_1,a_2)\in\bR^2|\; a_1\geq \delta-1\}$,
we have
\begin{equation}\label{eq.case1-3}
\alphaob(f',z',(u_1,\phi')) \geq \alphaob(f',y',(u_1,\phi'))
= \alpha(f',y',(u_1,\phi'))=\delta-1.
\end{equation}
By \eqref{eq.claim5-3} we have
$$
\alpha(f',y',(u_1,\phi'))= \alpha(f',z',(u_1,\phi'))\geq
\alphaob(f',z',(u_1,\phi')),
$$
where the inequality holds since
$\Delta(f',z',(u_1,\phi'))\subset \Deltaob(f',z',(u_1,\phi'))$.
Thus we get
\begin{equation}\label{eq.case1-4}
\begin{aligned}
&\alpha(f',z',(u_1,\phi'))=\alphaob(f',z',(u_1,\phi')),\\
&\betaob(f',z',(u_1,\phi'))\leq \beta(f',z',(u_1,\phi')).
\end{aligned}
\end{equation}
On the other hand, from Lemma \ref{lem.claim1} and \eqref{eq.2star.on.9-4},
we get
$$
\beta'':=\beta(f',z',(u_1,\phi'))<
\frac{\gamma^+(f,y,u)}{d} +1 \leq \frac{\gamma^+(f,y,u)}{2} +1.
$$
If $\gamma^+:=\gamma^+(f,y,u)\geq 2$, then $\beta''<\gamma^+$,
which shows the desired inequality \eqref{eq.betabeta}
thanks to \eqref{eq.case1-1} and \eqref{eq.case1-4}.
If $\gamma^+<2$, then Lemma \ref{lem.claim1} implies
$\beta':=\beta(f',y',(u_1,\phi'))\leq \gamma^+/2 < 1$ so that $\beta'=0$
since $\beta'\in\bZ$ as noted in \eqref{betadeltaZ}.
Hence \eqref{eq.case1-4} and \eqref{eq.2star.on.9-4} implies
$$
\betaob(f',z',(u_1,\phi'))\leq \beta(f',z',(u_1,\phi')) <\beta'+1=1.
$$
On the other hand we have
\begin{equation}\label{eq.case1-6}
\betaob(f,y,u) \geq \deltaob(f,y,u) - \alphaob(f,y,u) > \deltaob(f,y,u)-1=
\delta(f,y,u)-1 \geq 1.
\end{equation}
Here the first inequality holds in general, the second inequality holds
since $\alphaob(f,y,u)<1$ by Lemma \ref{case1.lem3} and the assumption (a),
and the last inequality follows from \eqref{betadeltaZ}.
This proves the desired inequality \eqref{eq.betabeta} in Case (1).
\medbreak

\textbf{Assume we are in Case (2).}
Write $\deltaob=\deltaob(f,y,u)$. The assumption implies
$$
\gammapob:=\gammapob(f,y,u)=\gammapob(y,u).
\quad\text{(cf. Definition \ref{loginv})}
$$
By \eqref{eq.alphabeta} and \eqref{eq.alphaob} it also implies
\begin{equation*}\label{eq.case2-0.2}
\alphaob(f',y',(u_1,\phi'))=\deltaob-1 < \alpha(f',y',(u_1,\phi'))=\delta-1.
\end{equation*}
Since the coordinate transformation $y'\to z'$ in \eqref{star.on.6up} affects
only those vertices lying in $\{(a_1,a_2)\in\bR^2|\; a_1\geq \delta-1\}$
and $\deltaob<\delta$, this implies
\begin{equation*}\label{eq.case2-0.3}
\alphaob(f',z',(u_1,\phi'))=\deltaob-1  < \alpha(f',z',(u_1,\phi')).
\end{equation*}
and hence
\begin{equation}\label{eq.case2-0.4}
\betaob(f',z',(u_1,\phi'))=\betaob(z',(u_1,\phi')).
\end{equation}
Recall that the $\delta$-face of $\Deltaob(f,y,u)$ is
\begin{equation*}\label{eq.edgeob}
\Deltaob(f,y,u)\cap \{A\in\bR^2|\; \Lob(A)=1\}\qwith
\Lob :\bR^2\to \bR; \; (a_1,a_2) \to \frac{a_1+a_2}{\deltaob}.
\end{equation*}
For $B\in \OBx$, the initial form of $l_B$ along this face is written as:
\begin{equation}\label{eq.case2-1}
\inLob(l_B)=L_B(Y) + \Phi(U)^{s_B} \Gamma_B(U),
\end{equation}
where $L_B(Y)\in k[Y]$ is a linear form, $s_B\in \bZ_{\geq 0}$, and
$\Gamma_B(U) \in k[U]$ is either $0$ or homogeneous of degree
$\deltaob-d \cdot s_B$ and not divisible by $\Phi$.
Note that $\Gamma_B(U)\not\equiv 0$ for some $B\in \OBx$.
From this we compute
\begin{equation}\label{eq.case2-2}
\begin{aligned}
\gammapob=\gammapob(y,u)
&=\sup \left\{\deg_{U_2}(\Phi(1,U_2)^{s_B}\Gamma_B(1,U_2))
\;|\; B\in \OBx, \; \Gamma_B(U)\not\equiv 0 \right\} \\
&\geq  d\cdot\left\{ s_B  \;|\;
B\in \OBx,\; \Gamma_B(U)\not\equiv 0 \right\}.\\
\end{aligned}
\end{equation}
Choose lifts $\tilde{L}_B(Y)\in R[Y]$ and $\tilde{\Gamma_B}(U)\in R[U]$
of $L_B(Y)\in k[Y]$ and $\Gamma_B(U)\in k[U]$, respectively, and set
$$
\gamma_B=\Gamma_B(u_1,u_2)\in R,\qaq \gamma_B'=\gamma_B/u_1^{\deg \Gamma_B}.
$$
Note $\gamma_B'\in \fm'$ if and only if $\Gamma_B(U)\equiv 0$.
\eqref{eq.case2-1} implies
\begin{equation*}\label{eq.case2-3}
l_B=\tilde{L}_B(y)+ \phi^{s_B}\gamma_B + \epsilon
\qwith v_{\Lob}(\epsilon)_{(y,u)}> v_{\Lob}(l_B)_{(y,u)}= 1.
\end{equation*}
By Lemma \ref{lem2.1} (1), this implies
\begin{equation*}\label{eq.case2-3.1}
l'_B=l_B/u_1 =\tilde{L}_B(y')+ u_1^{\deltaob-1}{\phi'}^{s_B} \gamma_B'
+ \epsilon' \qwith v_{\Lamdob}(\epsilon')_{(y',(u_1,\phi'))} >1,
\end{equation*}
$$\displaystyle{\Lamdob :\bR^2\to \bR
\;;\; (a_1,a_2) \to \frac{a_1}{\deltaob-1}}.$$
Substituting $y'=z'- u_1^{\delta-1} {\phi'}^{\beta'}\cdot a'$
(cf. \eqref{star.on.6up}), we get
\begin{equation*}\label{eq.case2-5}
l'_B=\tilde{L}_B(z')+ u_1^{\deltaob-1}{\phi'}^{s_B}  \gamma_B' +
u_1^{\delta-1}\cdot h + \epsilon' \quad (h\in R').
\end{equation*}
By Lemma \ref{lem1.1} (3), $v_{\Lamdob}(\epsilon')_{(y',(u_1,\phi'))}>1$
implies $v_{\Lamdob}(\epsilon')_{(z',(u_1,\phi'))}>1$ since
we have
$$
v_{\Lamdob}(z'_i-y'_i)_{(z',(u_1,\phi'))}=\frac{\delta-1}{\deltaob-1}>1
\qfor i=1,\dots,r
$$
by \eqref{star.on.6up}. Hence we get
\begin{equation*}
\betaob(f',z',(u_1,\phi'))=\betaob(z',(u_1,\phi'))
=\inf \left\{s_B \;|\; B\in \OBx,\; \Gamma_B(U)\not\equiv 0 \right\}
\end{equation*}
By \eqref{eq.case2-2} this implies
$$
\beta'':=\betaob(f',z',(u_1,\phi'))\leq \frac{\gammapob}{d}\leq
\frac{\gammapob}{2}.
$$
If $\gammapob\not=0$, this implies
$\beta''<\gammapob\leq \betaob(f,y,u)$ as desired.
If $\gammapob=0$, then $\beta''=0$.
On the other hand, as is seen in \eqref{eq.case1-6}, we have
$$
\betaob(f,y,u) \geq \deltaob(f,y,u) - \alphaob(f,y,u) > \deltaob(f,y,u)-1>0,
$$
where the last inequality was noted in \eqref{eq.eOB2}.
This proves the desired assertion \eqref{eq.betabeta} and the proof of
Lemma \ref{lem.claim5} is complete, up to the proof of Lemma \ref{lem.claim5-2}.
$\square$
\bigskip

\textbf{Proof of Lemma \ref{lem.claim5-2}}
Assume the contrary, i.e., that we have $R_{i,B}(U)\equiv 0$ in $k[U]$
 for all $(i,B)$ with $|B| < n_i$. Then, for all $(i,B)$ we have
\begin{equation}\label{eq.dot.on.9-7}
\Phi(1,U_2)^{b_i(B)}Q_{i,B}(1,U_2)=\sum\limits_{|D|=n_i-|B|}
\overline{K_{i,B,D}}A(1,U_2)^D\in k[U_2]
\end{equation}
Write
$$
\Gamma_j(U_1,U_2)=U_1^{\delta-(d \beta'+\alpha_j)}A_j(U_1,U_2)
\;\in\; k[U_1,U_2][U_1^{-1}]
$$
($\delta-(d\beta'+\alpha_j)$ may be negative).
Multiplying \eqref{eq.dot.on.9-7} by $U_1^{(n_i-|B|)(\delta-d \beta')}$,
we get
\begin{equation}\label{eq.2dot.on.9-7}
P_{i,B}(U_1,U_2) =  \sum\limits_{|D|=n_i-|B|}
\overline{K_{i,B,D}}\Gamma(U_1,U_2)^D
\end{equation}
where $\Gamma(U)^D = \prod\limits^r_{j=1}\Gamma_j(U)^{D_j}$.
In fact, for the left hand side we note that
$b_i(B)=e_i(B)-(n_i-|B|)\beta'$ and that
$Q_{i,B}(U_1,U_2)\Phi(U_1,U_2)^{e_i(B)}= P_{i,B}(U_1,U_2)$
is either $0$ or homogeneous of degree $\delta(n_i-|B|)$.
For the right hand side recall that
$A_j(U_1,U_2)=U_1^{\alpha_j}A_j(1,U_2)$ and that
$$
U_1^{|D|(\delta-d \beta')}=\prod\limits^r_{j=1}
U_1^{D_j(\delta-d \beta'-\alpha_j)}\cdot U_1^{D_j\cdot \alpha_j}\,.
$$
In view of \eqref{eq.2star.on.9-1},
equations \eqref{eq.2dot.on.9-7} and \eqref{eq.13-4} imply
\begin{equation}\label{eq.star.on.9-8}
\begin{aligned}
F_i(Y+\Gamma(U))  = & F_i(Y)+\sum\limits_{|B+D|=n_i}
(-1)^{|D|}\overline K_{i,B,D}Y^B\Gamma(U)^D\\
 = & F_i(Y)+\sum\limits_{|B|<n_i} P_{i,B}(U) =  in_{\delta}(f_i)
\end{aligned}
\end{equation}
Now we claim:
\begin{equation}\label{claim.13-1}
\delta\geq d \beta'+\alpha_j \text{ for $j=1,\dots,r$,}\quad\text{i.e.,}\quad
\Gamma_j(U)\in k[U] \subset k[U][U_1^{-1}].
\end{equation}
Admitting this claim, \eqref{eq.star.on.9-8} implies that one can dissolve
all the vertices of $\Delta(f,y,u)$ on the line
$\{A\in \bR^2|\; L(A)=1\}$, which contradicts assumption (b).
This completes the proof of Lemma \ref{lem.claim5-2}.

\medbreak\medbreak
It remains to show claim \eqref{claim.13-1}.
We show that $\Gamma_j(U)\in k[U_1.U_2] \subset k[U_1,U_2,U_1^{-1}]$.
Recall from \eqref{eq.star.on.9-8} that in any case
\begin{equation}\label{eq.Gamma}
F_i(Y + \Gamma(U)) = in_{\delta}(f_i) \in k[Y,U]\,.
\end{equation}
Denote the variables $Y_1,\ldots,Y_r,U_1,U_2$ by $X_1,\ldots,X_n$
(so that $n=r+2$),
and define derivations $D_A: k[X]\to k$ for $A\in \bZ^n_{\geq 0}$ by
$$
G(X+\rho) = \sum_A D_AG \cdot \rho^A\quad (G(X)\in k[X]),
$$
where $\rho=(\rho_1,\ldots,\rho_n)$ are new variables.
We now apply \cite{H5} (1.2) and \cite{Gi3} Lemmas 1.7 and 3.3.4
(the assumptions of the lemmas are satisfied by Lemma \ref{lem.N} and
the normalizedness of $(F_1,\dots,F_N)$ implied by Proposition \ref{pr.13}(c)).
According to these results,
after possibly changing the ordering of $X_1,\dots,X_n$, we can find:
\begin{itemize}
\item
$f$, an integer with $0 \leq f \leq n$,
\item
$P_j$ for $1\leq j \leq f$, homogeneous polynomials in the variables
$X_{A,i}$ indexed by $A\in \bZ^n_{\geq 0}$ and $1\leq i\leq N$,
with coefficients in $k$,
\item
$q_1 \leq q_2 \leq \ldots q_f$, numbers which are powers of the
exponential characteristic of $k$ (so that $q_j=1$ for $j=1,\dots,f$
if $\Char(k) = 0$),
\item
$\psi_j=c_{j,j+1} X_{j+1}^{q_j} + \ldots c_{j,r} X_n^{q_j}$
($c_{j,\nu} \in k,\; 1\leq j \leq f,\; j\leq \nu\leq n$),
additive polynomials homogeneous of degree $q_j$,
\end{itemize}
such that for $i=1,\dots, N$, we have
$$
\begin{matrix}
P_1(D_AF_i) & = & X_1^{q_1} & + &  \psi_1(X_2,\ldots,X_n) \\
P_2(D_AF_i) & = & X_2^{q_2} & + &  \psi_2(X_3,\ldots,X_n) \\
& \cdot & \\
& \cdot & \\
P_f(D_AF_i) & = & X_f^{q_f} & + & \psi_f(X_{f+1},\ldots,X_n)\,.
\end{matrix}
$$
Moreover the equations on the right hand side define the so-called ridge
(fa\^{\i}te in French) $F(C_x(X))=F(C(R/J))$ of the tangent cone
$C_x(X)=C(R/J)=\Spec(\grm(R/J))$, i.e., the biggest group subscheme
of $T_x(Z)=\Spec(\grm(R))=\Spec(k[X_1,\ldots, X_n])$ which respects
$C(R/J)$ with respect to the additive structure of $T_x(Z)$.
Since $\Dir(R/J) \subseteq F(C(R/J))$ and $e(R/J)=\eb(R/J)=2$
(cf. Definition \ref{def.dir3}) by the assumption,
all these schemes have dimension $2$. Hence we must have $f\geq n-2=r$.
Since $F_i(Y)\in k[Y]$, the variables $U_1$ and $U_2$ do not occur in
the above equations so that $f\leq r=n-2$. Thus we get $f=r$.
Hence, after a permutation of the variables $Y_1,\dots,Y_r$,
the equations become
$$
\begin{matrix}
P_1(D_AF_i(Y)) & = & Y_1^{q_1} & + &
\psi_1(Y_2,\ldots,Y_r) \\
P_2(D_AF_i(Y)) & = & Y_2^{q_2} & + &
\psi_2(Y_3,\ldots,Y_r) \\
& \vdots & \\
P_{r-1}(D_AF_i(Y)) & = & Y_{r-1}^{q_{r-1}} & + &
\psi_{r-1}(Y_r) \\
P_r(D_AF_i(Y)) & = & Y_r^{q_r} \,,
\end{matrix}
$$
This implies
$$
P_r(D_AF_i(Y + \Gamma(U))) = Y_r^{q_r} + \Gamma_r(U)^{q_r}\,.
$$
By \eqref{eq.Gamma} this gives
$\Gamma_r(U)^{q_r} \in k[U]$ and hence $\Gamma_r(U) \in k[U]$.
Noting
$\psi_j=c_{j,j+1} Y_{j+1}^{q_j} + \ldots + c_{j,r} Y_r^{q_j}$,
we easily conclude inductively from \eqref{eq.Gamma} that
$\Gamma_{r-1}(U), \ldots, \Gamma_1(U) \in k[U]$.
This completes the proof of claim \eqref{claim.13-1}.
$\square$

\newpage
\section{Non-existence of maximal contact in dimension $2$}
\bigskip

Let $Z$ be an excellent regular scheme and let $X\subset Z$ be a closed subscheme.

\begin{definition}\label{def.max.contact}
A closed subscheme $W\subset Z$ is said to have maximal contact with $X$ at
$x\in X$ if the following conditions are satisfied:
\begin{itemize}
\item[(1)]
$x\in W$.
\item[(2)]
Take any sequence of permissible blowups (cf. \eqref{seqBpermbu1}):
\begin{equation*}
\begin{array}{ccccccccccccccc}
  Z=&Z_0 & \lmapo{\pi_1} &   Z_1 & \lmapo{\pi_2} & Z_2 &
 \leftarrow\ldots \leftarrow & Z_{n-1} & \lmapo{\pi_n} & Z_{n}&& \gets \cdots\\
  &\cup & & \cup && \cup & & \cup && \cup \\
 X=&X_0 & \lmapo{\pi_1} &   X_1 & \lmapo{\pi_2} & X_2 &
 \leftarrow\ldots \leftarrow & X_{n-1} & \lmapo{\pi_n} & X_{n}&& \gets \cdots \\
\end{array}
\end{equation*}
where for any $n\geq 0$
$$
\begin{matrix}
Z_{n+1}&=& \Bl_{D_{n}}(Z_{n}) &\rmapo{\pi_{n+1}}& Z_{n} \\
&& \cup  && \cup\\
X_{n+1}&=& \Bl_{D_{n}}(X_{n}) &\rmapo{\pi_{n+1}}& X_{n} \\
\end{matrix}
$$
and $D_{n}\subset X_{n}$ is permissible. Assume that there exists a sequence of
points $x_n\in D_n$ ($n=0,1,\dots$) such that $x_0=x$ and $x_n$ is near to
$x_{n-1}$ for all $n\geq 1$. Then $D_n\subset W_n$ for all $n\geq 0$,
where $W_n$ is strict transform $W$ in $Z_n$.
\end{itemize}
\end{definition}

\begin{remark}\label{rem.max.contact}
The above definition is much weaker than Hironaka's original definition
(see \cite{AHV}).
\end{remark}

In this section we prove the following:

\begin{theorem}\label{thm.max.contact}
Let $k$ be a field of characteristic $p>0$ and let $y,u_1,u_2$ be three variables
over $k$. Consider $Z=\bA^3_k=\Spec(k[y,u_1,u_2])$ and let $X\subset Z$
be the hypersurface defined by the equation:
\begin{equation}\label{eq0.max.contact}
f=y^p+ y u_1^N u_2^N + u_1^a u_2^b(u_1+u_2)^{pA}
\end{equation}
where $A,a,b,N$ are integers satisfying the condition:
\begin{equation}\label{eq0-1.max.contact}
0<a,b<p,\; a+b=p,\; A>p,\; p\not|A,\; N \geq p^2 A.
\end{equation}
Let $x$ be the origin of $Z=\bA^3_k$ and let $U\subset Z$ be an open neighborhood of
$x$. Then there is no smooth hypersurface $W\subset U$ which has maximal contact
with $X\times_Z U$ at $x$.
\end{theorem}
\medbreak

We observe the following consequence of maximal contact, from which we will derive a contradiction.
Let $x\in X\subset Z$ be as in Definition \ref{def.max.contact}.
Let $R=\cO_{Z,x}$ with maximal ideal $\fm$ and residue field $k=R/\fm$, and
let $J\subset R$ be the ideal defining $X\times_Z\Spec(R)$.

\begin{lemma}\label{lem1.max.contact}
Assume $Z=\Spec(R)$. Set $e=e_x(X)$ and $r=\dim(R)-e$, and assume that $e\geq 1$.
Let $t_1,\dots,t_r$ be a part of a system of regular parameters of $R$ and
assume $W:=\Spec(R/\langle t_1,\dots,t_r\rangle)$ has maximal contact with $X$
at $x$. Let $(f,y,u)$ be a $\delta$-prepared label of $(X,Z)$ at $x$ (see
Definition \ref{def.label}). Assume $\delta(f,y,u)>2$, and
let $m\geq 2$ be the integer such that
\begin{equation}\label{eq1.max.contact}
m< \delta(f,y,u)\leq m+1.
\end{equation}
Then the following holds.
\begin{itemize}
\item[(1)]
$\langle t_1\dots,t_r\rangle = \langle y_1\dots,y_r\rangle$ in $R/\fm^2$.
In particular $(t,u)=(t_1,\dots,t_r,u_1,\dots,u_e)$ is a system of regular
parameters of $R$.
\item[(2)]
$m \leq \delta(f,t,u)\leq m+1$ so that $|\delta(f,y,u)-\delta(f,t,u)|\leq 1$.
\end{itemize}
\end{lemma}
\textbf{Proof }
Consider the fundamental sequence of permissible blowups as in Definition \ref{Def.fupb}
with $\cB=\emptyset$:
\begin{equation*}
\begin{array}{ccccccccccccccc}
Z= & Z_0 & \lmapo{\pi_1} &   Z_1 & \lmapo{\pi_2} & Z_2 &
\leftarrow\ldots \leftarrow & Z_{m-1} & \lmapo{\pi_m} & Z_m\\
& \cup & & \cup && \cup & & \cup && \cup \\
X= & X_0 & \lmapo{\pi_1} &   X_1 & \lmapo{\pi_2} & X_2 &
\leftarrow\ldots \leftarrow & X_{m-1} & \lmapo{\pi_m} & X_m\\
& \uparrow & & \cup && \cup & & \cup && \cup \\
& x & \leftarrow & C_1 & \stackrel{\sim}{\leftarrow} & C_2 &
\stackrel{\sim}{\leftarrow} \ldots \stackrel{\sim}{\leftarrow} & C_{m-1}
&\stackrel{\sim}{\leftarrow} & C_{m} & \\
\end{array}
\end{equation*}
By Corollary \ref{cor.fu} the integer $m$ in \eqref{eq1.max.contact} coincides
with the length of the above sequence. Let $W_q$ is the strict transform of $W$
in $Z_q$. By Definition \ref{def.max.contact} we have
\begin{equation}\label{eq2.max.contact}
C_q\subset W_q \;\text{ so that }\; W_{q+1}\simeq \Bl_{C_q}(W_q)
\text{  for all $q=1,\dots, m-1$}.
\end{equation}
For $1\leq q\leq m$, let $\eta_q$ be the generic point of $C_q$ and let
$\Retaq=\cO_{Z_q,\eta_q}$ with maximal ideal $\fmetaq$.
Let $\Jetaq\subset\Retaq$ be the ideal defining $X_q\times_{Z_q}\Spec(\Retaq)$. Write
$$
f_i^{(q)}=f_i/u_1^{q n_i},\; y_i^{(q)}=y_i/u_1^{q},\;
u_i'=u_i/u_1\; (2\leq i\leq e),\; t_i^{(q)}=t_i/u_1^{q},\; .
$$
By Claim \ref{fu.claim1}, we know
\begin{equation}\label{eq2-1.max.contact}
\text{$(f^{(q)},y^{(q)},u_1)$ is a $\delta$-prepared label of of $(X_q,Z_q)$
at $\eta_q$ for $q\leq m-1$.}
\end{equation}

\begin{claim}\label{claim1.max.contact}
$(t^{(q)},u_1)=(t^{(q)}_1,\dots, t^{(q)}_r,u_1)$ is a system of regular parameters
of $\Retaq$.
\end{claim}
\textbf{Proof }
Condition \eqref{eq2.max.contact} implies (by induction on $q$) that the strict transform
of $W$ in $\Spec(\Retaq)$ is defined by
$\langle t^{(q)}_1,\dots, t^{(q)}_r\rangle $.
It also implies that $W_q$ is transversal with the exception divisor of
$Z_q\to Z_{q-1}$ which is defined by $\langle u_1\rangle $ in $\Spec(\Retaq)$.
Thus the claim follows.
$\square$

\bigbreak

By claim \eqref{claim1.max.contact} we have
$\langle y^{(q)},u_1\rangle = \fmetaq=\langle t^{(q)},u_1\rangle\subset \Retaq$.
Lemma \ref{lem1.max.contact} (1) follows easily using this fact for $q=1$.
By Definition \ref{def.label} (2),
\eqref{eq2-1.max.contact} implies $\delta(f^{(q)},y^{(q)},u_1)>1$ so that
$$
f^{(q)}_j\in \fmetaq^{n_j}=\langle t^{(q)},u_1\rangle ^{n_j}\subset \Retaq
\quad (j=1,\dots, N)
$$
By Definition \ref{def.intial} (3) this implies
\begin{equation}\label{eq3.max.contact}
\delta(f^{(q)},t^{(q)},u_1) \geq 1 \quad\text{ if } 1\leq q\leq m-1.
\end{equation}
On the other hand, we have
\begin{equation}\label{eq4.max.contact}
\delta(f^{(q)},t^{(q)},u_1)=\delta(f,t,u)-q.
\end{equation}
Indeed, noting that $(t,u)$ is a system of regular parameters of $R$ by Lemma
\ref{lem1.max.contact} (1), we can write, as in \eqref{eq1.2.1}:
\begin{equation*}
f_i = \sum\limits_{(A, B)} C_{i, A, B} \; t^B u^A \quad \mbox{with}
\quad C_{i, A, B} \in R^{\times} \cup \{ 0 \}
\end{equation*}
We compute
\begin{equation*}
f^{(q)}_i=f_i/u_1^{n_j} = \sum\limits_{(A, B)} \big(C_{i, A, B}{u'}^A\big)  \;
(t^{(q)})^B u_1^{|A| + q(|B|-n_i)} ,
\end{equation*}
$$
{u'}^A={u'_2}^{a_2}\cdots {u'_e}^{a_e} \qfor A=(a_1,a_2,\dots,a_e).
$$
Noting that $k(\eta_q)=\Retaq/\fmetaq \simeq k(u')=k(u_1',\dots,u_e')$
and that $R\to \Retaq\to k(\eta_q)$ factors through $R\to k$
(see Claim \ref{fu.claim1}), we see
$$
\begin{aligned}
\delta(f^{(q)},t^{(q)},u_1)
&=\min\{\frac{|A|+q(|B|-n_i)}{n_i-|B|}\;|\; |B|<n_i,\; C_{i,A,B}\not=0 \}\\
&=\delta(f,y,u)-q.
\end{aligned}
$$

Now \eqref{eq3.max.contact} and \eqref{eq4.max.contact} give us
$\delta(f,t,u)\geq  m$. As $(f,y,u)$ is $\delta$-prepared,
Theorem \ref{thm.wellprepared2} implies $\delta(f,y,u)\geq \delta(f,t,u)$.
Thus we get $m \leq  \delta(f,t,u)\leq \delta(f,y,u)\leq m+1$ and
Lemma \ref{lem1.max.contact} (2) is shown.
$\square$
\bigskip

Now we start the proof of Theorem \ref{thm.max.contact}.
Let $x$ be the origin of $Z=\bA^3_k$ and let $R=\cO_{Z,x}$ with the maximal ideal
$\fm$. Write $u = (u_1,u_2)$.

\begin{claim}\label{claim2.max.contact}
$(f,y,u)$ is a $\delta$-prepared label of $(X,Z)$ at $x$ and
$\delta(f,y,u)=A+1$.
\end{claim}
\textbf{Proof }
From \eqref{eq0.max.contact} and \eqref{eq0-1.max.contact} one easily checks that $I\Dir_x(X)=\langle Y\rangle$
with $Y=\inm(y)$ and that $\Delta(f,y,u)\subset \bR^2$ is the polyhedron with the two vertices
$$
(A+\frac{a}{p},\frac{b}{p}),\quad (\frac{a}{p},A+\frac{b}{p})
$$
Since these vertices are not integral points, the polyhedron is $\delta$-prepared.
We have $\delta(f,y,u)=A+1>1$ as $a+b=p$ and the claim follows from
Definition \ref{def.label} (2).
$\square$
\medbreak

We will use three sequences of permissible blowups:
\medbreak\noindent
\textbf{Sequence I}: First consider $\pi:Z'=Bl_x(Z) \to Z$ and the strict
transform $X'\subset Z'$ of $X$. Look at the point $x'\in Z'$ of parameters
$$
\displaystyle{(y',u_1,v)=(\frac{y}{u_1},u_1,1+\frac{u_2}{u_1}}).
$$
Let $R'=\cO_{Z',x'}$. Then $X'\times_{Z'}\Spec(R')$ is defined by the equation:
$$
f':=\frac{f}{u_1^p}={y'}^p + y' u_1^{2N-p+1} (v-1)^N + u_1^{pA} v^{pA}(v-1)^b.
$$
Using \eqref{eq0-1.max.contact}, one can check that $\Delta(f',y',(u_1,v))$
has the unique vertex $(A,A)$ so that $\delta(f',y',(u_1,v))=2A>1$.
Hence $x'$ is very near to $x$ by Theorem \ref{deltanearcriteria}.
By Theorems \ref{prop2.1} and \ref{thm3.1}, $f$ is a $(u_1,v)$-standard base
and $(u_1,v)$  is admissible for $J=\langle f\rangle\subset R$.
On the other hand, $(f',y',(u_1,v))$ is not prepared at the vertex and we
dissolve it by the coordinate change
$$
z: = y' +\epsilon u_1^A v^A.\quad (\epsilon^p=(-1)^b)
$$
Setting $\lambda=(v-1)^N\in {R'}^\times$ and $\mu=v^{-1}\big((v-1)^b-\epsilon^p\big)\in {R'}^\times$, we compute
\begin{equation}\label{eq5.max.contact}
f'=z^p + \lambda z u_1^{2N-p+1}  + \mu u_1^{pA} v^{pA+1}
-\lambda \epsilon u_1^{2N-p+1+A} v^A .
\end{equation}
Using \eqref{eq0-1.max.contact}, we see that $\Delta(f',z,(u_1,v))$ is the polyhedron with three vertices
$$
(A,A+\frac{1}{p}),\quad (\frac{2N-p+1+A}{p},\frac{A}{p}),
\quad (\frac{2N}{p-1}-1,0)\,,
$$
and that the $\delta$-face of $\Delta(f',z,(u_1,v))$ is the first vertex.
In fact, we have
$$\frac{2N-p+1+A}{p}+\frac{A}{p} > A+A+ \frac{1}{p} \qaq
\frac{2N}{p-1}-1 > A+A+ \frac{1}{p}.$$
Since this first vertex is not an integral point, the polyhedron is $\delta$-prepared.
We have $\delta(f',z,(u_1,v))=2A+\frac{1}{p}>1$. Hence $(f',z,(u_1,v))$ is a
$\delta$-prepared label of $(X',Z')$ at $x'$ by Definition \ref{def.label} (2).
\medbreak

We now extend $\pi:Z'\to Z$ to the following sequence of permissible blowups:
\begin{equation*}
\begin{array}{ccccccccccccccc}
  Z& \lmapo{\pi}  &   Z' &  \lmapo{\pi_1} &  Z_1 & \lmapo{\pi_2} & Z_2 &
 \leftarrow\ldots \lmapo{\pi_p} & Z_{p}\\
\cup & & \cup && \cup & & \cup && \cup \\
 X& \lmapo{\pi}  &   X' &  \lmapo{\pi_1} &  X_1 & \lmapo{\pi_2} & X_2 &
 \leftarrow\ldots \lmapo{\pi_p} & X_{p}\\
 \uparrow & & \uparrow && \uparrow & & \uparrow && \uparrow \\
 x& \gets  &   x' &  \gets &  x_1 & \gets & x_2 &  \leftarrow\ldots\gets & x_{p}\\
\end{array}
\end{equation*}
Here $Z_1=\Bl_{x'}(Z')$ and $Z_{q+1}=\Bl_{x_q}(Z_q)$ ($q=1,\dots,p-1$)
where $x_q\in Z_q$ is the unique point lying on the strict transform of
$\{z=v=0\}\subset \Spec(R')$, and $X_q\subset Z_q$ is the strict transform of $X$.
For $q=1,\dots,p$, write $R_q=\cO_{Z_q,x_q}$ and
$$
f^{(q)}=f'/u_1^{pq}=f/u_1^{p(q+1)},\; z^{(q)}=z/u_1^{q},\;
v_i^{(q)}=v/u_1^{q} .
$$
By convention we write $x_0=x'$.

\begin{claim}\label{claim3.max.contact}
For $q=1,\dots,p$, $x_q$ is very near to $x_{q-1}$,
$(f^{(q)},z^{(q)},(u_1,v^{(q)}))$ is a $\delta$-prepared label of $(X_q,Z_q)$
at $x_q$, and the $\delta$-face of $\Delta(f^{(q)},z^{(q)},(u_1,v^{(q)}))$ is
the vertex
$$
(A+q(A-1+\frac{1}{p}),A+\frac{1}{p}).
$$
\end{claim}

\textbf{Proof }
By induction on $q$, one easily shows that $(z^{(q)},u_1,v^{(q)})$ is a system
of regular parameters of $R_q$. From \eqref{eq5.max.contact}, one computes
\begin{multline}\label{eq5-1.max.contact}
f^{(q)}=(z^{(q)})^p + \lambda z^{(q)} u_1^{2N-p+1+q(1-p)}  +
\mu u_1^{pA+q(pA+1-p)} (v^{(q)})^{pA+1} \\
-\lambda \epsilon u_1^{2N-p+1+A+q(A-p)} (v^{(q)})^A .
\end{multline}
This implies $f^{(q)}\in R_q$ and
$X_q\times_{Z_q}\Spec(R_q)=\Spec(R_q/\langle f^{(q)}\rangle)$.
It also implies the last assertion of the claim
by noting \eqref{eq0-1.max.contact}.
In fact, for the vertex coming from the second term we have
$$\frac{2N}{p-1}-1-q  >  A+q(A-1+\frac{1}{p})+A+\frac{1}{p}  \quad\mbox{for}\quad 0\leq q \leq p$$
(see the worst case $q=p$), and the case of the last term is even simpler.
Since the vertex in the claim is not an integral point, $(f^{(q)},z^{(q)},(u_1,v^{(q)}))$ is $\delta$-prepared. Noting $A>p$, we get
$$
\delta(f^{(q)},y^{(q)},(u_1^{(q)},u_2))=1+A+q(\frac{a}{p}-1)>1+a.
$$
Now the claim follows from Theorems \ref{deltanearcriteria}, \ref{prop2.1},
\ref{thm3.1} and Corollary \ref{cor.wellprepared2} as before.
$\square$

\medbreak\noindent
\textbf{Sequence II}: We consider the following sequence of permissible blowups:
\begin{equation*}
\begin{array}{ccccccccccccccc}
  Z&=Z_0& \lmapo{\pi_1}  &   Z_1 &  \lmapo{\pi_2} &  Z_2 &
 \leftarrow\ldots \lmapo{\pi_p} & Z_{p}\\
\cup && & \cup && \cup & &  \cup\\
 X&=X_0& \lmapo{\pi}  &   X_1 &  \lmapo{\pi_2} &  X_2 &
 \leftarrow\ldots \lmapo{\pi_p} & X_{p}\\
 \uparrow && & \uparrow && \uparrow & & \uparrow  \\
 x&=x_0& \gets  &   x_1 &  \gets &  x_2 & \leftarrow\ldots\gets & x_{p}\\
\end{array}
\end{equation*}
Here $Z_{q+1}=\Bl_{x_q}(Z_q)$ ($q=0,\dots,p-1$)
where $x_q\in Z_q$ is the unique point lying on the strict transform of
$\{y=u_2=0\}\subset \Spec(R)$, and $X_q\subset Z_q$ is the strict transform of $X$.
For $q=0,\dots,p$, write $R_q=\cO_{Z_q,x_q}$ and
$$
f^{(q)}=f/u_2^{pq},\; y^{(q)}=y/u_2^{q},\;
u_1^{(q)}=u_1/u_2^{q} .
$$

\begin{claim}\label{claim4.max.contact}
For $q=1,\dots,p$, $x_q$ is very near to $x_{q-1}$, and $(f^{(q)},y^{(q)},(u_1^{(q)},u_2))$ is
a $\delta$-prepared label of $(X_q,Z_q)$ at $x_q$, and
$\Delta(f^{(q)},y^{(q)},(u_1^{(q)},u_2))$ is the polyhedron
with the unique vertex
$$
(\frac{a}{p},\frac{b}{p}+A+q(\frac{a}{p}-1)).
$$
\end{claim}
\textbf{Proof }
By induction on $q$, one easily shows that $(y^{(q)},(u_1^{(q)},u_2))$ is a system
of regular parameters of $R_q$. From \eqref{eq0.max.contact}, one computes
\begin{equation*}
f^{(q)}=(y^{(q)})^p + y^{(q)} (u^{(q)}_1)^N u_2^{N+q(N+1-p)}  +
(u^{(q)}_1)^a u_2^{b+pA+q(a-p)} \big(u^{(q)}_1 u_2^{q-1}+1)^{pA}.
\end{equation*}
This shows $f^{(q)}\in R_q$ and it is an equation for $X_q\times_{Z_q}\Spec(R_q)$.
The last assertion easily follows by noting \eqref{eq0-1.max.contact}.
As $p\not|a$, the polyhedron is $\delta$-prepared. Then we get
$$
\delta(f^{(q)},y^{(q)},(u_1^{(q)},u_2))=1+A+q(\frac{a}{p}-1)>1+a.
$$
by noting $A>p\geq q$ and $p=a+b$.
Now the claim follows from Theorems \ref{deltanearcriteria}, \ref{prop2.1},
\ref{thm3.1} and Corollary \ref{cor.wellprepared2} as before.
$\square$

\medbreak\noindent
\textbf{Sequence III}: This is the sequence of permissible blowups,
which looks the same as the sequence II, except that now $x_q\in Z_q$ is the unique point
lying on the strict transform of $\{y=u_1=0\}\subset \Spec(R)$.

\medskip\medskip
Now assume that there exists $s\in \fm-\fm^2$ such that
\begin{enumerate}
\item[$(*)$]
$W:=\Spec(R/\langle s \rangle)\subset \Spec(R)$ has maximal contact with
$X\times_Z\Spec(R)$ at $x$.
\end{enumerate}
Then we want to deduce a contradiction.
Let $\hat{R}=k[[y,u_1,u_2]]$ be the completion of $R$.
It is easy to see that the assumption $(*)$ and Claims \ref{claim2.max.contact},
\ref{claim3.max.contact} and \ref{claim4.max.contact} hold even after
replacing $R$ with $\hat{R}$. Thus we may work with
$\hat{X}:=\Spec(\hat{R}/\langle f\rangle)$ and
$\hat{W}=\Spec(\hat{R}/\langle s \rangle)\subset \hat{Z}=\Spec(\hat{R})$
which has maximal contact with $\hat{X}$
at the closed point $x\in \Spec(\hat{R})$.

\begin{claim}\label{claim.ty}
There exists $t\in \hat{R}$ such that
$\langle s\rangle =\langle t\rangle \subset \hat{R}$ and
$$
t=y+\gamma\;\in \widehat{R} \quad\text{ where  }
\gamma\in \fn^{A+1}.
$$
Here $\fn=\langle u_1,u_2\rangle\subset k[[u]]=k[[u_1,u_2]]\subset \hat{R}$.
\end{claim}
\textbf{Proof }
By Lemma \ref{lem1.max.contact} (1) and Claim \ref{claim2.max.contact},
we can write
$\epsilon s = y+ a$ with $\epsilon\in k^\times$ and $a\in \fm^2$.
Noting $\hat{R}=k[[y,u_1,u_2]]=k[[t,u_1,u_2]]$, we can write
$$
a=\gamma + b \qwith \gamma\in \fn^2, \; b\in s\hat{R}\cap \fm^2\hat{R}.
$$
Setting $t=\epsilon s - b$, this implies $y=t+\gamma$ and
$\langle s\rangle =\langle t\rangle \in \hat{R}$.
From \eqref{eq0.max.contact} we get
\begin{equation}\label{eq6.max.contact}
f=t^p+ t u_1^N u_2^N + u_1^a u_2^b(u_1+u_2)^{pA}-\gamma^p - \gamma u_1^N u_2^N .
\end{equation}
From Definition \ref{def.intial} (3), this implies
$$
\delta(f,t,u)\leq c=\max\{i|\; \gamma\in \fn^i\}.
$$
Since $\delta(f,y,u)=A+1$ by Claim \ref{claim2.max.contact},
Lemma \ref{lem1.max.contact} (2) implies $A<\delta(f,t,u)$ and hence
$c\geq A+1$. This completes the proof of the claim.
$\square$

\bigbreak
To ease the notation, we write $R,X,Z$ for $\hat{R},\hat{X},\hat{Z}$ in what follows. Write
$$
\gamma=\Gamma +\theta \qwith \theta\in \fn^{A+2} \text{  and $\Gamma=0$ or
$\Gamma\in k[u_1,u_2]$ homogeneous of degree $A+1$,}
$$
and let $C=\max\{m\in \bN\cup\{\infty\} \;|\; (u_1+u_2)^m|\Gamma\}$.
There are two cases:

\medbreak\noindent
Case $C\not=A$: In this case \eqref{eq6.max.contact} implies
$$
f=t^p+ t u_1^N u_2^N + (u_1+u_2)^{pB}H(u_1,u_2) +\phi,
$$
where $B=\min\{C,A\}$, and $H(u_1,u_2)\in k[u_1,u_2]$ is homogeneous of degree $p(A+1-B)$ which is not divisible by $u_1+u_2$, and
$\phi\in \fn^{p(A+2)}$.
\medbreak

Now we consider the sequence I. Let $t'=t/u_1$ and $t^{(q)}=t'/u_1^q$.
By the same argument as in the proof of Claim \ref{claim1.max.contact},
we can show that $(t',u_1,v)$ (resp. $(t^{(q)},u_1,v^{(q)})$)
is a system of regular parameters for $R'=\cO_{Z',x'}$ (resp. $R_q=\cO_{Z_q,x_q}$).
We compute
$$
f'={t'}^p + \lambda t' u_1^{2N-p+1}+ \mu' u_1^{pA}v^{pB} +u_1^{p(A+1)}\phi',
$$
where $\lambda=(v-1)^N\in {R'}^\times$, and
$\mu'=H(1,u_2/u_1)\in k[v]\subset R'$, and
$\phi'=\phi/u_1^{p(A+2)}\in k[[u_1,u_2]][v]\subset \widehat{R'}$.
Note $\mu'\in {R'}^\times$ since $H(u_1,u_2)$ is not divisible by $u_1+u_2$.
From this we compute
$$
f^{(q)}= (t^{(q)})^p + \lambda t^{(q)} u_1^{2N-p+1+q(1-p)}+
\mu' u_1^{p(A+q(B-1))} (v^{(q)})^{pB} +u_1^{p(A+1-q)}\phi'.
$$
By Definition \ref{def.intial} (3), this implies
$$
\delta(f^{(q)},t^{(q)},(u_1,v^{(q)}))\leq A+q(B-1) + B.
$$
Setting $\sigma_q=
\delta(f^{(q)},z^{(q)},(u_1,v^{(q)}))-\delta(f^{(q)},t^{(q)},(u_1,v^{(q)}))$,
and noting $A\geq B$, we get
$$
\sigma_q \geq \big(A+q(A-1+\frac{1}{p})+A+\frac{1}{p}\big)-\big(A+q(B-1) + B\big)
\geq \frac{q+1}{p},
$$
from Claim \ref{claim3.max.contact}. Hence $\sigma_q>1$ for $q=p$, which contradicts
Lemma \ref{lem1.max.contact}(2) since $(f^{(q)},z^{(q)},(u_1,v^{(q)}))$ is a
$\delta$-prepared label of $(X_q,Z_q)$ at $x_q$ and the strict transform
$W_q\subset \Spec(R_q)$ of $W$ has maximal contact with
$X_q\times_{Z_q}\Spec(R_q)$ at $x_q$ by Definition \ref{def.max.contact}.

\begin{figure}[h]
\begin{center}
\scalebox{1} {\setlength{\unitlength}{0.00077865in}
\begingroup\makeatletter\ifx\SetFigFontNFSS\undefined%
\gdef\SetFigFontNFSS#1#2#3#4#5{%
  \reset@font\fontsize{#1}{#2pt}%
  \fontfamily{#3}\fontseries{#4}\fontshape{#5}%
  \selectfont}%
\fi\endgroup%
{\renewcommand{\dashlinestretch}{30}
\begin{picture}(7677,4407)(0,-10)
\path(465,4335)(465,330)(7665,330)
\path(420,1680)(510,1680)
\path(1815,375)(1815,285)
\path(420,2580)(510,2580)
\path(1815,4380)(1815,1680)
\path(1815,1680)(2265,1680)
\path(2145.000,1650.000)(2265.000,1680.000)(2145.000,1710.000)
\path(1815,2580)(3210,2580)
\path(3090.000,2550.000)(3210.000,2580.000)(3090.000,2610.000)
\path(1815,2580)(5865,1680)(7215,1680)
\path(7095.000,1650.000)(7215.000,1680.000)(7095.000,1710.000)
\path(5100,2445)(4830,1905)
\path(1995,1095)(1860,1545)
\path(1923.217,1438.681)(1860.000,1545.000)(1865.747,1421.440)
\dashline{60.000}(1815,1680)(2670,1140)(3345,1095)
\dottedline{45}(3300,1095)(3930,1095)
\put(1680,780){\makebox(0,0)[lb]{\smash{{\SetFigFontNFSS{11}{13.2}{\rmdefault}{\mddefault}{\updefault}``phantom polyhedron''$\Delta(f,t,(u_1,v))$}}}}
\put(60,1635){\makebox(0,0)[lb]{\smash{{\SetFigFontNFSS{11}{13.2}{\rmdefault}{\mddefault}{\updefault}$B$}}}}
\put(-90,2535){\makebox(0,0)[lb]{\smash{{\SetFigFontNFSS{11}{13.2}{\rmdefault}{\mddefault}{\updefault}$A+\frac 1 p$}}}}
\put(1725,15){\makebox(0,0)[lb]{\smash{{\SetFigFontNFSS{11}{13.2}{\rmdefault}{\mddefault}{\updefault}$A$}}}}
\put(2400,1635){\makebox(0,0)[lb]{\smash{{\SetFigFontNFSS{11}{13.2}{\rmdefault}{\mddefault}{\updefault}slower}}}}
\put(2220,2715){\makebox(0,0)[lb]{\smash{{\SetFigFontNFSS{11}{13.2}{\rmdefault}{\mddefault}{\updefault}faster under blow-ups}}}}
\put(4740,2580){\makebox(0,0)[lb]{\smash{{\SetFigFontNFSS{11}{13.2}{\rmdefault}{\mddefault}{\updefault}true polyhedron $\Delta(f,z,(u_1,v))$}}}}
\end{picture}
}}
\end{center}
\end{figure}

\bigskip\noindent
Case $C=A$: There exist
$c_1,c_2\in k$ with $c_1\not=c_2$ and $(c_1,c_2)\not=(0,0)$ such that
$$
\gamma=(u_1+u_2)^A(c_1 u_1+c_1 u_2) + \theta
\qwith \theta\in \fn^{A+2}\subset k[[u_1,u_2]].
$$
From \eqref{eq6.max.contact} we get
$$
f=t^p+ t u_1^N u_2^N + (u_1+u_2)^{pA}
\big(u_1^a u_2^b - (c_1u_1+c_2u_2)^p\big) + \phi,
$$
where $\phi=-\theta^p-\gamma u_1^Nu_2^N \in \fn^{p(A+2)}$.
\medbreak

Assume $c_2\not=0$ and consider the sequence II. Let $t^{(q)}=t/u_2^q$.
By the same argument as in the proof of Claim \ref{claim1.max.contact},
we can show that $(t^{(q)},u^{(q)}_1,u_2)$ is a system of regular parameters
of $R_q=\cO_{Z_q,x_q}$. We compute
$$
f^{(q)}= (t^{(q)})^p + t^{(q)} (u_1^{(q)})^N u_2^{N+q(N+1-p)} +
\xi u_2^{p(A+1-q)} + u_1^{p(A+2-q)}\phi',
$$
where $\phi'=\phi/u_2^{p(A+2)}\in R_q$ and
$$
\xi =(u_1^{(q)} u_2^{q-1} +1)^{pA}
\big((u_1^{(q)})^a u_2^{a(q-1)} - (c_1 u_1^{(q)} u_2^{q-1} +c_2)^p\big)\;\in
R_q^\times\;.
$$
Here we have used $a+b=p$. By Definition \ref{def.intial} (3), this implies
$$
\delta(f^{(q)},t^{(q)},(u_1^{(q)},u_2))\leq A+1-q.
$$
Set $\sigma_q=
\delta(f^{(q)},z^{(q)},(u_1,v^{(q)}))-\delta(f^{(q)},t^{(q)},(u_1,v^{(q)}))$.
Then Claim \ref{claim4.max.contact} implies
$$
\sigma_q \geq \big(1+ A + q(\frac{a}{p}-1)\big)-\big(A+1-q\big)
\geq \frac{qa}{p}.
$$
Hence $\sigma_q>1$ for $q=p$, which contradicts Lemma \ref{lem1.max.contact} (2)
by the same reason as in the previous case.
\medbreak

It remains to treat the case $c_1\not=0$. By symmetry, the proof is given
by the same argument using the sequence III instead of II.
This completes the proof of Theorem \ref{thm.max.contact}.

\newpage
\section{An alternative proof of Theorem \ref{thm.max-elim}}\label{sec:max-elim}

\bigskip

We give another proof of Theorem \ref{thm.max-elim} which uses more classical tools in
algebraic geometry.

\medbreak
Recall that the Hilbert polynomial $P=P(A)$ of a graded $k$-algebra $A$ (for a field $k$) is the unique
polynomial in $\bQ[T]$ with
$$
P(n) = H(A)(n) \quad\quad\mbox{  for  }\quad\quad n \gg 0\,.
$$
It is known that the degree of $P(A)$ is $\dim(A)-1$ where, by definition, the zero polynomial
has degree $-1$. The following property is well-known (see \cite{AP} (Corollary 3.2 and following
remarks for a proof).

\begin{lemma}\label{lem.AP}
A polynomial $P(T)\in \bQ[T]$ is the Hilbert polynomial of a standard graded $k$-algebra $A$ if and
only if
\begin{equation}\label{eq.AP}
P(T) = \binom{T+a_1}{a_1} + \binom{T+a_2-1}{a_2} + \ldots + \binom{T+a_i-i+1}{a_i} + \ldots
+ \binom{T+a_s-s+1}{a_s}
\end{equation}
for certain integers $a_1\geq a_2 \geq \ldots \geq a_s \geq 0$ with $s\geq 0$. Moreover,
one has $\deg(P) = a_1$ for $s\neq 0$, and the family $a(P) := (a_1,\ldots,a_s)$ is uniquely determined by $P$.
\end{lemma}

The equality $a_1 = \dim(A)-1$ follows from the fact that
$$
\binom{T+a_i-i+1}{a_i} = \frac{(T+a_i-i+1)(T+a_i-i)\cdots (T-i+2)}{a_i!}
$$
has degree $a_i$. The case $P(T)=0$ is obtained as the case where $s=0$ and $a(P)$ is the empty family.
The lemma also shows that the set $HP$ of all Hilbert polynomials of standard graded $k$-algebras
does not depend on $k$.

Recall that $HP$ is totally ordered with respect to the ordering
\begin{equation}\label{eq.ordHP}
P(T) \geq P'(T) \quad \Longleftrightarrow \quad P(n) \geq P'(n) \quad\quad\mbox{  for  }\quad\quad n \gg 0\,.
\end{equation}
We shall need the following description of this ordering.

\begin{lemma}\label{lem.ordHP} For two Hilbert polynomials $P, P' \in HP$ one has
$P \geq P'$ if and only if $a(P) \geq a(P')$ in the lexicographic ordering (where we
formally fill up the shorter family with entries $-\infty$ at the right until it has the
same length as the longer family).
\end{lemma}

\medbreak\noindent \textbf{Proof } Since $P = P' \Longleftrightarrow a(P) = a(P')$, and
since both orderings are total orderings, it suffices to show that $P > P'$ implies $a(P) > a(P')$.
Let $a(P) = (a_1,\ldots,a_s)$ and $a(P') = (b_1,\ldots,b_t)$, and assume $P > P'$.
We proceed by induction on min($s,t$). If min($s,t)=0$, then we must have
$t=0$ and $P \neq 0$, so that $a(P) > a(P')$. If $s,t\geq 1$ and $\deg(P) > \deg(P')$,
then $a_1 > b_1$ and hence $a(P) > a(P')$. If $a_1 =b_1$, then for
$$
P_1(T) = P(T) - \binom{T+a_1}{a_1} \quad, \quad P'_1(T) = P'(T) - \binom{T+a_1}{a_1}
$$
the polynomials $Q(T) = P_1(T+1)$ and $Q'(T) = P'_1(T+1)$ are again Hilbert polynomials
by Lemma \ref{lem.AP}, with associated invariants $a(Q) = (a_2,\ldots,a_s)$ and $a(Q')=(b_2,\ldots,b_t)$,
respectively, and $P>P'$ implies $Q > Q'$, so that $a(Q) > a(Q')$ in the lexicographic ordering
by induction. Together with $a_1 = b_1$ this implies the claim for $P$ and $P'$.

\bigbreak\noindent
Lemma \ref{lem.ordHP} immediately implies:

\begin{theorem}\label{thm.ordHP} In the ordered set $HP$, every
strictly descending sequence $P_1 > P_2 > P_3 > \ldots$ is finite.
\end{theorem}

Now let $HF \subset \bNN$ be the set of all Hilbert functions of standard graded $k$-algebras.
For $\nu\in HF$, write $P_\nu$ for the associated Hilbert polynomial. We trivially have
\begin{equation}\label{eq.ord.HF-HP}
\nu \geq \nu' \quad \Longrightarrow \quad P_\nu \geq P_{\nu'}\,.
\end{equation}

Let $X$ be an excellent scheme. For $x\in X$, let $P_X(x)$ be the Hilbert polynomial
associated to $\HSfX x$ and let
$$
\SigmaXP=\{\nu\in \Sigma_X|\; P_\nu=P_X^{\max}\}\qaq
\HSP X=\underset{\nu\in \SigmaXP}{\cup}\HSa X,
$$
where
$P_X^{\max} =\max\{P_X(x)\;|\; x\in X\}$.
Theorem \ref{HSf.usc} implies:

\begin{lemma}
$\HSP X$ is closed in $X$.
\end{lemma}

In fact, if $\nu \in \SigmaXP$ and $\mu \geq \nu$, then $\mu \in \SigmaXP$ by \eqref{eq.ord.HF-HP},
so that $\HSP X$ is the union of the finitely many closed sets $\HSgeqa X$ for $\nu \in \SigmaXP$.
Theorem \ref{thm.pbu.inv} (1) and \eqref{eq.ord.HF-HP} imply

\begin{lemma}\label{lem.pbu.HP}
If $\pi: X' \longrightarrow X$ is a permissible blow-up, then $P_{X'}^{\max} \leq P_X^{\max}$.
\end{lemma}

\textbf{Another proof of Theorem \ref{thm.max-elim}}
Suppose there exists an infinite sequence $X = X_0 \leftarrow X_1 \leftarrow \ldots $
of $\Sigmamax$-eliminations such that no $X_n$ is locally equisingular.
For each $n\geq 0$ let $X_n^0 \subseteq X_n$ be the union of those connected
components of $X_n$ which are not equisingular, and let $Y_n=\HSP {(X_n^0)}$.
Then $X_{n+1}^0 \subseteq \pi_n^{-1}(X_n^0)$,
so by Remark \ref{rem.equising1} (a) and Lemma \ref{lem.pbu.HP} we have
$P_{X_{n+1}^0}^{max} \leq P_{X_n^0}^{max}$.
By Theorem \ref{thm.ordHP} we may assume that
$P^{\max}_{X_n^0}=P^{\max}_{X_{n+1}^0}$ for all $n\geq 0$.
By Theorem \ref{thm.pbu.inv} this implies $\pi_n(Y_{n+1})\subset Y_n$, and we get
an infinite sequence of proper morphisms $Y=Y_0\gets Y_1\gets \cdots.$
By Lemma \ref{isp} below one can choose a sequence of points $x_n\in Y_n$
for $n=0,1,\dots$ such that $x_n=\pi_n(x_{n+1})$. By [H2] Theorem (1.B) and
Theorem \ref{thm.pbu.inv} (3) we may assume
$\HSf {X_n} {x_n}=\HSf {X_{n+1}} {x_{n+1}}$ for all $n\geq 0$, so that there
exists a $\nu_0\in HF$ such that
$\nu_0\in\Sigma_{X_n^0}^P$ for all $n\geq 0$. We claim that $\nu_0\in \Sigma_{X_n^0}^{max}$
for some $n$. Then \ref{def.max-elim1} (ME2) contradicts $\nu_0\in\Sigma_{X_{n+1}^0}$. Let
$S_n=\{\nu\in \Sigma_{X_n^0}^max\;|\; \nu>\nu_0\}$. We want to show $S_n=\emptyset$
for some $n$. Let $\nu_1,\dots,\nu_r$ be the elements of $S_0$ and put
$$
\Lambda=\{\mu\in HF\;|\; \nu_0 \leq \mu<\nu_i\text{ for some }i\in [1,r]\}.
$$
Thus the claim follows from the following:
\begin{lemma}
$\Lambda$ is finite, $S_n\subset \Lambda$, and
$S_n\cap S_m=\emptyset$ if $n\not=m\geq 0$.
\end{lemma}
\textbf{Proof }
The first claim follows from the assumption that $P_{\nu_i}=P_\nu$ so that
there exists $N>0$ such that $\nu_0(n)=\nu_i(n)$ for all $n\geq N$ and all
$i=1,\dots,r$. The second and third claim follow from \ref{def.max-elim1}
(ME3) and (ME2), respectively.
$\square$
\bigskip

\begin{lemma}\label{isp}
Let $Z_0 \lmapo{\pi_0}Z_1\lmapo{\pi_1} Z_2\lmapo{\pi_2}\cdots$ be an infinite
sequence of proper morphisms of non-empty noetherian schemes.
Then there exists a sequence of points
$x_n\in Z_n$ ($n=0,1,2,\dots$) such that $x_n=\pi_n(x_{n+1})$.
\end{lemma}
\textbf{Proof }
For $m>n\geq 0$, let $\pi_{n,m}:Z_m\to Z_n$ be the composite of $\pi_i$ for
$i=n,\dots,m-1$ and put $Z_{n,m}=\pi_{n,m}(Z_m)$.
By the properness, $Z_{n,m}$ is non-empty and closed in $Z_n$.
Clearly we have $Z_n\supset Z_{n,n+1}\supset Z_{n,n+2}\supset \cdots$ and
\begin{equation}\label{Ynml}
Z_{n,l}=\pi_{n,m}(Z_{m,l})\qfor l>m>n.
\end{equation}
Put
$\displaystyle{Z_{n,\infty}=\underset{m>n}{\cap} Z_{n,m}}$.
By the Noetherian condition, there exists $N(n)>n$ such that $Z_{n,\infty}=Z_{n,N(n)}$
so that $Z_{n,\infty}\not=\emptyset$. Then \eqref{Ynml} implies
$Z_{n,\infty}=\pi_{n,m}(Z_{m,\infty})$ for $m>n$. Thus we get an infinite
sequence of proper surjective morphisms
$Z_{0,\infty} \gets Z_{1,\infty} \gets Z_{2,\infty} \gets \cdots$,
which inductively gives the desired claim.
$\square$

\newpage
\section{Functoriality, locally noetherian schemes, algebraic spaces and stacks}\label{sec:funct}

In this section we reformulate the obtained functoriality for our resolution
and apply this to obtain resolution of singularities for (two-dimensional)
locally noetherian excellent schemes, algebraic spaces and algebraic stacks.

\begin{definition}\label{def.res.morphism}
Let $X$ be a noetherian excellent scheme of dimension at most two.
The resolution morphism
$$
r = r(X): X' = R(X) \longrightarrow X
$$
is defined as the composition of the morphisms $X' = X_n \rightarrow \ldots X_1 \rightarrow X_0 = X$
in the canonical resolution sequence of Theorem \ref{thm.intro.1} (or, what amounts to the same,
of the morphisms in the canonical resolution sequence $S_0(X,\emptyset)$ of Theorem \ref{thm.main.1} (b)).
\end{definition}

Then we can state the functoriality in the mentioned theorems as follows.

\begin{theorem}\label{thm.funct.1}
(a) If $\alpha: U \rightarrow X$ is an \'etale morphism, then there is a unique morphism
$\alpha' = R(\alpha): R(U) \rightarrow R(X)$ making the diagram
$$
\begin{CD}
R(U) @>{R(\alpha)}>> R(X) \\
@V{r(U)}VV @VV{r(X)}V \\
U @>{\alpha}>> X
\end{CD}
$$
commutative. Moreover, the diagram is cartesian.

(b) The same functoriality holds for a boundary $\cB$ on $X$ and the analogous resolution morphism
$r(X,\cB): X'' = R(X,\cB) \rightarrow X$ obtained from Theorem \ref{thm.main.1} by composing the
morphisms in $S(X,\cB)$, and for an embedded situation $(X,Z,\cB)$ as in the Theorem \ref{thm.main.3}
and the resolution morphism $r(X,Z,\cB): (X''',Z''',\cB''') \rightarrow (X,Z,\cB)$ obtained by
composing the morphisms in the sequence $S(X,Z,\cB)$ of that theorem.

(c) More generally, let $X$ be of arbitrary dimension, let $\cB$ be a boundary on $X$, let
$\alpha: X^\ast \rightarrow X$ be a flat morphism with geometrically regular fibers (e.g., a smooth morphism), and let $\cB^\ast$
be the pull-back of $\cB$ to $X^\ast$. Assume that the canonical resolution
sequence constructed by Remark \ref{rem.strategy} and Corollary \ref{cor.Omax-elim1} for $(X,\cB)$ is finite.
Then the corresponding sequence for $(X^\ast,\cB^\ast)$ is finite as well, and there is a unique morphism
$\alpha'=R(\alpha): R(X^\ast,\cB^\ast) \rightarrow R(X,\cB)$ making the diagram
$$
\begin{CD}
(X^\ast)'=R(X^\ast,\cB^\ast) @>{\alpha'=R(\alpha)}>> R(X,\cB) = X' \\
@V{r^\ast := r(X^\ast,\cB^\ast)}VV @VV{r := r(X,\cB)}V \\
X^\ast @>{\alpha}>> X
\end{CD}
$$
commutative, where the vertical morphisms are the compositions of the morphisms in the respective
resolution sequences. Moreover, the diagrams are cartesian, and the boundary $(\cB^\ast)'$ on $(X^\ast)'$
is the pull-back of the boundary $\cB'$ on $X'$.
\end{theorem}

\textbf{Proof } The finiteness of the resolution sequence for $(X^\ast,\cB^\ast)$, the existence
of $\alpha'$ and the fact that the diagram is cartesian follow from Proposition \ref{prop.funct.strat} (ii).
Since $r(X,\cB)$ is the composition of a sequence
$X'=X_n \rightarrow X_{n-1} \rightarrow \ldots \rightarrow X_1 \rightarrow X_0 = X$
of blow-ups, the uniqueness follows from the following Lemma.

\begin{lemma}\label{lem.funct}
Consider a cartesian diagram of schemes
$$
\begin{CD}
(Y)' @>{\alpha'}>> X' \\
@V{\rho}VV @VV{r}V  \\
Y @>{\alpha}>> X \;,
\end{CD}
$$
where $\alpha$ is flat and $r$ is a composition of blow-ups.
Then $\alpha'$ is the unique morphism making the diagram commutative.
\end{lemma}

\textbf{Proof } Let $\beta$ be another morphism making the diagram
commutative. Write $\pi$ as the composition of blow-ups
$X' = X_n \rightarrow X_{n-1} \rightarrow \ldots \rightarrow X_1 \rightarrow X_0 = X$. Then we get a
commutative diagram with cartesian squares
$$
\begin{array}{rcccl}
Y' = & Y_n & \mathop{\longrightarrow}\limits^{\alpha_n=\alpha'} & X_n & = X' \\
& \downarrow && \downarrow &  \\
& \vdots && \vdots & \\
& Y_{m+1} & \mathop{\longrightarrow}\limits^{\alpha_{m+1}} & X_{m+1} \\
& \downarrow && \downarrow\rlap{$\pi_m$} & \\
& Y_{m} & \mathop{\longrightarrow}\limits^{\alpha_{m}} & X_{m} \\
& \downarrow && \downarrow & \\
& \vdots && \vdots & \\
& \downarrow && \downarrow & \\
Y = & Y_0 & \mathop{\longrightarrow}\limits^{\alpha_0=\alpha} & X_0 & = X'\;.
\end{array}
$$
Since $\alpha$ is flat, all $\alpha_m$ are flat as well. Hence, if
$\pi_m: X_{m+1} \rightarrow X_m$ is the blow-up in $D_n$, the morphism
$Y_{m+1} \rightarrow Y_m$ identifies with the blow-up in $Y_m\times_{X_m} D_m$.
Write $r_m: X_n \rightarrow X_m$ for the composition in the right column
and $\rho_m: Y_n \rightarrow Y_m$ for the composition in the left column.
We show inductively that
\begin{equation}\label{eq.alpha.beta.rho}
r_m \beta = r_m \alpha'\;,
\end{equation}
which gives the claim for $m=n$, where $r_n=id$. By assumption \eqref{eq.alpha.beta.rho} holds for $m=0$,
since $r_0 \beta = r\beta = \alpha \rho = r\alpha' = r_0 \alpha'$. If equation \eqref{eq.alpha.beta.rho} holds
for $m$, then we show it for $m+1$ as follows. If $J_m$ is the ideal sheaf of
$D_m$ in $X_m$, then $L_1 = \pi_m^{-1}J_m\cO_{X_{m+1}}$ is an invertible sheaf on $X_{m+1}$ by
the universal property of blow-ups. Since $\alpha_{m+1}$ is flat, $L_2 = \alpha_{m+1}^{-1}L_1\cO_{Y_{m+1}}$
is invertible on $Y_{m+1}$. Then $L_3 = \rho_m^{-1}L_2\cO_{Y_n}$ is invertible on $Y_n$.
In fact, for any blow-up morphism $f: Z' \rightarrow Z$ it is obvious that
$f^{-1}M\cO_{Z'}$ is invertible for any invertible ideal sheaf $M$ on $Z$. Consider the morphism
$$
\varphi = r_m \beta = r_m \alpha' = \alpha_m \rho_m = \pi_m\alpha_{m+1}\rho_{m+1}: Y_n \longrightarrow X_m\,.
$$
Since we have shown that $\varphi^{-1}J_m\cO_{Y_n}$ is invertible, it
follows by the universal property of the blow-up $\pi_m: X_{m+1} \rightarrow X_m$
that there is a unique morphism $\psi: Y_n \rightarrow X_{m+1}$ with $\pi_m \psi = \varphi$.
Since
$$
\pi_m r_{m+1} \beta = r_m \beta = \varphi = r_m \alpha' = \pi_m r_{m+1} \alpha'\,,
$$
we get $r_{m+1} \beta = r_{m+1} \alpha'$ as wanted.

\bigskip
First we apply this to obtain resolution of singularities for excellent schemes
of dimension $\leq 2$ which are only locally noetherian, and not necessarily noetherian.

\begin{theorem}\label{thm.loc.noeth.schemes}
Let $X$ be an excellent scheme of dimension $\leq 2$. Then there exists a canonical
resolution morphisms
$$
r = r(X): X' = R(X) \longrightarrow X
$$
for $X$, and
$$
r(X,\cB): X'' = R(X,\cB) \rightarrow X
$$
for each boundary $\cB$ on $X$, and
$$
r(X,Z,\cB): (X''',Z''',\cB''') \rightarrow (X,Z,\cB)
$$
for each embedded situation $(X,Z,\cB)$ with regular excellent $Z$ and normal crossings divisor $\cB$
on $Z$, which satisfy the same properties as the corresponding morphisms in Theorem \ref{thm.funct.1}.
\end{theorem}

\textbf{Proof} This follows from Theorem \ref{thm.funct.1} by looking at an open cover $(U_i)_{i\in I}$
of $X$ by noetherian (excellent) schemes and gluing the corresponding morphisms for the $U_i$ by the
uniqueness on the intersections. By the same argument, the uniqueness and functorial properties follow.

\bigskip
Finally we apply this to get resolution for algebraic stacks of dimension $\leq 2$.

\begin{theorem}\label{thm.stacks}
Let $\cX$ be an excellent algebraic stack for a Grothendieck topology whose coverings are
flat morphisms with geometrically regular fibers (e.g., a Deligne-Mumford stack or an Artin stack).
Assume that there is a representable covering morphism $X \rightarrow \cX$
in the chosen topology such that $X$ is an excellent scheme with $\dim(X) \leq 2$. Then there exists
a proper surjective morphism $\cX' \rightarrow \cX$ of stacks such that $\cX'$ is regular.
\end{theorem}

\textbf{Proof } This follows from theorem \ref{thm.loc.noeth.schemes} and the functoriality obtained
in there, in the same way as in \cite{Te} section 5.2:
By the assumption one may represent $\cX$ by a groupoid scheme
$$
R \mathop{\rightrightarrows}\limits^{s}\limits_{t} X
$$
(with the usual further structural morphisms $m, e$ and $i$) such that $R$ and $X$ are excellent schemes
and the morphisms $s$ and $t$ are flat with regular fibers. By the functorial resolution for $X$
(which then gives finiteness for the canonical resolution strategy for $R$ as well,
by \ref{thm.funct.1} (c) and using $s$, say) one defines the structural morphism $s'$ by the
cartesian diagram
$$
\begin{CD}
R' @>{s'}>>  X \\
@V{r(R)}VV    @VV{r(X)}V \\
R  @>{s}>>  X
\end{CD}
$$
with the vertical resolution morphisms, in the same way one defines $t'$ from $t$, and
similarly one defines the `composition' structural morphism $m'$ by the cartesian diagram
$$
\begin{CD}
R'\times_{t' X' s'} R' = (R \times_{t X s} R)'   @>{m'}>>  R' \\
@V{r(R\times_{tXs}R)}VV            @VV{r(R)}V    \\
 R \times_{t X s} R   @>{m}>>    R
\end{CD}
$$
with vertical resolution morphisms, since $m$ is flat with regular fibers as well.
The further structural morphisms $e$, $i$ are obtained via base change, and the rules
for $s',t',m',e'$ and $i'$ follow by the functoriality from those for $s,t,m,e$ and $i$,
because the equalities are checked on flat morphisms with regular fibers.

\bigskip
In the special case of an algebraic space one sees that one gets a resolution by
an algebraic space again.

\end{document}